\pdfoutput=1
\RequirePackage{ifpdf}
\ifpdf 
\documentclass[pdftex]{sigma}
\else
\documentclass{sigma}
\fi

\numberwithin{equation}{section}

\usepackage{tikz-cd,enumitem}

\newtheorem{thm}{Theorem}[section]
\newtheorem{defn}[thm]{Definition}
\newtheorem{prop}[thm]{Proposition}
\newtheorem{cor}[thm]{Corollary}
\newtheorem{lem}[thm]{Lemma}
\newtheorem{conj}[thm]{Conjecture}

\theoremstyle{definition}
\newtheorem{rema}[thm]{Remark}

\newtheorem{assum}[thm]{Assumption}

\newcommand{\nn}{\nonumber \\}

 \newcommand{\res}{\operatorname{Res}}

\newcommand{\wt}{\mbox{\rm wt}\,}

\newcommand{\mbar}{\big\vert}
\newcommand{\lbar}{\big\vert}

\newcommand{\Y}{\mathcal{Y}}
\newcommand{\C}{\mathbb{C}}
\newcommand{\Z}{\mathbb{Z}}
\newcommand{\R}{\mathbb{R}}

\newcommand{\N}{\mathbb{N}}

\newcommand{\I}{\mathbb{I}}

\renewcommand{\i}{\mathrm{i}}
\newcommand{\PB}{{\rm PB}}

 \makeatletter
\newlength{\@pxlwd} \newlength{\@rulewd} \newlength{\@pxlht}
\catcode`.=\active \catcode`B=\active \catcode`:=\active \catcode`|=\active
\def\sprite#1(#2,#3)[#4,#5]{
 \edef\@sprbox{\expandafter\@cdr\string#1\@nil @box}
 \expandafter\newsavebox\csname\@sprbox\endcsname
 \edef#1{\expandafter\usebox\csname\@sprbox\endcsname}
 \expandafter\setbox\csname\@sprbox\endcsname =\hbox\bgroup
 \vbox\bgroup
 \catcode`.=\active\catcode`B=\active\catcode`:=\active\catcode`|=\active
 \@pxlwd=#4 \divide\@pxlwd by #3 \@rulewd=\@pxlwd
 \@pxlht=#5 \divide\@pxlht by #2
 \def .{\hskip \@pxlwd \ignorespaces}
 \def B{\@ifnextchar B{\advance\@rulewd by \@pxlwd}{\vrule
 height \@pxlht width \@rulewd depth 0 pt \@rulewd=\@pxlwd}}
 \def :{\hbox\bgroup\vrule height \@pxlht width 0pt depth
0pt\ignorespaces}
 \def |{\vrule height \@pxlht width 0pt depth 0pt\egroup
 \prevdepth= -1000 pt}
 }
\def\endsprite{\egroup\egroup}
\catcode`.=12 \catcode`B=11 \catcode`:=12 \catcode`|=12\relax
\makeatother

\makeatletter
\newcommand{\raisemath}[1]{\mathpalette{\raisem@th{#1}}}
\newcommand{\raisem@th}[3]{\raisebox{#1}{$#2#3$}}
\makeatother

\makeatletter
\newcommand{\subalign}[1]{%
	\vcenter{%
		\Let@ \restore@math@cr \default@tag
		\baselineskip\fontdimen10 \scriptfont\tw@
		\advance\baselineskip\fontdimen12 \scriptfont\tw@
		\lineskip\thr@@\fontdimen8 \scriptfont\thr@@
		\lineskiplimit\lineskip
		\ialign{\hfil$\m@th\scriptstyle##$&$\m@th\scriptstyle{}##$\hfil\crcr
			#1\crcr
		}%
	}%
}
\makeatother

\def\hboxtr{\FormOfHboxtr} 
\sprite{\FormOfHboxtr}(25,25)[0.5 em, 1.2 ex] 

:BBBBBBBBBBBBBBBBBBBBBBBBB |
:BB......................B |
:B.B.....................B |
:B..B....................B |
:B...B...................B |
:B....B..................B |
:B.....B.................B |
:B......B................B |
:B.......B...............B |
:B........B..............B |
:B.........B.............B |
:B..........B............B |
:B...........B...........B |
:B............B..........B |
:B.............B.........B |
:B..............B........B |
:B...............B.......B |
:B................B......B |
:B.................B.....B |
:B..................B....B |
:B...................B...B |
:B....................B..B |
:B.....................B.B |
:B......................BB |
:BBBBBBBBBBBBBBBBBBBBBBBBB |

\endsprite

\sprite{\FormOfShboxtr}(25,25)[0.3 em, 0.72 ex] 

:BBBBBBBBBBBBBBBBBBBBBBBBB |
:BB......................B |
:B.B.....................B |
:B..B....................B |
:B...B...................B |
:B....B..................B |
:B.....B.................B |
:B......B................B |
:B.......B...............B |
:B........B..............B |
:B.........B.............B |
:B..........B............B |
:B...........B...........B |
:B............B..........B |
:B.............B.........B |
:B..............B........B |
:B...............B.......B |
:B................B......B |
:B.................B.....B |
:B..................B....B |
:B...................B...B |
:B....................B..B |
:B.....................B.B |
:B......................BB |
:BBBBBBBBBBBBBBBBBBBBBBBBB |

\endsprite

\begin{document}

\allowdisplaybreaks

\newcommand{\arXivNumber}{2501.15003}

\renewcommand{\thefootnote}{}

\renewcommand{\PaperNumber}{086}

\FirstPageHeading

\ShortArticleName{Twisted Intertwining Operators and Tensor Products}

\ArticleName{Twisted Intertwining Operators and Tensor Products\\ of (Generalized) Twisted Modules\footnote{This paper is a~contribution to the Special Issue on Recent Advances in Vertex Operator Algebras in honor of James Lepowsky. The~full collection is available at \href{https://sigma-journal.com/Lepowsky.html}{https://sigma-journal.com/Lepowsky.html}}}

\Author{Jishen DU and Yi-Zhi HUANG}

\AuthorNameForHeading{J.~Du and Y.-Z.~Huang}

\Address{Department of Mathematics, Rutgers University,\\
110 Frelinghuysen Rd., Piscataway, NJ 08854-8019, USA}
\Email{\mail{jd1422@scarletmail.rutgers.edu}, \mail{yzhuang@math.rutgers.edu}}
\URLaddress{\url{https://jishendumath.github.io/}, \url{https://sites.math.rutgers.edu/~yzhuang/}}

\ArticleDates{Received July 09, 2025, in final form July 30, 2026; Published online September 02, 2026}

\Abstract{We study the general twisted intertwining operators (intertwining operators among twisted modules) for a~vertex operator algebra $V$. We give the skew-symmetry and contragredient isomorphisms between spaces of twisted intertwining operators and also prove some other properties of twisted intertwining operators. Using twisted intertwining operators,
we introduce a~notion of $P(z)$-tensor product of two objects for $z\in \C^{\times}$ in a~category of suitable $g$-twisted $V$-modules for $g$ in a~group of automorphisms of $V$ and give a~construction of such a~$P(z)$-tensor product under suitable assumptions. We also construct $G$-crossed commutativity isomorphisms and $G$-crossed braiding isomorphisms. We formulate a~$P(z)$-compatibility condition and a~$P(z)$-grading-restriction condition and use these conditions to give another construction of the $P(z)$-tensor product.}


\Classification{17B69; 18M15; 81T40}

\renewcommand{\thefootnote}{\arabic{footnote}}
\setcounter{footnote}{0}

\section{Introduction}

Modular tensor categories associated to conformal field theories
were discovered first in physics by
Moore and Seiberg~\cite{MS}. In~\cite{T1}, Turaev formulated a~precise notion of modular tensor category based on his joint
work~\cite{RT} with Reshetikhin on the construction of quantum invariants
of three manifolds using representations of quantum groups.
In~\cite{H-rigidity}, the second author proved the following theorem.

\begin{thm}\label{mod-tensor-cat}
Let $V$ be a~simple vertex operator algebra satisfying the following conditions:
\begin{enumerate}[label=$(\arabic*)$]\itemsep=0pt
\item For $n<0$, $V_{(n)}=0$ and $V_{(0)}=\mathbb{C}\mathbf{1}$
and as a~$V$-module, $V$ is equivalent to its
contragredient $V$-module $V'$ $($or equivalently, there exists a~nondegenerate invariant bilinear form on~$V)$.
\item Every lower-bounded $($generalized$)$ $V$-module is completely reducible.
\item $V$ is $C_{2}$-cofinite.
\end{enumerate}
Then the category of $V$-modules has a~natural structure of
modular tensor category in the sense of Turaev~{\rm~\cite{T1}}.
\end{thm}
The proof of Theorem~\ref{mod-tensor-cat} was
based on the results obtained by Lepowsky and the second author in
\cite{tensor1,tensor2,tensor3} and the results
obtained by the second author in
\cite{tensor4,H-diff-eqn,H-modular,H-verlinde-conj}.

It is natural to expect that Theorem~\ref{mod-tensor-cat} has generalizations
in two-dimensional orbifold conformal field theory.
Two-dimensional orbifold conformal field theories are two-dimensional
conformal field theories constructed from
known theories and their automorphisms. The first example of
two-dimensional orbifold conformal field theories is
the moonshine module constructed by
Frenkel, Lepowsky and Meurman~\cite{FLM1,FLM2,FLM3}
in mathematics. In string theory, the systematic study of two-dimensional
orbifold conformal field theories was started by Dixon,
Harvey, Vafa and Witten
\cite{DHVW1,DHVW2}. See~\cite{H-orbifold} for
an exposition on general results, conjectures and open problems
in the construction of two-dimensional orbifold conformal field
theories using the
approach of the representation theory of vertex operator algebras.

In~\cite{K3}, Kirillov Jr.\ stated that the category of $g$-twisted modules
for a~vertex operator algebra $V$ for $g$ in a~finite subgroup $G$
of the automorphism group of $V$ is a~$G$-equivariant fusion category
($G$-crossed braided (tensor) category in the sense of Turaev~\cite{T2}).
For general $V$, this is certainly not true. The vertex operator
algebra $V$ must satisfy certain conditions.
Here is a~precise conjecture formulated by the second author in~\cite{H-problems}.

\begin{conj}\label{orb-conj}
Let $V$ be a~vertex operator algebra satisfying the three conditions
in Theorem~{\rm\ref{mod-tensor-cat}}
and let $G$ be a~finite group of automorphisms of $V$. Then the category
of $g$-twisted $V$-modules for
all $g\in G$ is a~$G$-crossed braided tensor category.
\end{conj}

We also conjecture that the category of $g$-twisted
$V$-modules for all $g\in G$ is a~$G$-crossed modular
tensor category in a~suitable sense.
Since the definitions of $G$-crossed modular
tensor category
in~\cite{K3,T2} are different, more work
needs to be done to find out which definition is the correct
one for the category of twisted modules for a~vertex operator algebra.
But we do believe that this stronger
$G$-crossed modular tensor category conjecture should be true in a~suitable sense.

In the case that $G$ is trivial (the group containing only the identity),
Conjecture~\ref{orb-conj} and even
the stronger $G$-crossed modular tensor category
conjecture is true by Theorem~\ref{mod-tensor-cat}.
Thus the $G$-crossed modular tensor category
conjecture is a~natural generalization of
Theorem~\ref{mod-tensor-cat} to
the category of $g$-twisted $V$-modules for $g\in G$.

In the case that the fixed point subalgebra $V^{G}$ of $V$ under $G$
satisfies the conditions in
Theorem~\ref{mod-tensor-cat} above, the category of $V^{G}$-modules
is a~modular tensor category. In this case,
Conjecture~\ref{orb-conj} can be proved using the modular tensor category
structure on the category of $V^{G}$-modules and the results
on tensor categories by Kirillov Jr. \cite{K1,K2,K3}
and M\"{u}ger~\cite{Mu1,Mu2}.
In the special case that $G$ is a~finite cyclic group and
$V$ satisfies the conditions in Theorem~\ref{mod-tensor-cat},
Carnahan--Miyamoto~\cite{CM}
proved that $V^{G}$ also satisfies the conditions in
Theorem~\ref{mod-tensor-cat}. In the case that $G$ is a~finite cyclic group
and $V$ is in addition a~holomorphic
vertex operator algebra (meaning that the only irreducible $V$-module
is $V$ itself),
Conjecture~\ref{orb-conj} can be obtained as a~consequence
of the results of
van Ekeren--M\"{o}ller--Scheithauer~\cite{EMS} and M\"{o}ller~\cite{Mo}
on the modular tensor category
of $V^{G}$-modules. Assuming that the category of grading-restricted $V^{G}$-modules
has a~natural structure of vertex tensor category in the sense of~\cite{HL}, McRae~\cite{Mc} constructed a~nonsemisimple $G$-crossed
braided tensor category structure on the category of grading-restricted
(generalized) $g$-twisted $V$-modules.

For general finite group $G$, the
conjecture that the fixed point subalgebra $V^{G}$ of $V$ under $G$
also satisfies the conditions in
Theorem~\ref{mod-tensor-cat} is still open and seems to be a~difficult problem. On the other hand, using twisted modules
and twisted intertwining operators to construct $G$-crossed braided
tensor categories seems to be a~more conceptual
and direct approach. If this approach works, we expect that
the category of $V^{G}$-modules can also be studied
using the $G$-crossed braided
tensor category structure on the category of twisted $V$-modules.

In the case that the vertex operator algebra $V$ does not satisfy
the three conditions in Theorem~\ref{mod-tensor-cat} and/or the group $G$ is
not finite, it is not even clear what should be the precise conjecture.
This was proposed as an open problem in~\cite{H-problems}.

In the present paper, we prove some initial results in a~long term program
to prove the conjecture and to solve the open problem above.
We introduce a~more general notion
of twisted intertwining operator
than the one introduced by the second author in~\cite{H-twisted-int}.
In~\cite{H-twisted-int}, the correlation functions obtained from
the products and iterates of
twisted intertwining operators and twisted vertex operators
are required to be of a~special explicit form. But for a~twisted intertwining
operator in this paper, such correlation functions are not required
to have such an explicit form.

As in~\cite{H-twisted-int}, we prove some
basic properties and construct the skew-symmetry
and contragredient isomorphisms for our general twisted intertwining
operators.
Using such general twisted intertwining operators,
we introduce a~notion of $P(z)$-tensor product of
two twisted modules for~${z\in \C^{\times}}$
and give a~construction of such a~$P(z)$-tensor product
under suitable assumptions. We also prove a~result showing
that under suitable conditions, these assumptions are satisfied.\looseness=1

We need $P(z)$-tensor products
for $z\in \C^{\times}$ because we would like to construct $G$-crossed
vertex tensor categories in the future,
not just $G$-crossed braided tensor categories.
Also note that to give the correct
notion of $P(z)$-tensor product of twisted modules, we need
to use the most general twisted intertwining operators.
If we use only a~certain special
set of twisted intertwining operators as in~\cite{H-twisted-int}
to define and construct the $P(z)$-tensor products,
we would obtain submodules of the correct $P(z)$-tensor products.

We note that in the untwisted case, a~$P(z)$-compatibility
condition and a~$P(z)$-grading-restriction condition
(see~\cite{tensor3,HLZ4}) play an important role in the
proof of associativity (operator product expansion)
of intertwining operators and in the construction of the
associativity isomorphisms for the vertex tensor category
structure (see~\cite{tensor4,HLZ6}). In this paper, we also formulate a~$P(z)$-compatibility
condition and a~$P(z)$-grading-restriction condition and use
these conditions to give another construction of
the $P(z)$-tensor product. In the untwisted case,
the $P(z)$-compatibility condition
is formulated using a~formula obtained from the Jacobi identity
in the definition of intertwining operators (see~\cite{tensor3,HLZ4}). But since in general we do not have a~Jacobi identity
that can be used as the main axioms for the definition of
twisted intertwining operators, our formulation of this condition
and the construction of tensor products using this condition are
complex analytic and are very different from the formulation and
construction in~\cite{tensor3,HLZ4}. We expect that these two conditions
will play the same important role in the future proof of
the conjectured associativity of twisted intertwining operators
formulated in~\cite{H-orbifold} (where twisted intertwining operators
should be
replaced by the most general twisted intertwining operators introduced
in this paper).

This paper is organized as follows: In Section \ref{sec2}, we recall
the definitions of (generalized) twisted module, lower-bounded
(generalized) twisted module and grading-restricted (generalized)
twisted module. We then introduce the general notion
of twisted intertwining operator mentioned above.
In Section \ref{sec3}, we give the skew-symmetry and contragredient isomorphisms
for these general twisted intertwining operators.
We introduce the notion of $P(z)$-tensor product and give a~construction
under suitable assumptions in Section \ref{sec4}, we prove a~result showing
that under suitable conditions, these assumptions are satisfied.
We also construct
$G$-crossed commutativity isomorphisms and
$G$-crossed braiding isomorphisms in this section.
We give the $P(z)$-compatibility
condition and $P(z)$-grading-restriction condition and give another construction
of the $P(z)$-tensor product in Section \ref{sec5}.

\section{Twisted modules and twisted intertwining operators}\label{sec2}

We first recall in this section the notion of (generalized) twisted module
from~\cite{H-log-twisted-mod}. We then introduce a~notion
of twisted intertwining
operators more general than the one in
\cite{H-twisted-int}. We also give some basic results on such
twisted intertwining operators.

For $z\in \C^{\times}$ and $p\in \Z$, we shall use the notation
$l_{p}(z)=\log |z|+{\rm i}\arg z+2\pi p {\rm i}$, where $0\le \arg z<2\pi$.
We shall also use the notation $\log z=l_{0}(z)=\log |z|+{\rm i}\arg z$.
For a~vector space~$U$,~${p\in \Z}$ and a~formal series
\[f(x)=\sum_{k=0}^{K}\sum_{n\in \C}a_{n, k}x^{n}(\log x)^{k},\]
where $a_{n, k}\in U$,
the series
\[f^{p}(z)=\sum_{k=0}^{K}\sum_{n\in \C}
a_{n, k}{\rm e}^{nl_{p}(z)}(l_{p}(z))^{k}\]
is called the $p$-th analytic branch of $f(x)$. We also denote
$f^{0}(z)$ simply by $f(z)$.

Let $V$ be a~vertex operator algebra and $g$ an automorphism of $V$.
We recall the definition of generalized $g$-twisted
$V$-module first introduced in~\cite{H-log-twisted-mod}.
Note that as in~\cite{H-log-twisted-mod,H-twisted-int} in this paper, the vertex operator map for a~generalized $g$-twisted $V$-module
in general contains the logarithm of the variable and the operator
$L(0)$ in general does not have to act semisimply.

\begin{defn}\label{twisted-mod}
A~{\it generalized $g$-twisted
$V$-module without a~$g$-action} is a~$\C$-graded
vector space~${W = \coprod_{n \in \C} W_{[n]}}$ {\rm(}graded by {\it weights}$)$
equipped with a~linear map
\begin{gather*}
Y_{W}^g\colon \ V\otimes W \to W\{x\}[\log x],\\
v \otimes w \mapsto Y_{W}^g(v, x)w=\sum_{k=0}^{K}\sum_{n\in \C}(Y_{W}^g)_{n, k}(v)wx^{-n-1}(\log x)^{k}
\end{gather*}
satisfying the following conditions:
\begin{enumerate}[label=$(\arabic*)$]\itemsep=0pt
\item The {\it equivariance property}: For $p \in \mathbb{Z}$, $z
\in \mathbb{C}^{\times}$, $v \in V$ and $w \in W$,
$\bigl(Y^{g}_{W}\bigr)^{p + 1}(gv,
z)w = \bigl(Y^{g}_{W}\bigr)^{p}(v, z)w$,
where for $p \in \mathbb{Z}$, $\bigl(Y^{g}_{W}\bigr)^{ p}(v, z)$
is the $p$-th analytic branch of $Y_{W}^g(v, x)$.

\item The {\it identity property}: For $w \in W$, $Y^g_{W}({\bf 1}, x)w
= w$.

\item The {\it duality property}: For
any $u, v \in V$, $w \in W$ and $w' \in W'$, there exists a~maximally-extended
multivalued analytic function with preferred branch of the form
\[
f(z_1, z_2) = \sum_{i,
j, k, l = 0}^N a_{ijkl}z_1^{m_i}z_2^{n_j}(\operatorname{log}z_1)^k({\rm
log}z_2)^l(z_1 - z_2)^{-t}
\]
for $N \in \mathbb{N}$, $m_1, \dots,
m_N$, $n_1, \dots, n_N \in \mathbb{C}$ and $t \in \mathbb{Z}_{+}$,
such that the series
\begin{gather*}
\big\langle w', \bigl(Y^{g}_{W}\bigr)^{ p}(u, z_1)\bigl(Y^{g}_{W}\bigr)^{p}(v,
z_2)w\big\rangle = \sum_{n \in \mathbb{C}}\big\langle w', \bigl(Y^{g}_{W}\bigr)^{ p}(u,
z_1)\pi_n\bigl(Y^{g}_{W}\bigr)^{ p}(v, z_2)w\big\rangle,
\\
\big\langle w', \bigl(Y^{g}_{W}\bigr)^{ p}(v,
z_2)\bigl(Y^{g}_{W}\bigr)^{p}(u, z_1)w\big\rangle = \sum_{n \in \mathbb{C}}\big\langle w',
\bigl(Y^{g}_{W}\bigr)^{ p}(v, z_2)\pi_n\bigl(Y^{g}_{W}\bigr)^{ p}(u, z_1)w\big\rangle,
\\
\big\langle w', \bigl(Y^{g}_{W}\bigr)^{ p}(Y_{V}(u, z_1 - z_2)v,
z_2)w\big\rangle = \sum_{n \in \mathbb{C}}\big\langle w', \bigl(Y^{g}_{W}\bigr)^{
p}(\pi_nY_{V}(u, z_1 - z_2)v, z_2)w\big\rangle
\end{gather*}
are absolutely convergent on
the regions $|z_1| > |z_2| > 0$, $|z_2| > |z_1| > 0$, $|z_2| > |z_1
- z_2| > 0$, respectively, and their sums are equal to the branch
\[
f^{p,p}(z_{1}, z_{2})= \sum_{i, j, k, l = 0}^N
a_{ijkl}{\rm e}^{m_il_p(z_1)}{\rm e}^{n_jl_p(z_2)}l_p(z_1)^kl_p(z_2)^l(z_1 -
z_2)^{-t}
\]
of $f(z_1, z_2)$ on the region $|z_1| > |z_2| > 0$, the region $|z_2| > |z_1| > 0$,
the region given by $|z_2| > |z_1- z_2| > 0$ and $|{\arg z_{1}-\arg z_{2}}|<\frac{\pi}{2}$, respectively.

\item The {\it $L(0)$-grading condition}:
Let $L_{W}^{g}(0)=\res_{x}xY_{W}^g(\omega, x)$. Then for $n\in \C$,
$w \in W_{[n]}$,
there exists $K\in \Z_{+}$ such that \smash{$\bigl(L_{W}^g (0)-n\bigr)^{K} w
=0$}.

\item The $L(-1)$-{\it derivative property}: For $v \in V$,
\[
\frac{\rm d}{{\rm d}x}Y^g_{W}(v, x) = Y^g_{W}(L_{V}(-1)v, x).
\]

\end{enumerate} A~{\it lower-bounded generalized $g$-twisted $V$-module}
is a~generalized $g$-twisted
$V$-module $W$ such that
for each ${n\in \C}$, $W_{[n + l]} = 0$ for
sufficiently negative real numbers $l$. A~generalized $g$-twisted $V$-module $W$
is said to be {\it grading-restricted} if it is lower
bounded and for each ${n \in \mathbb{C}}$, \smash{$\dim W_{[n]}<\infty$}. A~{\it generalized $g$-twisted $V$-module with a~$g$-action} is a~${\C}\times \C/\Z$-graded
vector space \smash{$W = \coprod_{n \in \C, \alpha\in \C/\Z} W_{[n]}^{[\alpha]}$} {\rm(}graded by {\it weights} and {\it $g$-weights}{\rm)}
equipped with a~linear map \smash{$Y_{W}^g\colon V\otimes W \to W\{x\}[\log x]$} as above and an action of $g$
such that $W$ equipped with $Y_{W}^g$ is a~generalized $g$-twisted $V$-module
without a~$g$-action satisfying in addition the following {\it $g$-grading condition}:
There exists $\Lambda\in \Z_{+}$ such that for $\alpha\in \C/\Z$,
\smash{$w \in W^{[\alpha]}=\coprod_{n\in \N}W_{[n]}^{[\alpha]}$}, \smash{$\bigl(g-{\rm e}^{2\pi \alpha{\rm i}}\bigr)^{\Lambda}w=0$}. Moreover,
$gY_{W}^{g}(u,x)w=Y_{W}^{g}(gu,x)gw$ for $u\in V$ and $w\in W$.
\end{defn}

Module maps between $g$-twisted $V$-modules are defined in the obvious way.
For two different automorphisms $g_{1}$ and $g_{2}$ of $V$, there might exist
linear maps from a~$g_{1}$-twisted $V$-module to a~$g_{2}$-twisted $V$-module commuting with the twisted vertex operator maps when
$V$ is not simple. But we will not consider such maps in this paper. So in the
categories of twisted $V$-modules considered in this paper,
we always choose the space of
 module maps between a~$g_{1}$-twisted $V$-module and a~$g_{2}$-twisted $V$-module
to be $0$ if $g_{1}\ne g_{2}$.

In this paper, generalized $g$-twisted $V$-modules with or without $g$-actions
are all called {\it generalized $g$-twisted $V$-modules}. Only when the $g$-actions are specifically discussed,
we shall add the words ``with a~$g$-action'' or ``without a~$g$-action''. Later, we shall also need
generalized $g$-twisted $V$-modules with $g$-actions equipped with an $h$-action for another
automorphism $h$ of $V$.
Also, for simplicity, we shall sometimes omit the subscript
$W$ to denote the twisted vertex operator map~$Y_{W}^{g}$
by $Y^{g}$.

\begin{rema}
Note that a~lower-bounded generalized
$g$-twisted $V$-module defined above also satisfies
the lower-truncation property: For $v\in V$ and $w\in W$, $n\in \C$ and $k=0, \dots, K$,
\smash{$\bigl(Y_{W}^g\bigr)_{n+l, k}(w_{1})w_{2}=0$} for $l\in \N$ sufficiently large.
This in fact follows from the $L(0)$-commutator formula,
which is a~consequence of the duality property and the $L(0)$-grading condition.
\end{rema}

\begin{rema}
In this paper, we need to study generalized $g$-twisted $V$-modules with and without
$g$-actions in a~different way because the skew-symmetry and contragredient isomorphisms in Section \ref{sec3} are different in these two types
of twisted modules. In fact, there is also a~notion of generalized $g$-twisted module with a~$G$ action for a~group $G$ of automorphisms of $V$. See~\cite{H-cofinite-P(z)} for a~definition and
the construction $P(z)$-tensor product of such twisted modules.
\end{rema}

In this paper, we shall consider only lower-bounded generalized $g$-twisted $V$-modules and grading-restricted
generalized $V$-modules. For simplicity, we shall often call lower-bounded generalized
$g$-twisted $V$-modules simply $g$-twisted $V$-modules and
grading-restricted
generalized $V$-modules simply grading-restricted $g$-twisted $V$-modules.

Let $\bigl(W, Y^{g}_{W}\bigr)$ be a~$g$-twisted
$V$-module without a~$g$-action. Let $h$ be an automorphism of~$V$. We recall
the $hgh^{-1}$-twisted
$V$-module $\bigl(W, \phi_{h}\bigl(Y^{g}_{W}\bigr)\bigr)$ without a~$hgh^{-1}$-action (see, for example,~\cite{H-twisted-int}).
Let
\begin{gather*}
\phi_{h}\bigl(Y^{g}_{W}\bigr)\colon\ V\times W\to W\{x\}[\operatorname{log} x],\qquad
v \otimes w \mapsto \phi_{h}\bigl(Y^{g}_{W}\bigr)(v, x)w
\end{gather*}
be the linear map defined by{\samepage
\[\phi_{h}\bigl(Y^{g}_{W}\bigr)(v, x)w=Y^{g}_{W}\bigl(h^{-1}v, x\bigr)w.\]
Then the pair $\bigl(W, \phi_{h}\bigl(Y^{g}_{W}\bigr)\bigr)$ is an $hgh^{-1}$-twisted
$V$-module without a~$hgh^{-1}$-action.}

In the case that $W$ is a~$g$-twisted $V$-module
with a~$g$-action and also has an action
of $h$ in the sense that $W$ is a~module for the group generated by $g$ and $h$, we define
\begin{gather*}
\varphi_{h}\bigl(Y^{g}_{W}\bigr)\colon\ V\times W\to W\{x\}[\operatorname{log} x],\qquad
v \otimes w \mapsto \varphi_{h}\bigl(Y^{g}_{W}\bigr)(v, x)w
\end{gather*}
by
\[\varphi_{h}\bigl(Y^{g}_{W}\bigr)(v, x)w=hY^{g}_{W}\bigl(h^{-1}v, x\bigr)h^{-1}w.\]
Then $\bigl(W, \varphi_{h}\bigl(Y^{g}_{W}\bigr)\bigr)$ is also an $hgh^{-1}$-twisted
$V$-module without a~$hgh^{-1}$-action.
Since the conformal element $\omega$ is fixed under the action of
$g^{-1}$, the coefficients
\[L_{\varphi_{h}(W)}(n)=\res_{x}x^{n+1}\varphi_{h}\bigl(Y^{g}_{W}\bigr)(\omega, x)
=h\res_{x}x^{n+1}Y_{W}^{g}(\omega, x)h^{-1}=hL_{W}^{g}(n)h^{-1}\]
 of the vertex operator $\varphi_{h}\bigl(Y^{g}_{W}\bigr)(\omega, x)$ for $n\in \N$
also give Virasoro operators on $W$ with the same central charge.
Since $h$ is an isomorphism of the vector space $W$,
\[hW=\coprod_{n\in \C, \alpha\in \C/\Z}hW_{[n]}^{[\alpha]}.\]
We denote the graded space $hW$ by $\varphi_{h}(W)$ and
\smash{$hW_{[n]}^{[\alpha]}$} by \smash{$\varphi_{h}(W)_{[n]}^{[\alpha]}$}
$n\in \C, \alpha\in \C/\Z$. Then we have
\[\varphi_{h}(W)=\coprod_{n\in \C, \alpha\in \C/\Z}\varphi_{h}(W)_{[n]}^{[\alpha]}.\]
Note that as a~vector space $\varphi_{h}(W)=hW$ is the same as $W$, but as graded spaces,
they are different.
For \smash{$w\in \varphi_{h}(W)_{[n]}=\coprod_{\alpha\in \C/\Z}\varphi_{h}(W)_{[n]}^{[\alpha]}$},
\smash{$h^{-1}w\in W_{[n]}=\coprod_{\alpha\in \C/\Z}W_{[n]}^{[\alpha]}$}. So there exists ${K\in \Z_{+}}$ such that \smash{$\bigl(L_{W}^{g}-n\bigr)^{K}h^{-1}w=0$}. Then
\[
(L_{\varphi_{h}(W)}(0)-n)^{K}w=h\bigl(L_{W}^{g}-n\bigr)^{K}h^{-1}w=0.
\]
For \smash{$w\in \varphi_{h}(W)^{[\alpha]}=\coprod_{n\in \Z}\varphi_{h}(W)_{[n]}^{[\alpha]}$},
\smash{$h^{-1}w\in W^{[\alpha]}=\coprod_{n\in \Z}W_{[n]}^{[\alpha]}$}. So there exists
$\Lambda\in \Z_{+}$ such that \smash{$\bigl(g-{\rm e}^{2\pi {\rm i}\alpha}\bigr)^{\Lambda}h^{-1}w=0$}.
Then \smash{$\bigl(hgh^{-1}-{\rm e}^{2\pi {\rm i}\alpha}\bigr)^{\Lambda}w=h\bigl(g-{\rm e}^{2\pi {\rm i}\alpha}\bigr)^{\Lambda}h^{-1}w=0$}.
Also we have
\begin{align*}
hgh^{-1}\varphi_{h}\bigl(Y^{g}_{W}\bigr)(v, x)w&=hgY^{g}_{W}\bigl(h^{-1}v, x\bigr)h^{-1}w=hY^{g}_{W}\bigl(h^{-1}hgh^{-1}v, x\bigr)h^{-1}hgh^{-1}w\\
&=\varphi_{h}\bigl(Y^{g}_{W}\bigr)\bigl(hgh^{-1}v, x\bigr)hgh^{-1}w.
\end{align*}
So $(\varphi_{h}(W), \varphi_{h}\bigl(Y^{g}_{W}\bigr))$ is a~$hgh^{-1}$-twisted $V$-module with an
$hgh^{-1}$-action.
We shall denote this $hgh^{-1}$-twisted
$V$-module with an
$hgh^{-1}$-action by $\varphi_{h}(W)$.

Note that when $h=g$, $\varphi_{h}(W)=W$ and
$\phi_{g}(W)$ with the twisted vertex operator
given by
\[\phi_{g}\bigl(Y^{g}_{W}\bigr)(v, x)w=Y^{g}_{W}\bigl(g^{-1}v, x\bigr)w\]
is equivalent to the original
$g$-twisted module $W$ with the equivalence $g^{-1}\colon W\to W$. By the equivariance property,
we have \smash{$\bigl(Y^{g}_{W}\bigr)^{p}\bigl(g^{-1}v, x\bigr)=\bigl(Y^{g}_{W}\bigr)^{p+1}(v, x)$}
and, if
\[Y^{g}_{W}(v, x)w=\sum_{k=0}^{K}\sum_{n\in \C}
\bigl(Y^{g}_{W}\bigr)_{n, k}(v)x^{-n-1}(\log x)^{k}\]
for $v\in V$ and $w\in W$,
we have
\[Y^{g}_{W}\bigl(g^{-1}v, x\bigr)w=\sum_{k=0}^{K}\sum_{n\in \C}
\bigl(Y^{g}_{W}\bigr)_{n, k}(v){\rm e}^{-2\pi {\rm i} n}x^{-n-1}(\log x+2\pi {\rm i})^{k}.\]

\begin{rema}
 In the definitions of $\bigl(\varphi_{h}(W), \varphi_{h}\bigl(Y_{W}^{g}\bigr)\bigr)$,
we assume that there is an action of~$h$ on~$W$. One immediate question is whether
there exist such $g$-twisted $V$-modules with $h$-actions. In fact, the same construction
in Section \ref{sec5} of~\cite{H-const-twisted-mod} can be used to construct such a~$g$-twisted $V$-module with an $h$-action
or even a~$g$-twisted $V$-module with a~$G$-action for a~group~$G$ of automorphisms of $V$ satisfying
the additional condition $hY_{W}^{g}(v, x)w=Y_{W}^{g}(hv, x)hw$ for~${h\in G}$. In this case, though
$\varphi_{h}(W)=W$ as vector spaces and
$\varphi_{h}\bigl(Y_{W}^{g}\bigr)=Y_{W}^{g}$ as linear maps from~${V\otimes W}$ to $W\{x\}[\log x]$ when we view
$\varphi_{h}(W)$ as a~vector space equal to $W$, $\varphi_{h}(W)$ is graded by generalized eigenspaces of $hgh^{-1}$
but $W$is graded by the generalized eigenspaces of $g$.
In particular, $\bigl(\varphi_{h}(W), \varphi_{h}\bigl(Y_{W}^{g}\bigr)\bigr)$ and $\bigl(W, Y_{W}^{g}\bigr)$ are indeed different twisted $V$-modules.
\end{rema}

We also recall contragredient twisted $V$-modules (see, for example, \cite{H-twisted-int}).
Let $\bigl(W, Y^{g}_{W}\bigr)$ be a~$g$-twisted $V$-module relative to $G$.
Let $W'$ be the graded dual of $W$. Define a~linear map
\begin{gather*}
\bigl(Y_{W}^{g}\bigr)'\colon\ V\otimes W' \to W'\{x\}[\operatorname{log} x],\qquad
v \otimes w' \mapsto \bigl(Y^g_{W}\bigr)'(v, x)w'
\end{gather*}
by
\[\big\langle \bigl(Y^g_{W}\bigr)'(v, x)w', w\big\rangle=\big\langle w', Y^{g}_{W}\bigl({\rm e}^{xL(1)}\bigl(-x^{-2}\bigr)^{L(0)}v, x^{-1}\bigr)w\big\rangle\]
for $v\in V$, $w\in W$ and $w'\in W'$.
Then the pair $\bigl(W', \bigl(Y^{g}_{W}\bigr)'\bigr)$ is a~$g^{-1}$-twisted $V$-module.

\begin{figure}[!b]
\centering
\includegraphics[scale=0.6]{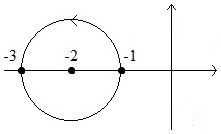}
\caption{The loop for $b_{12}$.}\label{b12}
\end{figure}

\begin{figure}[!t]
\centering
\includegraphics[scale=0.6]{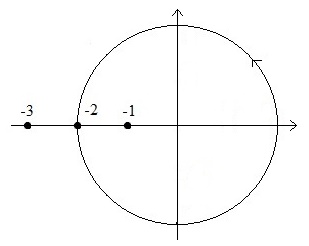}
\caption{The loop for $b_{23}$.}\label{b23}
\end{figure}

\begin{figure}[!t]
\centering
\includegraphics[scale=0.6]{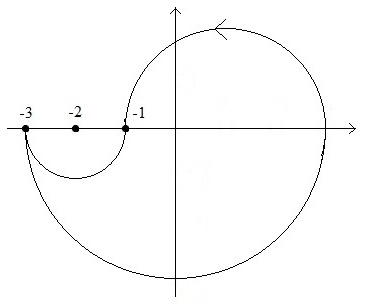}
\caption{The loop for $b_{13}$.}\label{b13}
\end{figure}

Let $M^{2}=\big\{(z_{1}, z_{2})\in \C^{2}\mid z_{1}\ne 0,\, z_{2}\ne 0,\,
z_{1}\ne z_{2}\big\}$.
Let $f(z_{1}, z_{2})$ be a~maximally extended
multivalued analytic function
on $M^{2}$
with a~preferred single-valued branch $f^{e}(z_{1}, z_{2})$
on the simply-connected
region $M^{2}_{0}$ given by
cutting $M^{2}$ along the positive real lines in the $z_{1}$-, $z_{2}$- and
$(z_{1}-z_{2})$-planes, that is, the sets
\begin{gather*}
\big\{(z_{1}, z_{2})\in M^{2}\mid z_{1}\in \R_{+}\big\},\qquad
\big\{(z_{1}, z_{2})\in M^{2}\mid z_{2}\in \R_{+}\big\},\\
\big\{(z_{1}, z_{2})\in M^{2}\mid z_{1}-z_{2}\in \R_{+}\big\},
\end{gather*}
with these sets attached to the upper half $z_{1}$-, $z_{2}$- and
$(z_{1}-z_{2})$-planes, that is, the points in these sets are
viewed as limits in $M_{0}^{2}$ of sequences of points approaching them
from the upper $z_{1}$-, $z_{2}$- and
$(z_{1}-z_{2})$-planes. Note that given any point
$\bigl(z^{0}_{1}, z^{0}_{2}\bigr)\in M^{2}$
and any loop $\gamma$ in $M^{2}$ based at $\bigl(z^{0}_{1}, z^{0}_{2}\bigr)$, we obtain from
the preferred branch of $f(z_{1}, z_{2})$ another single-valued branch
by going around $\gamma$. The resulting single-valued branch
depends only on the homotopy class of the loop and is independent of
the choices of the base point. Thus we obtain a~right action of the
fundamental group of $M^{2}$ on the set of single-valued
branches of $f(z_{1}, z_{2})$. Note that $M^{2}$ is
homotopically equivalent to the configuration space
\[F_{3}(\C)=\big\{(z_{1}, z_{2}, z_{3})\in \C^{3}\mid
z_{i}\ne z_{j},\, i\ne j\big\}.\]
So the fundamental group of $M^{2}$
is in fact the pure braid group $\PB_{3}$ and has three generators~$b_{12}$,~$b_{13}$ and $b_{23}$, which are given as follows:
Choose the base point
to be $(-3, -2)$. Then the generator $b_{12}$
is the homotopy class of the loop given by letting $z_{1}$ go
counterclockwise around the circle of radius $1$
centered at $-2$ (see Figure~\ref{b12}).
The generator $b_{23}$ is the homotopy class of
the loop given by letting $z_{2}$ go counterclockwise
around the circle of radius $2$ centered at $0$ (see Figure~\ref{b23}).
The generator $b_{13}$
is the homotopy class of the loop given by letting
$z_{1}$ go around first the lower half circle
of radius $3$ centered at $0$ counterclockwise, then the upper half circle
of radius $2$ centered at $1$ counterclockwise and finally the
lower half circle of radius $1$ centered at $-2$ clockwise (see Figure~\ref{b13}).
We know that the pure braid group $\PB_3$ is isomorphic to the
group generated by $b_{12}$, $b_{13}$, $b_{23}$ with the relations
$
b_{13}b_{12}b_{23}=b_{12}b_{23}b_{13}=b_{23}b_{13}b_{12}$.
See
\cite{D} for more detailed discussions on $M^{2}$, its fundamental group,
the pure braid groups and
multivalued functions on $M^{2}$.

For $a\in \PB_{3}$,
we denote the single-valued branch of $f(z_{1}, z_{2})$
obtained by applying
$a$ to $f^{e}(z_{1}, z_{2})$
by $f^{a}(z_{1}, z_{2})$.

Let $\delta=(\delta_1,\dots,\delta_n)$, where $\delta_i\in\{0,\infty\}$,
	$i=1,\dots,n$. Suppose $f(z_1,\dots,z_n)$ is a~multi-valued analytic
	function defined on an open region $\Omega$ of $\C^n$. We say
	$(z_1,\dots,z_n)=\delta$ is a~{\it component-isolated singularity
	of $f(z_1,\dots,z_n)$} if there exists $r=(r_{1}, \dots, r_{n})\in \R_{+}^n$ such that
	$\Delta^\times(\delta,r)\subset\Omega$, where
$\Delta^\times(\delta,r)=\{(z_{1}, \dots, z_{n})\in \C^{n}\mid |z_{i}|>r_{i} \text{ if } \delta_{i}=\infty,
0< |z_{i}|<r_{i} \text{ if } \delta_{i}=0\}$.
	Let~${A\in\operatorname{GL}(n,\C)}$
	and $\beta\in\C^n$ (written as a~row vector).
Then $\zeta_{1}, \dots, \zeta_{n}$ given by
\[(\zeta_1,\dots,\zeta_n)=(z_{1}, \dots, z_{n})A-\beta\]
are also independent variables. Define
	\begin{equation}\label{isolated}
	g(\zeta_1,\dots,\zeta_n)=f\left((\zeta_1,\dots,\zeta_n)
	A^{-1}+\beta A^{-1}\right).
	\end{equation}
For $\delta\in (\C\cup \{\infty\})^{n}$, we say that
$(\zeta_{1}, \dots, \zeta_{n})=\delta$ is a~{\it component-isolated singularity
of the function $f(z_1,\dots,z_n)$}
	if $(\zeta_{1}, \dots, \zeta_{n})=\delta$ is a~component-isolated singularity of
	$g(\zeta_1,\dots,\zeta_n)$.
	
\begin{rema}
{\rm Notice that $(\zeta_1,\dots,\zeta_n)(=(z_{1}, \dots, z_{n})A-\beta)
=\delta$ being a~component-isolated singularity
of a~function is not equivalent to
$(z_1,\dots,z_n)=\delta A^{-1}+\beta A^{-1}$ being a~component-isolated singularity of the same function.
This is because we have different sets of independent variables. For example, consider
$f(z_1,z_2)=\frac{1}{z_1-z_2}$. Since $(\zeta_1,\zeta_2)=(0,0)$ is a~component-isolated
singularity of the function $g(\zeta_1,\zeta_2)=1/\zeta_1$, we also say that
$(z_1-z_2,z_2)=(0,0)$ is a~component-isolated singularity of $f(z_1,z_2)$.
In this case, $z_{1}-z_{2}$ and $z_{2}$ are independent variables.
However, $(z_1,z_2)=(0,0)$ is clearly not a~component-isolated singularity of $f(z_1,z_2)$.
In this case, the independent variables are $z_{1}$ and $z_{2}$. }
\end{rema}

\begin{defn}
{\rm Let $f(z_1,\dots,z_n)$ be a~multi-valued analytic function defined on
an open region of $\C^n$. Let $\delta=(\delta_{1},\dots,\delta_n)
\in \{0,\infty\}^n$. Suppose $(z_1,\dots,z_n)=\delta$ is a~component-isolated
 singularity of $f(z_1,\dots,z_n)$. Let $\{f^b(z_1,\dots,z_n)\}_{b\in B}$
 be the set of all single valued branches of $f(z_1,\dots,z_n)$ near
 $\delta$ with cuts at $z_i\in\R_{+}$. Then for each $b\in B$,
 there exists $r_b\in \R_{+}^n$ such that $f^b(z_1,\dots,z_n)$
 is analytic on $\Delta^\times(\delta,r_b)$. We say that $(z_1,\dots,z_n)
 =\delta$ is a~\textit{regular singularity} of $f(z_1,\dots,z_n)$
 if there exists $K\in\N$,
 \smash{$D_i=\cup_{j=1}^{N_i}r_j^{(i)}+\N$}
 \big(or \smash{$D_i=\cup_{j=1}^{N_i}r_j^{(i)}-\N$}\big)
 where~\smash{$r_1^{(i)},\dots,r_{N_i}^{(i)}\in\C$} for $\delta_i=0$ (or
 $\delta_i=\infty$), and \smash{$\alpha^{(b)}_{a_1,j_1;\dots;a_n,j_n}\in\C$},
 such that on the region~$\Delta^\times(\delta,r_b)$, the right-hand-side
 of the following equation is absolutely convergent, and
\begin{align}\label{reg-sing-pt-expansion}
f^b(z_1,\dots,z_n)=\sum_{i=1}^{n}\sum_{a_i\in D_i}
\sum_{j_1,\dots,j_n=0}^{K}\alpha^{(b)}_{a_1,j_2;
\dots;a_n,j_n}z_1^{a_1}(\log z_1)^{j_1}\cdots z_n^{a_n}(\log z_n)^{j_n}.
\end{align}
If $K=0$ and $D_{i}=-n^{(i)}+\N$
in the case $\delta_{i}=0$
and $D_{i}=n^{(i)}-\N$ in the case $\delta_{i}=\infty$ for $i=1, \dots, n$
in (\eqref{reg-sing-pt-expansion}), where $n^{(i)}\in \Z_{+}$
for $i\in I\subset \{1, \dots, n\}$
and $n^{(i)}=0$ for $i\ne I$, we say that
$z_i=\delta_{i}$ for $i\in I$ are {\it poles of
$f(z_1,\dots,z_n)$}. If $n^{(i)}\in \Z_{+}$
for $i\in I\subset \{1, \dots, n\}$, $\delta_{i}=0$ are the smallest
and	such that~\eqref{reg-sing-pt-expansion} holds, we call
$n_{i}$ the orders of the poles $z_i=\delta_{i}$, respectively,
for $i\in \I$.
Let $A\in\operatorname{GL}(n,\C)$
and $\beta\in\C^n$ be the same as above. We say that
 $(\zeta_1,\dots,\zeta_n)=\delta$ is a~\textit{regular singularity}
 of $f(z_1,\dots,z_n)$ if $(\zeta_1,\dots,\zeta_n)=\delta$
 is a~regular
 singularity of $g(\zeta_1,\dots,\zeta_n)$, where $g(\zeta_1,
 \dots,\zeta_n)$ is given by~\eqref{isolated}. We say that
 $\zeta_{i}=\delta_{i}$ for $i\in I$ are
poles of $f(z_1,\dots,z_n)$ with orders $n_{i}$, respectively, if
$\zeta_{i}=\delta_{i}$ for $i\in I$ are poles of $g(\zeta_1,\dots,\zeta_n)$ with
orders $n_{i}$, respectively.}
\end{defn}

\begin{rema}
{\rm Let $f(z_1,\dots,z_n)$ and $g(z_1,\dots,z_n)$ be multivalued analytic functions
with preferred branches defined on an open region.
Then $\lambda f(z_1,\dots,z_n)+\mu g(z_1,\dots,z_n)$ for
$\lambda, \mu\in \C$ and
$f(z_1,\dots,z_n) g(z_1,\dots,z_n)$ are well defined using the preferred
branches and are also multivalued analytic functions on the same region
with preferred branches.
If $(\zeta_1,\dots,\zeta_n)
=((z_1,\dots,z_n)A-\beta=)\delta$ is a~regular singular point of both $f(z_1,\dots,z_n)$ and $g(z_1,\dots,z_n)$,
then it is also a~regular singular point for $\lambda f(z_1,\dots,z_n)+
\mu g(z_1,\dots,z_n)$
and $f(z_1,\dots,z_n) g(z_1,\dots,z_n)$.
Therefore, the set of multivalued analytic functions
on the same region with a~preferred branch
such that $(\zeta_1,\dots,\zeta_n)=\delta$ is a~regular
singular point form a~commutative associative algebra over $\C$.}
\end{rema}

We also need the region
$M^{n}=\{(z_{1}, \dots, z_{n})\in \C^{n}\mid
z_{i}\ne 0, z_{i}\ne z_{j}, i\ne j \}$
for $n\in \Z_{+}$.

\begin{defn}\label{def-tw-int-op}
 Let $g_{1}$, $g_{2}$, $g_{3}$ be automorphisms of $V$ and let $W_{1}$, $W_{2}$ and $W_{3}$ be $g_{1}$-, $g_{2}$-
and $g_{3}$-twisted
$V$-modules $($with or without actions of the corresponding automorphisms of $V${\rm)}, respectively. A~{\it twisted intertwining operator of type $\binom{W_{3}}{W_{1}W_{2}}$} is a~linear map
\begin{gather*}
\mathcal{Y}\colon\ W_{1}\otimes W_{2}\to W_{3}\{x\}[\log x],\\
w_{1}\otimes w_{2}\mapsto \mathcal{Y}(w_{1}, x)w_{2}=\sum_{k=0}^{K}\sum_{n\in \C}
\mathcal{Y}_{n, k}(w_{1})w_{2}x^{-n-1}(\log x)^{k}
\end{gather*}
satisfying the following conditions:
\begin{enumerate}[label=$(\arabic*)$]\itemsep=0pt
\item {\it The lower truncation property}: For $w_{1}\in W_{1}$ and $w_{2}\in W_{2}$, $n\in \C$ and $k=0, \dots, K$,
$\mathcal{Y}_{n+l, k}(w_{1})w_{2}=0$ for $l\in \N$ sufficiently large.

\item The {\it duality property}:
For $u\in V$, $w_{1}\in W_{1}$, $w_{2}\in W_{2}$
and $w_{3}'\in W_{3}'$, there exists a~maximally extended
multivalued analytic function
$f(z_1, z_2; u, w_{1}, w_{2}, w_{3}')$
on $M^{2}$
with a~preferred single-valued branch
$f^e(z_1, z_2; u, w_{1}, w_{2}, w_{3}')$ on $M_{0}^{2}$
such that the series
\begin{gather}
\big\langle w'_{3}, Y_{W_{3}}
^{g_{3}}(u, z_{1})\mathcal{Y}(w_{1}, z_{2})w_{2}\big\rangle
=\sum_{n\in \C}\big\langle w'_{3}, Y_{W_{3}}^{g_{3}}(u, z_{1})\pi_{n}
\mathcal{Y}(w_{1}, z_{2})w_{2}\big\rangle,
\label{int-prod}\\
\big\langle w'_{3}, \mathcal{Y}(w_{1}, z_{2})
Y_{W_{2}}^{g_{2}}(u, z_{1})w_{2}\big\rangle
=\sum_{n\in \C}\big\langle w'_{3}, \mathcal{Y}(w_{1}, z_{2})
\pi_{n}Y_{W_{2}}^{g_{2}}(u, z_{1})w_{2}\big\rangle,
\label{int-rev-prod}\\
\big\langle w'_{3}, \mathcal{Y}\bigl(
Y_{W_{1}}^{g_{1}}(u, z_{1}-z_{2})w_{1}, z_{2}\bigr)w_{2}\big\rangle
=\sum_{n\in \C}\big\langle w'_{3}, \mathcal{Y}(\pi_{n}
Y_{W_{1}}^{g_{1}}(u, z_{1}-z_{2})w_{1}, z_{2})w_{2}\big\rangle
\label{int-iter}
\end{gather}
are absolutely convergent on the regions
$|z_{1}|>|z_{2}|>0$,
$|z_{2}|>|z_{1}|>0$, $|z_{2}|>{|z_{1}-z_{2}|>0}$, respectively. Moreover, their sums are equal
to $f^e(z_1, z_2; u, w_{1}, w_{2}, w_{3}')$ on the region given by
$|z_{1}|>|z_{2}|>0$ and $|{\arg (z_{1}-z_{2})-\arg z_{1}}|<\frac{\pi}{2}$,
the region given by ${|z_{2}|>|z_{1}|>0}$ and~$-\frac{3\pi}{2}< \arg (z_{1}-z_{2})-\arg z_{2}<-\frac{\pi}{2}$,
the region given by ${|z_{2}|>|z_{1}-z_{2}|>0}$ and $|{\arg z_{1}- \arg z_{2}}|<\frac{\pi}{2}$, respectively.

\item The {\it convergence and analytic extension
for products with more than one twisted vertex
operators}:
For $k\in \N+3$, $u_{1}, \dots, u_{k-1} \in V$,
$w_{1}\in W_{1}$, $w_{2}\in W_{2}$ and $w_{3}'\in W_{3}'$,
the series
\begin{gather*}
\big\langle w'_{3}, Y_{W_{3}}^{g_{3}}(u_{1}, z_{1})
\cdots Y_{W_{3}}^{g_{3}}(u_{k-1}, z_{k-1})
\mathcal{Y}(w_{1}, z_{k})w_{2}\big\rangle\\
\qquad=\sum_{n_{1}, \dots, n_{k-1}\in \C}\big\langle w'_{3}, Y_{W_{3}}
^{g_{3}}(u_{1}, z_{1})\pi_{n_{1}}
\cdots \pi_{n_{k-2}}Y_{W_{3}}^{g_{3}}(u_{k-1}, z_{k-1})
\pi_{k-1}\mathcal{Y}(w_{1}, z_{k})w_{2}\big\rangle
\end{gather*}
is absolutely convergent on the region
$|z_{1}|>\cdots >|z_{k}|>0$ and can be maximally extended to a~multivalued analytic function $f(z_{1}, \dots, z_{k})$ on the region
$M^{k}$. Moreover, for any fixed $z_{k}\in \C^{\times}$ and $M_{ij}\in \Z_{+}$ such that $x^{M_{ij}}Y_{V}(u_{i}, x)u_{j}\in V[[x]]$
for $i, j=1, \dots, {k-1}$, any component-isolated singularity of the form
$(z_1,\dots, z_{k-1})-\beta=\delta$ for $\beta\in \C^{k-1}$
and~${\delta\in\{0,\infty\}^{k-1}}$ of the multivalued analytic function of $z_{1}, \dots, z_{k-1}$
obtained as a~value at $\zeta=z_{k}$ of
$\prod_{1\le i<j\le k-1}(z_{i}-z_{j})^{M_{ij}}f(z_{1}, \dots, \zeta)$
 is regular.

\item The {\it $L(-1)$-derivative property}
$\frac{\rm d}{{\rm d}x}\mathcal{Y}(w_{1}, x)=\mathcal{Y}(L(-1)w_{1}, x)$.
\end{enumerate}
\end{defn}

\begin{rema} Note that the series in the definition of twisted intertwining operator
$\Y$ above are all summed over unique expansion sets (see \cite[Definition 7.5]{HLZ5}). In fact, by the duality
property and the fact that $\omega$ is fixed under $g_{1}$, $g_{2}$ and $g_{3}$,
the $L(0)$-commutator formula for $\Y$ holds. Then for
homogeneous $v\in V$, $w_{1}\in W_{1}$,
$w_{2}\in W_{2}$ and $w_{3}'\in W'_{3}$, the series above are in fact summed
over sets of the form $(\{r_{1}, \dots, r_{m}\}+\N)\times \{0, \dots, K\}$, which is a~unique expansion set (see~\cite{H-tensor-applicability}).
\end{rema}

\begin{rema}\label{general-branch}
 For simplicity, in the duality property in
Definition~\ref{def-tw-int-op}, we
use only the preferred branches of the twisted intertwining operator and
the twisted vertex operators. One can derive what the products and iterates
of other branches of a~twisted intertwining operator
converge to using the actions of the elements of $\PB_{3}$
on the single-valued branch
$f^e(z_1, z_2; u, w_{1}, w_{2}, w_{3}')$ of the multivalued function
$f(z_1, z_2; u, w_{1}, w_{2}, w_{3}')$ in the definition.
Let $\mathcal{Y}$ be a~twisted intertwining operator of type $\binom{W_{3}}{W_{1}W_{2}}$.
For any $p_1,p_2,p_{12}\in\Z$, the series
\begin{gather*}
\big\langle w'_{3}, \bigl(Y_{W_{3}}^{g_{3}}\bigr)^{p_{1}}(u, z_{1})
\mathcal{Y}^{p_{2}}(w_{1}, z_{2})w_{2}\big\rangle
=\sum_{n\in \C}\big\langle w'_{3},
\bigl(Y_{W_{3}}^{g_{3}}\bigr)^{p_{1}}(u, z_{1})\pi_{n}
\mathcal{Y}^{p_{2}}(w_{1}, z_{2})w_{2}\big\rangle,
\\
\big\langle w'_{3}, \mathcal{Y}^{p_{2}}(w_{1}, z_{2})
\bigl(Y_{W_{2}}^{g_{2}}\bigr)^{p_{1}}(u, z_{1})w_{2}\big\rangle
=\sum_{n\in \C}\big\langle w'_{3},
\mathcal{Y}^{p_{2}}(w_{1}, z_{2})
\pi_{n}\bigl(Y_{W_{2}}^{g_{2}}\bigr)^{p_{1}}(u, z_{1})w_{2}\big\rangle,
\\
\big\langle w'_{3}, \mathcal{Y}^{p_{2}}(
\bigl(Y_{W_{1}}^{g_{1}}\bigr)^{p_{12}}(u, z_{1}-z_{2})w_{1}, z_{2})
w_{2}\rangle
=\sum_{n\in \C}\big\langle w'_{3}, \mathcal{Y}^{ p_{2}}(\pi_{n}
\bigl(Y_{W_{1}}^{g_{1}}\bigr)^{p_{12}}(u, z_{1}-z_{2})w_{1}, z_{2})
w_{2}\big\rangle
\end{gather*}
are absolutely convergent on the regions
$|z_{1}|>|z_{2}|>0$,
$|z_{2}|>|z_{1}|>0$, $|z_{2}|>|z_{1}-z_{2}|>0$, respectively.
Moreover, their sums are equal to the branches
\begin{gather*}
f^{(b_{13}b_{12})^{p_{1}}b_{23}^{p_{2}}}
(z_{1}, z_{2}; u, w_{1}, w_{2}, w_{3}'),\qquad
f^{(b_{12}b_{23})^{p_{2}}b_{13}^{p_{1}}}
(z_{1}, z_{2}; u, w_{1}, w_{2}, w_{3}'),\\
f^{(b_{23}b_{13})^{p_{2}}b_{12}^{p_{12}}}
(z_{1}, z_{2}; u, w_{1}, w_{2}, w_{3}'),
\end{gather*}
respectively, of $f(z_1, z_2; u, w_{1}, w_{2}, w_{3}')$
on the region given by
$|z_{1}|>|z_{2}|>0$ and $|{\arg (z_{1}-z_{2})}-\arg z_{1} |
<\frac{\pi}{2}$,
the region given by $|z_{2}|>|z_{1}|>0$ and $-\frac{3\pi}{2}
< \arg (z_{1}-z_{2})-\arg z_{2}<-\frac{\pi}{2}$,
the region given by $|z_{2}|>|z_{1}-z_{2}|>0$ and
$|{\arg z_{1}- \arg z_{2}}|<\frac{\pi}{2}$, respectively.
See~\cite{D} for more details.
\end{rema}

\begin{rema}\label{fixed-pt-subalg-int}
 From Definition~\ref{def-tw-int-op}, we see that a~twisted intertwining operator
of type $\binom{W_{3}}{W_{1}W_{2}}$ is in fact
an intertwining operator of type \smash{$\binom{W_{3}}{W_{1}W_{2}}$} when $W_{1}$, $W_{2}$, and $W_{3}$ are
viewed as (untwisted) modules for the fixed-point subalgebra of $V$ under
$g_{1}$ and $g_{2}$. In particular, all the usual properties for intertwining operators
holds for the vertex operators
associated to elements of the fixed-point subalgebra and the twisted intertwining operator.
For example, the Jacobi identity and the commutator formula for these operators are
satisfied.
\end{rema}

\begin{prop}
Let $g_{1}$, $g_{2}$, $g_{3}$ be automorphisms of $V$,
$W_{1}$, $W_{2}$, $W_{3}$ $g_{1}$-, $g_{2}$-, $g_{3}$-twisted
generalized $V$-modules and
 $\Y$ a~twisted intertwining operator of type
 $\binom{W_{3}}{W_{1}W_{2}}$.
 Assume that the map $u\mapsto Y_{W_{3}}^{g_{3}}(u,x)$ is injective
and $\Y$ is surjective in the sense that the coefficients of
the series~${\mathcal{Y}(w_1,x)w_2}$ for $w_1\in W_{1}$, $w_2\in
W_{2}$ span $W_3$.
Then $g_{3}=g_{1}g_{2}$.
\end{prop}
\begin{proof}
By the definition of twisted intertwining
operator, for $u\in V$, $w_{1}\in W_{1}$, $w_{2}\in W_{2}$
and $w_{3}'\in W_{3}'$, there exists a~multivalued analytic function
$f(z_1,z_2; u,w_1,w_2,w_3')$ on $M^{2}$ with a~preferred
single-valued branch $f^{e}(z_1,z_2; u,w_1,w_2,w_3')$
on $M^{2}_{0}$ such that
\eqref{int-prod}, \eqref{int-rev-prod} and~\eqref{int-iter}
are absolutely convergent to $f^{e}(z_1,z_2; u,w_1,w_2,w_3')$
on the corresponding regions in given in Definition~\ref{def-tw-int-op}.
In particular, on the region
$|z_{1}|>|z_{2}|>0$, $|{\arg (z_{1}-z_{2})-\arg z_{1}}|<\frac{\pi}{2}$
the sum of the series
$\big\langle w'_{3}, Y_{W_{3}}^{g_{3}}(g_3u, z_{1})
\mathcal{Y}(w_{1}, z_{2})w_{2}\big\rangle$ is equal to
$f^{e}(z_1,z_2;g_3u,w_1,w_2,w_3')$.
By the equivariance property for $W_{3}$,
\begin{equation}\label{prod-g_3}
\big\langle w'_{3}, Y_{W_{3}}^{g_{3}}(g_3u, z_{1})
\mathcal{Y}(w_{1}, z_{2})w_{2}\big\rangle
=\big\langle w'_{3}, \bigl(Y_{W_{3}}^{g_{3}}\bigr)^{-1}(u, z_{1})
\mathcal{Y}(w_{1}, z_{2})w_{2}\big\rangle,
\end{equation}
where as above $\bigl(Y_{W_{3}}^{g_{3}}\bigr)^{-1}$ is the $(-1)$-th branch of
the twisted vertex operator map $Y_{W_{3}}^{g_{3}}$.
But the right-hand side of~\eqref{prod-g_3} can be obtained
from $\langle w'_{3}, Y_{W_{3}}^{g_{3}}(u, z_{1})
\mathcal{Y}(w_{1}, z_{2})w_{2}\rangle$
on the region $|z_{1}|>|z_{2}|>0$, $|{\arg (z_{1}-z_{2})-\arg z_{1}}|<\frac{\pi}{2}$
by letting $z_{1}$ go around clockwise a~circle of radius larger than $|z_{2}|$. The homotopy
class of such a~circle is equal to $(b_{13}b_{12})^{-1}$.
So~\eqref{prod-g_3} gives
\[f^{e}(z_1,z_2;g_3u,w_1,w_2,w_3')
=f^{(b_{13}b_{12})^{-1}}(z_1,z_2;u,w_1,w_2,w_3')\]
or equivalently
\begin{equation}\label{prod-g_3-1}
f^{b_{13}b_{12}}(z_1,z_2;g_3u,w_1,w_2,w_3')=f^e(z_1,z_2;u,w_1,w_2,w_3').
\end{equation}

Similarly, on the region
$|z_{2}|>|z_{1}|>0$, $-\frac{3\pi}{2}<\arg (z_{1}-z_{2})
-\arg z_{2}<-\frac{\pi}{2}$,
the sum of the series
$\langle w'_{3}, \mathcal{Y}(w_{1}, z_{2})
Y_{W_{3}}^{g_{2}}(g_2u, z_{1})w_{2}\rangle$ is equal to
$f^{e}(z_1,z_2;g_2u,w_1,w_2,w_3')$.
By the equivariance property for $W_{2}$,
\begin{equation}\label{prod-g_2}
\big\langle w'_{3}, \mathcal{Y}(w_{1}, z_{2})
Y_{W_{3}}^{g_{2}}(g_2u, z_{1})w_{2}\big\rangle
=\big\langle w'_{3}, \mathcal{Y}(w_{1}, z_{2})
\bigl(Y_{W_{2}}^{g_{2}}\bigr)^{-1}(u, z_{1})w_{2}\big\rangle,
\end{equation}
where \smash{$\bigl(Y_{W_{2}}^{g_{2}}\bigr)^{-1}$} is the $(-1)$-th branch of
the twisted vertex operator map \smash{$Y_{W_{2}}^{g_{2}}$}.
The right-hand side of~\eqref{prod-g_2} can be obtained
from \smash{$\big\langle w'_{3}, \mathcal{Y}(w_{1}, z_{2})
Y_{W_{3}}^{g_{2}}(u, z_{1})w_{2}\big\rangle$}
on the region
$|z_{2}|>|z_{1}|>0$, $-\frac{3\pi}{2}<\arg (z_{1}-z_{2})
-\arg z_{2}<-\frac{\pi}{2}$
by letting $z_{1}$ go around clockwise a~circle of radius less than $|z_{2}|$. Such a~circle as a~loop must have a~base point in the region
$|z_{2}|>|z_{1}|>0$, $-\frac{3\pi}{2}<\arg (z_{1}-z_{2})
-\arg z_{2}<-\frac{\pi}{2}$. But there is a~canonical isomorphism
between the fundamental group of $M^{2}$ with such a~base point
and $\PB_{3}$ which has a~base point $(-3, -2)$.
To see how the loop given by the circle above
acts on the single-valued branches of
$f(z_1,z_2;u,w_1,w_2,w_3')$, we need to find the element
of $\PB_{3}$ corresponding to this loop.
We choose the following loop $\gamma$ based at
$(-3, -2)$:
The first part $\gamma_{1}$
of $\gamma$ is the lower half circle centered at $-2$
with radius $1$ from $-3$ to $-1$ for $z_{1}$ and trivial for $z_{2}$
(meaning $z_{2}$ is always equal to $-2$).
The second part $\gamma_{2}$
is the counterclockwise circle centered at $0$ with
radius $1$ based at $-1$ for $z_{1}$ and trivial for $z_{2}$.
The third part $\gamma_{3}
=\gamma_{1}^{-1}$ is also
the lower half circle centered at $-1$
with radius $1$ but from $-1$ to $-3$ for $z_{1}$
and trivial for $z_{2}$. It is clear that $\gamma$ is
homotopically equivalent to the loop given in Figure~\ref{b13}.
Then we have $b_{13}=[\gamma]=
[\gamma_{1}][\gamma_{2}]
[\gamma_{1}]^{-1}$, where $[\gamma]$ for a~path $\gamma$
means its homotopy class. Equivalently, we have
$[\gamma_{2}]=[\gamma_{1}]^{-1} b_{13}
[\gamma_{1}]$.
When we let $z_{1}$ go from $-1$ to $-3$ along~$\gamma_{1}^{-1}$,
since $\gamma_{1}^{-1}$ passes the cut along the positive real line
in the $z_{1}-z_{2}$-plane, the single-valued branch
$f^{e}(z_1,z_2;u,w_1,w_2,w_3')$ is changed to the single-valued
branch \smash{$f^{b_{12}^{-1}}(z_1,z_2;u,w_1,w_2,w_3')$}.
Similarly,
when we let $z_{1}$ go from $-1$ to $-3$ along $\gamma_{1}$,
$f^{b}(z_1,z_2;u,w_1,w_2,w_3')$ is changed to the single-valued
branch $f^{bb_{12}}(z_1,z_2;u,w_1,w_2,w_3')$ for any $b\in \PB_{3}$.
Thus when we
let $z_{1}$ go around the loop $\gamma_{2}$,
$f^{e}(z_1,z_2;u,w_1,w_2,w_3')$ is changed to
\smash{$f^{b_{12}^{-1}b_{13}b_{12}}(z_1,z_2;u,w_1,w_2,w_3')$}, that is,
\[f^{[\gamma_{2}]}(z_1,z_2;u,w_1,w_2,w_3')
=f^{b_{12}^{-1}b_{13}b_{12}}(z_1,z_2;u,w_1,w_2,w_3').\]

Note that the circle $\gamma_{2}^{-1}$ is exactly the
circle we use to obtain
the right-hand side of~\eqref{prod-g_2} from~${\big\langle w'_{3}, \mathcal{Y}(w_{1}, z_{2})
Y_{W_{3}}^{g_{2}}(u, z_{1})w_{2}\big\rangle}$
on the region
$|z_{2}|>|z_{1}|>0$, $-\frac{3\pi}{2}<\arg (z_{1}-z_{2})
-\arg z_{2}<-\frac{\pi}{2}$. So the right-hand side of
\eqref{prod-g_2} is equal to
\[f^{[\gamma_{2}^{-1}]}(z_1,z_2;u,w_1,w_2,w_3')
=f^{b_{12}^{-1}b_{13}^{-1}b_{12}}(z_1,z_2;u,w_1,w_2,w_3').\]
But the sum of the left-hand side of~\eqref{prod-g_2}
is equal to $f^{e}(z_1,z_2;g_{2}u,w_1,w_2,w_3')$.
So we obtain
\[f^{e}(z_1,z_2;g_2u,w_1,w_2,w_3')
=f^{b_{12}^{-1}b_{13}^{-1}b_{12}}(z_1,z_2;u,w_1,w_2,w_3')\]
or equivalently
\begin{equation}\label{prod-g_2-1}
f^{b_{12}^{-1}b_{13}b_{12}}(z_1,z_2;g_2u,w_1,w_2,w_3')=f^e(z_1,z_2;u,w_1,w_2,w_3').
\end{equation}

Similarly we also have
\begin{equation}\label{prod-g_1-1}
f^{b_{12}}(z_1,z_2;g_1u,w_1,w_2,w_3')=f^e(z_1,z_2;u,w_1,w_2,w_3').
\end{equation}
Using the right
action of $\PB_{3}$ on the set of single-valued branches of
the multivalued analytic function
$f(z_1,z_2;u,w_1,w_2,w_3')$, for
$b\in \PB_{3}$, we obtain from~\eqref{prod-g_1-1}
\begin{equation}\label{prod-g_1-2}
f^{b_{12}b}(z_1,z_2;g_2u,w_1,w_2,w_3')
=f^b(z_1,z_2;u,w_1,w_2,w_3').
\end{equation}

From~\eqref{prod-g_3-1}, \eqref{prod-g_2-1} and
\eqref{prod-g_1-2}, we have
\begin{align}
f^{b_{13}b_{12}}(z_1,z_2;g_3u,w_1,w_2,w_3')&=f^e(z_1,z_2;u,w_1,w_2,w_3')=f^{b_{12}^{-1}b_{13}b_{12}}(z_1,z_2;g_2u,w_1,w_2,w_3')\notag\\
&=f^{b_{12}(b_{12}^{-1}b_{13}b_{12})}
(z_1,z_2;g_1g_2u,w_1,w_2,w_3')\nonumber\\
&=f^{b_{13}b_{12}}
(z_1,z_2;g_1g_2u,w_1,w_2,w_3')\label{2.27}
\end{align}
When $|z_1|>|z_2|>0$ and $|{\arg(z_1-z_2)-\arg z_1}|<\frac{\pi}{2}$,
the left- and right-hand sides of~\eqref{2.27} are equal to
the sum of \smash{$\big\langle w'_{3}, \bigl(Y_{W_{3}}^{g_{3}}\bigr)^{1}(g_3u, z_{1})
\mathcal{Y}(w_{1}, z_{2})w_{2}\big\rangle$} and
\smash{$\big\langle w'_{3}, \bigl(Y_{W_{3}}^{g_{3}}\bigr)^{1}(g_1g_2u, z_{1})
\mathcal{Y}(w_{1}, z_{2})w_{2}\big\rangle$}, respectively.
Therefore, we obtain
\begin{equation}\label{4.18-in-H}
\big\langle w'_{3}, \bigl(Y_{W_{3}}
^{g_{3}}\bigr)^{1}(g_1g_2u-g_3u, z_{1})
\mathcal{Y}(w_{1}, z_{2})w_{2}\big\rangle=0
\end{equation}
for $w_1\in W_1$, $w_2\in W_2$, $w_3\in W_3'$.
Since $\Y$ is surjective in the sense above and $w'_{3}$
is arbitrary, \eqref{4.18-in-H}
implies \smash{$\bigl(Y_{W_{3}}
^{g_{3}}\bigr)^{1}(g_1g_2u-g_3u, z_{1})=0$}.
But the map given by
$v\to Y_{W_{3}}
^{g_{3}}(v, x)$ is injective, we obtain $g_1g_2u-g_3u=0$
for $u\in V$. So we have
$g_{3}=g_{1}g_{2}$.
\end{proof}

\section{Skew-symmetry and contragredient isomorphisms}\label{sec3}

Let $g_{1}$, $g_{2}$ be automorphisms of $V$, $W_{1}$, $W_{2}$ and $W_{3}$ $g_{1}$-, $g_{2}$-
and $g_{1}g_{2}$-twisted
$V$-modules without $g_{1}$-, $g_{2}$- and $g_{1}g_{2}$-actions, respectively
and $\mathcal{Y}$ a~twisted intertwining operator
of type $\binom{W_{3}}{W_{1}W_{2}}$.
We define linear maps
$
\Omega_{\pm}(\Y)\colon W_{2}\otimes W_{1}\to W_{3}\{x\}[\log x]$, $
w_{2}\otimes w_{1}\mapsto \Omega_{\pm}(\Y)(w_{2}, x)w_{1}
$
by
\begin{equation}\label{omega1}
\Omega_{\pm}(\Y)(w_{2}, x)w_{1}={\rm e}^{xL(-1)}\Y(w_{1}, y)w_{2}\lbar_{y^{n}={\rm e}^{\pm n\pi \i}x^{n},
\log y=\log x\pm \pi \i}
\end{equation}
for $w_{1}\in W_{1}$ and $w_{2}\in W_{2}$.
Note that we can also define $\Omega_{p}$ for $p\in \Z$ by
changing $\pm $ in the right-hand side of
 \eqref{omega1} to $+(2p+1)$. But we will not discuss
$\Omega_{p}$ in this paper.

From the definition~\eqref{omega1}, for $p\in \Z$, $w_{1}\in W_{1}$, $w_{2}\in W_{2}$ and $z\in \C^{\times}$,
\begin{align*}
\Omega_{\pm}(\Y)^{p}(w_{2}, z)w_{1}&=\Omega_{\pm}(\Y)(w_{2}, x)w_{1}\lbar_{x^{n}={\rm e}^{nl_{p}(z)}, \log x=l_{p}(z)}\\
&=\big({\rm e}^{xL(-1)}\Y(w_{1}, y)w_{2}\lbar_{y^{n}={\rm e}^{\pm n\pi \i}x^{n},
\log y=\log x\pm \pi \i}\big)\lbar_{x^{n}={\rm e}^{nl_{p}(z)}, \log x=l_{p}(z)}\\
&={\rm e}^{zL(-1)}\Y(w_{1}, y)w_{2}\lbar_{y^{n}={\rm e}^{n(l_{p}(z)\pm \pi \i)}, \log y=l_{p}(z)\pm \pi \i}.
\end{align*}
When $\arg z<\pi$ and $\arg z\ge \pi$, $\arg (-z)=\arg z+\pi$ and $\arg (-z)=\arg z-\pi$, respectively.
Hence
\[{\rm e}^{zL(-1)}\Y(w_{1}, y)w_{2}\lbar_{y^{n}={\rm e}^{n(l_{p}(z)+ \pi \i)}, \log y=l_{p}(z)+ \pi \i}=
{\rm e}^{zL(-1)}\Y^{p}(w_{1}, -z)w_{2}
\]
when $\arg z<\pi$ and
\[{\rm e}^{zL(-1)}\Y(w_{1}, y)w_{2}\lbar_{y^{n}={\rm e}^{n(l_{p}(z)- \pi \i)}, \log y=l_{p}(z)- \pi \i}=
{\rm e}^{zL(-1)}\Y^{p}(w_{1}, -z)w_{2}\]
when $\arg z\ge \pi$. In particular, for $w_{1}\in W_{1}$, $w_{2}\in W_{2}$ and $z\in \C^{\times}$ satisfying
$\arg z<\pi$ and $\arg z \ge \pi$, we have
\begin{equation}\label{omega2}
\Omega_{+}(\Y)^{p}(w_{2}, z)w_{1}={\rm e}^{zL(-1)}\Y^{p}(w_{1}, -z)w_{2}
\end{equation}
and
\begin{equation}\label{omega3}
\Omega_{-}(\Y)^{p}(w_{2}, z)w_{1}={\rm e}^{zL(-1)}\Y^{p}(w_{1}, -z)w_{2},
\end{equation}
respectively.

Now assume that $W_{1}$, $W_{2}$, and $W_{3}$ are $g_{1}$-, $g_{2}$-, and $g_{1}g_{2}$-twisted $V$-modules with
$g_{1}$-, $g_{2}$- and $g_{1}g_{2}$-actions, respectively.
In the case that $W_{1}$ also has an action of $g_{2}^{-1}$,
we have the $g_{2}^{-1}g_{1}g_{2}$-twisted $V$-module $\varphi_{g_{2}^{-1}}(W_{1})$.
In this case, we define
\begin{gather*}
\Omega_{+}^{g_{2}^{-1}}(\Y)\colon\ W_{2}\otimes W_{1}\to W_{3}\{x\}[\log x],\qquad
w_{2}\otimes w_{1}\mapsto \Omega_{+}^{g_{2}^{-1}}(\Y)(w_{2}, x)w_{1}
\end{gather*}
by
\smash{$
\Omega_{+}^{g_{2}^{-1}}(\Y)(w_{2}, x)w_{1}={\rm e}^{xL(-1)}\Y(g_{2}w_{1}, y)w_{2}
\big|_{y^{n}={\rm e}^{n\pi \i}x^{n},
\log y=\log x+ \pi \i}
$}
for $w_{1}\in W_{1}$ and $w_{2}\in W_{2}$.
In the case that $W_{2}$ has an action of $g_{1}$,
we have the $g_{1}g_{2}g_{1}^{-1}$-twisted $V$-module $\varphi_{g_{1}}(W_{2})$.
In this case, we define
\begin{gather*}
\Omega_{-}^{g_{1}}(\Y)\colon\ W_{2}\otimes W_{1}\to W_{3}\{x\}[\log x],\qquad
w_{2}\otimes w_{1}\mapsto \Omega_{-}^{g_{1}}(\Y)(w_{2}, x)w_{1}
\end{gather*}
by
\smash{$
\Omega_{-}^{g_{1}}(\Y)(w_{2}, x)w_{1}={\rm e}^{xL(-1)}\Y(w_{1}, y)g_{1}^{-1}w_{2}\big|_{y^{n}={\rm e}^{- n\pi \i}x^{n},
\log y=\log x- \pi \i}
$}
for $w_{1}\in W_{1}$ and $w_{2}\in W_{2}$.

For $w_{1}\in W_{1}$, $w_{2}\in W_{2}$ and $z\in \C^{\times}$ satisfying
$\arg z<\pi$ and $\arg z \ge \pi$, we have
\begin{equation}\label{omega4}
\Omega_{+}^{g_{2}^{-1}}(\Y)^{p}(w_{2}, z)w_{1}={\rm e}^{zL(-1)}g_{2}\Y^{p}(g_{2}w_{1}, -z)g_{2}^{-1}w_{2}
\end{equation}
and
\begin{equation}\label{omega5}
\Omega_{-}^{g_{1}}(\Y)^{p}(w_{2}, z)w_{1}={\rm e}^{zL(-1)}\Y^{p}(w_{1}, -z)g_{1}^{-1}w_{2},
\end{equation}
respectively.

\begin{thm}\label{skew-sym}
The linear maps $\Omega_{+}(\Y)$, $\Omega_{-}(\Y)$, \smash{$\Omega_{+}^{g_{2}^{-1}}(\Y)$} and
$\Omega_{-}^{g_{1}}(\Y)$ are twisted intertwining operators
of types
\[
\binom{W_{3}}{W_{2} \phi_{g_{2}^{-1}}(W_{1})} , \quad \binom{W_{3}}{\phi_{g_{1}}(W_{2}) W_{1}},\quad
 \binom{W_{3}}{W_{2} \varphi_{g_{2}^{-1}}(W_{1})}, \quad \binom{W_{3}}{\varphi_{g_{1}}(W_{2}) W_{1}},
 \]
respectively $($recall the~definition of $\phi_{g}$ and $\varphi_{g}$ for an automorphism $g$ of $V$ in Section~{\rm\ref{sec2})}.
\end{thm}
\begin{proof}
We need only prove the results for $\Omega_{+}(\Y)$ and $\Omega_{-}(\Y)$. Using
the definitions of $\varphi_{g_{1}}$ and~$\varphi_{g_{2}^{-1}}$ and
\eqref{omega4} and~\eqref{omega5}, we see that the proofs for
\smash{$\Omega_{+}^{g_{2}^{-1}}(\Y)$} and $\Omega_{-}^{g_{1}}(\Y)$ are reduced to
 the proofs for $\Omega_{+}(\Y)$ and $\Omega_{-}(\Y)$.

The proofs of the lower-truncation property and the $L(-1)$-derivative property of
$\Omega_{+}(\Y)$ and~$\Omega_{-}(\Y)$ are easy to verify. As in the proof of
\cite[Theorem 5.1]{H-twisted-int}, we omit their proofs here.
The main difference between the proof here and the proof of~\cite[Theorem 5.1]{H-twisted-int} is that here we cannot use the explicit
form of the correlation functions in~\cite{H-twisted-int}
obtained from the products and iterates
of twisted intertwining operators and twisted vertex operators.
So our proof here is much more complicated and involves
some technical convergence and analytic extension results, even though the
idea is the same as in the proof of \cite[Theorem 5.1]{H-twisted-int}.

Let $u\in V$, $w_{1}\in W_{1}$, $w_{2}\in W_{2}$
and $w_{3}'\in W_{3}'$. As in the proof of \cite[Theorem 5.1]{H-twisted-int},
we use~${f(z_1, z_2; u, w_{1}, w_{2}, w_{3}')}$ to denote
the multivalued analytic function
in the duality property for the twisted intertwining operator $\Y$ with
the preferred branch $f^{e}(z_1, z_2; u, w_{1}, w_{2}, w_{3}')$.
Note that~${f(z_1, z_2; u, w_{1}, w_{2}, w_{3}')}$
in~\cite{H-twisted-int} is of the particular
form in the definition of twisted intertwining operator there.
But in this proof, it is a~maximally-extended multivalued
analytic function on~$M^{2}$ with regular singular points
at $z_{1}, z_{2}=0$ and $z_{1}-z_{2}=0$ and in general
might not have the
special form in~\cite{H-twisted-int}.

Define
\[
g_{\pm}(z_{1}, z_{2}; u, w_{2}, w_{1}, w_{3}')
 = f\bigl(z_1-z_{2}, -z_2; u, w_{1}, w_{2}, {\rm e}^{z_{2}L'(1)}w_{3}'\bigr)
\]
and choose the preferred branch
$g_{\pm}^{e}(z_{1}, z_{2}; u, w_{2}, w_{1}, w_{3}')$ of
$g_{\pm}(z_{1}, z_{2}; u, w_{2}, w_{1}, w_{3}')$ as follows:
On the subregion $|z_{1}|>|z_{2}|>0$, 
$|{\arg (z_{1}-z_{2})-\arg z_{1}}|<\frac{\pi}{2}$ and
$\arg z_{2}< \pi$ (for $\Omega_{+}$) or $\arg z_{2}\ge \pi$
(for $\Omega_{-}$) of $M^{2}_{0}$, let
\[
g_{\pm}^{e}(z_{1}, z_{2}; u, w_{2}, w_{1}, w_{3}')
 = f^{e}\bigl(z_1-z_{2}, -z_2; u, w_{1}, w_{2}, {\rm e}^{z_{2}L'(1)}w_{3}'\bigr).
\]
For general $(z_{1}, z_{2})\in M^{2}_{0}$, we define
$g_{\pm}^{e}(z_{1}, z_{2}; u, w_{2}, w_{1}, w_{3}')$
to be the unique analytic extension on $M^{2}_{0}$.

When $|z_{1}|>|z_{2}|>0$, $|{\arg (z_{1}-z_{2})-\arg z_{1}}|<\frac{\pi}{2}$
and $\arg z_{2}< \pi$ (for $\Omega_{+}$) or $\arg z_{2}\ge \pi$
(for~$\Omega_{-}$),
from~\eqref{omega2} and~\eqref{omega3} and the $L(-1)$-derivative
property for $Y_{W_{3}}
^{g_{3}}$, we have
\begin{gather}
\big\langle w'_{3}, Y_{W_{3}}
^{g_{3}}(u, z_{1})\Omega_{\pm}(\mathcal{Y})(w_{2}, z_{2})w_{1}\big\rangle\nonumber\\
\qquad=\big\langle w'_{3}, Y_{W_{3}}
^{g_{3}}(u, x_{1})\Omega_{\pm}(\mathcal{Y})(w_{2}, x_{2})w_{1}\big\rangle\mbar_{x_{1}^{n}={\rm e}^{n\log z_{1}}, \log x_{1}
=\log z_{1},
x_{2}^{n}={\rm e}^{n\log z_{2}}, \log x_{2}=\log z_{2}}\nonumber\\
\qquad=\big\langle w'_{3}, Y_{W_{3}}
^{g_{3}}(u, x_{1}){\rm e}^{-yL(-1)}\mathcal{Y}(w_{1},
y)w_{2}\big\rangle\mbar_{x_{1}^{n}={\rm e}^{n\log z_{1}}, \log x_{1}=\log z_{1},
y^{n}={\rm e}^{n\log (-z_{2})}, \log y=\log (-z_{2})}\nonumber\\
\qquad=\big\langle {\rm e}^{z_{2}L'(1)}w'_{3}, (Y_{W_{3}}
^{g_{3}})(u, x_{1}+y)\nonumber\\
\phantom{\qquad=}{}\times\mathcal{Y}(w_{1}, y)w_{2}\big\rangle\mbar_{x_{1}^{n}={\rm e}^{n\log z_{1}}, \log x_{1}=\log z_{1},
y^{n}={\rm e}^{n\log (-z_{2})}, \log y=\log (-z_{2})}.\label{skew-sym-0}
\end{gather}
We first prove that the right-hand side of formula~\eqref{skew-sym-0}
is absolutely convergent on the region $|z_{1}|>|z_{2}|>0$ and is
convergent to $g_{\pm}^{e}(z_1, z_2; u, w_{1}, w_{2}, w_{3}')$ on the region
$|z_{1}|>|z_{2}|>0$, $|{\arg (z_{1}-z_{2})-\arg z_{1}}|<\frac{\pi}{2}$.
The proof is in fact the same as the proof that the right-hand side of~\cite[equation~(9.170)]{HLZ6} is absolutely convergent in the region
$|z_{2}|>|z_{0}|>0$. Here we give a~slightly different proof.

We can always take $u\in V$, $w_{1}\in W_{1}$, $w_{2}\in W_{2}$
and ${\rm e}^{z_{2}L'(1)}w'_{3}\in W_{3}'$
to be homogeneous.
Let~${\Delta=-\wt {\rm e}^{z_{2}L'(1)}w_3'+\wt u+\wt w_1+\wt w_2}$.
Let
\begin{align*}
 D=\big\{n\in\C\mid \text{there exist } i, j\in \N, \,
\big\langle
\bigl(Y_{W_3}^{g_3}\bigr)_{\Delta-n-2, j}(u)
\mathcal{Y}_{n, i}(w_{1})w_{2}
\big\rangle \neq0\big\}
\end{align*}
and $M, N\in \N$ such that
\smash{$\bigl(Y_{W_3}^{g_3}\bigr)_{m, j}(u)=0$} for $m\in \C$, $j> M$
and $\Y_{n, i}(w_{1})=0$ for $i>N$.
Then by the lower truncation property of $\Y$, the fact that
$u$ is a~finite sum of generalized eigenvectors of $g_{3}$
and the equivariance property of the $g_{3}$-twisted module $W_{3}$,
we know that there exist a~finite subset $A$ of $\C/\Z$
and $R_\mu\in\mu$ for each $\mu\in A$ such that
\[
D\subset \bigcup_{\mu\in A}=(R_\mu-\N).
\]
Let
\[a_{n, j, i}=\big\langle {\rm e}^{z_{2}L'(1)}w_{3}^{\prime},
\bigl(Y_{W_3}^{g_3}\bigr)_{\Delta-n-2, j}(u)
(\mathcal{Y})_{n, i}(w_{1})w_{2}\big\rangle
\in \mathbb{C}\]
for $n\in D$, $0\le j\le M$ and $0\le n\le N$.
Then, by the convergence of
\eqref{int-prod}, the $L(-1)$-derivative properties for $Y_{W_3}^{g_3}$
and $\Y$, and \cite[Proposition 7.9]{HLZ5}, we know that the triple series
\begin{align}
&\sum_{n\in D}\sum_{j=0}^{M}\sum_{i=0}^{N}a_{n,j,i}
{\rm e}^{(-\Delta+n+1)\log z_1}(\log z_1)^j{\rm e}^{(-n-1)
\log z_2}(\log z_2)^i\label{kkl-1}
\end{align}
is absolutely convergent on the region given by
$|z_{1}|>|z_{2}|>0$ and is convergent on the region
$|z_{1}|>|z_{2}|>0$, $|{\arg (z_{1}-z_{2})-\arg z_{1}}|<\frac{\pi}{2}$ to
\[
f^{e}\bigl(z_{1}, z_{2}; u, w_{1}, w_{2}, {\rm e}^{z_{2}L'(1)}w_3'\bigr)
=\big\langle {\rm e}^{z_{2}L'(1)}w_3',Y_{W_3}^{g_3}(u, z_1)
\mathcal{Y}(w_1, z_2) w_2\big\rangle.
\]

For $n\in D$, $j=0,\dots,M$,
$k\in\Z_{\geq0}$, $s=0,\dots,j$, let
$b_{n,j,k,s}\in\C$ be the numbers defined
in~\eqref{b-defi}. Then
\begin{gather*}
\sum_{m\in D-\N}\sum_{s=0}^M
\sum_{i=0}^{N}
 \Bigg(\sum_{j=s}^{M}\sum_{\substack{n-k=m\\n\in D,\, k\in\N}}
a_{n,j,i}b_{n,j,k,s}\Bigg){\rm e}^{(-\Delta+m+1)
\log z_1}
(\log z_1)^s\\
\qquad\times
{\rm e}^{(-m-1)\log (-z_2)}(\log (-z_2))^i
\end{gather*}
 is equal to the right-hand side of~\eqref{skew-sym-0}
 and, by Lemma~\ref{convergence}, is absolutely convergent
on the region $|z_{1}|>|z_{2}|>0$ and is convergent to
\smash{$f^{e}\bigl(z_{1}-z_{2}, -z_{2};
u, w_{1}, w_{2}, {\rm e}^{z_{2}L'(1)}w_3'\bigr)$} on the region~${|z_1|>|z_2|>0}$,
$|{\arg z_1-\arg (z_1-z_{2})}|<\frac{\pi}{2}$.

Now it is easy to see that the left-hand side of~\eqref{skew-sym-0}
is absolutely convergent on the region~${|z_{1}|>|z_{2}|>0}$ and its sum is equal to
$g^{e}_{\pm}(z_{1}, z_{2};
u, w_{1}, w_{2}, w_3')$ on the region
$|z_1|>|z_2|>0$,
$|{\arg z_1-\arg (z_1-z_{2})}|<\frac{\pi}{2}$. In fact,
we know that the left-hand side of~\eqref{skew-sym-0}
as a~series is equal to the right-hand side of~\eqref{skew-sym-0}
on the region $|z_{1}|>|z_{2}|>0$,
$|{\arg (z_{1}-z_{2})-\arg z_{1}}|<\frac{\pi}{2}$
and $\arg z_{2}< \pi$ (for $\Omega_{+}$) or $\arg z_{2}\ge \pi$
(for $\Omega_{-}$). But we have just proved that
the right-hand side of~\eqref{skew-sym-0} is absolutely convergent
on the larger region $|z_{1}|>|z_{2}|>0$. The left-hand side of
\eqref{skew-sym-0} is also a~series of the same form
as~\eqref{kkl-1}. In particular, its absolute
convergence on the
smaller region $|z_{1}|>|z_{2}|>0$,
$|{\arg (z_{1}-z_{2})-\arg z_{1}}|<\frac{\pi}{2}$
and $\arg z_{2}< \pi$ (for $\Omega_{+}$) or $\arg z_{2}\ge \pi$
(for $\Omega_{-}$) implies its absolute
convergence on the
larger region $|z_{1}|>|z_{2}|>0$.

On the region $|z_{1}|>|z_{2}|>0$,
$|{\arg (z_{1}-z_{2})-\arg z_{1}}|<\frac{\pi}{2}$
and $\arg z_{2}< \pi$ (for $\Omega_{+}$) or $\arg z_{2}\ge \pi$
(for $\Omega_{-}$), by~\eqref{skew-sym-0} and the discussion
above, the left-hand side of~\eqref{skew-sym-0} is convergent to
\smash{$f^{e}\bigl(z_{1}-z_{2}, -z_{2};
u, w_{1}, w_{2}, {\rm e}^{z_{2}L'(1)}w_3'\bigr)$}, which in turn is
by definition equal to $g^{e}_{\pm}(z_{1}, z_{2};
u, w_{1}, w_{2},\allowbreak w_3')$ on the same region. We know that
the left-hand side of
\eqref{skew-sym-0} on the region given by
$|z_{1}|>|z_{2}|>0$ and $|{\arg (z_{1}-z_{2})-\arg z_{1}}|<\frac{\pi}{2}$
with cuts along the positive real lines
on the $z_{1}$- and $z_{2}$-planes is convergent to the analytic extension of
the sum of the left-hand side of~\eqref{skew-sym-0} on the smaller region
$|z_{1}|>|z_{2}|>0$,
$|{\arg (z_{1}-z_{2})-\arg z_{1}}|<\frac{\pi}{2}$
and $\arg z_{2}< \pi$ (for~$\Omega_{+}$) or $\arg z_{2}\ge \pi$
(for~$\Omega_{-}$). Also, by definition, $g^{e}_{\pm}(z_{1}, z_{2};
u, w_{1}, w_{2}, w_3')$ on $M_{0}^{2}$ is obtained by analytically
extending $g^{e}_{\pm}(z_{1}, z_{2};
u, w_{1}, w_{2}, w_3')$ on the smaller region
$|z_{1}|>|z_{2}|>0$,
$|{\arg (z_{1}-z_{2})-\arg z_{1}}|<\frac{\pi}{2}$
and $\arg z_{2}< \pi$ (for $\Omega_{+}$) or $\arg z_{2}\ge \pi$
(for $\Omega_{-}$). Thus the left-hand side of
\eqref{skew-sym-0} on the region given by~${|z_{1}|>|z_{2}|>0}$ and $|{\arg (z_{1}-z_{2})-\arg z_{1}}|<\frac{\pi}{2}$
is absolutely convergent to~${g^{e}_{\pm}(z_{1}, z_{2};
u, w_{1}, w_{2}, w_3')}$ on $M_{0}^{2}$.

Generalizing the proof of the convergence of~\eqref{skew-sym-0}
above, we can prove that
\begin{equation}\label{skew-sym-k-prod-1}
\big\langle w'_{3}, Y_{W_{3}}^{g_{3}}(u_{1}, z_{1})
\cdots Y_{W_{3}}^{g_{3}}(u_{k-1}, z_{k-1})
\Omega_{\pm}(\mathcal{Y})(w_{2}, z_{k})w_{1}\big\rangle
\end{equation}
is absolutely convergent on the region
$|z_{1}|>\cdots >|z_{k}|>0$ and can be maximally extended to a~multivalued analytic function on the region
$M^{k}$ for $k\in \N+3$, $u_{1}, \dots, u_{k-1} \in V$,
$w_{1}\in W_{1}$, $w_{2}\in W_{2}$ and $w_{3}'\in W_{3}'$.
In fact, generalizing~\eqref{skew-sym-0}, we see that
\eqref{skew-sym-k-prod-1} is equal to
\begin{gather}
\big\langle {\rm e}^{z_{k}L'(1)}w'_{3}, \bigl(Y_{W_{3}}
^{g_{3}}\bigr)(u_{1}, x_{1}+y)\cdots \bigl(Y_{W_{3}}
^{g_{3}}\bigr)(u_{k-1}, x_{k-1}+y)\nonumber\\
\qquad
\times \mathcal{Y}(w_{1}, y)w_{2}\big\rangle
\mbar_{x_{i}^{n}={\rm e}^{n\log z_{i}}, \log x_{i}=\log z_{i}, i=1, \dots, k-1,
y^{n}={\rm e}^{n\log (-z_{k})}, \log y=\log (-z_{k})}.\label{skew-sym-k-prod-2}
\end{gather}
From Definition~\ref{twisted-mod}, we see that the convergence
and analytic extension of~\eqref{skew-sym-k-prod-2} on the region~${|z_{1}|>\cdots>|z_{k}|>0}$ is equivalent to the convergence
and analytic extension of
\begin{gather}
\prod_{1\le i<j\le k-1}(z_{i}-z_{j})^{M_{ij}}
 \big\langle {\rm e}^{z_{k}L'(1)}w'_{3}, \bigl(Y_{W_{3}}
^{g_{3}}\bigr)(u_{1}, x_{1}+y)\cdots \bigl(Y_{W_{3}}
^{g_{3}}\bigr)(u_{k-1}, x_{k-1}+y)\nonumber\\
\qquad
\times \mathcal{Y}(w_{1}, y)w_{2}\big\rangle
\mbar_{x_{i}^{n}={\rm e}^{n\log z_{i}}, \log x_{i}=\log z_{i}, i=1, \dots, k-1,
y^{n}={\rm e}^{n\log (-z_{k})}, \log y=\log (-z_{k})}\label{skew-sym-k-prod-3}
\end{gather}
on the region $|z_{i}|>|z_{k}|>0$ for $i=1, \dots, k-1$, $z_{i}\ne z_{j}$
for $i\ne j$, where $M_{ij}\in \Z_{+}$ for~${i\ne j}$ satisfy~${x^{M_{ij}}Y_{V}(u_{i}, x)u_{j}\in V[[x]]}$.
From Remark~\ref{convergence-multi}, we know that Lemma~\ref{convergence} can be generalized to
the case of more than two variables for a~series of the form
\eqref{skew-sym-k-prod-3}. It is important to note that
the polynomial $\prod_{1\le i<j\le k-1}(z_{i}-z_{j})^{M_{ij}}$ in
\eqref{skew-sym-k-prod-3} is multiplied to cancel the poles at $z_{i}=z_{j}$ for~${1\le i<j\le k-1}$
so that
the series~\eqref{skew-sym-k-prod-3} is a~different iterated sum of a~convergent iterated series of the form~\eqref{gen-88}.
Using the convergence of products with more than one twisted vertex
operator
for the
twisted intertwining operator $\Y$ and this generalization of
Lemma~\ref{convergence}, we see that~\eqref{skew-sym-k-prod-3}
is absolutely convergent
on the region $|z_{i}|>|z_{k}|>0$ for $i=1, \dots, k-1$, $z_{i}\ne z_{j}$
for~${i\ne j}$ and its sum has analytic extension on the region
$M^{k}$. Thus~\eqref{skew-sym-k-prod-1}
is absolutely convergent
on the region $|z_{1}|>\cdots >|z_{k}|>0$ and its sum has
maximal analytic extension on the region~$M^{k}$. Moreover, from the discussion
above, we have
\[g(z_{1}, \dots, z_{k})=\sum_{j=0}^{J}\frac{z_{k}^{j}}{j!}f_{j}(z_{1}-z_{k}, \dots, z_{k-1}- z_{k}, -z_{k}),\]
where $g(z_{1}, \dots, z_{k})$ is the maximal analytic extension
of the sum of~\eqref{skew-sym-k-prod-1} on the region $M^{k}$, $J\in \Z_{+}$ such that
$L'(1)^{J+1} w_{3}'=0$, and
$f_{j}(z_{1}, \dots, z_{k})$ for $j=0, \dots, J$
are the maximal analytic extensions
of the sum of
\[\big\langle L'(1)^{j}w'_{3}, Y_{W_{3}}
^{g_{3}}(u_{1}, z_{1})\cdots Y_{W_{3}}
^{g_{3}}(u_{k-1}, z_{k-1}) \mathcal{Y}(w_{1}, z_{k})w_{2}\big\rangle.\]
For fixed $z_{k}\in \C^{\times}$, let $(z_1,\dots, z_{k-1})-\beta=\delta$ for
$\beta=\C^{k-1}$
and $\delta\in\{0,\infty\}^{k-1}$ be a~component-isolated singularity of the multivalued analytic function of
of $z_{1}, \dots, z_{k-1}$ obtained as a~value $\prod_{1\le i<j\le k-1}(z_{i}-z_{j})^{M_{ij}}g(z_{1}, \dots, -z_{k})$
at $\zeta=-z_{k}$ of
$\prod_{1\le i<j\le k-1}(z_{i}-z_{j})^{M_{ij}}g(z_{1}, \dots, \zeta)$. Then this singularity, which can also be
written as $(z_1-z_{k},\dots, z_{k-1}-z_{k})-(\beta+(z_{k},
\dots,z_{k}))=\delta$,
 is also a~component-isolated singularity
of the corresponding values
\[
\prod_{1\le i<j\le k-1}(z_{i}-z_{j})^{M_{ij}}f_{j}(z_{1}-z_{k}, \dots, z_{k-1}-z_{k}, -z_{k})
\]
 at $\zeta=-z_{k}$
of
\[\prod_{1\le i<j\le k-1}(z_{i}-z_{j})^{M_{ij}}f_{j}(z_{1}-z_{k}, \dots, z_{k-1}-z_{k}, \zeta)
\]
 for $j=0, \dots, J$.
Thus this singularity of
\[
\prod_{1\le i<j\le k-1}(z_{i}-z_{j})^{M_{ij}}f_{j}(z_{1}-z_{k}, \dots, z_{k-1}-z_{k}, -z_{k})
\]
corresponds to the component-isolated singularity $(z_1,\dots, z_{k-1})-(\beta+(z_{k}, \dots, z_{k}))=\delta$ of ${\prod_{1\le i<j\le k-1}(z_{i}-z_{j})^{M_{ij}}f_{j}(z_{1}, \dots, z_{k})}$.
But we know that a~component-isolated singularity
of \smash{$\prod_{1\le i<j\le k-1}(z_{i}-z_{j})^{M_{ij}}f_{j}(z_{1}, \dots, z_{k-1}, -z_{k})$} of such a~form is regular.
So this component-isolated singularity $(z_1,\dots, z_{k-1})-(\beta+(z_{k}, \dots, z_{k}))=\delta$ is regular.
Equivalently, the component-isolated singularity $(z_1,\dots, z_{k-1})-\beta=\delta$ of
\[
\prod_{1\le i<j\le k-1}(z_{i}-z_{j})^{M_{ij}}f_{j}(z_{1}-z_{k}, \dots, z_{k-1}-z_{k}, -z_{k})
\]
 for $j=0, \dots, J$
is regular and thus the component-isolated singularity of $\smash{\prod_{1\le i<j\le k-1}}(z_{i}-z_{j})^{M_{ij}}g(z_{1}, \dots, z_{k-1},
-z_{k})$
is regular. But $z_{k}$ is arbitrary, $-z_{k}$ is also arbitrary. We see that
for any fixed $z_{k}\in \C^{\times}$, $(z_1,\dots, z_{k-1})-\beta=\delta$ is a~regular singularity of any value
at $\zeta=z_{k}$ of~\smash{$\prod_{1\le i<j\le k-1}(z_{i}-z_{j})^{M_{ij}}g(z_{1}, \dots, \zeta)$}.

When $|z_{2}|>|z_{1}|>0$ and
$\arg z_{2}\ge \pi$,
\begin{gather}
\big\langle w'_{3}, \Omega_{-}(\mathcal{Y})(w_{2}, z_{2})Y_{W_{1}}
^{g_{1}}(u, z_{1})w_{1}\big\rangle\nn
\qquad=\big\langle w'_{3}, {\rm e}^{z_{2}L(-1)}\mathcal{Y}\bigl(Y_{W_{1}}
^{g_{1}}(u, z_{1})w_{1}, -z_{2}\bigr)w_{2}\big\rangle\nn
\qquad=\big\langle {\rm e}^{z_{2}L'(1)}w'_{3}, \mathcal{Y}\bigl(Y_{W_{1}}
^{g_{1}}(u, (z_{1}-z_{2})-(-z_{2}))w_{1}, -z_{2}\bigr)w_{2}\big\rangle\label{skew-sym-4}
\end{gather}
converges absolutely and if in addition, $|{\arg (z_{1}-z_{2})-\arg (-z_{2})}|<\frac{\pi }{2}$, its sum is equal to
\begin{equation}\label{f-ch-var}
f^{e}\bigl(z_1-z_{2}, -z_2; u, w_{1}, w_{2}, {\rm e}^{z_{2}L'(1)}w_{3}'\bigr).
\end{equation}
Note that by definition, \eqref{f-ch-var} is a~single-valued
analytic function on the set $\widetilde{M}_{0}^{2}$
given by cutting~$M^{2}$ along the
positive real line in the $z_{1}$- and $(z_{1}-z_{2})$-planes
and along the negative real line in the $z_2$-plane, with these
positive real lines in the $z_{1}$- and $(z_{1}-z_{2})$-planes
attached to the upper half $z_{1}$- and $(z_{1}-z_{2})$-planes
and the negative real line in the $z_2$-plane attached to the lower
half $z_{2}$-plane. Then the subset
\begin{equation}\label{subset-M-2-0}
\left\{(z_{1}, z_{2})\in M^{2}_{0} \mid |z_{2}|>|z_{1}|>0,\,
\arg z_{2}\ge \pi,\,
|{\arg (z_{1}-z_{2})-\arg (-z_{2})}|<\frac{\pi }{2}\right\}
\end{equation}
of $\widetilde{M}_{0}^{2}$ is also a~subset of $M_{0}^{2}$.
But on the subset of $M^{2}_{0}$ given by
$|z_{1}|>|z_{2}|>0$, $\arg z_{2}\ge \pi$, by definition,
\eqref{f-ch-var} is equal to $g^{e}_{-}(z_{1}, z_{2}; u, w_{2}, w_{1}, w_{3}')$.
Since $g^{e}_{-}(z_{1}, z_{2}; u, w_{2}, w_{1}, w_{3}')$ on
$M^{2}_{0}$ is obtained by analytic extension,
we see that~\eqref{f-ch-var} is equal to
$g^{e}_{-}(z_{1}, z_{2}; u, w_{2}, w_{1}, w_{3}')$ also
on the subset~\eqref{subset-M-2-0}. Thus the left-hand side
of~\eqref{skew-sym-4} is absolutely convergent to
$g^{e}_{-}(z_{1}, z_{2}; u, w_{2}, w_{1}, w_{3}')$ when
$|z_{2}|>|z_{1}|>0$,
 $|{\arg (z_{1}-z_{2})-\arg (-z_{2})}|<\frac{\pi }{2}$ and
 $\arg z_{2}\ge \pi$. Since
 $g^{e}_{-}(z_{1}, z_{2}; u, w_{2}, w_{1}, w_{3}')$
 and the sum of~\eqref{skew-sym-4} are both analytic extensions of
 their restrictions on the subset given by~${|z_{2}|>|z_{1}|>0}$,
 $|{\arg (z_{1}-z_{2})-\arg (-z_{2})}|<\frac{\pi }{2}$ and
 $\arg z_{2}\ge \pi$, we see that the left-hand side
of~\eqref{skew-sym-4} is absolutely convergent to
$g^{e}_{-}(z_{1}, z_{2}; u, w_{2}, w_{1}, w_{3}')$ when
$|z_{2}|>|z_{1}|>0$ and $|{\arg (z_{1}-z_{2})-\arg (-z_{2})}|<\frac{\pi }{2}$.
 But when
$\arg z_{2}\ge \pi$, $\arg (-z_{2})=\arg z_{2}-\pi$.
Hence in this case, the
inequality~$|{\arg (z_{1}-z_{2})-\arg (-z_{2})}|<\frac{\pi }{2}$ becomes
$-\frac{3\pi}{2}<\arg (z_{1}-z_{2})-\arg z_{2}<-\frac{\pi}{2}$.
Also both the left-hand side of
\eqref{skew-sym-4} and $g_{-}^{e}(z_{1}, z_{2}; u, w_{2}, w_{1}, w_{3}')$
 are single
valued analytic functions in $z_{1}$ and $z_{2}$ with cuts at $z_{1}\in \R_{+}$ and $z_{2}\in \R_{+}$. Thus
when $|z_{2}|>|z_{1}|>0$ and
$-\frac{3\pi}{2}<\arg (z_{1}-z_{2})-\arg z_{2}<-\frac{\pi}{2}$, the
left-hand side of~\eqref{skew-sym-4} is absolutely convergent to
$g_{-}^{e}(z_{1}, z_{2}; u, w_{2}, w_{1}, w_{3}')$.

Next we discuss the iterate of $\Omega_{-}(\Y)$ and the twisted vertex operator
map $\phi_{g_{1}}\bigl(Y^{g_{2}}_{W_{2}}\bigr)$. When $|z_{2}|>|z_{1}-z_{2}|>0$ and $\arg z_{2}\ge \pi$,
\begin{gather}
\lefteqn{\bigl\langle w'_{3}, \Omega_{-}(\mathcal{Y})\bigl(\phi_{g_{1}}\bigl(Y^{g_{2}}_{W_{2}}\bigr)(u, z_{1}-z_{2})w_{2}, z_{2}\bigr)w_{1}\bigr\rangle}\nn
\qquad=\big\langle w'_{3}, \Omega_{-}(\mathcal{Y})\bigl(\bigl(Y_{W_{2}}
^{g_{2}}\bigr)\bigl(g_{1}^{-1}u, z_{1}-z_{2}\bigr)w_{2}, z_{2}\bigr)w_{1}\big\rangle\nn
\qquad=\big\langle w'_{3}, {\rm e}^{z_{2}L(-1)}\mathcal{Y}(w_{1}, -z_{2})
\bigl(Y_{W_{2}}
^{g_{2}}\bigr)\bigl(g_{1}^{-1}u, z_{1}-z_{2}\bigr)w_{2}\big\rangle\nn
\qquad=\big\langle {\rm e}^{z_{2}L'(1)}w'_{3}, \mathcal{Y}(w_{1}, -z_{2})
\bigl(Y_{W_{2}}^{g_{2}}\bigr)\bigl(g_{1}^{-1}u, z_{1}-z_{2}\bigr)w_{2}\big\rangle\label{skew-sym-7}
\end{gather}
converges absolutely and if in addition,
$-\frac{3\pi}{2}<\arg z_{1}-\arg (-z_{2})<-\frac{\pi}{2}$, its sum
is equal to $f^{e}(z_1-z_{2}, -z_2;\smash{ g_{1}^{-1}u, w_{1}, w_{2},
{\rm e}^{z_{2}L'(1)}w_{3}}')$.
Note that the proofs of \cite[Lemmas 4.5 and 4.6]{H-twisted-int} do not use the explicit form of the
 multivalued analytic functions in the duality property of the
 twisted intertwining operators introduced in~\cite{H-twisted-int}.
Then the same proof of \cite[Lemma 4.5]{H-twisted-int} shows that
the sum of the right-hand side of~\eqref{skew-sym-7} is equal to
\smash{$f^{b_{12}^{-1}}\bigl(z_1-z_{2}, -z_2; g_{1}^{-1}u, w_{1},
w_{2}, {\rm e}^{z_{2}L'(1)}w_{3}'\bigr)$} when $|z_{2}|>|z_{1}-z_{2}|>0$,
$\arg z_{2}\ge \pi$ and
$\frac{\pi}{2}<\arg z_{1}-\arg (-z_{2})<\frac{3\pi}{2}$.
The same proof of (4.4) in~\cite[Lemma 4.6]{H-twisted-int} shows that
\[f^{b_{12}^{-1}}\bigl(z_1-z_{2}, -z_2; g_{1}^{-1}u, w_{1},
w_{2}, {\rm e}^{z_{2}L'(1)}w_{3}'\bigr)
=f^{e}\bigl(z_1-z_{2}, -z_2; u, w_{1},
w_{2}, {\rm e}^{z_{2}L'(1)}w_{3}'\bigr).\]
Since when $\arg z_{2}\ge \pi$, $\arg (-z_{2})=\arg z_{2}-\pi$
and $\frac{\pi}{2}<\arg z_{1}-\arg (-z_{2})<\frac{3\pi}{2}$ becomes
$|{\arg z_{1}-\arg z_{2}}|<\frac{\pi}{2}$,
we see that the sum of the right-hand side of~\eqref{skew-sym-7}
is equal to~\eqref{f-ch-var} when~${|z_{2}|>|z_{1}-z_{2}|>0}$,
$\arg z_{2}\ge \pi$ and
$|{\arg z_{1}-\arg z_{2}}|<\frac{\pi}{2}$.

We now use the same argument as above to finish the proof in this case.
The subset
\begin{equation}\label{another-subset-M-2-0}
\left\{(z_{1}, z_{2})\in M^{2}_{0} \mid |z_{2}|>|z_{1}-z_{2}|>0,\,
\arg z_{2}\ge \pi,\,
|{\arg (z_{1})-\arg z_{2}}|<\frac{\pi }{2}\right\}
\end{equation}
of $\widetilde{M}_{0}^{2}$ is also a~subset of $M_{0}^{2}$.
By definition, on the subset of $M^{2}_{0}$ given by
$|z_{1}|>|z_{2}|>0$, ${\arg z_{2}\ge \pi}$,
\eqref{f-ch-var} is equal to $g^{e}_{-}(z_{1}, z_{2}; u, w_{2}, w_{1}, w_{3}')$.
Since $g^{e}_{-}(z_{1}, z_{2}; u, w_{2}, w_{1}, w_{3}')$ on
$M^{2}_{0}$ is obtained by analytic extension,
we see that~\eqref{f-ch-var} is equal to
$g^{e}_{-}(z_{1}, z_{2}; u, w_{2}, w_{1}, w_{3}')$ also
on the subset~\eqref{another-subset-M-2-0}. Thus the left-hand side
of~\eqref{skew-sym-7} is absolutely convergent to
$g^{e}_{-}(z_{1}, z_{2}; u, w_{2}, w_{1}, w_{3}')$ when
$|z_{2}|>|z_{1}-z_{2}|>0$,
 $|{\arg z_{1}-\arg z_{2}}|<\frac{\pi }{2}$ and
 $\arg z_{2}\ge \pi$. Since
 $g^{e}_{-}(z_{1}, z_{2}; u, w_{2}, w_{1}, w_{3}')$
 and the sum of~\eqref{skew-sym-7} are both analytic extensions of
 their restrictions on the subset given by $|z_{2}|>|z_{1}-z_{2}|>0$,
 $|{\arg z_{1}-\arg z_{2}}|<\frac{\pi }{2}$ and
 $\arg z_{2}\ge \pi$, we see that the left-hand side
of~\eqref{skew-sym-7} is absolutely convergent to
$g^{e}_{-}(z_{1}, z_{2}; u, w_{2}, w_{1}, w_{3}')$ when
$|z_{2}|>|z_{1}-z_{2}|>0$ and~$|{\arg z_{1}-\arg z_{2}}|<\frac{\pi }{2}$.

We still need to prove the two other cases for $\Omega_{+}(\Y)$.
When $|z_{2}|>|z_{1}|>0$ and
$\arg z_{2}< \pi$,
\begin{gather}
\big\langle w'_{3}, \Omega_{+}(\mathcal{Y})(w_{2}, z_{2})\phi_{g_{2}^{-1}}\bigl(Y_{W_{1}}
^{g_{1}}\bigr)(u, z_{1})w_{1}\big\rangle\nn
\qquad=\big\langle w'_{3}, {\rm e}^{z_{2}L(-1)}\mathcal{Y}\bigl(\phi_{g_{2}^{-1}}\bigl(Y_{W_{1}}
^{g_{1}}\bigr)(u, z_{1})w_{1}, -z_{2}\bigr)w_{2}\big\rangle\nn
\qquad=\big\langle {\rm e}^{z_{2}L'(1)}w'_{3}, \mathcal{Y}\bigl(\bigl(Y_{W_{1}}
^{g_{1}}\bigr)(g_{2}u, z_{1})w_{1}, -z_{2}\bigr)w_{2}\big\rangle\label{skew-sym-4-1}
\end{gather}
converges absolutely and if in addition, $|{\arg (z_{1}-z_{2})-\arg (-z_{2})}|<\frac{\pi}{2}$,
its sum is equal to
\begin{equation}\label{f-g-2-ch-var}
f^{e}\bigl(z_1-z_{2}, -z_2; g_{2}u, w_{1}, w_{2},
{\rm e}^{z_{2}L'(1)}w_{3}'\bigr).
\end{equation}
The subset
\begin{equation}\label{3rd-subset-M-2-0}
\left\{(z_{1}, z_{2})\in M^{2}_{0} \mid |z_{2}|>|z_{1}|>0,\,
\arg z_{2}< \pi,\,
|{\arg (z_{1}-z_{2})-\arg (-z_{2})}|<\frac{\pi }{2}\right\}
\end{equation}
of $\widetilde{M}_{0}^{2}$ is also a~subset of $M_{0}^{2}$.
By definition, on the subset of $M^{2}_{0}$ given by
$|z_{1}|>|z_{2}|>0$, ${\arg z_{2}< \pi}$,
\eqref{f-g-2-ch-var}
is equal to $g^{e}_{+}(z_{1}, z_{2}; g_{2}u, w_{2}, w_{1}, w_{3}')$.
Since $g^{e}_{+}(z_{1}, z_{2}; g_{2}u, w_{2}, w_{1}, w_{3}')$ on
$M^{2}_{0}$ is obtained by analytic extension,
\eqref{f-g-2-ch-var} is equal to
$g^{e}_{+}(z_{1}, z_{2}; g_{2}u, w_{2}, w_{1}, w_{3}')$ also
on the subset~\eqref{3rd-subset-M-2-0}. Thus the left-hand side
of~\eqref{skew-sym-4-1} is absolutely convergent to
$g^{e}_{+}(z_{1}, z_{2}; g_{2}u, w_{2}, w_{1}, w_{3}')$ when~${|z_{2}|>|z_{1}|>0}$,
 $|{\arg (z_{1}-z_{2})-\arg (-z_{2})}|<\frac{\pi }{2}$ and
 $\arg z_{2}<\pi$. Since
 $g^{e}_{+}(z_{1}, z_{2}; g_{2}u, w_{2},\allowbreak w_{1}, w_{3}')$
 and the sum of~\eqref{skew-sym-4-1} are both analytic extensions of
 their restrictions on the subset given~by ${|z_{2}|>|z_{1}|>0}$,
 $|{\arg (z_{1}-z_{2})-\arg (-z_{2})}|<\frac{\pi }{2}$ and
 $\arg z_{2}<\pi$, the left-hand side
of~\eqref{skew-sym-4-1} is absolutely convergent to
$g^{e}_{+}(z_{1}, z_{2}; g_{2}u, w_{2}, w_{1}, w_{3}')$ when
$|z_{2}|>|z_{1}|>0$,
 $|{\arg (z_{1}-z_{2})}-\arg (-z_{2}) |<\frac{\pi }{2} $.

The same proof as that of (4.5) in~\cite[Lemma 4.6]{H-twisted-int} gives
$
g^{e}_{+}(z_1, z_2; g_{2}u, w_{1}, w_{2}, w_{3}')
=\smash{g^{b_{13}^{-1}}_{+}}(z_1,\allowbreak z_2; u, w_{1}, w_{2}, w_{3}')$,
since $\arg z_{1}< \pi$, $\arg (-z_{2})=\arg z_{2}+\pi$ and the
inequality $|{\arg (z_{1}-z_{2})}-\arg (-z_{2})|<\frac{\pi}{2}$
becomes $\frac{\pi}{2}<\arg (z_{1}-z_{2})-\arg z_{2}<\frac{3\pi}{2}$.
So the sum of the left-hand side~of \eqref{skew-sym-4-1} is absolutely convergent to
\smash{$g^{b_{13}^{-1}}_{+}(z_1, z_2; u, w_{1}, w_{2}, w_{3}')$} when
$|z_{2}|>|z_{1}|>0$ and ${\frac{\pi}{2}}<\arg (z_{1}-z_{2})-\arg z_{2}<\frac{3\pi}{2}$.
Then the same proof as that of \cite[Lemma~4.5]{H-twisted-int} shows that this is equivalent to
the sum of left-hand side of~\eqref{skew-sym-4-1} being equal to
$g_{+}^{e}(z_{1}, z_{2}; u, w_{2}, w_{1}, w_{3}')$
when $|z_{2}|>|z_{1}|>0$ and
$-\frac{3\pi}{2}<\arg (z_{1}-z_{2})-\arg z_{2}<-\frac{\pi}{2}$.

Finally, we discuss the iterate of $\Omega_{+}(\Y)$ and the twisted
vertex operator
map $Y^{g_{2}}_{W_{2}}$. When $|z_{2}|>|z_{1}-z_{2}|>0$ and $\arg z_{2}< \pi$,
\begin{gather}
\lefteqn{\big\langle w'_{3}, \Omega_{+}(\mathcal{Y})\bigl(Y_{W_{2}}
^{g_{2}}(u, z_{1}-z_{2})w_{2}, z_{2}\bigr)w_{1}\big\rangle}\nn
\qquad=\big\langle w'_{3}, {\rm e}^{z_{2}L(-1)}\mathcal{Y}(w_{1}, -z_{2})
Y_{W_{2}}
^{g_{2}}(u, z_{1}-z_{2})w_{2}\big\rangle\nn
\qquad=\big\langle {\rm e}^{z_{2}L'(1)}w'_{3}, \mathcal{Y}(w_{1}, -z_{2})
Y_{W_{2}}
^{g_{2}}(u, z_{1}-z_{2})w_{2}\big\rangle\label{skew-sym-7-1}
\end{gather}
converges absolutely and if in addition,
$-\frac{3\pi}{2}<\arg z_{1}-\arg (-z_{2})<-\frac{\pi}{2}$, its sum is equal to~\eqref{f-ch-var}.
The subset
\begin{equation}\label{4th-subset-M-2-0}
\left\{(z_{1}, z_{2})\in M^{2}_{0} \mid |z_{2}|>|z_{1}-z_{2}|>0,\,
\arg z_{2}< \pi,\,
-\frac{3\pi}{2}<\arg z_{1}-\arg (-z_{2})<-\frac{\pi}{2}\right\}
\end{equation}
of $\widetilde{M}_{0}^{2}$ is also a~subset of $M_{0}^{2}$.
As we have discussed above, on the subset of $M^{2}_{0}$ given by
$|z_{1}|>|z_{2}|>0$, $\arg z_{2}< \pi$,
\eqref{f-ch-var}
is equal to $g^{e}_{+}(z_{1}, z_{2}; u, w_{2}, w_{1}, w_{3}')$.
Since $g^{e}_{+}(z_{1}, z_{2}; u, w_{2}, w_{1},\allowbreak w_{3}')$ on~$M^{2}_{0}$ is obtained by analytic extension,
\eqref{f-ch-var} is equal to
$g^{e}_{+}(z_{1}, z_{2}; u, w_{2}, w_{1}, w_{3}')$ also
on the subset~\eqref{4th-subset-M-2-0}. Thus the left-hand side
of~\eqref{skew-sym-7-1} is absolutely convergent to
$g^{e}_{+}(z_{1}, z_{2}; u, w_{2},\allowbreak w_{1}, w_{3}')$ when
$|z_{2}|>|z_{1}-z_{2}|>0$,
$-\frac{3\pi}{2}<\arg z_{1}-\arg (-z_{2})<-\frac{\pi}{2}$ and
 $\arg z_{2}<\pi$. When $\arg z_{2}<\pi$, $\arg (-z_{2})=\pi+\arg z_{2}$.
So in this case, $-\frac{3\pi}{2}<\arg z_{1}-\arg (-z_{2})<-\frac{\pi}{2}$
is the same as $|{\arg z_{1}-\arg z_{2}}|<\frac{\pi}{2}$.
Since
 $g^{e}_{+}(z_{1}, z_{2}; u, w_{2}, w_{1}, w_{3}')$
 and the sum of~\eqref{skew-sym-7-1} are both analytic extensions of
 their restrictions on the subset given by $|z_{2}|>|z_{1}-z_{2}|>0$,
$|{\arg z_{1}-\arg z_{2}}|<\frac{\pi}{2}$ and
 $\arg z_{2}<\pi$, the left-hand side
of~\eqref{skew-sym-7-1} is absolutely convergent to
$g^{e}_{+}(z_{1}, z_{2}; u, w_{2}, w_{1}, w_{3}')$ when
$|z_{2}|>|z_{1}-z_{2}|>0$ and
$|{\arg z_{1}-\arg z_{2}}|<\frac{\pi}{2}$.
\end{proof}

Let \smash{$\mathcal{V}_{W_{1}W_{2}}^{W_{3}}$} be the space of twisted intertwining operators
of type $\binom{W_{3}}{W_{1}W_{2}}$. Then we have the following.

\begin{cor}
The maps
\begin{alignat*}{3}
& \Omega_{+}\colon\ \mathcal{V}_{W_{1}W_{2}}^{W_{3}}\to \mathcal{V}_{W_{2}\phi_{g_{2}^{-1}}(W_{1})}^{W_{3}},\qquad&&
\Omega_{-}\colon \ \mathcal{V}_{W_{1}W_{2}}^{W_{3}}\to \mathcal{V}_{\phi_{g_{1}}(W_{2})W_{1}}^{W_{3}},&\\
&\Omega_{+}^{g_{2}^{-1}}\colon\ \mathcal{V}_{W_{1}W_{2}}^{W_{3}}\to \mathcal{V}_{W_{2}\varphi_{g_{2}^{-1}}(W_{1})}^{W_{3}},\qquad&& \Omega_{-}^{g_{1}}\colon \ \mathcal{V}_{W_{1}W_{2}}^{W_{3}}\to \mathcal{V}_{\varphi_{g_{1}}(W_{2})W_{1}}^{W_{3}}&
\end{alignat*}
are linear isomorphisms. In particular,
\begin{gather*}
\mathcal{V}_{W_{1}W_{2}}^{W_{3}}, \quad \mathcal{V}_{\phi_{g_{1}}(W_{2})W_{1}}^{W_{3}},\quad
\mathcal{V}_{W_{2}\phi_{g_{2}^{-1}}(W_{1})}^{W_{3}}, \quad \mathcal{V}_{\varphi_{g_{1}}(W_{2})W_{1}}^{W_{3}}\quad
\text{and}\quad \mathcal{V}_{W_{2}\varphi_{g_{2}^{-1}}(W_{1})}^{W_{3}}
\end{gather*}
are linearly isomorphic.
\end{cor}
\begin{proof}
It is clear that $\Omega_{+}$ and $\Omega_{-}$ are inverse of each other. Since \smash{$\varphi_{g_{1}}(W_{2})$}
and \smash{$\varphi_{g_{2}^{-1}}(W_{1})$} are equivalent to \smash{$\phi_{g_{1}}(W_{2})$}
and \smash{$\phi_{g_{2}^{-1}}(W_{1})$}, respectively, \smash{$\Omega_{+}^{g_{2}^{-1}}$} and \smash{$\Omega_{-}^{g_{1}}$}
are also linear isomorphisms.
\end{proof}

The linear isomorphisms $\Omega_{+}$, $\Omega_{-}$, \smash{$\Omega_{+}^{g_{2}^{-1}}$} and \smash{$\Omega_{-}^{g_{1}}$} are
called {\it skew-symmetry isomorphisms}.

Let $g_{1}$, $g_{2}$ be automorphisms of $V$, $W_{1}$, $W_{2}$ and $W_{3}$ $g_{1}$-, $g_{2}$-
and $g_{1}g_{2}$-twisted
$V$-modules without $g_{1}$-, $g_{2}$- and $g_{1}g_{2}$-actions
and $\mathcal{Y}$ a~twisted intertwining operator
of type $\binom{W_{3}}{W_{1}W_{2}}$.
We define linear maps
\smash{$
A_{\pm}(\Y)\colon W_{1}\otimes W_{3}'\to W_{2}'\{x\}[\log x]$}, $
w_{1}\otimes w_{3}'\mapsto A_{\pm}(\Y)(w_{1}, x)w_{3}'
$
by
\begin{align}\label{A1}
\langle A_{\pm}(\Y)(w_{1}, x)w_{3}', w_{2}\rangle
=\big\langle w_{3}', \Y\bigl({\rm e}^{xL(1)}{\rm e}^{\pm \pi \i L(0)}\bigl(x^{-L(0)}\bigr)^{2}w_{1}, x^{-1}\bigr)w_{2}\big\rangle
\end{align}
for $w_{1}\in W_{1}$ and $w_{2}\in W_{2}$ and $w_{3}'\in W_{3}'$.
Similarly to $\Omega_{p}$ for~${p\in \Z}$, we can also define
$A_{p}$ for~${p\in \Z}$ by replacing $\pm$ in~\eqref{A1} by
$(2p+1)$. But we will also not discuss $A_{p}$ in this paper.

Let $L_{W_{1}}^{s}(0)$ be the semisimple part of $L_{W_{1}}(0)$.
From the definition~\eqref{A1}, for $p\in \Z$, $w_{1}\in W_{1}$, $w_{2}\in W_{2}$, $w_{3}'\in W_{3}'$
and $z\in \C^{\times}$, we have
\begin{gather}
\big\langle A_{\pm}(\Y)^{p}(w_{1}, z)w_{3}', w_{2}\big\rangle\nn
\qquad=\big\langle A_{\pm}(\Y)^{p}(w_{1}, x)w_{3}', w_{2}\big\rangle\lbar_{x^{n}={\rm e}^{nl_{p}(z)}, \log x=l_{p}(z)}\nn
\qquad=\big\langle w_{3}', \Y\bigl({\rm e}^{xL_{W_{1}}(1)}{\rm e}^{\pm \pi \i L_{W_{1}}(0)}\bigl(x^{-L_{W_{1}}(0)}\bigr)^{2}w_{1}, x^{-1}\bigr)w_{2}\big\rangle
\lbar_{x^{n}={\rm e}^{nl_{p}(z)}, \log x=l_{p}(z)}\nn
\qquad=\big\langle w_{3}', \Y\bigl({\rm e}^{xL_{W_{1}}(1)}{\rm e}^{\pm \pi \i L_{W_{1}}(0)}\bigl(x^{-L_{W_{1}}^{s}(0)}\bigr)^{2}\cdot\nn
\phantom{  \qquad =}{}\times \bigl({\rm e}^{-(L_{W_{1}}(0)-L_{W_{1}}^{s}(0))\log x}\bigr)^{2}w_{1}, x^{-1}\bigr)w_{2}\big\rangle
\lbar_{x^{n}={\rm e}^{nl_{p}(z)}, \log x=l_{p}(z)}\nn
\qquad=\big\langle w_{3}', \Y\bigl({\rm e}^{zL_{W_{1}}(1)}{\rm e}^{\pm \pi \i L_{W_{1}}(0)}\bigl({\rm e}^{-l_{p}(z) L_{W_{1}}^{s}(0)}\bigr)^{2}\cdot\nn
\phantom{  \qquad =}{} \times\bigl({\rm e}^{-(L_{W_{1}}(0)-L_{W_{1}}^{s}(0))l_{p}(z)}\bigr)^{2}w_{1}, y)w_{2}\big\rangle
\lbar_{y^{n}={\rm e}^{-nl_{p}(z)}, \log y=-l_{p}(z)}\nn
\qquad=\big\langle w_{3}', \Y\bigl({\rm e}^{zL_{W_{1}}(1)}{\rm e}^{\pm \pi \i L_{W_{1}}(0)}{\rm e}^{-2l_{p}(z) L_{W_{1}}(0)}w_{1}, y\bigr)w_{2}\big\rangle
\lbar_{y^{n}={\rm e}^{-nl_{p}(z)}, \log y=-l_{p}(z)}.\label{A2}
\end{gather}
When $\arg z=0$, $\arg z^{-1}=\arg z=0$ and $-l_{p}(z)=l_{-p}(z^{-1})$.
When $\arg z\ne 0$, $\arg z^{-1}=-\arg z+2\pi$ and $-l_{p}(z)=l_{-p-1}(z^{-1})$.
Hence when $\arg z=0$, the right-hand side of~\eqref{A2} is equal to
\begin{gather}
\big\langle w_{3}', \Y\bigl({\rm e}^{zL_{W_{1}}(1)}{\rm e}^{\pm \pi \i L_{W_{1}}(0)}{\rm e}^{2l_{-p}(z^{-1}) L_{W_{1}}(0)}w_{1}, y\bigr)w_{2}\big\rangle
\lbar_{y^{n}={\rm e}^{nl_{-p}(z^{-1})}, \log y=l_{-p}(z^{-1})}\nn
\qquad=\big\langle w_{3}', \Y^{-p}\bigl({\rm e}^{zL_{W_{1}}(1)}{\rm e}^{\pm \pi \i L_{W_{1}}(0)}{\rm e}^{2l_{-p}(z^{-1})L_{W_{1}}(0)}
w_{1}, z^{-1}\bigr)w_{2}\big\rangle\label{A3}
\end{gather}
and when $\arg z\ne 0$, it is equal to
\begin{gather}
\big\langle w_{3}', \Y\bigl({\rm e}^{zL_{W_{1}}(1)}{\rm e}^{\pm \pi \i L_{W_{1}}(0)}{\rm e}^{2l_{-p-1}(z^{-1}) L_{W_{1}}(0)}
w_{1}, y\bigr)w_{2}\big\rangle
\lbar_{y^{n}={\rm e}^{nl_{-p-1}(z^{-1})}, \log y=l_{-p-1}(z^{-1})}\nn
\qquad=\big\langle w_{3}', \Y^{-p-1}\bigl({\rm e}^{zL_{W_{1}}(1)}{\rm e}^{\pm \pi \i L_{W_{1}}(0)}{\rm e}^{2l_{-p-1}(z^{-1}) L_{W_{1}}(0)}
w_{1}, z^{-1}\bigr)w_{2}\big\rangle.\label{A4}
\end{gather}
{\samepage From~\eqref{A2}--\eqref{A4}, for $w_{1}\in W_{1}$, $w_{2}\in W_{2}$, $w_{3}'\in W_{3}'$
and $z\in \C^{\times}$,
we have
\begin{equation}\label{A5}
\big\langle A_{\pm}(\Y)^{p}(w_{1}, z)w_{3}', w_{2}\big\rangle
=\big\langle w_{3}', \Y^{-p}\bigl({\rm e}^{zL_{W_{1}}(1)}{\rm e}^{\pm \pi \i L_{W_{1}}(0)}{\rm e}^{2l_{-p}(z^{-1}) L_{W_{1}}(0)}
w_{1}, z^{-1}\bigr)w_{2}\big\rangle
\end{equation}
when $\arg z=0$ and
\begin{gather}
\big\langle A_{\pm}(\Y)^{p}(w_{1}, z)w_{3}', w_{2}\big\rangle\nn
\qquad
=\big\langle w_{3}', \Y^{-p-1}\bigl({\rm e}^{zL_{W_{1}}(1)}{\rm e}^{\pm \pi \i L_{W_{1}}(0)}{\rm e}^{2l_{-p-1}(z^{-1}) L_{W_{1}}(0)}
w_{1}, z^{-1}\bigr)w_{2}\big\rangle\label{A6}
\end{gather}
when $\arg z\ne 0$.}

Now assume that $W_{1}$, $W_{2}$, and $W_{3}$ are $g_{1}$-, $g_{2}$- and $g_{1}g_{2}$-twisted
$V$-modules with
$g_{1}$-, $g_{2}$- and $g_{1}g_{2}$-actions, respectively.
In the case that $W_{2}'$ has an action of $g_{1}$ (which induces a~$g_{1}^{-1}$-action on $W_{2}$),
we have the $g_{1}g_{2}^{-1}g_{1}^{-1}$-twisted $V$-module $\varphi_{g_{1}}(W_{2}')$.
In this case, we define
\begin{gather*}
A_{+}^{g_{1}}(\Y)\colon\ W_{1}\otimes W_{3}'\to W_{2}'\{x\}[\log x],\qquad
w_{1}\otimes w_{3}'\mapsto A_{+}^{g_{1}}(\Y)(w_{1}, x)w_{3}'
\end{gather*}
by
\[
\big\langle A_{+}^{g_{1}}(\Y)(w_{1}, x)w_{3}', w_{2}\big\rangle
=\big\langle w_{3}', \Y\bigl({\rm e}^{xL(1)}{\rm e}^{\pi \i L(0)}\bigl(x^{-L(0)}\bigr)^{2}w_{1}, x^{-1}\bigr)g_{1}^{-1}w_{2}\big\rangle
\]
for $w_{1}\in W_{1}$, $w_{2}\in W_{2}$ and $w_{3}'\in W_{3}'$.
In the case that $W_{3}'$ has an action of $g_{1}^{-1}$ (which induces a~$g_{1}$-action on $W_{3}$),
we have the $g_{1}^{-1}g_{2}^{-1}$-twisted $V$-module $\varphi_{g_{1}^{-1}}(W_{3}')$.
In this case, we define
\begin{gather*}
A_{-}^{g_{1}^{-1}}(\Y)\colon\ W_{1}\otimes W_{3}'\to W_{2}'\{x\}[\log x],\qquad
w_{1}\otimes w_{3}'\mapsto A_{-}^{g_{1}^{-1}}(\Y)(w_{1}, x)w_{3}'
\end{gather*}
by
\[
\big\langle A_{-}^{g_{1}^{-1}}(\Y)(w_{1}, x)w_{3}', w_{2}\big\rangle
=\big\langle g_{1}w_{3}', \Y\big({\rm e}^{xL(1)}{\rm e}^{-\pi \i L(0)}\bigl(x^{-L(0)}\bigr)^{2}w_{1}, x^{-1}\big)w_{2}\big\rangle
\]
for $w_{1}\in W_{1}$, $w_{2}\in W_{2}$ and $w_{3}'\in W_{3}'$.

For $w_{1}\in W_{1}$, $w_{2}\in W_{2}$, $w_{3}'\in W_{3}'$
and $z\in \C^{\times}$,
we have from~\eqref{A5} and~\eqref{A6}
\begin{gather}
\big\langle A_{+}^{g_{1}}(\Y)^{p}(w_{1}, z)w_{3}', w_{2}\big\rangle\nn
\qquad
=\big\langle w_{3}', \Y^{-p}\big({\rm e}^{zL_{W_{1}}(1)}{\rm e}^{ \pi \i L_{W_{1}}(0)}{\rm e}^{2l_{-p}(z^{-1}) L_{W_{1}}(0)}
w_{1}, z^{-1}\big)g_{1}^{-1}w_{2}\big\rangle\label{A9}
\end{gather}
and
\begin{gather}
\big\langle A_{-}^{g_{1}^{-1}}(\Y)^{p}(w_{1}, z)w_{3}', w_{2}\big\rangle\nonumber\\
\qquad =\big\langle g_{1}w_{3}', \Y^{-p}({\rm e}^{zL_{W_{1}}(1)}{\rm e}^{- \pi \i L_{W_{1}}(0)}{\rm e}^{2l_{-p}(z^{-1}) L_{W_{1}}(0)}
w_{1}, z^{-1})w_{2}\big\rangle\label{A10}
\end{gather}
when $\arg z=0$ and
\begin{gather}
\big\langle A_{+}^{g_{1}}(\Y)^{p}(w_{1}, z)w_{3}', w_{2}\big\rangle\nonumber\\
\qquad =\big\langle w_{3}', \Y^{-p-1}\bigl({\rm e}^{zL_{W_{1}}(1)}{\rm e}^{\pi \i L_{W_{1}}(0)}{\rm e}^{2l_{-p-1}(z^{-1}) L_{W_{1}}(0)}
w_{1}, z^{-1}\bigr)g_{1}^{-1}w_{2}\big\rangle\label{A11}
\end{gather}
and
\begin{gather}
\big\langle A_{-}^{g_{1}^{-1}}(\Y)^{p}(w_{1}, z)w_{3}', w_{2}\big\rangle\nonumber\\
\qquad
=\big\langle g_{1}w_{3}', \Y^{-p-1}\bigl({\rm e}^{zL_{W_{1}}(1)}{\rm e}^{- \pi \i L_{W_{1}}(0)}{\rm e}^{2l_{-p-1}(z^{-1}) L_{W_{1}}(0)}
w_{1}, z^{-1}\bigr)w_{2}\big\rangle\label{A12}
\end{gather}
when $\arg z\ne 0$.

Let $\bigl(W, Y_{W}^{g}\bigr)$ be a~$g$-twisted $V$-module. When $W_{1}=V$, $W_{2}=W_{3}=W$ and
$\Y=Y_{W}^{g}$, by definition, \smash{$A_{+}\bigl(Y_{W}^{g}\bigr)=A_{-}\bigl(Y_{W}^{g}\bigr)
=A_{+}^{1_{V}}\bigl(Y_{W}^{g}\bigr)=A_{-}^{1_{V}}\bigl(Y_{W}^{g}\bigr)=\bigl(Y_{W}^{g}\bigr)'$} (see Section \ref{sec4}).

\begin{thm}\label{A}
The linear maps $A_{+}(\Y)$, $A_{-}(\Y)$, $A_{+}^{g_{1}}(\Y)$ and \smash{$A_{-}^{g_{1}^{-1}}(\Y)$}
are twisted intertwining operators
of types
\[
\binom{\phi_{g_{1}}(W_{2}')}{W_{1}W_{3}'} ,\quad
 \binom{W_{2}'}{W_{1}\phi_{g_{1}^{-1}}(W_{3}')} , \quad \binom{\varphi_{g_{1}}(W_{2}')}{W_{1}W_{3}'} \quad \text{and} \quad
 \binom{W_{2}'}{W_{1}\varphi_{g_{1}^{-1}}(W_{3}')},
 \] respectively.
\end{thm}
\begin{proof}
As in the proof of Theorem~\ref{skew-sym}, we need only prove the results for
$A_{+}(\Y)$ and $A_{-}(\Y)$. The proofs for $A_{+}^{g_{1}}(\Y)$ and \smash{$A_{-}^{g_{1}^{-1}}(\Y)$}
are reduced to the proofs for $A_{+}(\Y)$ and $A_{-}(\Y)$ using the definitions of $\varphi_{g_{1}}$
and $\varphi_{g_{1}^{-1}}$ and~\eqref{A9}--\eqref{A12}.

The proof for $A_{+}(\Y)$ and $A_{-}(\Y)$ is essentially the same as the proof of~\cite[Theorem 6.1]{H-twisted-int}. But the proof here is much more complicated because the correlation functions involved are
not of the explicit form as in~\cite{H-twisted-int}.
As in~\cite{H-twisted-int},
we prove only the duality property since other properties are essentially the same as in the untwisted case.

We first give the multivalued analytic functions with preferred branches
in the duality property. We
shall denote these multivalued analytic functions
for $A_{+}(\Y)$ and $A_{-}(\Y)$
by
$h_{+}(z_{1}, z_{2}; u, w_{1}, w_{3}', w_{2})$ and
$h_{-}(z_{1}, z_{2}; u, w_{1}, w_{3}', w_{2})$, respectively. Let
$f(z_1, z_2; u, w_{1}, w_{3}', w_{2})$
be the multivalued analytic function with the preferred branch
$f^{e}(z_1, z_2; u, w_{1}, w_{2}, w_{3}')$
in the duality property for the twisted intertwining operator $\Y$.
Define
\begin{gather}
 h_{\pm}(z_{1}, z_{2}; u, w_{1}, w_{3}', w_{2})\nn
\qquad = f\bigl(z_1^{-1}, z_2^{-1}; {\rm e}^{z_{1}L_{V}(1)}\bigl(-z_{1}^{2}\bigr)^{-L_{V}(0)}u,
 {\rm e}^{z_{2}L_{W_{1}}(1)}{\rm e}^{\pm \pi \i L_{W_{1}}(0)}\bigl(z_{2}^{2}\bigr)^{-L_{W_{1}}(0)}
w_{1}, w_{2}, w_{3}'\bigr)\label{A-1}
\end{gather}
and choose the preferred branch
$h_{\pm}^{e}(z_{1}, z_{2}; u, w_{1}, w_{3}', w_{2})$ of
$h_{\pm}(z_{1}, z_{2}; u, w_{1}, w_{3}', w_{2})$ as follows:
On the subregion $|z_{1}|>|z_{2}|>0$,
$|{\arg (z_{1}-z_{2})-\arg z_{1}}|<\frac{\pi}{2}$ and
$\arg z_{1}, \arg z_{2}\ne 0$ of $M^{2}_{0}$, let
\begin{align}
h_{+}^{e}(z_{1}, z_{2}; u, w_{1}, w_{3}', w_{2})={}&f^{^{b_{13}^{-1}b_{23}^{-1}}}
\bigl(z_1^{-1}, z_2^{-1}; {\rm e}^{z^{-1}L_{V}(1)}\bigl(-z_{1}^{2}\bigr)^{-L_{V}(0)}u,\nonumber\\
& {\rm e}^{z_{2}L_{W_{1}}(1)}{\rm e}^{ \pi \i L_{W_{1}}(0)}{\rm e}^{-2(\log z_{2}) L_{W_{1}}(0)}
w_{1}, w_{2}, w_{3}'\bigr).\label{A-1.1}
\end{align}
and
\begin{align*}
h_{-}^{e}(z_{1}, z_{2}; u, w_{1}, w_{3}', w_{2})
={}&f^{^{b_{13}^{-1}b_{23}^{-1}b_{12}^{-1}}}
\bigl(z_1^{-1}, z_2^{-1}; {\rm e}^{z^{-1}L_{V}(1)}\bigl(-z_{1}^{2}\bigr)^{-L_{V}(0)}u,\\
& {\rm e}^{z_{2}L_{W_{1}}(1)}{\rm e}^{- \pi \i L_{W_{1}}(0)}{\rm e}^{-2(\log z_{2}) L_{W_{1}}(0)}
w_{1}, w_{2}, w_{3}'\bigr).
\end{align*}
For general $(z_{1}, z_{2})\in M^{2}_{0}$, we define
$h_{\pm}^{e}(z_{1}, z_{2}; u, w_{1}, w_{3}', w_{2})$
to be the unique analytic extensions on $M^{2}_{0}$.

Let $u\in V$, $w_{1}\in W_{1}$, $w_{2}\in W_{2}$
and $w_{3}'\in W_{3}'$. We consider $z_{1}, z_{2}\in \C$ satisfying
$\big|z_{2}^{-1}\big|>\big|z_{1}^{-1}\big|>0$ (or equivalently $|z_{1}|>|z_{2}|>0$) and
 $\arg z_{1}, \arg z_{2}\ne 0$.
Since \smash{$\big|z_{2}^{-1}\big|>\big|z_{1}^{-1}\big|>0$} and~${\arg z_{1}, \arg z_{2}\ne 0}$,
from~\eqref{A6}, \smash{$\bigl(Y_{W_{2}}
^{g_{2}}\bigr)'=A_{+}\bigl(Y_{W_{2}}
^{g_{2}}\bigr)$} and the duality property for $\Y$,
we know that
\begin{gather}
\big\langle w'_{3}, \mathcal{Y}\bigl({\rm e}^{z_{2}L_{W_{1}}(1)}{\rm e}^{ \pi \i L_{W_{1}}(0)}{\rm e}^{-2(\log z_{2}) L_{W_{1}}(0)}
w_{1}, z_{2}^{-1}\bigr)\nonumber\\
\qquad\times
Y_{W_{2}}
^{g_{2}}\bigl({\rm e}^{z_{1}L_{V}(1)}\bigl(-z_{1}^{-2}\bigr)^{L_{V}(0)}g_{1}^{-1}u, z_{1}^{-1}\bigr)w_{2}\big\rangle\label{A-1.2}
\end{gather}
is absolutely convergent and if in addition, $-\frac{3\pi}{2}<\arg \bigl(z_{1}^{-1}-z_{2}^{-1}\bigr)-\arg z_{2}^{-1}<-\frac{\pi}{2}$,
its sum is equal to
\begin{gather*}
f^{e}
\bigl(z_{1}^{-1}, z_{2}^{-1}; {\rm e}^{z_{1}L_{V}(1)}\bigl(-z_{1}^{-2}\bigr)^{L_{V}(0)}g_{1}^{-1}u,
 {\rm e}^{z_{2}L_{W_{1}}(1)}{\rm e}^{ \pi \i L_{W_{1}}(0)}{\rm e}^{-2(\log z_{2})
 L_{W_{1}}(0)}
w_{1}, w_{2}, w_{3}'\bigr)\nn
\qquad=f^{e}
\bigl(z_{1}^{-1}, z_{2}^{-1};
 g_{1}^{-1}{\rm e}^{z_{1}L_{V}(1)}\bigl(-z_{1}^{-2}\bigr)^{L_{V}(0)}u,
 {\rm e}^{z_{2}L_{W_{1}}(1)}{\rm e}^{ \pi \i L_{W_{1}}(0)}{\rm e}^{-2(\log z_{2})
 L_{W_{1}}(0)}
w_{1}, w_{2}, w_{3}'\bigr).
\end{gather*}
We know that
\begin{gather}
\langle \phi_{g_{1}}\bigl(\bigl(Y_{W_{2}}^{g_{2}}\bigr)'\bigr)(u, z_{1})A_{+}(\mathcal{Y})(w_{1}, z_{2})w_{3}', w_{2}\big\rangle\nn
\qquad=\big\langle \bigl(Y_{W_{2}}
^{g_{2}}\bigr)'\bigl(g_{1}^{-1}u, z_{1}\bigr)A_{+}(\mathcal{Y})(w_{1}, z_{2})w_{3}', w_{2}\big\rangle\nn
\qquad=\big\langle w'_{3}, \mathcal{Y}^{-1}\bigl({\rm e}^{z_{2}L_{W_{1}}(1)}{\rm e}^{ \pi \i L_{W_{1}}(0)}{\rm e}^{-2(\log z_{2}) L_{W_{1}}(0)}
w_{1}, z_{2}^{-1}\bigr)\nn
\qquad \phantom{=}{}\times
\bigl(Y_{W_{2}}
^{g_{2}}\bigr)^{-1}\bigl({\rm e}^{z_{1}L_{V}(1)}\bigl(-z_{1}^{-2}\bigr)^{L_{V}(0)}g_{1}^{-1}u, z_{1}^{-1}\bigr)w_{2}\big\rangle\label{A-0}
\end{gather}
can be obtained using the multivalued analytic function~\eqref{A-1}
on the region $\big|z_{2}^{-1}\big|>\big|z_{1}^{-1}\big|>0$
starting from the value given by~\eqref{A-1.2} by
letting $z_{1}^{-1}$
go around $0$ clockwise once \big(corresponding to $b_{13}^{-1}$\big),
then letting $z_{2}^{-1}$ go around
$0$ clockwise once \big(corresponding to $b_{23}^{-1}b_{12}^{-1}$\big).
Then~\eqref{A-0} also
converges absolutely on the region $\big|z_{2}^{-1}\big|>\big|z_{1}^{-1}\big|>0$
and if in addition,
$-\frac{3\pi}{2}<\arg \bigl(z_{1}^{-1}-z_{2}^{-1}\bigr)-\arg z_{2}^{-1}
<-\frac{\pi}{2}$,
its sum is equal to
\begin{gather}
 f^{b_{13}^{-1}b_{23}^{-1}b_{12}^{-1}}
\bigl(z_{1}^{-1}, z_{2}^{-1};
 g_{1}^{-1}{\rm e}^{z_{1}L_{V}(1)}\bigl(-z_{1}^{-2}\bigr)^{L_{V}(0)}u,
 {\rm e}^{z_{2}L_{W_{1}}(1)}{\rm e}^{ \pi \i L_{W_{1}}(0)}{\rm e}^{-2(\log z_{2})
 L_{W_{1}}(0)}
w_{1}, w_{2}, w_{3}'\bigr)\nn
\qquad =f^{b_{12}^{-1}b_{13}^{-1}b_{23}^{-1}}
\bigl(z_{1}^{-1}, z_{2}^{-1};
 g_{1}^{-1}{\rm e}^{z_{1}L_{V}(1)}\bigl(-z_{1}^{-2}\bigr)^{L_{V}(0)}u,\nn
 \phantom{\qquad =}{}
 {\rm e}^{z_{2}L_{W_{1}}(1)}{\rm e}^{ \pi \i L_{W_{1}}(0)}{\rm e}^{-2(\log z_{2})
 L_{W_{1}}(0)}
w_{1}, w_{2}, w_{3}'\bigr).\label{A-2}
\end{gather}
Using the equivariance property for $W_{1}$ and the convergence of
\eqref{int-iter} to $f^{e}(z_{1}, z_{2}; u,
w_{1}, w_{2}, w_{3}')$ on the region
$|z_{2}|>|z_{1}-z_{2}|>0$, $|{\arg z_{1}-\arg z_{2}}|<\frac{\pi}{2}$,
we have
\[f^{e}\bigl(z_{1}, z_{2}; g_{1}^{-1}u,
w_{1}, w_{2}, w_{3}'\bigr)=f^{b_{12}}(z_{1}, z_{2}; u,
w_{1}, w_{2}, w_{3}').\]
Applying $b\in \PB_{3}$ to both sides of this equality, we obtain
\begin{equation}\label{g_1=b-12}
f^{b}\bigl(z_{1}, z_{2}; g_{1}^{-1}u,
w_{1}, w_{2}, w_{3}'\bigr)=f^{b_{12}b}(z_{1}, z_{2}; u,
w_{1}, w_{2}, w_{3}')
\end{equation}
for $b\in \PB_{3}$. By~\eqref{g_1=b-12} with $b=
b_{12}^{-1}b_{13}^{-1}b_{23}^{-1}$, we see that
\eqref{A-2} is equal to
\[
f^{b_{13}^{-1}b_{23}^{-1}}\bigl(z_{1}^{-1}, z_{2}^{-1}; {\rm e}^{z_{1}L_{V}(1)}\bigl(-z_{1}^{-2}\bigr)^{L_{V}(0)}u,
 {\rm e}^{z_{2}L_{W_{1}}(1)}{\rm e}^{ \pi \i L_{W_{1}}(0)}{\rm e}^{-2(\log z_{2}) L_{W_{1}}(0)}
w_{1}, w_{2}, w_{3}'\bigr),
\]
which by definition, is equal to
$h_{+}^{e}(z_{1}, z_{2}; u, w_{2}, w_{1}, w_{3}')$ when
$|z_{1}|>|z_{2}|>0$,
$|{\arg (z_{1}-z_{2})}-\arg z_{1}|<\frac{\pi}{2}$ and
$\arg z_{1}, \arg z_{2}\ne 0$.

What we need to prove is that when
$|z_{1}|>|z_{2}|>0$ and $|{\arg (z_{1}-z_{2})-\arg z_{1}}|<\frac{\pi}{2}$,
the sum of the left-hand side of~\eqref{A-0} is equal to
$h_{+}^{e}(z_{1}, z_{2}; u, w_{2}, w_{1}, w_{3}')$.
We choose ${z_{1}, z_{2}\in -\R_{+}}$ such that $|z_{1}|=-z_{1}>-z_{2}={|z_{2}|>0}$.
Then $\big|z_{2}^{-1}\big|=-z_{2}^{-1}>-z_{1}^{-1}=\big|z_{1}^{-1}\big|>0$. In~this case,
$\arg (z_{1}-z_{2})=\pi$,
$\arg \bigl(z_{1}^{-1}-z_{2}^{-1}\bigr)=0$, and $\arg z_{1}=\arg z_{2}=\arg z_{1}^{-1}={\arg z_{2}^{-1}=\pi}$. Hence
$|{\arg (z_{1}-z_{2})-\arg z_{1}}|=0<\frac{\pi}{2}$ and
$-\frac{3\pi}{2}<0-\pi=\arg \bigl(z_{1}^{-1}-z_{2}^{-1}\bigr)-\arg z_{2}^{-1}<-\frac{\pi}{2}$.
From the discussion above, for such $z_{1}$ and $z_{2}$,
the sum of the left-hand side of~\eqref{A-0} is equal to
$h_{+}^{e}(z_{1}, z_{2}; u, w_{2}, w_{1}, w_{3}')$. On the other hand,
such $z_{1}$ and $z_{2}$ satisfy $|z_{1}|>|z_{2}|>0$ and~$|{\arg (z_{1}-z_{2})-\arg z_{1}}|<\frac{\pi}{2}$.
Since both $h_{+}^{e}(z_{1}, z_{2}; u, w_{2}, w_{1}, w_{3}')$ and the sum of the left-hand side of~\eqref{A-0}
are analytic extensions of their restrictions on the subset
given by $(z_{1}, z_{2})$ for~${z_{1}, z_{2}\in -\R_{+}}$ such that $|z_{1}|=-z_{1}>-z_{2}=|z_{2}|>0$,
the sum of the left-hand side of~\eqref{A-0} is equal to
$h_{+}^{e}(z_{1}, z_{2}; u, w_{2}, w_{1}, w_{3}')$ when
$|z_{1}|>|z_{2}|>0$ and $|{\arg (z_{1}-z_{2})-\arg z_{1}}|<\frac{\pi}{2}$.

Generalizing the convergence and analytic extension of~\eqref{A-0}
above, we can prove that
\begin{equation}\label{contrag-k-prod-1}
\big\langle \phi_{g_{1}}\bigl(\bigl(Y_{W_{2}}^{g_{2}}\bigr)'\bigr)(u_{1}, z_{1})
\cdots \phi_{g_{1}}\bigl(\bigl(Y_{W_{2}}^{g_{2}}\bigr)'\bigr)(u_{k-1}, z_{k-1})
A_{+}(\mathcal{Y})(w_{1}, z_{k})w_{3}', w_{2}\big\rangle
\end{equation}
is absolutely convergent on the region
$|z_{1}|>\cdots >|z_{k}|>0$ and its sum can be maximally extended to a~multivalued analytic function on the region
$M^{k}$ for $k\in \Z_{+}+3$, $u_{1}, \dots, u_{k-1} \in V$,
$w_{1}\in W_{1}$, $w_{2}\in W_{2}$ and $w_{3}'\in W_{3}'$.
In fact, the same calculations as in~\eqref{A-0} shows that
\eqref{contrag-k-prod-1} is equal to
\begin{gather}
\big\langle w'_{3}, \mathcal{Y}^{-1}\bigl({\rm e}^{z_{k}L_{W_{1}}(1)}{\rm e}^{ \pi \i L_{W_{1}}(0)}{\rm e}^{-2(\log z_{k}) L_{W_{1}}(0)}
w_{1}, z_{k}^{-1}\bigr)\nn
\qquad  \times
\bigl(Y_{W_{2}}
^{g_{2}}\bigr)^{-1}\bigl({\rm e}^{z_{k-1}L_{V}(1)}\bigl(-z_{k-1}^{-2}\bigr)^{L_{V}(0)}
g_{1}^{-1}u_{k-1}, z_{k-1}^{-1}\bigr)\cdots\nn
 \qquad {}\times
\bigl(Y_{W_{2}}
^{g_{2}}\bigr)^{-1}\bigl({\rm e}^{z_{1}L_{V}(1)}\bigl(-z_{1}^{-2}\bigr)^{L_{V}(0)}
g_{1}^{-1}u_{1}, z_{1}^{-1}\bigr)w_{2}\big\rangle\label{contrag-k-prod-2}
\end{gather}
when $\arg z_{1}, \dots, \arg z_{k}\ne 0$.
From the duality properties in
Definitions~\ref{twisted-mod} and~\ref{def-tw-int-op}, we know that
\begin{gather}
\big\langle w'_{3}, \bigl(Y_{W_{2}}
^{g_{2}}\bigr)^{-1}\bigl({\rm e}^{z_{k-1}L_{V}(1)}\bigl(-z_{k-1}^{-2}\bigr)^{L_{V}(0)}
g_{1}^{-1}u_{k-1}, z_{k-1}^{-1}\bigr)\cdots\nn
\qquad{}\times
\bigl(Y_{W_{2}}
^{g_{2}}\bigr)^{-1}\bigl({\rm e}^{z_{1}L_{V}(1)}\bigl(-z_{1}^{-2}\bigr)^{L_{V}(0)}
g_{1}^{-1}u_{1}, z_{1}^{-1}\bigr)\nn
\qquad {} \times
\mathcal{Y}^{-1}\bigl({\rm e}^{z_{k}L_{W_{1}}(1)}{\rm e}^{ \pi \i L_{W_{1}}(0)}{\rm e}^{-2(\log z_{k}) L_{W_{1}}(0)}
w_{1}, z_{k}^{-1}\bigr)
w_{2}\big\rangle\label{contrag-k-prod-2.5}
\end{gather}
is absolutely convergent on the region $\big|z_{k-1}^{-1}\big|>\cdots >\big|z_{1}^{-1}\big|>\big|z_{k}^{-1}\big|>0$.
Then for $M_{ij}\in \Z_{+}$ ($1\le j< i\le k-1$) such that
\[
x^{M_{ij}}Y_{V}\bigl({\rm e}^{z_{i}L_{V}(1)}\bigl(-z_{i}^{-2}\bigr)^{L_{V}(0)}
g_{1}^{-1}u_{i}, x\bigr){\rm e}^{z_{j}L_{V}(1)}\bigl(-z_{j}^{-2}\bigr)^{L_{V}(0)}
g_{1}^{-1}u_{j}\in V[[x]],\]
we have
\begin{gather}
\prod_{1\le j<i\le k-1}\bigl(z_{i}^{-1}-z_{j}^{-1}\bigr)^{M_{ij}}
\langle w'_{3}, \bigl(Y_{W_{2}}
^{g_{2}}\bigr)^{-1}\bigl({\rm e}^{z_{k-1}L_{V}(1)}\bigl(-z_{k-1}^{-2}\bigr)^{L_{V}(0)}
g_{1}^{-1}u_{k-1}, z_{k-1}^{-1}\bigr)\cdots\nn
 \qquad\times
\bigl(Y_{W_{2}}
^{g_{2}}\bigr)^{-1}\bigl({\rm e}^{z_{1}L_{V}(1)}\bigl(-z_{1}^{-2}\bigr)^{L_{V}(0)}
g_{1}^{-1}u_{1}, z_{1}^{-1}\bigr)\nn
\qquad \times
\mathcal{Y}^{-1}\bigl({\rm e}^{z_{k}L_{W_{1}}(1)}{\rm e}^{ \pi \i L_{W_{1}}(0)}{\rm e}^{-2(\log z_{k}) L_{W_{1}}(0)}
w_{1}, z_{k}^{-1}\bigr)
w_{2}\big\rangle\label{contrag-k-prod-3}
\end{gather}
is absolutely convergent on the region $\big|z_{i}^{-1}\big|
>\big|z_{k}^{-1}\big|>0$ for $i=1, \dots, k-1$,
for $i\ne j$ because~${z_{i}^{-1}=z_{j}^{-1}}$ for $i\ne j$ are poles of~\eqref{contrag-k-prod-2.5} and
are canceled by $\prod_{1\le j<i\le k-1}(z_{i}^{-1}-z_{j}^{-1})^{M_{ij}}$ in~\eqref{contrag-k-prod-3}.
Moreover, the sum of~\eqref{contrag-k-prod-3} can be analytically
extended to a~maximal multivalued analytic function on
$\big\{(z_{1}, \dots, z_{k})\in \C^{k}\mid
z_{i}\ne z_{k},\, i=1, \dots, k-1\big\}$.
Using the duality property in Definition~\ref{def-tw-int-op} repeatedly,
we see that
\begin{gather*}
\prod_{1\le i<j\le k-1}\bigl(z_{i}^{-1}-z_{j}^{-1}\bigr)^{M_{ij}}
\big\langle w'_{3}, \mathcal{Y}^{-1}\bigl({\rm e}^{z_{k}L_{W_{1}}(1)}{\rm e}^{ \pi \i L_{W_{1}}(0)}{\rm e}^{-2(\log z_{k}) L_{W_{1}}(0)}
w_{1}, z_{k}^{-1}\bigr)\\
 \qquad \times
\bigl(Y_{W_{2}}
^{g_{2}}\bigr)^{-1}\bigl({\rm e}^{z_{k-1}L_{V}(1)}\bigl(-z_{k-1}^{-2}\bigr)^{L_{V}(0)}
g_{1}^{-1}u_{k-1}, z_{k-1}^{-1}\bigr)\cdots\\
 \qquad \times
\bigl(Y_{W_{2}}
^{g_{2}}\bigr)^{-1}\bigl({\rm e}^{z_{1}L_{V}(1)}\bigl(-z_{1}^{-2}\bigr)^{L_{V}(0)}
g_{1}^{-1}u_{1}, z_{1}^{-1}\bigr)w_{2}\big\rangle
\end{gather*}
is absolutely convergent on the region \smash{$\big|z_{i}^{-1}\big|
<\big|z_{k}^{-1}\big|$} for $i=1, \dots, k-1$, \smash{$z_{i}^{-1}\allowbreak\ne z_{j}^{-1}$}
for~${i\ne j}$, and its sum can be analytically
extended to a~maximal multivalued analytic function on
$\big\{(z_{1}, \dots, z_{k})\allowbreak\in \C^{k}\mid
z_{i}\ne z_{k},\, i=1, \dots, k-1\big\}$.
Thus~\eqref{contrag-k-prod-2} is absolutely convergent
on the region $|z_{1}|>\cdots >|z_{k}|>0$ and consequently
\eqref{contrag-k-prod-1} is absolutely convergent
on the region $|z_{1}|>\cdots >|z_{k}|>0$, $\arg z_{1}, \dots, \arg z_{k}\ne 0$.
But the absolute convergence of~\eqref{contrag-k-prod-1}
on the region $|z_{1}|>\cdots >|z_{k}|>0$, $\arg z_{1}, \dots, \arg z_{k}\ne 0$
implies the absolute convergence of~\eqref{contrag-k-prod-1}
on the region $|z_{1}|>\cdots >{|z_{k}|>0}$. Moreover,
its sum has a~maximal analytic extension on the region $M^{k}$ and from the discussion
above, we have
\[
h(z_{1}, \dots, z_{k})=\sum_{j_{1}, \dots, j_{k-1}=J'}^{J}\sum_{i=1}^{I}
\sum_{m=0}^{M}z_{1}^{j_{1}}\cdots z_{k-1}^{j_{k-1}}z_{k}^{r_{i}}(\log z_{k})^{m}
f_{j_{1}, \dots, j_{k-1}; i, m}\bigl(z_{1}^{-1}, \dots,, z_{k}^{-1}\bigr),
\]
where $h(z_{1}, \dots, z_{k})$ is the maximal analytic extension
of the sum of~\eqref{contrag-k-prod-1} on the region $M^{k}$,
$J', J\in \Z$, $r_{1}, \dots, r_{I}\in \C$, $M\in \N$ such that
\begin{gather*}
{\rm e}^{z_{l}L_{V}(1)}\bigl(-z_{l}^{-2}\bigr)^{L_{V}(0)}
g_{1}^{-1}u_{l}=\sum_{j_{l}=J'}^{J}z_{l}^{j_{l}}u_{l}^{j_{l}},\qquad l=1, \dots, k-1,\\
{\rm e}^{z_{k}L_{W_{1}}(1)}{\rm e}^{ \pi \i L_{W_{1}}(0)}{\rm e}^{-2(\log z_{k}) L_{W_{1}}(0)}
w_{1}=\sum_{i=1}^{I}
\sum_{m=0}^{M}z_{k}^{r_{i}}(\log z_{k})^{m}w_{1}^{i, m}
\end{gather*}
for $u_{l}^{j_{l}}\in V$ and $w_{1}^{i, m}\in W_{1}$, and
$f_{j_{1}, \dots, j_{k-1}; i, m}(z_{1}, \dots, z_{k})$ for $J', J\in \Z$, $r_{1}, \dots, r_{I}\in \C$, $M\in \N$
are the maximal analytic extensions
of the sums of
$
\smash{\big\langle w'_{3}, Y_{W_{2}}
^{g_{2}}\bigl(u_{1}^{j_{1}}, z_{1}\bigr)\cdots
Y_{W_{2}}
^{g_{2}}\bigl(u_{k-1}^{j_{k-1}}, z_{k-1}\bigr)
\mathcal{Y}\bigl(w_{1}^{i, m}},\allowbreak z_{k}\bigr)
w_{2}\big\rangle$.
We will show at the end of this proof that a~component-isolated singularity
of the form $(z_{1}, \dots,\allowbreak z_{k-1})-\beta=\delta$ for $\beta\in \C^{k-1}$
and $\delta\in \{0, \infty\}^{k}$ of any value
at
$\zeta=z_{k}$ of $\smash{\prod_{1\le i<\le j\le k-1}}(z_{i}-z_{j})^{M_{ij}}h(z_{1}, \dots, z_{k-1}, \zeta)$
 is regular.

Next we consider the product of $A_{+}(\Y)$ and the twisted vertex operator
$\bigl(Y_{W_{3}}^{g_{2}}\bigr)'$. Let $u$, $w_{1}$, $w_{2}$
and $w_{3}'$ be the same as above.
When \smash{$\big|z_{1}^{-1}\big|>\big|z_{2}^{-1}\big|>0$} (or equivalently $|z_{2}|>|z_{1}|>0$),
\begin{gather}
\big\langle w'_{3}, Y_{W_{2}}
^{g_{2}}\bigl({\rm e}^{z_{1}L_{V}(1)}\bigl(-z_{1}^{-2}\bigr)^{L_{V}(0)}u, z_{1}^{-1}\bigr)\nn
\qquad\times
\mathcal{Y}\bigl({\rm e}^{z_{2}L_{W_{1}}(1)}{\rm e}^{ \pi \i L_{W_{1}}(0)}{\rm e}^{-2(\log z_{2})
 L_{W_{1}}(0)}
w_{1}, z_{2}^{-1}\bigr) w_{2}\big\rangle\label{A-5.3}
\end{gather}
converges absolutely and if in addition, \smash{$\big|{\arg \bigl(z_{1}^{-1}-z_{2}^{-1}\bigr)-\arg z_{1}^{-1}}\big|<\frac{\pi}{2}$},
its sum is equal to
\begin{gather*}
f^{e}\bigl(z_{1}^{-1}, z_{2}^{-1};
{\rm e}^{z_{1}L_{V}(1)}\bigl(-z_{1}^{-2}\bigr)^{L_{V}(0)}u,
 {\rm e}^{z_{2}L_{W_{1}}(1)}{\rm e}^{ \pi \i L_{W_{1}}(0)}{\rm e}^{-2(\log z_{2})
 L_{W_{1}}(0)}
w_{1}, w_{2}, w_{3}'\bigr).
\end{gather*}
We know that when $\arg z_{1}, \arg z_{2}\ne 0$,
\begin{gather}
\big\langle A_{+}(\mathcal{Y})(w_{1}, z_{2})
\bigl(Y_{W_{2}}^{g_{2}}\bigr)'(u, z_{1})w_{3}', w_{2}\big\rangle\nn
\qquad=\big\langle w'_{3}, (Y_{W_{2}}
^{g_{2}})^{-1}\bigl({\rm e}^{z_{1}L_{V}(1)}\bigl(-z_{1}^{-2}\bigr)^{L_{V}(0)}u, z_{1}^{-1}\bigr)\nn
\phantom{\qquad=}{} \times
\mathcal{Y}^{-1}\bigl({\rm e}^{z_{2}L_{W_{1}}(1)}{\rm e}^{ \pi \i L_{W_{1}}(0)}{\rm e}^{-2(\log z_{2}) L_{W_{1}}(0)}
w_{1}, z_{2}^{-1}\bigr) w_{2}\big\rangle\label{A-6}
\end{gather}
can be obtained using the multivalued analytic function~\eqref{A-1}
on the region $\big|z_{1}^{-1}\big|>\big|z_{2}^{-1}\big|>0$
starting from the value given by
\eqref{A-5.3} by letting $z_{2}^{-1}$
go around $0$ clockwise once (corresponding to $b_{23}^{-1}$),
then letting $z_{1}^{-1}$ go around
$0$ clockwise once (corresponding to $b_{12}^{-1}b_{13}^{-1}$).
Then~\eqref{A-6} also
converges absolutely and if in addition,
\smash{$\big|{\arg \bigl(z_{1}^{-1}-z_{2}^{-1}\bigr)-\arg z_{1}^{-1}}\big|<\frac{\pi}{2}$},
its sum is equal~to
\begin{gather}
f^{b_{23}^{-1}b_{12}^{-1}b_{13}^{-1}}\bigl(z_{1}^{-1}, z_{2}^{-1};
{\rm e}^{z_{1}L_{V}(1)}\bigl(-z_{1}^{-2}\bigr)^{L_{V}(0)}u,
 {\rm e}^{z_{2}L_{W_{1}}(1)}{\rm e}^{ \pi \i L_{W_{1}}(0)}{\rm e}^{-2(\log z_{2})
 L_{W_{1}}(0)}
w_{1}, w_{2}, w_{3}'\bigr)\nn
\qquad= f^{b_{13}^{-1}b_{23}^{-1}b_{12}^{-1}}\bigl(z_{1}^{-1}, z_{2}^{-1};
{\rm e}^{z_{1}L_{V}(1)}\bigl(-z_{1}^{-2}\bigr)^{L_{V}(0)}u,\nn
\phantom{\qquad=}{}
 {\rm e}^{z_{2}L_{W_{1}}(1)}{\rm e}^{ \pi \i L_{W_{1}}(0)}{\rm e}^{-2(\log z_{2})
 L_{W_{1}}(0)}
w_{1}, w_{2}, w_{3}'\bigr).\label{A-7}
\end{gather}

We now show that the sum of~\eqref{A-6} is equal to~\eqref{A-7}.
We choose $z_{1}, z_{2}\in -\R_{+}$ such that $|z_{2}|=-z_{2}>-z_{1}=|z_{1}|>0$.
Then $\big|z_{1}^{-1}\big|=-z_{1}^{-1}>-z_{2}^{-1}=\big|z_{2}^{-1}\big|>0$. In this case,
$\arg (z_{1}-z_{2})=0$,
$\arg \bigl(z_{1}^{-1}-z_{2}^{-1}\bigr)=\pi$, and $\arg z_{1}=\arg z_{2}=\arg z_{1}^{-1}=\arg z_{2}^{-1}=\pi$. Hence
$-\frac{3\pi}{2}<0-\pi= \arg (z_{1}-z_{2})-\arg z_{2}<-\frac{\pi}{2}$ and
$|{\arg \bigl(z_{1}^{-1}-z_{2}^{-1}\bigr)-\arg z_{1}^{-1}}|=0<\frac{\pi}{2}$.
From the discussion above, for such $z_{1}$ and $z_{2}$,
the sum of~\eqref{A-6} is equal to~\eqref{A-7}.
Using analytic extension, we see that when $|z_{2}|>|z_{1}|>0$ and
$-\frac{3\pi}{2}<\arg (z_{1}-z_{2})-\arg z_{2}<-\frac{\pi}{2}$,
the sum of~\eqref{A-6} is equal to~\eqref{A-7}.

On the region
$|z_{2}|>|z_{1}|>0$
, \eqref{A-7}
is in fact equal to $h_{+}^{e}(z_{1}, z_{2}; u, w_{2}, w_{1}, w_{3}')$.
This can be proved as follows: On
the intersection of the region $|z_{2}|>|z_{1}|>0$ and $M^{2}_{0}$, the branch~${h_{+}^{e}(z_{1}, z_{2}; u, w_{2}, w_{1}, w_{3}')}$
is obtained by analytically extending
\eqref{A-1.1} defined on the intersection of the region
$|z_{1}|>|z_{2}|>0$, $|{\arg (z_{1}-z_{2})-\arg z_{1}}|<\frac{\pi}{2}$, $\arg z_{1}, \arg z_{2}\ne 0$ and
$M^{2}_{0}$.
We need to find what is this analytic extension.
Let $\xi, \zeta\in -\R_{+}$ satisfying $\xi<\zeta<0$.
Then we have~${|\xi|>|\zeta|>0}$, $\arg (\xi-\zeta)=\arg \xi=\arg \zeta=\pi$ so that $|{\arg (\xi-\zeta)-\arg \xi}|=0<\frac{\pi}{2}$. By definition,
\begin{align}
h_{+}^{e}(\xi, \zeta; u, w_{2}, w_{1}, w_{3}')={}&f^{b_{13}^{-1}b_{23}^{-1}}
\bigl(\xi^{-1}, \zeta^{-1}; {\rm e}^{\xi^{-1}L_{V}(1)}\bigl(-\xi^{2}\bigr)^{-L_{V}(0)}u,\nn
& {\rm e}^{\zeta L_{W_{1}}(1)}{\rm e}^{ \pi \i L_{W_{1}}(0)}{\rm e}^{-2(\log \zeta) L_{W_{1}}(0)}
w_{1}, w_{2}, w_{3}'\bigr).\label{A-7.5}
\end{align}
Let $\gamma=(\gamma_{1}, \gamma_{2})$ be the path from
$(\xi, \zeta)$ to $(\zeta, \xi)$
given by the upper half circle $\gamma_{1}$ centered at~\smash{$\frac{\xi+\zeta}{2}$}
with radius $\frac{-\xi+\zeta}{2}$ from $\xi$ to $\zeta$
and the lower half circle $\gamma_{2}$ centered at $\frac{\xi+\zeta}{2}$
with radius~\smash{$\frac{-\xi+\zeta}{2}$} from $\zeta$ to $\xi$. Since by definition,
$(\xi, \zeta)$, $(\zeta, \xi)$ are in $M^{2}_{0}$
and the path $\gamma$ does not cross the boundary of $M^{2}_{0}$,
$\gamma$ is a~continuous path in $M^{2}_{0}$.
So by definition, the value at $(z_{1}, z_{2})=(\zeta, \xi)$ of the multivalued analytic function
$h_{+}(z_{1}, z_{2}; u, w_{2}, w_{1}, w_{3}')$ obtained
 by starting with the
value~${h_{+}^{e}(\xi, \zeta; u, w_{2}, w_{1}, w_{3}')}$ and analytically extending
along $\gamma$ is $h_{+}^{e}(\zeta, \xi; u, w_{2}, w_{1}, w_{3}')$. On the other hand, we need to find out what
the right-hand side of~\eqref{A-7.5} is analytically extended to along the same path $\gamma$. We have
the path $\gamma'=\bigl(\gamma_{1}^{-1}, \gamma_{2}^{-1}\bigr)$
from $\bigl(\xi^{-1}, \zeta^{-1}\bigr)$ to $\bigl(\zeta^{-1}, \xi^{-1}\bigr)$, where $\gamma_{1}^{-1}$
and $\gamma_{2}^{-1}$ are the paths on the complex plane from $\xi^{-1}$ to $\zeta^{-1}$ and from
$\zeta^{-1}$ to $\xi^{-1}$, respectively, given by $\gamma_{1}^{-1}(t)=\gamma_{1}(t)^{-1}$
and $\gamma_{2}^{-1}(t)=\gamma_{1}(t)^{-1}$.
But $\gamma'$ is not a~continuous path
in~$M_{0}^{2}$ because when $t$ goes from $0$ to $1$, $\bigl(\gamma_{1}^{-1}(t), \gamma_{2}^{-1}(t)\bigr)$
goes from $\bigl(\xi^{-1}, \zeta^{-1}\bigr)$ to $\bigl(\zeta^{-1}, \xi^{-1}\bigr)$ along the path
$\gamma'$ but $\gamma_{1}^{-1}(t)-\gamma_{2}^{-1}(t)$ when $t$ goes from $0$ to $1$
crosses the positive real line clockwise.
Crossing the positive real line clockwise corresponds to changing the branch
by an action of $b_{12}^{-1}$. So we see that the value at $\bigl(z_{1}^{-1}, z_{2}^{-1}\bigr)=\bigl(\zeta^{-1}, \xi^{-1}\bigr)$
of the multivalued
analytic function
\begin{equation}\label{A-7.6}
f\bigl(z_{1}^{-1}, z_{2}^{-1};
{\rm e}^{z_{1}L_{V}(1)}\bigl(-z_{1}^{-2}\bigr)^{L_{V}(0)}u,
 {\rm e}^{z_{2}L_{W_{1}}(1)}{\rm e}^{ \pi \i L_{W_{1}}(0)}{\rm e}^{-2(\log z_{2})
 L_{W_{1}}(0)}
w_{1}, w_{2}, w_{3}'\bigr)
\end{equation}
of $z_{1}^{-1}$ and $z_{2}^{-1}$
obtained by starting with the value
\[f^{b_{13}^{-1}b_{23}^{-1}}
\bigl(\xi^{-1}, \zeta^{-1}; {\rm e}^{\xi^{-1}L_{V}(1)}\bigl(-\xi^{2}\bigr)^{-L_{V}(0)}u,
{\rm e}^{\zeta L_{W_{1}}(1)}{\rm e}^{ \pi \i L_{W_{1}}(0)}{\rm e}^{-2(\log \zeta) L_{W_{1}}(0)}
w_{1}, w_{2}, w_{3}'\bigr)\]
and going along the path $\gamma'$ is
\[f^{b_{13}^{-1}b_{23}^{-1}b_{12}^{-1}}
\bigl(\zeta^{-1}, \xi^{-1}; {\rm e}^{\zeta^{-1}L_{V}(1)}\bigl(-\zeta^{2}\bigr)^{-L_{V}(0)}u,
{\rm e}^{\xi L_{W_{1}}(1)}{\rm e}^{ \pi \i L_{W_{1}}(0)}{\rm e}^{-2(\log \xi) L_{W_{1}}(0)}
w_{1}, w_{2}, w_{3}'\bigr).\]
Since we have~\eqref{A-7.5} and the path $\gamma$ is mapped to
$\gamma'$ under the map $(z_{1}, z_{2})\mapsto \bigl(z_{1}^{-1}, z_{2}^{-1}\bigr)$,
the value $h_{+}^{e}(\zeta, \xi; u, w_{2}, w_{1}, w_{3}')$ of the multivalued
analytic function $h_{+}(z_{1}, z_{2}; u, w_{2}, w_{1}, w_{3}')$ of~$z_{1}$ and~$z_{2}$
at the end point $(z_{1}, z_{2})=(\zeta, \xi)$ of $\gamma$ and the value of the multivalued analytic function~\eqref{A-7.6} of $z_{1}^{-1}$ and $z_{2}^{-1}$ at the end point $\bigl(z_{1}^{-1}, z_{2}^{-1}\bigr)=\bigl(\zeta^{-1}, \xi^{-1}\bigr)$ of $\gamma'$
must be the same, that~is,
\begin{align*}
h_{+}^{e}(\zeta, \xi; u, w_{2}, w_{1}, w_{3}')={}&f^{b_{13}^{-1}b_{23}^{-1}b_{12}^{-1}}
\bigl(\zeta^{-1}, \xi^{-1}; {\rm e}^{\zeta^{-1}L_{V}(1)}\bigl(-\zeta^{2}\bigr)^{-L_{V}(0)}u,\\
&
{\rm e}^{\xi L_{W_{1}}(1)}{\rm e}^{ \pi \i L_{W_{1}}(0)}{\rm e}^{-2(\log \xi) L_{W_{1}}(0)}
w_{1}, w_{2}, w_{3}'\bigr).
\end{align*}
Since a~single-valued branch on a~subregion of $M_{0}^{2}$ of a~multivalued analytic function
on $M^{2}$ is determined uniquely by the value at one point in this subregion of
$M^{2}_{0}$, we see that on the intersection of the region
$|z_{2}|>|z_{1}|>0$ and $M_{0}^{2}$,
\begin{align}
h_{+}^{e}(z_{1}, z_{2}; u, w_{2}, w_{1}, w_{3}')={}&f^{b_{13}^{-1}b_{23}^{-1}b_{12}^{-1}}
\bigl(z_{1}^{-1}, z_{2}^{-1}; {\rm e}^{z_{1}^{-1}L_{V}(1)}\bigl(-z_{1}^{2}\bigr)^{-L_{V}(0)}u,\nonumber\\
& {\rm e}^{z_{2} L_{W_{1}}(1)}{\rm e}^{\pm \pi \i L_{W_{1}}(0)}{\rm e}^{-2(\log z_{2})
 L_{W_{1}}(0)}
w_{1}, w_{2}, w_{3}'\bigr).\label{A-8}
\end{align}
Since on the region $|z_{2}|>|z_{1}|>0$
and $\bigl|{\arg \bigl(z_{1}^{-1}-z_{2}^{-1}\bigr)-\arg z_{1}^{-1}}\bigr|<\frac{\pi}{2}$,
the sum of~\eqref{A-6} is equal to~\eqref{A-7} (that is, the right-hand side of
\eqref{A-8}), the sum of~\eqref{A-6} is indeed equal to~${h_{+}^{e}(z_{1}, z_{2}; u, w_{2}, w_{1}, w_{3}')}$ on the same region.

We can prove the absolute convergence of
\[
\big\langle \bigl(Y_{W_{2}}^{g_{2}}\bigr)'(u, z_{1})
A_{-}(\mathcal{Y})(w_{1}, z_{2})w_{3}', w_{2}\big\rangle,\qquad\big\langle A_{-}(\mathcal{Y})(w_{1}, z_{2})\phi_{g_{1}^{-1}}\bigl(\bigl(Y_{W_{3}}^{g_{2}}\bigr)'\bigr)(u, z_{1})w_{3}', w_{2}\big\rangle\]
in the corresponding regions to $h_{-}^{e}(z_{1}, z_{2};
u, w_{2}, w_{1}, w_{3}')$ similarly by generalizing
the proofs in~\cite{H-twisted-int} using the same method above for
$A_{-}(\mathcal{Y})$. Here we omit the details.

Now we study the iterate of $A_{\pm}(\Y)$ and the twisted vertex operator
$Y_{W_{1}}^{g_{1}}$. When
$\arg z_{2}\ne 0$, from~\eqref{A6}, we have
\begin{gather}
\big\langle A_{\pm}(\mathcal{Y})\bigl(Y_{W_{1}}^{g_{1}}(u, z_{1}-z_{2})w_{1}, z_{2}\bigr)w_{3}', w_{2}\big\rangle\nn
\qquad=\big\langle w'_{3},
\mathcal{Y}^{-1}\bigl({\rm e}^{z_{2}L_{W_{1}}(1)}{\rm e}^{\pm \pi \i L_{W_{1}}(0)}
{\rm e}^{-2(\log z_{2}) L_{W_{1}}(0)}
 Y_{W_{1}}^{g_{1}}(u, z_{1}-z_{2})w_{1}, z_{2}^{-1}\bigr) w_{2}\big\rangle.\label{A-11}
\end{gather}
As in~\cite{H-twisted-int}, we have in the region $|z_{2}|>|z_{1}-z_{2}|>0$
\begin{gather*}
e^{z_{2}L_{W_{1}}(1)}{\rm e}^{\pm \pi \i L_{W_{1}}(0)}{\rm e}^{-2(\log z_{2}) L_{W_{1}}(0)}
 Y_{W_{1}}^{g_{1}}(u, z_{1}-z_{2})w_{1}\nn
\qquad=Y_{W_{1}}^{g_{1}}\left({\rm e}^{z_{1}L_{V}(1)}\bigl(-z_{1}^{-2}\bigr)^{L_{V}(0)}u,
\frac{xx_0}{(x_2+x_0)x_2}\right) \lbar_{\substack{x_{0}^{n}={\rm e}^{n\log (z_{1}-z_{2})},
\log x_{0}=\log (z_{1}-z_{2}), x_{2}^{n}={\rm e}^{n\log z_{2}}\\ \log x_{2}=\log z_{2},
 x^{n}={\rm e}^{\pm n\pi \i}, \log x=\pm \pi \i}} \nn
\phantom{\qquad=}{}\times {\rm e}^{z_{2}L_{W_{1}}(1)}{\rm e}^{\pm \pi \i L_{W_{1}}(0)}{\rm e}^{-2(\log z_{2}) L_{W_{1}}(0)}
w_{1}.
\end{gather*}
As in the proof of the first part of Theorem~\ref{skew-sym}, we see that
\begin{gather*}
\mathcal{Y}^{-1}\left(Y_{W_{1}}^{g_{1}}\left({\rm e}^{z_{1}L_{V}(1)}\bigl(-z_{1}^{-2}\bigr)^{L_{V}(0)}u,
\frac{xx_0}{(x_2+x_0)x_2}\right) \lbar_{\substack{x_{0}^{n}={\rm e}^{n\log (z_{1}-z_{2})},
\log x_{0}=\log(z_{1}-z_{2}), x_{2}^{n}={\rm e}^{n\log z_{2}}\\ \log x_{2}=\log z_{2},
 x^{n}={\rm e}^{\pm n\pi \i}, \log x=\pm \pi \i}}\right. \nn
\left. \qquad \times {\rm e}^{z_{2}L_{W_{1}}(1)}{\rm e}^{\pm \pi \i L_{W_{1}}(0)}{\rm e}^{-2(\log z_{2}) L_{W_{1}}(0)}
w_{1}, z_{2}^{-1}\right)
\end{gather*}
is an iterated series obtained by
expanding one variable inside the series obtained from
 the iterate of $\mathcal{Y}^{-1}$
and $Y_{W_{1}}^{g_{1}}$. Using the same method
as in the proof of the first part of Theorem~\ref{skew-sym} (including the use of Lemma
 \ref{convergence}),
we can prove that the corresponding multiple series
is absolutely convergent. Here we omit the detailed proof. In particular, we can calculate
the sum of the series using any of the iterated sums associated with
the multisum. Then the same calculations as those in the proof of \cite[Theorem 6.1]{H-twisted-int} shows that when
$|z_{1}|, |z_{2}|>|z_{1}-z_{2}|>0$ and~$|{\arg z_{1}-\arg z_{2}}|<\frac{\pi}{2}$,
the right-hand side of~\eqref{A-11} is equal to
\begin{gather}
\big\langle w'_{3},
\mathcal{Y}^{-1}\bigl(\bigl(Y_{W_{1}}^{g_{1}}\bigr)^{m+\frac{1\pm 1}{2}}\bigl({\rm e}^{z_{1}L_{V}(1)}\bigl(-z_{1}^{-2}\bigr)^{L_{V}(0)}u, z_{1}^{-1}-z_{2}^{-1}\bigr) \nn
\qquad
 \times {\rm e}^{z_{2}L_{W_{1}}(1)}{\rm e}^{\pm \pi \i L_{W_{1}}(0)}{\rm e}^{-2(\log z_{2})
 L_{W_{1}}(0)}
w_{1}, z_{2}^{-1}\bigr) w_{2}\big\rangle,\label{A-12}
\end{gather}
where $m\in \Z$ is given by
\[\log (z_{1}-z_{2})-\log z_{1}-\log z_{2}+ \pi {\rm i}
=l_{m+1}\bigl(z_{1}^{-1}-z_{2}^{-1}\bigr).\]

We know that
\begin{gather}
\big\langle w'_{3}, \mathcal{Y}\bigl(Y_{W_{1}}^{g_{1}}
\bigl({\rm e}^{z_{1}L_{V}(1)}\bigl(-z_{1}^{-2}\bigr)^{L_{V}(0)}u, z_{1}^{-1}-z_{2}^{-1}\bigr)
 {\rm e}^{z_{2}L_{W_{1}}(1)}{\rm e}^{\pm \pi \i L_{W_{1}}(0)}\nn
 \qquad\times{\rm e}^{-2(\log z_{2})
 L_{W_{1}}(0)}
w_{1}, z_{2}^{-1}\bigr) w_{2}\big\rangle\label{A-13}
\end{gather}
is absolutely convergent
on the region $\big|z_{2}^{-1}\big|>\big|z_{1}^{-1}-z_{2}^{-1}\big|>0$ and, if in addition
$\big|{\arg z_{1}^{-1}}-\arg z_{2}^{-1}\big|<\frac{\pi}{2}$, the sum
is equal to
\[
f^{e}\bigl(z_{1}^{-1}, z_{2}^{-1};
{\rm e}^{z_{1}L_{V}(1)}\bigl(-z_{1}^{-2}\bigr)^{L_{V}(0)}u,
{\rm e}^{z_{2}L_{W_{1}}(1)}{\rm e}^{\pm \pi \i L_{W_{1}}(0)}{\rm e}^{-2(\log z_{2})
 L_{W_{1}}(0)}
w_{1}, w_{2}, w_{3}'\bigr).
\]
Since~\eqref{A-12} can be obtained from~\eqref{A-13}
using the multivalued analytic function~\eqref{A-1}
on the region $\big|z_{2}^{-1}\big|>\big|z_{1}^{-1}-z_{2}^{-1}\big|>0$
starting from the value given by
\eqref{A-13} by letting $z_{2}^{-1}$ go around~$0$ clockwise
once while $z_{1}^{-1}-z_{2}^{-1}$ is fixed
(corresponding to $b_{13}^{-1}b_{23}^{-1}$ since to keep~${z_{1}^{-1}-z_{2}^{-1}}$ fixed, $z_{1}^{-1}$ must also
go around $0$ clockwise once) and then letting
$z_{1}^{-1}$ go around $z_{2}^{-1}$ counterclockwise~\smash{$m+\frac{1\pm 1}{2}$} times \big(corresponding to
\smash{$b_{12}^{m+\frac{1\pm 1}{2}}$}\big), we see that
\eqref{A-12} is absolutely convergent on the region
$\big|z_{2}^{-1}\big|>\big|z_{1}^{-1}-z_{2}^{-1}\big|>0$
and if $\big|{\arg z_{1}^{-1}-\arg z_{2}^{-1}}\big|<\frac{\pi}{2}$, the sum
is equal to
\begin{gather*}
f^{b_{13}^{-1}b_{23}^{-1}b_{12}^{m+\frac{1\pm 1}{2}}}
\bigl(z_{1}^{-1}, z_{2}^{-1};
{\rm e}^{z_{1}L_{V}(1)}\bigl(-z_{1}^{-2}\bigr)^{L_{V}(0)}u,\\
\qquad
{\rm e}^{z_{2}L_{W_{1}}(1)}{\rm e}^{\pm \pi \i L_{W_{1}}(0)}{\rm e}^{-2(\log z_{2})
 L_{W_{1}}(0)}
w_{1}, w_{2}, w_{3}'\bigr).
\end{gather*}

In the case of $A_{+}$, we consider the subset $S_{+}$ of $(z_{1}, z_{2})\in \C^{2}$ such that $z_{1}, z_{2}\in -\R_{+}$,
 $z_{1}<z_{2}$, and $-z_{1}^{-1}>z_{1}^{-1}-z_{2}^{-1}>0$.
This set is certainly nonempty since for $z_{1}=-\frac{1}{2}$, $z_{2}=-\frac{1}{3}$, we have
$z_{1}<z_{2}$ and $-z_{1}^{-1}=2>1=z_{1}^{-1}-z_{2}^{-1}>0$.
Then for $(z_{1}, z_{2})\in S_{+}$, we have
\begin{gather*}
|z_{1}|>|z_{2}|>0, \qquad\arg z_{1}=\arg z_{2}=\arg z_{1}^{-1}=\arg z_{2}^{-1}
=\arg (z_{1}-z_{2})=\pi,\\
 \arg \bigl(z_{1}^{-1}-z_{2}^{-1}\bigr)=0,\qquad
|{\arg (z_{1}-z_{2})-\arg z_{1}}|=0<\frac{\pi}{2},\\
\big|z_{2}^{-1}\big|=-z_{2}^{-1}>-z_{1}^{-1}>z_{1}^{-1}-z_{2}^{-1}=\big|z_{1}^{-1}-z_{2}^{-1}\big|>0,\qquad
\big|{\arg z_{1}^{-1}-\arg z_{2}^{-1}}\big|=0<\frac{\pi}{2},
\end{gather*}
and
\begin{gather*}
\log (z_{1}-z_{2})-\log z_{1}-\log z_{2}+\pi \i\nn
 \qquad =\log |z_{1}-z_{2}|+\arg (z_{1}-z_{2}) \i
-\log |z_{1}|-\arg z_{1} \i-\log |z_{2}|-\arg z_{2} \i+\pi \i\nn
 \qquad =\log \left|\frac{z_{1}-z_{2}}{z_{1}z_{2}}\right|=\log \big|z_{1}^{-1}-z_{2}^{-1}\big|+\arg \bigl(z_{1}^{-1}-z_{2}^{-1}\bigr)\nn
\qquad=\log \bigl(z_{1}^{-1}-z_{2}^{-1}\bigr).
\end{gather*}
In particular, for $(z_{1}, z_{2})\in S_{+}$, we have $m+1=0$.
Then on $S_{+}$, \eqref{A-12} is absolutely convergent~to
\begin{gather*}
f^{b_{13}^{-1}b_{23}^{-1}}
\bigl(z_{1}^{-1}, z_{2}^{-1};
{\rm e}^{z_{1}L_{V}(1)}\bigl(-z_{1}^{-2}\bigr)^{L_{V}(0)}u,
{\rm e}^{z_{2}L_{W_{1}}(1)}{\rm e}^{ \pi \i L_{W_{1}}(0)}{\rm e}^{-2(\log z_{2})
 L_{W_{1}}(0)}
w_{1}, w_{2}, w_{3}'\bigr)\\
\qquad=h_{+}^{e}(z_{1}, z_{2}; u, w_{1}, w_{3}', w_{2}).
\end{gather*}
By analytic extension, we see that the sum of~\eqref{A-12} is equal to
$h_{+}^{e}(z_{1}, z_{2}; u, w_{1}, w_{3}', w_{2})$
on the subregion of $M_{0}^{2}$ given by $\big|z_{2}^{-1}\big|>\big|z_{1}^{-1}-z_{2}^{-1}\big|>0$,
 $\big|{\arg z_{1}^{-1}-\arg z_{2}^{-1}}\big|<\frac{\pi}{2}$.
By analytic extension again, we see that the sum of
\eqref{A-11} is equal to $h_{+}^{e}(z_{1}, z_{2}; u, w_{1}, w_{3}', w_{2})$
on the subregion of $M_{0}^{2}$ given by
$|z_{2}|>|z_{1}-z_{2}|>0$, $|{\arg z_{1}-\arg z_{2}}|<\frac{\pi}{2}$.

In the case of $A_{-}$, we consider the
subset $S_{-}$ of $(z_{1}, z_{2})\in \C^{2}$ such that $z_{1}, z_{2}\in -\R_{+}$,
 $z_{1}<z_{2}$, and $-z_{1}^{-1}>z_{1}^{-1}-z_{2}^{-1}>0$.
This set is nonempty since for $z_{1}=-\frac{1}{2}$, $z_{2}=-\frac{1}{3}$, we have
$z_{1}<z_{2}$ and $-z_{1}^{-1}=2>1=-2+3=z_{1}^{-1}-z_{2}^{-1}>0$. Then for
 for $(z_{1}, z_{2})\in S_{-}$, we have $|z_{1}|>|z_{2}|>0$, $|{\arg (z_{1}-z_{2})-\arg z_{1}}|=0<\frac{\pi}{2}$,
$\arg (z_{1}-z_{2})=\pi \ge \pi$, $\big|z_{2}^{-1}\big|=-z_{2}^{-1}>-z_{1}^{-1}>z_{1}^{-1}-z_{2}^{-1}=\bigl|z_{1}^{-1}-z_{2}^{-1}\bigr|>0$,
$\bigl|{\arg z_{1}^{-1}-\arg z_{2}^{-1}}\bigr|=0<\frac{\pi}{2}$, and
\[\log (z_{1}-z_{2})-\log z_{1}-\log z_{2}+\pi \i =\log \left|\frac{z_{1}-z_{2}}{z_{1}z_{2}}\right| \i
=\log \bigl(z_{1}^{-1}-z_{2}^{-1}\bigr).\]
In particular, for $(z_{1}, z_{2})\in S_{-}$, we also have $m+1=0$.
Then on $S_{-}$, \eqref{A-12} is absolutely convergent to
\begin{gather*}
f^{b_{13}^{-1}b_{23}^{-1}b_{12}^{-1}}
\bigl(z_{1}^{-1}, z_{2}^{-1};
{\rm e}^{z_{1}L_{V}(1)}\bigl(-z_{1}^{-2}\bigr)^{L_{V}(0)}u,
{\rm e}^{z_{2}L_{W_{1}}(1)}{\rm e}^{- \pi \i L_{W_{1}}(0)}{\rm e}^{-2(\log z_{2})
 L_{W_{1}}(0)}
w_{1}, w_{2}, w_{3}'\bigr)\\
\qquad=h_{-}^{e}(z_{1}, z_{2}; u, w_{1}, w_{3}', w_{2}).
\end{gather*}
By analytic extension, we see that the sum of~\eqref{A-12} is equal to
$h_{-}^{e}(z_{1}, z_{2}; u, w_{1}, w_{3}', w_{2})$
on the subregion of $M_{0}^{2}$ given by $\big|z_{2}^{-1}\big|>\big|z_{1}^{-1}-z_{2}^{-1}\big|>0$,
 $\big|{\arg z_{1}^{-1}-\arg z_{2}^{-1}}\big|<\frac{\pi}{2}$.
By analytic extension again, we see that the sum of
\eqref{A-11} is equal to $h_{-}^{e}(z_{1}, z_{2}; u, w_{1}, w_{3}', w_{2})$
on the subregion of $M_{0}^{2}$ given by
$|z_{2}|>|z_{1}-z_{2}|>0$, $|{\arg z_{1}-\arg z_{2}}|<\frac{\pi}{2}$.

Finally we prove that for any fixed $z_{k}\in \C^{\times}$, a~component-isolated singularity
of the form $(z_{1}, \dots, z_{k-1})-\beta=\delta$ for $\beta\in \C^{k-1}$
and $\delta\in \{0, \infty\}^{k}$ of any value
at $\zeta=z_{k}$ of
$
\smash{\prod_{1\le i<\le j\le k-1}}(z_{i}-z_{j})^{M_{ij}}h(z_{1}, \dots, z_{k-1}, \zeta)
$
 is regular. Generalizing the proof above, we can show that
\[
\big\langle \phi_{g_{1}}\bigl(\bigl(Y_{W_{2}}^{g_{2}}\bigr)'\bigr)(u_{1}, z_{1})
\cdots \phi_{g_{1}}\bigl(\bigl(Y_{W_{2}}^{g_{2}}\bigr)'\bigr)(u_{k-2}, z_{k-2})
A_{\pm}(\mathcal{Y})(Y_{W_{1}}^{g_{1}}(u, z_{k-1}-z_{k})w_{1}, z_{k})w_{3}', w_{2}\big\rangle
\]
is absolutely convergent on the region $|z_{1}|>\cdots >|z_{k-2}|>|z_{k}|>|z_{k-1}-z_{k}|>0$ and the maximal
analytic extension of its sum on $M^{k}$ is $h(z_{1}, \dots, z_{k})$. But we have proved above that~\eqref{contrag-k-prod-1}
is absolutely convergent on the region $|z_{1}|>\cdots >|z_{k}|>0$ and the maximum analytic extension
of its sum is also $h(z_{1}, \dots, z_{k})$. Thus
\smash{$\prod_{1\le i<j\le k-1}(z_{i}-z_{j})^{M_{ij}}h(z_{1}, \dots, \zeta)$}
can be expanded as
\begin{gather}
\prod_{1\le i<j\le k-1}(z_{i}-z_{j})^{M_{ij}}\big\langle \phi_{g_{1}}\bigl(\bigl(Y_{W_{2}}^{g_{2}}\bigr)'\bigr)(u_{1}, z_{1})
\cdots \phi_{g_{1}}\bigl(\bigl(Y_{W_{2}}^{g_{2}}\bigr)'\bigr)(u_{k-1}, z_{k-1}) \nn
\qquad
 \times A_{\pm}(\mathcal{Y})\bigl(Y_{W_{1}}^{g_{1}}(u, z_{k-1}-z_{k})w_{1}, z_{k}\bigr)w_{3}', w_{2}\big\rangle,\label{A-14}
\end{gather}
a~series in powers of $z_{1}, \dots, z_{k-2}, z_{k-1}-z_{k}$ and
finitely many integral powers of $\log z_{1}, \dots, \allowbreak\log z_{k-2}, \log (z_{k-1}-z_{k})$.
From the expansion~\eqref{A-14}, we see that the component-isolated singularity
$(z_{1}, \dots, z_{k-2}, z_{k-1}-z_{k})=(\infty, \dots, \infty, 0)$ of the value~\smash{$\prod_{1\le i<j\le k-1}(z_{i}-z_{j})^{M_{ij}}h(z_{1}, \dots, z_{k})$} at $\zeta=z_{k}$
is regular. Using the commutativity for $A_{\pm}(\mathcal{Y})$ (which is already proved), we see that
the component-isolated singularity
$(z_{1}, \dots, z_{k-2}, z_{k-1}-z_{k})=(0, \dots, 0)$ of the value
$\prod_{1\le i<j\le k-1}(z_{i}-z_{j})^{M_{ij}}h(z_{1}, \dots, z_{k})$ at $\zeta=z_{k}$ is also regular.
The regularity of the other component-isolated singularity of the form above
can be proved either similarly or easily by using the relation between $h(z_{1}, \dots, z_{k})$ and
$f_{j_{1}, \dots, j_{k-1}; i, m}(z_{1}, \dots, z_{k})$.
\end{proof}

\section{Tensor product bifunctors and some natural isomorphisms}\label{sec4}

In this section, we introduce the notion of twisted $P(z)$-intertwining map
and give a~definition and a~construction of $P(z)$-tensor product of a~$g_{1}$-twisted module and a~$g_{2}$-twisted
$V$-module for~$g_{1}$,~$g_{2}$ in a~group
$G$ of automorphisms of $V$ in a~category of
twisted $V$-modules
under suitable assumptions.
Using the skew-symmetry
isomorphism $\Omega_{+}$ given in the preceding section,
we construct
$G$-crossed commutativity isomorphisms. We also construct parallel
transport isomorphisms. Using $G$-crossed commutativity
isomorphisms and parallel
transport isomorphisms, we construct
$G$-crossed braiding isomorphisms.
The material in this section is essentially the same as the corresponding
material in~\cite{H-rigidity,tensor1,HLZ3} except that $V$-modules and
intertwining maps are replaced by
twisted $V$-modules and twisted
intertwining maps.

Let $G$ be a~group of automorphisms of $V$. Let $\mathcal{C}_{1}$ be a~category of $g$-twisted $V$-modules without $g$-actions for $g\in G$.
For $g\in G$, let $W$ be a~$g$-twisted $V$-module with a~$g$-action and also with
an extension of the $g$-action on $W$ to an action of $G$. Here by an action of $G$ or a~$G$-action
on $W$, we mean that $W$ is a~left $G$-module. We shall call $W$ simply a~$g$-twisted $V$-module with a~$G$-action. But note that the $g$-action on $W$ has
additional properties but the actions of other elements of $G$ are not required to have any
properties.
Let $\mathcal{C}_{2}$ be a~category of
$g$-twisted $V$-modules with actions of $G$. In this paper, we do not require that morphisms in
$\mathcal{C}_{2}$ commute with actions of $G$. But to prove substantial results on $\mathcal{C}_{2}$,
we do need to consider categories whose morphisms commute with the actions of $G$. See, for example, \cite{H-cofinite-P(z)}
for such results.
The category $\mathcal{C}_{1}$ can be the category of grading-restricted
$g$-twisted $V$-modules without $g$-actions for~${g\in G}$ and the category $\mathcal{C}_{2}$
can be the category of grading-restricted
$g$-twisted $V$-modules with actions of $G$.
We shall use $\mathcal{C}$ to denote either $\mathcal{C}_{1}$ or $\mathcal{C}_{2}$.
Since many of the constructions
in the present paper work for any category satisfying suitable conditions,
we shall work with a~general category $\mathcal{C}$.
But the constructions for $\mathcal{C}_{1}$ and $\mathcal{C}_{2}$, especially the construction of
the $G$-crossed commutativity isomorphisms and thus the construction of
the $G$-crossed braiding isomorphisms, are slightly different.

\begin{defn}
Let $g_{1}, g_{2}\in G$,
$W_{1}$, $W_{2}$, $W_{3}$ $g_{1}$-, $g_{2}$-,
$g_{1}g_{2}$-twisted $V$-modules, respectively,
in the category $\mathcal{C}$
and $z\in \C^{\times}$. A~{\it twisted $P(z)$-intertwining map of type
\smash{$\binom{W_{3}}{W_{1}W_{2}}$}} is a~linear
map~${I\colon W_{1}\otimes W_{2}\to \overline{W}_{3}}$
given by
$I(w_{1}\otimes w_{2})=\Y(w_{1}, z)w_{2}$ for
$w_{1}\in W_{1}$ and $w_{2}\in W_{2}$, where
$\Y$ is a~twisted intertwining operator of type~\smash{$\binom{W_{3}}{W_{1}W_{2}}$}.
\end{defn}

Using the notion of twisted $P(z)$-intertwining map, we now define
the notion of tensor product of two twisted modules in $\mathcal{C}$.

\begin{defn}
{\rm Let $W_{1}$ and $W_{2}$ be $g_{1}$- and $g_{2}$-twisted $V$-modules,
respectively, in $\mathcal{C}$. A~{\it $P(z)$-product of $W_{1}$
and $W_{2}$ in $\mathcal{C}$}
is a~pair $(W_{3}, I)$ consisting of a~$g_{1}g_{2}$-twisted $V$-module $W_{3}$
in $\mathcal{C}$ and a~twisted $P(z)$-intertwining map $I$ of type
$\binom{W_{3}}{W_{1}W_{2}}$. A~{\it $P(z)$-tensor product
of $W_{1}$ and $W_{2}$ in~$\mathcal{C}$} is a~$P(z)$-product $(W_{1}\boxtimes_{P(z)}W_{2},
\boxtimes_{P(z)})$ satisfying the following universal property:
For any $P(z)$-product $(W_{3}, I)$ of $W_{1}$ and $W_{2}$,
there exists a~unique module map $f\colon W_{1}\boxtimes_{P(z)}W_{2}\to W_{3}$
such that we have the commutative diagram
\[\begin{tikzcd}
W_{1}\otimes W_{2}\arrow[swap]{d}{\boxtimes_{P(z)}}
\arrow{r}{I}&\overline{W}_{3}\\
\overline{W_{1}\boxtimes_{P(z)}W_{2}}\arrow[swap]{ru}{\bar{f}},&
\end{tikzcd}\]
where $\bar{f}$ is the natural extension of $f$ to
$\overline{W_{1}\boxtimes_{P(z)}W_{2}}$.}
\end{defn}

We now give a~construction of
$(W_{1}\boxtimes_{P(z)}W_{2}, \boxtimes_{P(z)})$
under a~suitable assumption using the same method as in
\cite{tensor1,HLZ3}.

Given a~$P(z)$-product $(W_{3}, I)$ of $W_{1}$ and $W_{2}$ in $\mathcal{C}$,
for $w_{3}'\in W_{3}'$, we have an element $\lambda_{I, w_{3}'}\in
(W_{1}\otimes W_{2})^{*}$ defined by
$\lambda_{I, w_{3}'}(w_{1}\otimes w_{2})
=\langle w_{3}', I(w_{1}\otimes w_{2})\rangle$
for $w_{1}\in W_{1}$ and $w_{2}\in W_{2}$.
Let~${W_{1}\hboxtr_{P(z)}W_{2}}$ be the subspace of
 $(W_{1}\otimes W_{2})^{*}$ spanned by
$\lambda_{I, w_{3}'}$ for all $P(z)$-products $(W_{3}, I)$
and~${w_{3}'\in W_{3}'}$. Note that $W_{1}\hboxtr_{P(z)}W_{2}$ in fact depends on the
category $\mathcal{C}$.
We define a~vertex operator map
\[Y^{(g_{1}g_{2})^{-1}}_{W_{1}\hboxtr_{P(z)}W_{2}}\colon\
V\otimes (W_{1}\hboxtr_{P(z)}W_{2})\to
(W_{1}\hboxtr_{P(z)}W_{2})\{x\}[\log x]\]
by
\[Y^{(g_{1}g_{2})^{-1}}_{W_{1}\hboxtr_{P(z)}W_{2}}(v, x)\lambda_{I, w_{3}'}
=\lambda_{I, Y_{W_{3}'}^{(g_{1}g_{2})^{-1}}
(v, x)w_{3}'}\]
for $v\in V$ and \smash{$\lambda_{I, w_{3}'}\in W_{1}\hboxtr_{P(z)}W_{2}$}.

This vertex operator map is indeed well defined. In fact, assuming that
\smash{$\lambda_{I, w_{3}'}=\lambda_{\tilde{I}, \tilde{w}_{3}'}$} for another $P(z)$-product
\smash{$(\widetilde{W}_{3}, \tilde{I})$} of $W_{1}$ and $W_{2}$ in $\mathcal{C}$ and \smash{$\tilde{w}_{3}'\in \widetilde{W}_{3}'$}.
Then for $v\in V$, $w_{1}\in W_{1}$ and~${w_{2}\in W_{2}}$,
\begin{align*}
\big\langle w_{3}', I\bigl(Y_{W_{1}}^{g_{1}}(v, x)w_{1}\otimes w_{2}\bigr)\big\rangle
&=\lambda_{I, w_{3}'}\bigl(Y_{W_{1}}^{g_{1}}(v, x)w_{1}\otimes w_{2}\bigr)=\lambda_{\tilde{I}, \tilde{w}_{3}'}\bigl(Y_{W_{1}}^{g_{1}}(v, x)w_{1}\otimes w_{2}\bigr)\nn
&=\big\langle \tilde{w}_{3}', \tilde{I}\bigl(Y_{W_{1}}^{g_{1}}(v, x)w_{1}\otimes w_{2}\bigr)\big\rangle.
\end{align*}
Let $\Y$ and $\widetilde{Y}$ be the twisted intertwining operators corresponding to $I$ and $\tilde{I}$.
We also need the opposite twisted vertex operator map \smash{$\bigl(Y_{W}^{g}\bigr)^{o}$} of a~$g$-twisted $V$-module $W$
defined by $\bigl(Y_{W}^{g}\bigr)^{o}(v, x)=\smash{Y_{W}^{g}\bigl({\rm e}^{xL(1)}\bigl(-x^{-2}\bigr)^{-L(0)}v, x^{-1}\bigr)}$ for $v\in V$.
Then on the region $|z_{1}|>|z|>|z_{1}-z|>0$, $|{\arg (z_{1}-z)-\arg z_{1}}|<\frac{\pi}{2}$,
$|{\arg z_{1}-\arg z_{2}}|<\frac{\pi}{2}$, we have
\begin{gather*}
\big\langle \bigl(Y_{W_{3}'}^{(g_{1}g_{2})^{-1}}\bigr)^{o}(v, z_{1})w_{3}', I(w_{1}\otimes w_{2})\big\rangle\nn
 \qquad =\big\langle w_{3}', Y_{W_{3}}^{(g_{1}g_{2})}(v, z)\Y(w_{1}, z)w_{2}\big\rangle =\langle w_{3}', \Y(Y_{W_{1}}(v, z_{1}-z)w_{1}, z)w_{2})\rangle\nn
 \qquad =\big\langle w_{3}', I\bigl(Y_{W_{1}}^{g_{1}}(v, z_{1}-z)w_{1}\otimes w_{2}\bigr)\big\rangle =\big\langle \tilde{w}_{3}', \tilde{I}\bigl(Y_{W_{1}}^{g_{1}}(v, z_{1}-z)w_{1}\otimes w_{2}\bigr)\big\rangle\nn
 \qquad =\langle \tilde{w}_{3}', \widetilde{\Y}(Y_{W_{1}}(v, z_{1}-z)w_{1}, z)w_{2})\rangle =\big\langle \tilde{w}_{3}', Y_{\widetilde{W}_{3}}^{(g_{1}g_{2})}(v, z)\Y(w_{1}, z)w_{2}\big\rangle\nn
 \qquad =\big\langle \bigl(Y_{\widetilde{W}_{3}'}^{(g_{1}g_{2})^{-1}}\bigr)^{o}(v, z_{1})\tilde{w}_{3}', I(w_{1}\otimes w_{2})\big\rangle.
\end{gather*}
Since the both sides are analytic on the region $|z_{1}|>|z|>0$, this equality
holds on this larger region. So we can replace the complex variable $z_{1}$ by a~formal variable $x$ to obtain
\[
\big\langle \bigl(Y_{W_{3}'}^{(g_{1}g_{2})^{-1}}\bigr)^{o}(v, x)w_{3}', I(w_{1}\otimes w_{2})\big\rangle
=\big\langle \bigl(Y_{\widetilde{W}_{3}'}^{(g_{1}g_{2})^{-1}}\bigr)^{o}(v, x)\tilde{w}_{3}', I(w_{1}\otimes w_{2})\big\rangle.\]
Then from the definition of opposite twisted vertex operator map, we obtain
\[
\big\langle Y_{W_{3}'}^{(g_{1}g_{2})^{-1}}(v, x)w_{3}', I(w_{1}\otimes w_{2})\big\rangle
=\big\langle Y_{\widetilde{W}_{3}'}^{(g_{1}g_{2})^{-1}}(v, x)\tilde{w}_{3}', I(w_{1}\otimes w_{2})\big\rangle.\]
Thus we have
\[
\lambda_{I, Y_{W_{3}'}^{(g_{1}g_{2})^{-1}}(v, x)w_{3}'}
=\lambda_{\tilde{I}, Y_{\widetilde{W}_{3}'}^{(g_{1}g_{2})^{-1}}(v, x)\tilde{w}_{3}'},\]
proving that \smash{$Y^{(g_{1}g_{2})^{-1}}_{W_{1}\hboxtr_{P(z)}W_{2}}(v, x)$} is indeed well defined.

\begin{prop}\label{hboxtr}
The pair \smash{$\bigl(W_{1}\hboxtr_{P(z)}W_{2},
Y^{(g_{1}g_{2})^{-1}}_{W_{1}\hboxtr_{P(z)}W_{2}}\bigr)$} is a~generalized
$(g_{1}g_{2})^{-1}$-twisted $V$-module without a~$(g_{1}g_{2})^{-1}$-action.
\end{prop}
\begin{proof}
Note that every element of \smash{$W_{1}\hboxtr_{P(z)}W_{2}$} is a~linear combination
of elements of the form~$\lambda_{I, w_{3}'}$
for a~$g_{1}g_{2}$-twisted $V$-module $W_{3}$
in $\mathcal{C}$, a~$P(z)$-intertwining map $I$ of type
of \smash{$\binom{W_{3}}{W_{1}W_{2}}$} and an element
$w_{3}'\in W_{3}'$. For fixed $W_{3}$ and $I$,
the space spanned by all $\lambda_{I, w_{3}'}$ for $w_{3}'\in W_{3}'$
is the image of $W_{3}'$ under the linear map from
$W_{3}'$ to \smash{$W_{1}\hboxtr_{P(z)}W_{2}$}
given by $w_{3}'\mapsto \lambda_{I, w_{3}'}$. This linear map
preserves the gradings. So the space of
spanned by all \smash{$\lambda_{I, w_{3}'}$} for $w_{3}'\in W_{3}'$
is a~generalized $(g_{1}g_{2})^{-1}$-twisted $V$-module without a~$g$-action.
Thus $W_{1}\hboxtr_{P(z)}W_{2}$ as a~sum of
generalized $(g_{1}g_{2})^{-1}$-twisted $V$-modules without a~$(g_{1}g_{2})^{-1}$-action is also a~generalized $(g_{1}g_{2})^{-1}$-twisted $V$-module without a~$(g_{1}g_{2})^{-1}$-action.
\end{proof}

\begin{assum}\label{assum}
We assume that the following conditions for $\mathcal{C}$
hold:
\begin{enumerate}[label=$(\arabic*)$]\itemsep=0pt
\item For objects $W_{1}$ and $W_{2}$
in $\mathcal{C}$, the contragredient of $W_{1}\hboxtr_{P(z)}W_{2}$
is also in $\mathcal{C}$. (In the case~${\mathcal{C}=\mathcal{C}_{2}}$, this means in particular
that for objects
$W_{1}$ and $W_{2}$ in $\mathcal{C}_{2}$,
there exists an action of $G$ on $W_{1}\hboxtr_{P(z)}W_{2}$ and a~decomposition of
$W_{1}\hboxtr_{P(z)}W_{2}$ as the direct sum of generalized eigenspaces of the $(g_{1}g_{2})$-action
such that the contragredient of $W_{1}\hboxtr_{P(z)}W_{2}$ with this action and this decomposition
is an object in $\mathcal{C}_{2}$.)
\item The double contragredient of an object in $\mathcal{C}$
is equivalent to the object itself.
\end{enumerate}
\end{assum}

See Theorem~\ref{grading-res}
for some conditions implying condition~(1) in Assumption~\ref{assum}.
Condition~(2) in Assumption~\ref{assum} can be replaced by the condition that
all the twisted $V$-modules in $\mathcal{C}$ is grading-restricted. In fact, the proof
in~\cite{FHL} in the untwisted case of the result that the double contragredient of a~grading-restricted
$V$-module is equivalent to the $V$-module itself still works for grading-restricted twisted $V$-modules.
 In the untwisted case, it has been shown by the second author in~\cite{H-C1-vtc} that Assumption~\ref{assum} is enough to construct $P(z)$-tensor product bifunctors.
In the case that $G$ is abelain,
the second author has shown in~\cite{H-cofinite-P(z)} that when $V$ is $C_{2}$-cofinite and
the category is the category of $C_{2}$-cofinite twisted modules, condition~(1) in Assumption~\ref{assum}
is satisfied and condition~(2) is a~consequence.

From Assumption~\ref{assum}\,(1), we see that
$(W_{1}\hboxtr_{P(z)}W_{2})'$ is in $\mathcal{C}$.
We take $W_{1}\boxtimes_{P(z)}W_{2}$
to be $(W_{1}\hboxtr_{P(z)}W_{2})'$. We still need to
give a~twisted $P(z)$-intertwining map $\boxtimes_{P(z)}$
of type~\smash{$\binom{W_{1}\boxtimes_{P(z)}W_{2}}{W_{1}W_{2}}$}
or equivalently,
an intertwining operator of the same type.

Let $W$ be a~$(g_{1}g_{2})^{-1}$-twisted $V$-module in $\mathcal{C}$
and \smash{$f\colon W\to W_{1}\hboxtr_{P(z)}W_{2}$} a~$V$-module map.
For~${w_{1}\in W_{1}}$, $w_{2}\in W_{2}$ and
$w\in W$,
we define an element $\Y_{f}(w_{1}, z)w_{2}\in W^{*}=\overline{W'}$ by
\begin{equation}\label{int-op-z}
\langle w, \Y_{f}(w_{1}, z)w_{2}\rangle
=(f(w))(w_{1}\otimes w_{2}).
\end{equation}
Then we define
\begin{align}
\Y_{f}(w_{1}, x)w_{2}={}&x^{L_{W'}(0)}
{\rm e}^{-(\log z) L_{W_{3}}(0)}
\Y_{f}\bigl(x^{-L_{W_{1}}(0)}{\rm e}^{(\log z) L_{W_{1}}(0)}w_{1}, z\bigr)\nn
&\times
x^{-L_{W_{2}}(0)}{\rm e}^{(\log z) L_{W_{2}}(0)}w_{2}\label{int-op-x}
\end{align}
for $w_{1}\in W_{1}$, $w_{2}\in W_{2}$.
We now have a~linear map
$\Y_{f}\colon W_{1}\otimes W_{2}\to W'\{x\}[\log x]$.

\begin{prop}\label{mod-map-int-op}
The linear map
$\Y_{f}\colon W_{1}\otimes W_{2}\to W'\{x\}[\log x]$
given by~\eqref{int-op-z} and~\eqref{int-op-x}
above
is a~twisted intertwining operator
of type \smash{$\binom{W'}{W_{1}W_{2}}$}. In particular,
in the case that $W=W_{1}\hboxtr_{P(z)}W_{2}$ and
\smash{$f=1_{W_{1}\hboxtr_{P(z)}W_{2}}\colon W\to W_{1}\hboxtr_{P(z)}W_{2}$}
is the identity map, we obtain a~twisted
intertwining operator \smash{$\Y_{1_{W_{1}\hboxtr_{P(z)}W_{2}}}$}
of type \smash{$\binom{W_{1}\boxtimes_{P(z)}W_{2}}{W_{1}W_{2}}$}.
\end{prop}
\begin{proof}
We first verify the $L(-1)$-derivative property,
\begin{align}
\frac{\rm d}{{\rm d}x}\Y_{f}(w_{1}, x)w_{2}={}&\frac{\rm d}{{\rm d}x}x^{L_{W'}(0)}{\rm e}^{-(\log z) L_{W'}(0)}
\Y_{f}\nn
&\times\bigl(x^{-L_{W'_{1}}(0)}{\rm e}^{(\log z) L_{W'_{1}}(0)}w_{1}, z\bigr)
x^{-L_{W_{2}}(0)}{\rm e}^{(\log z) L_{W_{2}}(0)}w_{2}\nn
={}&x^{L_{W'}(0)-1}{\rm e}^{-(\log z) L_{W'}(0)}L_{W'}(0)
\Y_{f}\bigl(x^{-L_{W_{1}}(0)}{\rm e}^{(\log z) L_{W_{1}}(0)}w_{1}, z\bigr)
x^{-L_{W_{2}}(0)}
\nn
&\times {\rm e}^{(\log z) L_{W_{2}}(0)}w_{2}-x^{L_{W'}(0)}{\rm e}^{-(\log z) L_{W'}(0)}\nn
&\times
\Y_{f}\bigl(L_{W_{1}}(0)x^{-L_{W_{1}}(0)-1}{\rm e}^{(\log z)
L_{W_{1}}(0)}w_{1}, z\bigr)
x^{-L_{W_{2}}(0)}{\rm e}^{(\log z) L_{W_{2}}(0)}w_{2} \nn
&\times -x^{L_{W'}(0)}{\rm e}^{-(\log z) L_{W'}(0)}
\Y_{f}\bigl(x^{-L_{W_{1}}(0)}{\rm e}^{(\log z)L_{W_{1}}(0)} w_{1},
 z\bigr)\nn
&\times L_{W_{2}}(0)
x^{-L_{W_{2}}(0)-1}{\rm e}^{(\log z) L_{W_{2}}(0)}w_{2}\nn
={}&x^{-1}x^{L_{W'}(0)}{\rm e}^{-(\log z) L_{W'}(0)}
\bigl(L_{W'}(0)\Y_{f}\bigl(x^{-L_{W_{1}}(0)}{\rm e}^{(\log z)
L_{W_{1}}(0)}w_{1}, z\bigr)
\nn
&
 -\Y_{f}\bigl(L_{W_{1}}(0)x^{-L_{W_{1}}(0)}{\rm e}^{(\log z)
 L_{W_{1}}(0)}w_{1}, z\bigr)\nn
&
-\Y_{f}\bigl(x^{-L_{W_{1}}(0)}{\rm e}^{(\log z)L_{W_{1}}(0)} w_{1},
 z\bigr)L_{W_{2}}(0)\bigr)
x^{-L_{W_{2}}(0)}{\rm e}^{(\log z) L_{W_{2}}(0)}w_{2}.\label{int-op-x-1}
\end{align}

Since $f$ is a~$V$-module map, $f(W)$ is a~submodule of
$W_{1}\hboxtr_{P(z)}W_{2}$.
By the definition of $W_{1}\hboxtr_{P(z)}W_{2}$, it is spanned by
elements of the form $\lambda_{I, w_{3}'}$ for a~$g_{1}g_{2}$-twisted
$V$-module $W_{3}$, a~$P(z)$-intertwining map $I$ of type
\smash{$\binom{W_{3}}{W_{1}W_{2}}$} and $w_{3}'\in W_{3}'$.
In particular, for $w\in W$,
\smash{$f(w)=\sum_{i=3}^{n}\lambda_{I^{i}, w_{i}'}$},
where for $i=3, \dots, n$, $w_{i}'$ is an element of the contragredient
module $W_{i}'$ of a~$g_{1}g_{2}$-twisted
$V$-module $W_{i}$, $I^{i}$ a~twisted $P(z)$-intertwining map of type
$\binom{W_{i}}{W_{1}W_{2}}$.

Let $\Y^{i}$ be the twisted intertwining operator of type
\smash{$\binom{W_{i}}{W_{1}W_{2}}$} such that $I^{i}
=\Y^{i}(\cdot, z)\cdot$.
Then we have
\begin{align*}
\begin{split}
\langle w, \Y_{f}(w_{1}, z)w_{2}\rangle
&=(f(w))(w_{1}\otimes w_{2})=\sum_{i=3}^{n}\lambda_{I^{i}, w_{i}'}(w_{1}\otimes w_{2})=\sum_{i=3}^{n}\langle w_{i}', I^{i}(w_{1}\otimes w_{2})\rangle\nn
&=\sum_{i=3}^{n}\langle w_{i}', \Y^{i}(w_{1}, z)w_{2})\rangle
\end{split}
\end{align*}
for $w_{1}\in W_{1}$ and $w_{2}\in W_{2}$.
Also,
\begin{align*}
(f(L_{W}(0)w))(w_{1}\otimes w_{2})
& =((L_{W_{1}\hboxtr_{P(z)}W_{2}}(0)f(w))(w_{1}\otimes w_{2})\nn
&=\sum_{i=3}^{n}(L_{W_{1}\hboxtr_{P(z)}W_{2}}(0)\lambda_{I^{i}, w_{i}'})
(w_{1}\otimes w_{2})\nn
&=\sum_{i=3}^{n}\res_{x}x\bigl(Y^{(g_{1}g_{2})^{-1}}_{W'}(\omega, x)
\lambda_{I^{i}, w_{i}'}\bigr)(w_{1}\otimes w_{2})\nn
&=\sum_{i=3}^{n}\res_{x}x\bigl(\lambda_{I^{i}, Y^{(g_{1}g_{2})^{-1}}_{W_{i}'}
(\omega, x)w_{i}'}\bigr)(w_{1}\otimes w_{2})\nn
&=\sum_{i=3}^{n}
\lambda_{I^{i}, L_{W_{i}'}(0)w_{i}'}(w_{1}\otimes w_{2})
\end{align*}
for $w_{1}\in W_{1}$ and $w_{2}\in W_{2}$.
So we have
\[
f(L_{W}(0)w)=\sum_{i=3}^{n}\lambda_{I^{i}, L_{W_{i}'}(0)w_{i}'}.
\]
Thus for $w_{1}\in W_{1}$, $w_{2}\in W_{2}$,
\begin{gather*}
\langle w, \left(L_{W'}(0)\Y_{f}(w_{1}, z)-
\Y_{f}(L_{W_{1}}(0)w_{1}, z)-\Y_{f}(w_{1}, z)L_{W_{2}}(0)\right)
w_{2}\rangle\nn
\qquad =\langle L_{W}(0)w, \Y_{f}(w_{1}, z)w_{2}\rangle
-\langle w, \Y_{f}(L_{W_{1}}(0)w_{1}, z)w_{2}\rangle
-\langle w, \Y_{f}(w_{1}, z)L_{W_{2}}(0)
w_{2}\rangle\nn
\qquad =(f(L_{W}(0)w))(w_{1}\otimes w_{2})
-(f(w))(L_{W_{1}}(0)w_{1}\otimes w_{2})
-(f(w))(w_{1}\otimes L_{W_{2}}(0)
w_{2})\nn
\qquad =\sum_{i=3}^{n}
\lambda_{I^{i}, L_{W_{i}'}(0)w_{i}'}(w_{1}\otimes w_{2})
-\sum_{i=3}^{n}
\lambda_{I, w_{i}'}(L_{W_{1}}(0)w_{1}\otimes w_{2})\\
\phantom{\qquad=}{}
 -\sum_{i=3}^{n}
\lambda_{I, w_{i}'}(w_{1}\otimes L_{W_{2}}(0)
w_{2})\nn
\qquad=\sum_{i=3}^{n}
\big\langle L_{W_{i}'}(0)w_{i}', \Y^{i}(w_{1}, z)w_{2}\big\rangle
-\sum_{i=3}^{n}
\big\langle w_{i}', \Y^{i}(L_{W_{1}}(0)w_{1}, z)w_{2}\big\rangle\\
\phantom{\qquad=}{}
-\sum_{i=3}^{n}
\big\langle w_{i}', \Y^{i}(w_{1}, z)L_{W_{2}}(0)
w_{2}\big\rangle\nn
\qquad=\sum_{i=3}^{n}
\big\langle w_{i}', \bigl(L_{W_{i}}(0)\Y^{i}(w_{1}, z)
-\Y^{i}(L_{W_{1}}(0)w_{1}, z)
-\Y^{i}(w_{1}, z)L_{W_{2}}(0)\bigr)
w_{2}\big\rangle\nn
\qquad =\sum_{i=3}^{n}
z\big\langle w_{i}', \Y^{i}(L_{W_{1}}(-1)w_{1}, z)
w_{2}\big\rangle =z\sum_{i=3}^{n}
\lambda_{I^{i}, w_{i}'}((L_{W_{1}}(-1)w_{1}\otimes
w_{2})\nn
\qquad =z(f(w))(L_{W_{1}}(-1)w_{1}\otimes
w_{2})=z\langle w, \Y_{f}(L_{W_{1}}(-1)w_{1}, z)
w_{2}\rangle,
\end{gather*}
where we have used the $L(0)$-commutator formula
for the twisted intertwining operators $\Y^{i}$ (see Remark~\ref{fixed-pt-subalg-int}).
Since $w\in W$ and $w_{2}\in W_{2}$ are arbitrary, we obtain
\begin{equation}\label{int-op-x-2}
L_{W}(0)\Y_{f}(w_{1}, z)-
\Y_{f}(L_{W_{1}}(0)w_{1}, z)-\Y_{f}(w_{1}, z)L_{W_{2}}(0)
=z\Y_{f}(L_{W_{1}}(-1)w_{1}, z)
\end{equation}
for $w_{1}\in W_{1}$.

Using~\eqref{int-op-x-2}, we see that the right-hand side of~\eqref{int-op-x-1}
is equal to
\begin{gather*}
x^{-1}x^{L_{W}(0)}{\rm e}^{-(\log z) L_{W}(0)}
z\Y_{f}\bigl(L_{W_{1}}(-1)x^{-L_{W_{1}}(0)}{\rm e}^{(\log z) L_{W_{1}}(0)}w_{1}, z\bigr)
x^{-L_{W_{2}}(0)}{\rm e}^{(\log z) L_{W_{2}}(0)}w_{2}\nn
 \qquad =x^{L_{W}(0)}{\rm e}^{-(\log z) L_{W}(0)}
\Y_{f}\bigl(x^{-L_{W_{1}}(0)}{\rm e}^{(\log z) L_{W_{1}}(0)}L_{W_{1}}(-1)w_{1}, z\bigr)
x^{-L_{W_{2}}(0)}{\rm e}^{(\log z) L_{W_{2}}(0)}w_{2}\nn
\qquad=\Y_{f}(L_{W_{1}}(-1)w_{1}, x)w_{2},
\end{gather*}
proving the $L(-1)$-derivative property.

For $v\in V$, $w_{1}\in W_{1}$ and $w_{2}\in W_{2}$,
we have
\begin{align*}
\bigl(f(Y_{W}^{(g_{1}g_{2})^{-1}}(v, x)w)\bigr)(w_{1}\otimes w_{2})
&=\bigl(Y_{W_{1}\hboxtr_{P(z)}W_{2}}^{(g_{1}g_{2})^{-1}}(v, x)f(w)\bigr)
(w_{1}\otimes w_{2})\nn
&=\sum_{i=3}^{n}
\bigl(Y_{W_{1}\hboxtr_{P(z)}W_{2}}^{(g_{1}g_{2})^{-1}}(v, x)
\lambda_{I^{i}, w_{i}'}\bigr)(w_{1}\otimes w_{2})\nn
&=\sum_{i=3}^{n}\lambda_{I^{i}, Y_{W_{i}}^{(g_{1}g_{2})^{-1}}
(v, x)w_{i}'}(w_{1}\otimes w_{2}).
\end{align*}
Then we obtain
\[f\bigl(Y_{W}^{(g_{1}g_{2})^{-1}}(v, x)w\bigr)=
\sum_{i=3}^{n}\lambda_{I^{i}, Y_{W_{i}}^{(g_{1}g_{2})^{-1}}
(v, x)w_{i}'}.\]

For $u\in V$, $w_{1}\in W_{1}$,
$w_{2}\in W_{2}$ and $w\in W$, we have
\begin{gather}
\big\langle w, Y_{W'}^{g_{1}g_{2}}(u, z_{1})\Y_{f}(w_{1}, z_{2})w_{2}\big\rangle\nn
 \qquad =\big\langle w, Y_{W'}^{g_{1}g_{2}}(u, z_{1})
{\rm e}^{(\log z_{2})L_{W'}(0)}
{\rm e}^{-(\log z) L_{W'}(0)}\nn
\qquad \phantom{=}{}\times
\Y_{f}\bigl({\rm e}^{-(\log z_{2})L_{W_{1}}(0)}
{\rm e}^{(\log z) L_{W_{1}}(0)}w_{1}, z\bigr)
{\rm e}^{-(\log z_{2})L_{W_{2}}(0)}
{\rm e}^{(\log z) L_{W_{2}}(0)}w_{2}\big\rangle\nn
\qquad =\big\langle {\rm e}^{(\log z_{2})L_{W}(0)}
{\rm e}^{-(\log z) L_{W}(0)}Y_{W}^{(g_{1}g_{2})^{-1}}\bigl({\rm e}^{z_{1}L_{V}(1)}
\bigl(-z_{1}^{-2}\bigr)^{L_{V}(0)}u, z_{1}^{-1}\bigr)w, \nn
\qquad \phantom{=}{}
\Y_{f}\bigl({\rm e}^{-(\log z_{2})L_{W_{1}}(0)}
{\rm e}^{(\log z) L_{W_{1}}(0)}w_{1}, z\bigr)
{\rm e}^{-(\log z_{2})L_{W_{2}}(0)}
{\rm e}^{(\log z) L_{W_{2}}(0)}w_{2}\big\rangle\nn
\qquad =\bigl(f({\rm e}^{(\log z_{2})L_{W}(0)}
{\rm e}^{-(\log z) L_{W}(0)}Y_{W}^{(g_{1}g_{2})^{-1}}({\rm e}^{z_{1}L_{V}(1)}
\bigl(-z_{1}^{-2}\bigr)^{L_{V}(0)}u, z_{1}^{-1})w)\bigr)\nn
\qquad \phantom{=}{}\times
\bigl({\rm e}^{-(\log z_{2})L_{W_{1}}(0)}
{\rm e}^{(\log z) L_{W_{1}}(0)}w_{1}\otimes {\rm e}^{-(\log z_{2})L_{W_{2}}(0)}
{\rm e}^{(\log z) L_{W_{2}}(0)}w_{2}\bigr)\nn
 \qquad =\sum_{i=3}^{n}\lambda_{I^{i}, {\rm e}^{(\log z_{2})L_{W_{i}'}(0)}
{\rm e}^{-(\log z) L_{W_{i}'}(0)}Y_{W_{i}'}^{(g_{1}g_{2})^{-1}}
\bigl({\rm e}^{z_{1}L_{V}(1)}
\bigl(-z_{1}^{-2}\bigr)^{L_{V}(0)}u, z_{1}^{-1}\bigr)w_{i}'}
\nn
\qquad \phantom{=}{}\times
\bigl({\rm e}^{-(\log z_{2})L_{W_{1}}(0)}
{\rm e}^{(\log z) L_{W_{1}}(0)}w_{1}\otimes {\rm e}^{-(\log z_{2})L_{W_{2}}(0)}
{\rm e}^{(\log z) L_{W_{2}}(0)}w_{2}\bigr)\nn
 \qquad =\sum_{i=3}^{n}
\big\langle {\rm e}^{(\log z_{2})L_{W_{i}'}(0)}
{\rm e}^{-(\log z) L_{W_{i}'}(0)}Y_{W_{i}'}^{(g_{1}g_{2})^{-1}}
\bigl({\rm e}^{z_{1}L_{V}(1)}
\bigl(-z_{1}^{-2}\bigr)^{L_{V}(0)}u, z_{1}^{-1}\bigr)w_{i}',
\nn
\qquad \phantom{=}{}
\Y^{i}\bigl({\rm e}^{-(\log z_{2})L_{W_{1}}(0)}
{\rm e}^{(\log z) L_{W_{1}}(0)}w_{1}, z\bigr)
{\rm e}^{-(\log z_{2})L_{W_{2}}(0)}
{\rm e}^{(\log z) L_{W_{2}}(0)}w_{2}\big\rangle\nn
 \qquad =\sum_{i=3}^{n}
\langle w_{i}', Y_{W_{i}}^{g_{1}g_{2}}(u, z_{1})
{\rm e}^{(\log z_{2})L_{W_{i}}(0)}
{\rm e}^{-(\log z) L_{W_{i}}(0)}
\nn
\qquad \phantom{=}{}\times
\Y^{i}\bigl({\rm e}^{-(\log z_{2})L_{W_{1}}(0)}
{\rm e}^{(\log z) L_{W_{1}}(0)}w_{1}, z\bigr)
{\rm e}^{-(\log z_{2})L_{W_{2}}(0)}
{\rm e}^{(\log z) L_{W_{2}}(0)}w_{2}\big\rangle\nn
\qquad =\sum_{i=3}^{n}
\big\langle w_{i}', Y_{W_{i}}^{g_{1}g_{2}}(u, z_{1})
\Y^{i}(w_{1}, z_{2}) w_{2}\big\rangle.\label{int-op-x-3}
\end{gather}
Similarly, we have
\begin{gather}
\big\langle w,
\Y_{f}(w_{1}, z_{2})Y_{W_{2}}^{g_{2}}(u, z_{1})w_{2}\big\rangle
=\sum_{i=3}^{n}\big\langle w_{i}',
\Y^{i}(w_{1}, z_{2})
Y_{W_{2}}^{g_{2}}(u, z_{1})
w_{2}\big\rangle\label{int-op-x-4}
\end{gather}
and
\begin{equation}\label{int-op-x-5}
\big\langle w, \Y_{f}\bigl(Y_{W_{2}}^{g_{1}}(u, z_{1}-z_{2})w_{1}, z_{2}\bigr)w_{2}\big\rangle=\sum_{i=3}^{n}\big\langle w_{i}',
\Y^{i}\bigl(Y_{W_{1}}^{g_{1}}(u, z_{1}-z_{2})w_{1}, z_{2}\bigr)
w_{2}\big\rangle.
\end{equation}
Since $\Y^{i}$ for $i=1, \dots, n$ are
twisted intertwining operators, the duality
property for $\Y$ follows from~\eqref{int-op-x-3},
\eqref{int-op-x-4}, \eqref{int-op-x-5} and the duality properties
for $\Y^{i}$.

The convergence for products of more than two operators
follows from the formula
\begin{gather*}
\big\langle w, Y_{W'}^{g_{1}g_{2}}(u_{1}, z_{1})
\cdots Y_{W'}^{g_{1}g_{2}}(u_{k-1}, z_{k-1})\Y_{f}(w_{1}, z_{k})w_{2}
\big\rangle\nn
\qquad =\sum_{i=3}^{n}
\big\langle w_{i}', Y_{W_{i}}^{g_{1}g_{2}}(u_{1}, z_{1})
\cdots Y_{W_{i}}^{g_{1}g_{2}}(u_{k-1}, z_{k-1})
\Y^{i}(w_{1}, z_{k}) w_{2}\big\rangle,
\end{gather*}
whose proof is the same as that of~\eqref{int-op-x-3}.
\end{proof}

Let \smash{$\boxtimes_{P(z)}=\Y_{1_{W_{1}\hboxtr_{P(z)}W_{2}}}
(\cdot, z)\cdot$}. Then $\boxtimes_{P(z)}$
is a~$P(z)$-intertwining map of type \smash{$\binom{W_{1}\boxtimes_{P(z)}W_{2}}
{W_{1}W_{2}}$}. Let
\smash{$w_{1}\boxtimes_{P(z)}w_{2}=\boxtimes_{P(z)}(w_{1}\otimes w_{2})
=\Y(w_{1}, z)w_{2}\in \overline{W_{1}\boxtimes_{P(z)}W_{2}}$}
for $w_{1}\in W_{1}$ and
$w_{2}\in W_{2}$. We call $w_{1}\boxtimes_{P(z)}w_{2}$ the
{\it tensor product}
of $w_{1}$ and $w_{2}$. By~\eqref{int-op-z}, we have
\begin{equation}\label{int-op-x-6}
\lambda(w_{1}\otimes w_{2})
=\langle \lambda, w_{1}\boxtimes_{P(z)}w_{2}\rangle
\end{equation}
for $\lambda\in W_{1}\hboxtr_{P(z)}W_{2}$, $w_{1}\in W_{1}$ and
$w_{2}\in W_{2}$.

\begin{thm}
The pair \smash{$(W_{1}\boxtimes_{P(z)}W_{2}, \boxtimes_{P(z)})$}
is a~$P(z)$-tensor product of $W_{1}$ and $W_{2}$ in $\mathcal{C}$.
\end{thm}
\begin{proof}
Let $(W_{3}, I)$ be a~$P(z)$-product of $W_{1}$ and $W_{2}$ in $\mathcal{C}$.
Then we have a~module map $g\colon W_{3}'\to W_{1}\hboxtr_{P(z)}W_{2}$
given by $g(w_{3}')=\lambda_{I, w_{3}'}$
for $w_{3}'\in W_{3}'$. By definition,
we have $(g(w_{3}'))(w_{1}\otimes w_{2})
=\lambda_{I, w_{3}'}(w_{1}\otimes w_{2})
=\langle w_{3}', I(w_{1}\otimes w_{2})\rangle$
for $w_{1}\in W_{1}$ and $w_{2}\in W_{2}$.
The adjoint of this module
map is a~module map from $W_{1}\boxtimes_{P(z)}W_{2}$ to $W_{3}''$.
By condition~(2) in Assumption~\ref{assum}, we have a~canonical module map
from $W_{3}''$ to $W_{3}$. The composition of these two module maps
gives a~module map $f\colon W_{1}\boxtimes_{P(z)}W_{2}\to W_{3}$.
 By definitions and~\eqref{int-op-x-6},
\begin{align*}
\big\langle w_{3}', (\bar{f}\circ \boxtimes_{P(z)})(w_{1}\otimes w_{2})\big\rangle
&=\big\langle w_{3}', \bar{f}(w_{1}\boxtimes_{P(z)}w_{2})\big\rangle=\big\langle g(w_{3}'), w_{1}\boxtimes_{P(z)}w_{2}\big\rangle\nn
&=(g(w_{3}'))(w_{1}\otimes w_{2})=\big\langle w_{3}', I(w_{1}\otimes w_{2})\big\rangle.
\end{align*}
So we obtain $\bar{f}\circ \boxtimes_{P(z)}=I$.
\end{proof}

We have assigned each object $(W_{1}, W_{2})$ in the
category $\mathcal{C}\times \mathcal{C}$ an object
$W_{1}\boxtimes_{P(z)}W_{2}$ in~$\mathcal{C}$. To obtain a~functor from
$\mathcal{C}\times \mathcal{C}$ to $\mathcal{C}$,
we still need to assign a~morphism $(f_{1}, f_{2})$
in $\mathcal{C}\times \mathcal{C}$ a~morphism~${f_{1}\boxtimes_{P(z)}f_{2}}$ in $\mathcal{C}$. Recall that in Section~\ref{sec2}, we define
module maps between a~$g_{1}$-twisted $V$-module and a~$g_{2}$-twisted $V$-module
to be $0$ if $g_{1}\ne g_{2}$. So we need only consider
module maps or morphisms in $\mathcal{C}$ from $g$-twisted $V$-modules
to $g$-twisted $V$-modules for $g\in G$.

Let $W_{1}$, $\widetilde{W}_{1}$ be $g_{1}$-twisted
$V$-modules in $\mathcal{C}$ and $W_{2}$, $\widetilde{W}_{2}$
$g_{2}$-twisted
$V$-modules in $\mathcal{C}$. Let
\[
f_{1}\colon \ W_{1}\to
\widetilde{W}_{1} \qquad \text{and}\qquad f_{2}\colon W_{2}\to
\widetilde{W}_{2}
\]
 be module maps.
Let $\widetilde{\Y}$ be the twisted intertwining operator~of type
\smash{$\binom{\widetilde{W}_{1}\boxtimes_{P(z)}\widetilde{W}_{2}}
{\widetilde{W}_{1}\widetilde{W}_{2}}$} such that
\[
\tilde{w}_{1}\boxtimes_{P(z)}\tilde{w}_{2}
=\widetilde{\Y}(\tilde{w}_{1}, z)\tilde{w}_{2}.
\]
Since $f_{1}$ and $f_{2}$ are module maps, $\Y=\widetilde{\Y}
\circ (f_{1}\otimes f_{2})$ is a~twisted intertwining operator of
type \smash{$\binom{\widetilde{W}_{1}\boxtimes_{P(z)}\widetilde{W}_{2}}
{W_{1}W_{2}}$}. Then $I=\Y(\cdot, z)\cdot$ is a~twisted $P(z)$-intertwining operator
of the same type. Hence we have a~$P(z)$-product
\smash{$\bigl(\widetilde{W}_{1}\boxtimes_{P(z)}\widetilde{W}_{2}, I\bigr)$}
of $W_{1}$ and $W_{2}$. By the universal property
of the tensor product \smash{$\bigl(W_{1}\boxtimes_{P(z)}W_{2}, \boxtimes_{P(z)}\bigr)$},
there exists a~unique module map $f\colon W_{1}\boxtimes_{P(z)}W_{2}
\to \widetilde{W}_{1}\boxtimes_{P(z)}\widetilde{W}_{2}$
such that $I=\bar{f}\circ \boxtimes_{P(z)}$. We define this module map
$f$ to be the $P(z)$-tensor product of~$f_{1}$ and~$f_{2}$ and denote it
by $f_{1}\boxtimes_{P(z)}f_{2}$.

\begin{thm}
The assignments given by $(W_{1}, W_{2})\mapsto
W_{1}\boxtimes_{P(z)}W_{2}$ and $(f_{1}, f_{2})\mapsto
f_{1}\boxtimes_{P(z)}f_{2}$ above is a~functor from
$\mathcal{C}\times \mathcal{C}$
to $\mathcal{C}$.
\end{thm}
\begin{proof}
It is easy to verify $1_{w_{1}}\boxtimes_{P(z)}1_{W_{2}}
=1_{W_{1}\boxtimes_{P(z)}W_{2}}$ and
$(f_{1}\boxtimes_{P(z)}f_{2})\circ (g_{1}\boxtimes_{P(z)}g_{2})
=(f_{1}g_{1})\boxtimes_{P(z)}(f_{2}g_{2})$ by using the
construction of the tensor products of module maps.
We omit the details of the proofs.
\end{proof}

 We call this functor the {\it $P(z)$-tensor product
bifunctor}.

We now construct $G$-crossed commutativity isomorphisms.
We assume that $\mathcal{C}_{1}$ is closed under $\phi_{h}$
for $h\in G$ and
$\mathcal{C}_{2}$ is closed under $\varphi_{h}$ for $h\in G$, that is,
for a~$g$-twisted module $W$ in $\mathcal{C}_{1}$ (or $\mathcal{C}_{2}$),
$\phi_{h}(W)$ (or $\varphi_{h}(W)$) for $h\in G$ is also in
$\mathcal{C}_{1}$ (or $\mathcal{C}_{2}$).
Note that the category of grading-restricted $g$-twisted $V$-modules without $g$-actions
(or the category of grading-restricted $g$-twisted $V$-modules with actions of $G$)
is indeed closed under $\phi_{h}$ (or $\varphi_{h}$) for $h\in G$.

Let $W_{1}$ and $W_{2}$ be objects of $\mathcal{C}$. Let
$\Y$ be the twisted intertwining operator of type
\smash{$\binom{\phi_{g_{1}}(W_{2})\boxtimes_{P(-z)}W_{1}}
{\phi_{g_{1}}(W_{2})W_{1}}$} when $\mathcal{C}=\mathcal{C}_{1}$
or the twisted intertwining operator of type
\smash{$\binom{\varphi_{g_{1}}(W_{2})\boxtimes_{P(-z)}W_{1}}
{\varphi_{g_{1}}(W_{2})W_{1}}$} when $\mathcal{C}=\mathcal{C}_{2}$
such that
\[
w_{2}\boxtimes_{P(-z)}w_{1}=\Y(w_{2}, -z)w_{1}
\]
for $w_{1}\in W_{1}$, $w_{2}\in \phi_{g_{1}}(W_{2})$ or
$w_{2}\in \varphi_{g_{1}}(W_{2})$, respectively.
Then by Theorem~\ref{skew-sym},
$\Omega_{+}(\Y)$ or \smash{$\Omega_{+}^{g_{1}^{-1}}(\Y)$} is a~twisted intertwining operator of
type
\[\binom{\phi_{g_{1}}(W_{2})\boxtimes_{P(-z)}W_{1}}
{W_{1} \phi_{g_{1}}^{-1}(\phi_{g_{1}}(W_{2}))}
=\binom{\phi_{g_{1}}(W_{2})\boxtimes_{P(-z)}W_{1}}
{W_{1} W_{2}},
\]
when $\mathcal{C}=\mathcal{C}_{1}$ or
\[\binom{\varphi_{g_{1}}(W_{2})\boxtimes_{P(-z)}W_{1}}
{W_{1} \varphi_{g_{1}}^{-1}(\varphi_{g_{1}}(W_{2}))}
=\binom{\varphi_{g_{1}}(W_{2})\boxtimes_{P(-z)}W_{1}}
{W_{1} W_{2}}\]
when $\mathcal{C}=\mathcal{C}_{2}$.
In particular, the pair \smash{$(\phi_{g_{1}}(W_{2})\boxtimes_{P(-z)}W_{1},
\Omega_{+}(\Y)(\cdot, z)\cdot)$} or $(\varphi_{g_{1}}(W_{2})\boxtimes_{P(-z)}W_{1},\smash{
\Omega_{+}^{g_{1}^{-1}}(\Y)(\cdot, z)\cdot)}$
is a~$P(z)$-product of $W_{1}$ and $W_{2}$ in $\mathcal{C}_{1}$ or $\mathcal{C}_{2}$,
respectively. By the universal property
of the tensor product $W_{1}\boxtimes_{P(z)}W_{2}$, there exists a~unique $g_{1}g_{2}$-twisted $V$-module map
\[\mathcal{R}_{P(z)}\colon\ W_{1}\boxtimes_{P(z)}W_{2}
\to \phi_{g_{1}}(W_{2})\boxtimes_{P(-z)}W_{1}\]
or
\[\mathcal{R}_{P(z)}\colon\ W_{1}\boxtimes_{P(z)}W_{2}
\to \varphi_{g_{1}}(W_{2})\boxtimes_{P(-z)}W_{1}\]
such that
\smash{$\Omega_{+}(\Y)(\cdot, z)\cdot=\overline{\mathcal{R}}_{P(z)}
\circ \boxtimes_{P(z)}$}
or
\smash{$\Omega_{+}^{g_{1}^{-1}}(\Y)(\cdot, z)\cdot=\overline{\mathcal{R}}_{P(z)}
\circ \boxtimes_{P(z)}$},
where $\overline{\mathcal{R}}_{P(z)}$ is the
natural extension of $\mathcal{R}_{P(z)}$. The
$g_{1}g_{2}$-twisted $V$-module map $\mathcal{R}_{P(z)}$
has an inverse
\[\mathcal{R}_{P(z)}^{-1}\colon\ \phi_{g_{1}}(W_{2})\boxtimes_{P(-z)}
W_{1}\to W_{1}\boxtimes_{P(z)}W_{2}\]
or
\[\mathcal{R}_{P(z)}^{-1}\colon\ \varphi_{g_{1}}(W_{2})\boxtimes_{P(-z)}
W_{1}\to W_{1}\boxtimes_{P(z)}W_{2}\]
constructed in the same way as above except that we use $\Omega_{-}$ or
$\Omega_{-}^{g_{1}}$
instead of $\Omega_{+}$ or \smash{$\Omega_{+}^{g_{1}^{-1}}$}, respectively. Then we obtain a~natural isomorphism
$\mathcal{R}_{P(z)}$ called the {\it $G$-crossed commutativity
isomorphism}.

As in~\cite{H-rigidity,HLZ3}, we also have parallel transport isomorphisms.
Let $z_{1}, z_{2}\in \C^{\times}$ and $\gamma$ a~path in~$\C^{\times}$
from $z_{1}$ to $z_{2}$. We denote the homotopy class of
$\gamma$ by $[\gamma]$. For the same $W_{1}$ and $W_{2}$,
let $\Y$ be the twisted intertwining operator of type
\smash{$\binom{W_{1}\boxtimes_{P(z_{2})}W_{2}}
{W_{1}W_{2}}$}
such that $w_{1}\boxtimes_{P(z_{2})}w_{2}=\Y(w_{1}, z_{2})w_{2}$
for~${w_{1}\in W_{1}}$, $w_{2}\in W_{2}$. Let
$l_{p}(z_{1})=\log |z_{1}|+{\rm i} \arg z_{1}+{\rm i} 2\pi p$
be the value of logarithm of $z_{1}$ obtained by analytically
extending the value $\log z_{2}$ along the path $\gamma$.
Then
\[
I=\Y(w_{1}, x)w_{2}\mbar_{x^{n}
={\rm e}^{nl_{p}(z_{1})}, \log x=l_{p}(z_{1})}
\]
is a~$P(z_{1})$-intertwining map
of type \smash{$\binom{W_{1}\boxtimes_{P(z_{2})} W_{2}}{W_{1}W_{2}}$}
and we have a~$P(z_{1})$-product
$(W_{1}\boxtimes_{P(z_{2})}W_{2}, I)$
of $W_{1}$ and $W_{2}$.
By the universal property of the $P(z_{1})$-tensor
product $W_{1}\boxtimes_{P(z_{1})} W_{2}$, there exists a~unique
$g_{1}g_{2}$-twisted $V$-module map
\[
\mathcal{T}_{[\gamma]}\colon\ W_{1}\boxtimes_{P(z_{1})} W_{2}\to
W_{1}\boxtimes_{P(z_{2})} W_{2}\]
such that $\overline{\mathcal{T}_{[\gamma]}}\circ \boxtimes_{P(z_{1})}
=I$. The $g_{1}g_{2}$-twisted $V$-module
map $\mathcal{T}_{[\gamma]}$
is invertible since the same construction
also gives a~$g_{1}g_{2}$-twisted $V$-module map
\[
\mathcal{T}_{[\gamma^{-1}]}\colon\ W_{1}\boxtimes_{P(z_{2})} W_{2}\to
W_{1}\boxtimes_{P(z_{1})} W_{2},
\]
which is clearly the inverse of $\mathcal{T}_{[\gamma]}$.
Thus the natural transformation $\mathcal{T}_{[\gamma]}$ is a~natural isomorphism called the {\it parallel transport isomorphism
from $z_{1}$ to $z_{2}$ along $[\gamma]$}.

Let $\gamma$ be a~path from $-1$ to $1$ in the closed
upper half plane with $0$ deleted. For the same $W_{1}$ and $W_{2}$,
we define the {\it $G$-crossed braiding isomorphism}
$\mathcal{R}\colon W_{1}\boxtimes_{P(1)}W_{2} \to
\phi_{g_{1}}(W_{2})\boxtimes_{P(1)} W_{1}$ or~${\mathcal{R}\colon W_{1}\boxtimes_{P(1)}W_{2} \to
\varphi_{g_{1}}(W_{2})\boxtimes_{P(1)} W_{1}}$
by
\smash{$\mathcal{R}=\mathcal{T}_{[\gamma]}\circ \mathcal{R}_{P(1)}$}.

We now give a~result on condition~(1) in Assumption~\ref{assum}.

\begin{thm}\label{grading-res}
Let $\mathcal{C}_{1}$ be the category of grading-restricted
$g$-twisted $V$-modules without $g$-actions for $g\in G$ and $\mathcal{C}_{2}$
the category of grading-restricted
$g$-twisted $V$-modules with actions of $G$.
Let~$\mathcal{C}$ be either $\mathcal{C}_{1}$ or $\mathcal{C}_{2}$,
Assume that the following
conditions are satisfied:
\begin{enumerate}[label=$(\arabic*)$]\itemsep=0pt
\item For $g\in G$,
there are only finitely many inequivalent irreducible
grading-restricted $g$-twisted $V$-modules in $\mathcal{C}$.
\item Every object in $\mathcal{C}$ is a~direct sum of irreducible objects in $\mathcal{C}$.

\item For $g_{1}, g_{2}\in G$ and $g_{1}$-, $g_{2}$-,
$g_{1}g_{2}$-twisted $V$-modules
$W_{1}$, $W_{2}$, $W_{3}$ in $\mathcal{C}$,
the fusion rule $N_{W_{1}W_{2}}^{W_{3}}
=\dim \mathcal{V}_{W_{1}W_{2}}^{W_{3}}$ is
finite.
\end{enumerate}
Then for objects $W_{1}$ and $W_{2}$ in $\mathcal{C}$, $W_{1}\hboxtr_{P(z)}W_{2}$
is in $\mathcal{C}$.
\end{thm}
\begin{proof}
Let $W_{1}$ and $W_{2}$ be $g_{1}$- and $g_{2}$-twisted $V$-modules, respectively,
in $\mathcal{C}$.
From the construction of $W_{1}\hboxtr_{P(z)}W_{2}$,
it is a~sum of objects in $\mathcal{C}$.
By condition~(2), $W_{1}\hboxtr_{P(z)}W_{2}$ must be a~direct
sum of irreducible
objects in $\mathcal{C}$.
But by condition~(1), there are only finitely many
irreducible
grading-restricted $(g_{1}g_{2})^{-1}$-twisted $V$-modules in $\mathcal{C}$.
If $W_{1}\hboxtr_{P(z)}W_{2}$ is an infinite direct sum of
irreducible
grading-restricted $(g_{1}g_{2})^{-1}$-twisted $V$-modules in $\mathcal{C}$,
at least one irreducible
grading-restricted $(g_{1}g_{2})^{-1}$-twisted $V$-module
$W_{3}$ has infinitely many copies in this decomposition of
$W_{1}\hboxtr_{P(z)}W_{2}$. But then we have infinitely many
linearly independent injective $V$-module maps from $W_{3}$
to the $W_{1}\hboxtr_{P(z)}W_{2}$. But by
Proposition~\ref{mod-map-int-op}, these infinite
injective $V$-module maps give linearly independent
twisted intertwining
operator of type \smash{$\binom{W_{3}'}{W_{1}W_{2}}$}.
Thus the fusion rule \smash{$N_{W_{1}W_{2}}^{W_{3}}$} is $\infty$.
By condition~(3), this is a~contradiction.
So $W_{1}\hboxtr_{P(z)}W_{2}$ must be a~finite direct sum of
irreducible
grading-restricted $(g_{1}g_{2})^{-1}$-twisted $V$-modules in $\mathcal{C}$.
In particular, it is grading restricted. In the case $\mathcal{C}=\mathcal{C}_{2}$,
the actions of $G$ on irreducible
grading-restricted $(g_{1}g_{2})^{-1}$-twisted $V$-modules in $\mathcal{C}$
and the decompositions of the finitely many irreducible
grading-restricted $(g_{1}g_{2})^{-1}$-twisted $V$-modules as direct sums of
eigenspaces of the actions of $(g_{1}g_{2})^{-1}$ give an action of $G$ and a~$(g_{1}g_{2})^{-1}$-grading on
$W_{1}\hboxtr_{P(z)}W_{2}$ satisfying the $(g_{1}g_{2})^{-1}$-grading condition.
\end{proof}

\begin{cor}
Under the assumptions in Theorem~{\rm\ref{grading-res}},
the categories $\mathcal{C}_{1}$ and $\mathcal{C}_{2}$ satisfy Assumption~{\rm\ref{assum}}.
\end{cor}
\begin{proof}
Note that $\mathcal{C}$ is closed under contragredient.
By Theorem~\ref{grading-res}, $W_{1}\hboxtr_{P(z)}W_{2}$ is in $\mathcal{C}$.
Thus the contragredient of
$W_{1}\hboxtr_{P(z)}W_{2}$ is also in $\mathcal{C}$.
So condition~(1) in Assumption~\ref{assum} holds.
Condition~(2) clearly holds for
grading-restricted twisted $V$-modules (see the discussions after Assumption~\ref{assum}).
\end{proof}

\section{Compatibility condition and grading-restriction condition}\label{sec5}

In this section, we introduce $P(z)$-compatibility condition and $P(z)$-local
grading-restriction condition and use these conditions to give another
construction of $W_{1}\hboxtr_{P(z)}W_{2}$ for two twisted $V$-modules
$W_{1}$ and $W_{2}$. As in Section 4, we take $\mathcal{C}$ to be either
$\mathcal{C}_{1}$ or $\mathcal{C}_{2}$ satisfying Assumption~\ref{assum}.
 In the untwisted case (the case that
$\mathcal{C}$ is the
category of (untwisted or $1_{V}$-twisted) $V$-modules), these conditions
and this second construction given in~\cite{tensor1,tensor2,HLZ4} play a~crucial role in the proof of the associativity of
intertwining operators and the construction of associativity isomorphisms
in~\cite{tensor4,HLZ6}. It is expected that they will play
the same crucial role in the proof of the associativity of
twisted intertwining operators and the construction of
associativity isomorphisms for the $P(z)$ tensor product bifunctors
on the category $\mathcal{C}$ of twisted $V$-modules.

In the untwisted case, The $P(z)$-compatibility condition
is formulated using a~formula corresponding to the Jacobi identity
for intertwining operators. Even though we can
obtain a~Jacobi identity for the rational coefficients of the expansions
in a~suitable basis of products and iterates of twisted intertwining operators
with twisted vertex operators as in~\cite{H-jacobi-int}, we do not have a~Jacobi identity for the products and iterates of twisted intertwining operators
with twisted vertex operators. Thus
we have to use the
analytic method to formulate the $P(z)$-compatibility condition
and prove the main results. In particular, the formulation and
proofs involving the $P(z)$-compatibility condition are completely
different from those in~\cite{tensor1,tensor2,HLZ4}.

For a~fixed $z\in \C_{\times}$,
we need to study multivalued analytic functions on the
region
\begin{align*}
&M^n(0,z)=\left\{(z_1,\dots,z_{n})\in\C^n\ \middle|\
\begin{aligned}
&z_i\neq0,\, z_i\neq z, \, i=1,\dots,n,\\
&z_i\neq z_j \text{ for } i,j=1,\dots,n,
\text{ and }i\neq j \text{if }n>1
\end{aligned}\right\}
\end{align*}
for $n\in \Z_{+}$
and its subregions
\begin{align*}
&\Omega_{m,k,l}(z)=
\left\{\begin{aligned}&(z_1,\dots,z_{m+k+l})\\
& \qquad \in\C^{m+k+l}
\end{aligned}\ \middle|\
\begin{aligned}
&|z| < |z_m| < \dots < |z_1|  \text{ if }m>0,\\
&0 < |z_{m+k}-z| < \cdots < |z_{m+1}-z| < |z|  \text{ if }k>0,\\
&0 < |z_{m+k+l}| < \cdots < |z_{m+k+1}| < |z|  \text{ if }l>0,\\
&|z_{m+1}-z|+|z_{m+k+1}|<|z|  \text{ if }k>0,\ l>0,\\
&|z_{m+1}-z|+|z_{1}|<|z| \text{ if }m>0,\ k>0
\end{aligned}\right\},\\
&\Omega_{m,k,l}^{(1)}(z)=
\left\{\begin{aligned}&(z_1,\dots,z_{m+k+l})\\
& \qquad \in\Omega_{m,k,l}(z)
\end{aligned}\ \middle| \
\begin{aligned}
&  |{\arg (z_{j}-z)-\arg (z_{j})}|<\frac{\pi}{2},\\
&  \qquad	 j=1,\dots,m,\text{if }m>0,\\
&  |{\arg(z_j)-\arg(z)}|<\frac{\pi}{2},\\
& \qquad j=m+1,\dots,m+k \text{ if }k>0,\\
&  - \frac{3\pi}{2} <\arg(z_j-z)-\arg(z) < -\frac{\pi}{2},\\
& \qquad j=m + k + 1,\dots,m + k + l \text{ if }l>0
\end{aligned}\right\},\\
&\Omega_{m,k,l}^{(2)}(z)=
\left\{\begin{aligned}&(z_1,\dots,z_{m+k+l})\\
& \qquad \in\Omega_{m,k,l}(z)
\end{aligned}\ \middle|
\begin{aligned}
&  |{\arg (z_{j}-z)-\arg (z_{j})}|<\frac{\pi}{2},\\
& \qquad j=1,\dots,m \text{ if }m>0,\\
& |{\arg(z_j)-\arg(z)}|<\frac{\pi}{2},\\
& \qquad j=m+1,\dots,m+k \text{ if }k>0,\\
&  \frac{\pi}{2}<\arg(z_j-z)-\arg(z)< \frac{3\pi}{2},\\
& \qquad j=m + k + 1,\dots,m + k + l  \text{ if }l>0
\end{aligned}\right\}
\end{align*}
for $m,k,l\in\N$.
Also, we define $M^n_0(0,z)$ to be the simply-connected region given by
cutting~$M^n(0,z)$ along the positive real lines in the $z_{i}$-planes,
$z_{i}-z_{j}$-planes and
$z_{i}-z$-planes, that is, the sets
\begin{align*}
&\{(z_{1},\dots,z_{n})\in M^n(0,z)\mid z_{i}\in \R_{+}\}, \qquad i=1,\dots,n,\\
&\{(z_{1},\dots,z_{n})\in M^n(0,z)\mid z_{i}-z_{j}\in \R_{+}\}, \qquad
i, j=1,\dots,n, i\ne j, \\
&\{(z_{1},\dots,z_{n})\in M^n(0,z)\mid z_{i}-z\in \R_{+}\}, \qquad i=1,\dots,n,
\end{align*}
with these sets attached to the upper half $z_{i}$-planes,
$z_{i}-z_{j}$-planes and
$z_{i}-z$-planes.

To formulate the $P(z)$-compatibility condition, we need the notion of
regular component-isolated singularity of a~multivalued analytic function of several variables
in Section \ref{sec2}.
Let $b=(b_1,\dots,b_n)\in
(\C\cup\{\infty\})^{n}$
and $r=(r_1,\dots,r_n)\in \R_{+}^{n}$. Let $I\subset\{1,\dots,n\}$.
We use the following notation for polydisks and polycircular domains
\begin{align*}	
&\Delta^I(b,r)=\left\{(z_1,\dots,z_n)\in(\C\cup \{\infty\})^n\ \middle| \
\begin{aligned}
& |z_i|>r_i \text{ or } z_{i}=\infty
\text{ if }b_i=\infty, \\
& \qquad |z_i-b_i|<r_i \text{ if }b_i\in\C,
\ \text{for }i\in I,\\
&  |z_j|>r_j\text{ if }b_j=\infty, \\
& \qquad 0<|z_j-b_j|<r_j\text{ if }b_j\in\C
 \ \text{for }j\notin I
\end{aligned}\right\},\\
&\Delta(b,r)=\Delta^{\{1,\dots,n\}}(b,r),\qquad\Delta^\times(b,r)=\Delta^\varnothing(b,r).
\end{align*}

We shall use $b_{z_{i}, z}$ to denote the homotopy class
of loops in $M^{l}(0, z)$ with $z_{i}$ going around $z$ counterclockwise
once with $z_{j}$ for $j\ne i$ fixed.

Let $g_{1}$ and $g_{2}$ be automorphisms of $V$,
$W_{1}$, $W_{2}$, $W_{3}$, $g_{1}$-, $g_{2}$-,
$g_{1}g_{2}$-twisted $V$-modules, respectively,
in the category $\mathcal{C}$
and $z\in \C^{\times}$. Let $I$ be a~twisted $P(z)$-intertwining
map of type $\binom{W_{3}}{W_{1}W_{2}}$ and $w_{3}'\in W_{3}'$.
Then we have an element
$\lambda_{I, w_{3}'}\in W_{1}\hboxtr_{P(z)}W_{2}$.

\begin{prop}\label{int-comp}
The element $\lambda_{I, w_{3}'}$ has the following property:
For $l\in\N$, $u_1,\dots,u_l\in V$, $w_1\in W_1$, and $w_2\in W_2$,
there exists a~multivalued analytic function
\begin{equation}\label{int-comp-cond-1}
f_{l}(z_1,\dots,z_l;u_1,\dots,u_l,w_{1}, w_{2};\lambda_{I, w_{3}'})
\end{equation}
on $M^l(0,z)$ with a~preferred branch
\begin{equation}\label{int-comp-cond-2}
f^e_{l}(z_1,\dots,z_l;u_1,\dots,u_l,w_{1}, w_{2};\lambda_{I, w_{3}'})
\end{equation}
on $M_0^l(0,z)$, satisfying the following:
\begin{enumerate}[label=$(\arabic*)$]\itemsep=0pt
\item \begin{enumerate}\itemsep=0pt
\item[$(a)$] For $i, j=1, \dots, l$, $i\ne j$,
$z_i-z_j=0$ are poles of~\eqref{int-comp-cond-1}.
In particular, there exist~${M_{ij}\in\Z_{+}}$ such that
\begin{align}
\left(\prod_{1\leq i<j\leq l}(z_i-z_j)^{M_{ij}}\right)
f_{l}(z_1,\dots,z_l;u_1,\dots,u_l,w_{1}, w_{2};
\lambda_{I, w_{3}'})\label{int-comp-cond-3}
\end{align}
can be analytically extended to a~multivalued analytic function on
\[\{(z_1,\dots,z_l)\in\C^l \mid z_i\neq0,\, z_i\neq z,\, i=1,\dots,l\}\]
	
\item[$(b)$] Any component-isolated singularity of the form
$(z_1,\dots, z_l)-\beta=\delta$ for $\beta\in \C^{l-1}$
and~${\delta\in\{0,\infty\}^{l-1}}$
of~\eqref{int-comp-cond-3} is
 regular.
\end{enumerate}
\item For $u_1, \dots, u_{l}\in V$, $w_1\in W_1$, and
$w_2\in W_2$,
\begin{enumerate}\itemsep=0pt
\item[$(a)$] The series
\begin{gather*}
\lambda_{I, w_{3}'}(Y^{g_1}(u_1,z_1-z)w_1\otimes
w_2)\nn
\qquad=\lambda_{I, w_{3}'}
(Y^{g_1}(u_1,x)w_1\otimes
w_2)\big|_{x^n={\rm e}^{nl_0(z_1-z)},
\log x=l_0(z_1-z)}
\end{gather*}
is absolutely convergent on the region $\Omega_{0,1,0}(z)$.
Moreover, it is convergent to
$f^e_{1}(z_1;u_1,\allowbreak w_{1}, w_{2};\lambda_{I, w_{3}'})$
on the region
$\Omega^{(1)}_{0,1,0}(z)=\Omega^{(2)}_{0,1,0}(z)$.
	
\item[$(b)$] For $l\in\N$, the multiple series
\begin{gather*}
\lambda_{I, w_{3}'}(w_1\otimes Y^{g_2}(u_1,z_1)\cdots
Y^{g_2}(u_l,z_l)w_2)\nn
 \qquad =\lambda_{I, w_{3}'}(w_1\otimes
Y^{g_2}(u_1,x_1)\cdots
Y^{g_2}(u_l,x_l)w_2)|_{x_i^n={\rm e}^{nl_0(z_1)},
 \log x_i=l_0(z_1), i=1,\dots,l}
\end{gather*}
in powers and logarithms of $z_1,\dots,z_l$
is absolutely convergent on the region $\Omega_{0,0,l}(z)$.
Moreover, it is absolutely convergent to
\eqref{int-comp-cond-2}
on the region \smash{$\Omega^{(1)}_{0,0,l}(z)$} and
to
\[f^{b_{z_{1},z}^{-1}\cdots
b_{z_{l},z}^{-1}}_{l}(z_1,\dots,z_l;u_1,\dots,u_l, w_{1},
w_{2};\lambda)\]
on the region \smash{$\Omega^{(2)}_{0,0,l}(z)$}.
		
\end{enumerate}
		
\end{enumerate}
\end{prop}
\begin{proof}
This result can be easily verified by using the definitions of
 $\lambda_{I, w_{3}'}$ and $P(z)$-intertwining maps and
 the properties of twisted intertwining operators. We omit the details. But see Remarks~\ref{remark-M-i-j},
\ref{Omega-regions}, and~\ref{remark_fundamentalgroup} for some subtle points.
\end{proof}

Let $g_{1}$ and $g_{2}$ be automorphisms of $V$,
$W_{1}$, $W_{2}$ $g_{1}$-, $g_{2}$-twisted $V$-modules, respectively,
in the category $\mathcal{C}$
and $z\in \C^{\times}$.
Motivated by Proposition~\ref{int-comp},
we formulate the following condition for
$\lambda\in (W_{1}\otimes W_{2})^{*}$:
\begin{description}\itemsep=0pt
	\item[$\boldsymbol{P(z)}$-compatibility condition.]
 An~element $\lambda\in(W_1\otimes W_2)^*$ is said to be
	$P(z)$\textit{-compatible} if for $l\in\N$, $u_1,\dots,u_l\in V$,
	$w_1\in W_1$, and $w_2\in W_2$, there exists a~multivalued analytic
	function%
	\begin{equation}\label{comp-cond-1}
	f_{l}(z_1,\dots,z_l;u_1,\dots,u_l,w_{1}, w_{2};\lambda)
	\end{equation}
	on $M^l(0,z)$ with a~preferred branch
	\begin{equation}\label{comp-cond-2}
	f^e_{l}(z_1,\dots,z_l;u_1,\dots,u_l,w_{1}, w_{2};\lambda)
	\end{equation}
	on $M_0^l(0,z)$, satisfying the following:
	\begin{enumerate}[label=$(\arabic*)$]\itemsep=0pt
\item
		\begin{enumerate}\itemsep=0pt
\item[$(a)$] For $i, j=1, \dots, l$, $i\ne j$, $z_{1}-z_{j}=0$
		are poles of~\eqref{comp-cond-1}. In particular,
		there exists~${M_{ij}\in\Z_{+}}$ such that
		\begin{align}\label{123}
		\prod_{1\leq i<j\leq n}(z_i-z_j)^{M_{ij}}
		f_{l}(z_1,\dots,z_n;u_1,\dots,u_n,w_{1}, w_{2};\lambda)
		\end{align}
		can be analytically extended to a~multivalued
		analytic function on
		\[\big\{(z_1,\dots,z_l)\in\C^l \mid z_i\neq0,\, z_i\neq z,\, i=1,\dots,l\big\}.\]

		\item[$(b)$] Any component-isolated singularity of the form
$(z_1,\dots, z_l)-\beta=\delta$ for $\beta\in \C^{l-1}$
and $\delta\in\{0,\infty\}^{l-1}$
of~\eqref{123} is
 regular.
		\end{enumerate}
\item For $u_1, \dots, u_{l}\in V$, $w_1\in W_1$, and
		$w_2\in W_2$,
\begin{enumerate}\itemsep=0pt
		\item[$(a)$] The series
		\begin{gather*}
		\lambda(Y^{g_1}(u_1,z_1-z)w_1\otimes w_2)\\
\qquad
		=\lambda(Y^{g_1}(u_1,x)w_1
		\otimes w_2)|_{x^n={\rm e}^{nl_0(z_1-z)},
		 \log x=l_0(z_1-z)}
		\end{gather*}
		is absolutely convergent on the region $\Omega_{0,1,0}(z)$.
		Moreover, it is convergent to~${f^e_{1}(z_1;u_1, w_{1},
		w_{2};\lambda)}$ on the region
		\smash{$\Omega^{(1)}_{0,1,0}(z)=\Omega^{(2)}_{0,1,0}(z)$}.
		
		\item[$(b)$] For $l\in\N$, the multiple series
		\begin{gather}
		\lambda(w_1\otimes Y^{g_2}(u_1,z_1)\cdots
		Y^{g_2}(u_l,z_l)w_2)\label{comp01}\\
		 \qquad =\lambda(w_1\otimes Y^{g_2}(u_1,x_1) \cdots
		Y^{g_2}(u_l,x_l)w_2)|_{x_i^n={\rm e}^{nl_0(z_1)},
		 \log x_i=l_0(z_1), i=1,\dots,l}\nonumber
		\end{gather}
		in powers and logarithms of $z_1,\dots,z_l$ is
		absolutely convergent on the region $\Omega_{0,0,l}(z)$. Moreover, it is absolutely convergent
		to~\eqref{comp-cond-2} on the region
		\smash{$\Omega^{(1)}_{0,0,l}(z)$} and
		to
		\[f^{b_{z_{1},z}^{-1}\cdots
		b_{z_{l},z}^{-1}}_{l}(z_1,\dots,z_l;u_1,\dots,u_l, w_{1},
		w_{2};\lambda)\]
		on the region \smash{$\Omega^{(2)}_{0,0,l}(z)$}.
	\end{enumerate}
	\end{enumerate}
	
	We denote the subspace of $P(z)$-compatible functionals in
	$(W_1\otimes W_2)^*$ as
	\[\operatorname{COMP}_{P(z)}((W_1\otimes W_2)^*)\]
	or $\operatorname{COMP}$ for short.
\end{description}

\begin{rema}\label{remark-M-i-j}
 Note that by 2\,(b) of Proposition~\ref{int-comp} or
2\,(b) of the $P(z)$-compatibility condition and the associativity of twisted vertex operators, $M_{ij}$
in Proposition~\ref{int-comp} or in the $P(z)$-compatibility condition
depend only on $u_{i}$ and $u_{j}$.
\end{rema}

\begin{rema}\label{Omega-regions}
 In 2\,(b) of Proposition~\ref{int-comp} and
2\,(b) of the $P(z)$-compatibility condition, the reason that we involve two
different regions \smash{$\Omega_{0,0,l}^{(1)}(z)$} and
\smash{$\Omega_{0,0,l}^{(2)}(z)$} is because either of these
two regions could be empty. Notice that here $z$ is a~\textit{fixed} nonzero complex number. Actually, when~${\arg z\in[0,\pi/2]}$, the region \smash{$\Omega_{0,0,l}^{(1)}(z)$} is
empty. When $\arg z\in[3\pi/2,2\pi)$, the region
\smash{$\Omega_{0,0,l}^{(2)}(z)$} is empty.	When
\smash{$\Omega_{0,0,l}^{(1)}(z)$} and \smash{$\Omega_{0,0,l}^{(2)}(z)$}
are both nonempty, the absolute convergences of~\eqref{comp01}
on these two regions are equivalent.
\end{rema}
	
\begin{rema}\label{remark_fundamentalgroup}
Because of the definition of the domain of~\eqref{comp-cond-2}, i.e., $M_0^l(0,z)$, and the fact that its singularities at $z_i=z_j$ for
$i\neq j$ are poles, branches of~\eqref{comp-cond-1}
can be indexed by elements in the fundamental group of the space
\[\{(z_1,\dots,z_n)\in\C^n \mid z_i\neq0,\, z_i\neq z,\, i=1,\dots,n\}
=\prod_{i=1}^{n}\{z_i\in\C \mid z_i\neq0,\, z_i\neq z\}.\] A~set of generators of this fundamental group can be chosen
to be $b_{z_{i},0}$, $b_{z_{i},z}$, $i=1,\dots,n$.
For each $i$, the elements
$b_{z_i,0}$ and $b_{z_i,z}$
correspond to $b_{13}$ and $b_{12}$
defined in Section~\ref{sec2}, and they freely generate
$\pi_1(\{z_i\in\C \mid z_i\neq0, z_i\neq z\})$.
Notice that
\[
\pi_{1}\left(\prod_{i=1}^{n}\{z_i\in\C \mid z_i\neq0,\, z_i\neq z\}\right)
=\prod_{i=1}^{n}\pi_{1}\bigl(\{z_i\in\C\mid z_i\neq0,\, z_i\neq z\}\bigr)
=\prod_{i=1}^{n}\langle b_{i,0},b_{i,z}\rangle.
\]
\end{rema}

For $\lambda\in\operatorname{COMP}$, we want to define \smash{$Y^{(g_{1}g_{2})^{-1}}_{P(z)}(u, x)\lambda
\in (W_{1}\otimes W_{2})^{*}
\{x\}[\log x]$}. We first define
\smash{$Y^{(g_{1}g_{2})^{-1}}_{P(z)}(u, x)\lambda
\in (W_{1}\otimes W_{2})^{*}
\{x\}[\log x]$} for $\lambda$ in a~larger subspace of
$(W_{1}\otimes W_{2})^{*}$ than $\operatorname{COMP}$.
Let $\operatorname{COM}_{P(z)}((W_{1}\otimes W_{2})^{*})$
 or simply
$\operatorname{COM}$ be the subspace of $ (W_{1}\otimes W_{2})^{*}$
consisting of $\lambda$ satisfying 1\,(a), (b) and 2\,(b) in the $P(z)$-compatibility condition.
By definition, $\operatorname{COMP}\subset \operatorname{COM}$.
To define \smash{$Y^{(g_{1}g_{2})^{-1}}_{P(z)}(u, x)\lambda\in
(W_{1}\otimes W_{2})^{*}
\{x\}[\log x]$} for $u\in V$ and
$\lambda\in \operatorname{COM}$
is equivalent to define
\[Y^{(g_{1}g_{2})^{-1}}_{P(z)}
\bigl({\rm e}^{xL(1)}\bigl(-x^{2}\bigr)^{-L(0)}u, x^{-1}\bigr)\lambda
=\bigl(Y^{(g_{1}g_{2})^{-1}}_{P(z)}\bigr)^{o}(u, x)\lambda
\in (W_{1}\otimes W_{2})^{*}
\{x\}[\log x].\]

As $z_1=\infty$ is a~regular singular point of
 $f_{1}(z_1;u_1, w_{1}, w_{2};\lambda)$, we know that there
exist unique $a_{i,n,j}(u_1,w_1,w_2;\lambda)\in\C$ and $r_i\in\C$,
for $i,j=0,\dots,K$ and $n\in\N$, such that on
\smash{$\Omega_{1,0,0}^{(1)}(z)=\Omega_{1,0,0}^{(2)}(z)$}
(i.e., the region given by $|z_1|>|z|$, $|{\arg(z_1-z)-\arg(z_1)}|
<\frac{\pi}{2}$),
\[
f^e_1(z_1;u, w_{1}, w_{2};\lambda)
=\sum_{i,j=0}^{K}\sum_{n\in\N}
a_{i,n,j}(u,w_1,w_2;\lambda)z_1^{r_i-n}(\log z_1)^j.
\]
For $i, j=0, \dots, K$, $n\in \N$ and $u\in V$, we define
\smash{$\bigl(Y^{(g_{1}g_{2})^{-1}}_{P(z)}\bigr)^{o}_{-r_{i}+n-1, j}(u)\lambda\in
(W_{1}\otimes W_{2})^{*}$} by
\[
\bigl(\bigl(Y^{(g_{1}g_{2})^{-1}}_{P(z)}\bigr)^{o}_{-r_{i}+n-1, k}(u)
\lambda\bigr)(w_{1}\otimes w_{2})
=a_{i,n,j}(u,w_1,w_2;\lambda)
\]
for $w_{1}\in W_{1}$ and $w_{2}\in W_{2}$. Then we define
\smash{$Y^{(g_{1}g_{2})^{-1}}_{P(z)}
\bigl({\rm e}^{xL(1)}\bigl(-x^{2}\bigr)^{L(0)}u, x^{-1}\bigr)\lambda$} to be
\[
\sum_{i,j=0}^{N}\sum_{n\in\N}
\bigl(Y^{(g_{1}g_{2})^{-1}}_{P(z)}\bigr)^{o}_{-r_{i}+n-1, j}(u)
\lambda x^{r_{i}-n}(\log x)^{j}\in (W_{1}\otimes W_{2})^{*}
\{x\}[\log x],\]
that is,
\begin{gather}
\bigl(Y^{(g_{1}g_{2})^{-1}}_{P(z)}
\bigl({\rm e}^{xL(1)}\bigl(-x^{2}\bigr)^{-L(0)}u, x^{-1}\bigr)\lambda\bigr)
(w_{1}\otimes w_{2})\nn
\qquad
=\sum_{i, j=0}^{N}\sum_{n\in \N}a_{i,n,j}(u,w_1,w_2;
\lambda)x^{r_{i}-n}(\log x)^{j}.\label{Pzaction}
\end{gather}
By definition,
on the region $|z_{1}|>|z|$ \smash{$
\bigl(\bigl(Y^{(g_{1}g_{2})^{-1}}_{P(z)}\bigr)^{o}
(u, z_{1})\lambda\bigr)
(w_{1}\otimes w_{2})
$}
 is absolutely convergent
and its sum on
\smash{$\Omega_{1,0,0}^{(1)}(z)=\Omega_{1,0,0}^{(2)}(z)$}
is equal to $f^e_1(z_1;u, w_{1}, w_{2};\lambda)$.
For simplicity, let
\[
\bigl(Y^{(g_{1}g_{2})^{-1}}_{P(z)}\bigr)^{o}_{m, k}(u)=0
\]
for $m\in \C$, $m\ne -r_{i}+n-1$ for $i=0, \dots, N$ and $n\in \N$
and $k=0, \dots, N$. Then we have
\[
\bigl(Y^{(g_{1}g_{2})^{-1}}_{P(z)}\bigr)^{o}(u, x)
=\sum_{m\in \C}\sum_{k=0}^{N}\bigl(Y^{(g_{1}g_{2})^{-1}}_{P(z)}\bigr)^{o}_{m, k}(u)
x^{-m-1}(\log x)^{k}.
\]

We have the following result.

\begin{prop}\label{COM-inv}
The space $\operatorname{COM}$ is invariant
under the action of the components of the twisted~vertex operators
\smash{$Y^{(g_{1}g_{2})^{-1}}_{P(z)}(u, x)$} for $u\in V$.
\end{prop}
\begin{proof}
We need to show that for $n_{1}\in \C$, $k=1, \dots, K$, $u_{1}\in V$ and
$\lambda\in \operatorname{COM}$,
\begin{equation}\label{com-0}
\bigl(Y^{(g_{1}g_{2})^{-1}}_{P(z)}\bigr)^{o}_{n_{1}, k_{1}}(u_{1})
\lambda\in \operatorname{COM}.
\end{equation}
For $u_{2}, \dots, u_{l}\in V$, we have
\begin{align*}
\begin{split}
Y^{g_2}(u_2,x_2)\cdots Y^{g_2}(u_{l},x_{l})
={}&\sum_{n_{2}, \dots, n_{l}\in \C}\sum_{k_{2}, \dots, k_{l}=0}^{K}
Y^{g_2}_{n_{2}, k_{2}}(u_2)\cdots Y^{g_2}_{n_{l}, k_{l}}(u_{l})
 \nn
&\times
 x_{2}^{-n_{2}-1}\cdots x_{l+1}^{-n_{l}-1}(\log x_{2})^{k_{2}}\cdots
(\log x_{l})^{k_{l}}.
\end{split}
\end{align*}
Since $\lambda\in \operatorname{COM}$,
for $u_{2}, \dots, u_{l}\in V$, $n_{2}, \dots, n_{l}\in \C$,
$k_{2}, \dots, k_{l}=0, \dots, K$,
$w_{1}\in W_{1}$ and $w_{2}\in W_{2}$,
\[
\lambda\bigl(w_{1}\otimes Y^{g_2}(u_{1}, z_{1})Y^{g_2}_{n_{2}, k_{2}}(u_2)
\cdots Y^{g_2}_{n_{l}, k_{l}}(u_{l})w_{2}\bigr)
\]
is absolutely convergent on the region $\Omega_{0,0,1}(z)$
and is absolutely convergent on the re\-gion \smash{$\Omega^{(1)}_{0,0,1}(z)$} and \smash{$\Omega^{(2)}_{0,0,1}(z)$}
to the preferred single-valued branches
\[
f_{1}^{e}\bigl(z_{1}; u_{1}, w_{1}, Y^{g_2}_{n_{2}, k_{2}}(u_2)
\cdots\allowbreak Y^{g_2}_{n_{l}, k_{l}}(u_{l})w_{2}\bigr),
\qquad f_{1}^{b_{z_{1}, z}^{-1}}\bigl(z_{1}; u_{1}, w_{1}, Y^{g_2}_{n_{2}, k_{2}}(u_2)
\cdots Y^{g_2}_{n_{l}, k_{l}}(u_{l})w_{2}\bigr),
\]
respectively, on $M^{1}_{0}(0, z)$ of the multivalued analytic function
\begin{equation}\label{com-1}
f_{1}\bigl(z_{1}; u_{1}, w_{1}, Y^{g_2}_{n_{2}, k_{2}}(u_2)
\cdots Y^{g_2}_{n_{l}, k_{l}}(u_{l})w_{2}\bigr)
\end{equation}
defined on $M^{1}(0, z)$. We discuss only the case that $\Omega^{(1)}_{0,0,1}(z)$
is nonempty (see Remark~\ref{Omega-regions}). By definition,
\smash{$\bigl(Y^{(g_{1}g_{2})^{-1}}_{P(z)}\bigr)^{o}_{-r_{i}+n-1, j}(u)
\lambda$} is obtained by expanding~\eqref{com-1} on the region
$\smash{\Omega_{1,0,0}^{(1)}(z)}=\smash{\Omega_{1,0,0}^{(2)}(z)}$
as a~series in suitable powers of $z_{1}$ and $\log z_{1}$
and then taking the corresponding coefficients.

We now consider the series
\begin{equation}\label{com-3}
f_{1}^{e}(z_{1}; u_{1}, w_{1}, Y^{g_2}(u_2, z_{2})
\cdots Y^{g_2}(u_{l}, z_{l})w_{2})
\end{equation}
in suitable powers of $z_{2}, \dots, z_{l}$ and
$\log z_{1}, \dots, \log z_{l}$. The series~\eqref{com-3}
on the region $\Omega^{(1)}_{0,0,1}(z)$ can be further expanded as the series
\begin{equation}\label{com-4}
\lambda(w_{1}\otimes
Y^{g_2}(u_{1}, z_{1})Y^{g_2}(u_2, z_{2})
\cdots Y^{g_2}(u_{l}, z_{l})w_{2}).
\end{equation}
But~\eqref{com-4} is absolutely convergent on the region
$\Omega_{0,0,l}(z)$ and is absolutely convergent
either to
\begin{gather}\label{com-5}
f^e_{l}(z_1, z_{2}, \dots, z_{l};u_1,\dots,u_l, w_{1}, w_{2};\lambda)
\end{gather}
on the region $\Omega^{(1)}_{0,0,l}(z)$ or to
\begin{gather}\label{com-6}
f^{b_{z_{1},z}^{-1}\cdots
b_{b_{l},z}^{-1}}_{l}(z_1, z_{2}, \dots,z_{l}; u_1,\dots,u_l, w_{1},
w_{2};\lambda)
\end{gather}
on the region $\Omega^{(2)}_{0,0,l}(z)$. Thus for $z_{1}$ satisfying $|z|>|z_{1}|>0$, the series
\eqref{com-3} as the sum of~\eqref{com-4} viewed as a~series in suitable powers
of $z_{1}$ and $\log z_{1}$
must be absolutely convergent on the region~${|z_{1}|>|z_{2}|>
\cdots >|z_{l}|>0}$
and if in addition, $-\frac{3\pi}{2}<\arg (z_{i}-z)<-\frac{\pi}{2}$ or
$\frac{\pi}{2}<\arg (z_{i}-z)<\frac{3\pi}{2}$ for $i=1, \dots, l+1$,
its sum must also be equal to~\eqref{com-5} or~\eqref{com-6}, respectively.
But
\eqref{com-5} can also be expanded
on the region $|z_{1}|>|z|>|z_{2}|>\cdots >|z_{l}|>0$
as a~series in suitable powers of
$z_{1}, \dots, z_{l}$ and $\log z_{1}, \dots,\log z_{l}$.
This fact can be seen as follows: By 1\,(a) in the $P(z)$-compatibility
condition, we know that
$\smash{\prod_{1\leq i<j\leq n}(z_i-z_j)^{M_{ij}}}
f_{l}(z_1,z_{2}, \dots,z_n;u_1,\dots,u_n,w_{1},\allowbreak w_{2};\lambda)$
can be analytically extended to an analytic function on
the region
\[\big\{(z_1,\dots,z_l)\in\C^l \mid z_i\neq0,\, z_i\neq z,\, i=1,\dots,l\big\}.
\]
Moreover, this analytic function has a~regular singularity at
$(\infty, 0, \dots, 0)$. In particular, this function can be
analytically expanded on the region $|z|<|z_{1}|$, $0<|z_{i}|<|z|$
for $i=2, \dots, l$ as a~series in suitable powers of
$z_{1}, \dots, z_{l}$ and $\log z_{1}, \dots,\log z_{l}$.
Thus~\eqref{com-5} can be expanded
on the region $|z_{1}|>|z|>|z_{2}|>\cdots >|z_{l}|>0$
as a~series in suitable powers of
$z_{1}, \dots, z_{l}$ and~${\log z_{1}, \dots,\log z_{l}}$.
By definition, this expansion is equal to
\begin{equation}\label{com-7}
\bigl(\bigl(Y^{(g_{1}g_{2})^{-1}}_{P(z)}\bigr)^{o}(u_{1}, z_{1})
\lambda\bigr)(w_{1}\otimes Y^{g_2}(u_2, z_{2})
\cdots Y^{g_2}(u_{l}, z_{l})w_{2}).
\end{equation}
The coefficients of
the powers of $z_{1}$ and $\log z_{1}$ in~\eqref{com-7} are
exactly
\[\bigl(\bigl(Y^{(g_{1}g_{2})^{-1}}_{P(z)}\bigr)^{o}_{-r_{i}+n-1, j}(u_{1})
\lambda\bigr)(w_{1}\otimes Y^{g_2}(u_2, z_{2})
\cdots Y^{g_2}(u_{l}, z_{l})w_{2})\]
for $i, j=0, \dots, N$, $n\in \N$. Since these are expansions of
single-valued branches of multivalued analytic functions on
$\Omega_{0, 0, l}$, they are
absolutely convergent to these multivalued analytic functions on the same region and
are absolutely convergent to their corresponding single-valued branches on~\smash{$\Omega_{0, 0, l}^{(1)}$} or \smash{$\Omega_{0, 0, l}^{(2)}$}. Moreover,
since $f_{l}(z_1,\dots,z_l; u_1,\dots,u_l, w_{1}, w_{2};\lambda)$
satisfies~1\,(a) and 1(b) in $P(z)$-compatibility condition,
the coefficients of its expansion on the region~${|z_{1}|>|z|>|z_{2}|>
\cdots >|z_{l}|>0}$ also satisfy these conditions.
This finishes the proof of~\eqref{com-0}.
\end{proof}

\begin{prop}\label{COM-comm}
For $u_{1}, \dots, u_{m+l}\in V$,
$w_{1}\in W_{1}$, $w_{2}\in W_{2}$
and $\lambda\in \operatorname{COM}$, the series
\begin{gather}
\bigl(\bigl(Y^{(g_{1}g_{2})^{-1}}_{P(z)}\bigr)^{o}(u_{m},
z_{m})\cdots \bigl(Y^{(g_{1}g_{2})^{-1}}_{P(z)}\bigr)^{o}(u_{1},
z_{1})\lambda\bigr)\nn
\qquad \times \bigl(w_{1}\otimes Y^{g_{2}}(u_{m+1}, z_{m+1})\cdots
Y^{g_{2}}(u_{m+l}, z_{m+l})w_{2}\bigr)\label{COM-inv-0}
\end{gather}
is absolutely convergent on the region $|z_{1}|>\cdots>|z_{m}|>|z|
>|z_{m+1}|>\cdots>|z_{l+m}|>0$ and its sum is equal to
\begin{equation}\label{COM-inv-1.5}
f_{m+l}^{e}(z_{1}, \dots, z_{m+l};
u_{1}, \dots, u_{m+l};
w_{1}, w_{2}; \lambda)
\end{equation}
on \smash{$\Omega_{m, 0, l}^{(1)}(z)$}
or to
\begin{equation}\label{COM-inv-1.7}
f_{m+l}^{b_{z_{m+1}, z}^{-1}\cdots b_{z_{m+l}, z}^{-1}}
(z_{1}, \dots, z_{m+l};
u_{1}, \dots, u_{m+l};
w_{1}, w_{2}; \lambda)
\end{equation}
on \smash{$\Omega_{m, 0, l}^{(2)}(z)$}.
Moreover,
we have the following commutativity for \smash{$Y^{(g_{1}g_{2})^{-1}}_{P(z)}$}:
For $u_{1}, \dots, u_{m}\in V$, $w_{1}\in W_{1}$, $w_{2}\in W_{2}$
and $\lambda\in \operatorname{COM}$, the series
\begin{equation}\label{COM-inv-0.1}
\bigl(Y^{(g_{1}g_{2})^{-1}}_{P(z)}(u_{1},
z_{1})\cdots Y^{(g_{1}g_{2})^{-1}}_{P(z)}(u_{m},
z_{m})\lambda\bigr)(w_{1}\otimes w_{2})
\end{equation}
is absolutely convergent on the region
$|z^{-1}|>|z_{1}|>\cdots>|z_{m}|>0$
and for $\sigma\in S_{m}$, the sums of~\eqref{COM-inv-0.1}
and
\begin{equation}\label{COM-inv-0.2}
\bigl(Y^{(g_{1}g_{2})^{-1}}_{P(z)}(u_{\sigma(1)},
z_{\sigma(1)})\cdots
Y^{(g_{1}g_{2})^{-1}}_{P(z)}(u_{\sigma(m)},
z_{\sigma(l)})\lambda\bigr)(w_{1}\otimes w_{2})
\end{equation}
are analytic extensions of each other. We also have
the weak commutativity for \smash{$Y^{(g_{1}g_{2})^{-1}}_{P(z)}$}:
For $u, v\in V$,
there exists $M\in \Z_{+}$ depending only on $u$ and $v$ such that
\begin{gather}
(x_{1}-x_{2})^{M}Y^{(g_{1}g_{2})^{-1}}_{P(z)}(u, x_{1})
Y^{(g_{1}g_{2})^{-1}}_{P(z)}(v, x_{2})\nn
\qquad
=(x_{1}-x_{2})^{M}Y^{(g_{1}g_{2})^{-1}}_{P(z)}(v, x_{2})
Y^{(g_{1}g_{2})^{-1}}_{P(z)}(u, x_{1}).\label{wk-comm}
\end{gather}
\end{prop}
\begin{proof}
To prove the convergence of~\eqref{COM-inv-0}, we use induction
on $m$. In the case $m=0$, the convergence of~\eqref{COM-inv-0} is given by condition~(2)\,(b).
Assume that the convergence of~\eqref{COM-inv-0} in the case that
$m$ is $m-1$ holds.
Then
\begin{gather}
\bigl(\bigl(Y^{(g_{1}g_{2})^{-1}}_{P(z)}\bigr)^{o}(u_{m-1},
z_{m-1})\cdots \bigl(Y^{(g_{1}g_{2})^{-1}}_{P(z)}\bigr)^{o}(u_{1},
z_{1})\lambda\bigr)\nn
\quad\times
\bigl(w_{1}\otimes Y^{g_{2}}(u_{m}, z_{m})\cdots
Y^{g_{2}}(u_{m+l}, z_{m+l})w_{2}\bigr)\nn
 \phantom{\quad\times}{} =\bigl(\bigl(Y^{(g_{1}g_{2})^{-1}}_{P(z)}\bigr)^{o}(u_{m-1},
x_{m-1})\cdots \bigl(Y^{(g_{1}g_{2})^{-1}}_{P(z)}\bigr)^{o}(u_{1},
x_{1})\lambda\bigr)\nn
 \phantom{\quad\times=}{}
\times \bigl(w_{1}\otimes Y^{g_{2}}(u_{m}, x_{m})\cdots
Y^{g_{2}}(u_{m+l}, x_{m+l})w_{2}\bigr)\lbar_{x_{i}^{n}={\rm e}^{n\log z_{i}},
\log x_{i}=\log z_{i}, i=1, \dots, m+l}\label{COM-inv-1}
\end{gather}
is absolutely convergent on the region $|z_{1}|>\cdots>|z_{m-1}|>|z|
>|z_{m}|>\cdots>|z_{m+l}|>0$ and its sum is equal to
\eqref{COM-inv-1.5}
on \smash{$\Omega_{m-1, 0, l+1}^{(1)}(z)$}
or to
\begin{equation}\label{COM-inv-1.8}
f_{l}^{b_{z_m, z}^{-1}\cdots b_{z_{m+l}, z}^{-1}}
(z_{1}, \dots, z_{m+l};
u_{1}, \dots, u_{m+l};
w_{1}, w_{2}; \lambda)
\end{equation}
on \smash{$\Omega_{m-1, 0, l+1}^{(2)}(z)$}.

Write
\begin{gather*}
\bigl(Y^{(g_{1}g_{2})^{-1}}_{P(z)}\bigr)^{o}(u_{m-1},
x_{m-1})\cdots \bigl(Y^{(g_{1}g_{2})^{-1}}_{P(z)}\bigr)^{o}(u_{1},
x_{1})\nn
 \qquad =\sum_{n_{m-1}, \dots, n_{1}\in \C}\sum_{k_{m-1},
\dots, k_{1}=0}^{K}
\bigl(Y^{(g_{1}g_{2})^{-1}}_{P(z)}\bigr)^{o}_{n_{m-1}, k_{m-1}}(u_{m-1})
\cdots \bigl(Y^{(g_{1}g_{2})^{-1}}_{P(z)}\bigr)^{o}_{n_{1}, k_{1}}(u_{1})
 \nn
\qquad \phantom{=}{}
 \times x_{m-1}^{-n_{m-1}-1}\cdots x_{1}^{-n_{1}-1}
 (\log x_{m-1})^{k_{m-1}}\cdots
(\log x_{1})^{k_{1}}
\end{gather*}
and
\begin{gather*}
Y^{g_2}(u_{m+1},x_{m+1})\cdots Y^{g_2}(u_{m+l},x_{m+l})\nn
\qquad=\sum_{n_{m+1}, \dots, n_{m+l}\in \C}\sum_{k_{m+1}, \dots, k_{m+l}=0}^{K}
Y^{g_2}_{n_{m+1}, k_{m+1}}(u_{m+1})\cdots Y^{g_2}_{n_{m+l}, k_{m+l}}(u_{m+l})
 \nn
\qquad \phantom{=}{}
\times x_{2}^{-n_{m+1}-1}\cdots x_{m+l}^{-n_{m+l}-1}
(\log x_{m+1})^{k_{m+1}}\cdots
(\log x_{m+l})^{k_{m+l}}.
\end{gather*}
By Proposition~\eqref{COM-inv},
\[
\bigl(Y^{(g_{1}g_{2})^{-1}}_{P(z)}\bigr)^{o}_{n_{m-1},
k_{m-1}}(u_{m-1})
\cdots \bigl(Y^{(g_{1}g_{2})^{-1}}_{P(z)}\bigr)^{o}_{n_{1}, k_{1}}
(u_{1})\lambda\in \operatorname{COM}.\]
Then
\begin{gather*}
\bigl(\bigl(Y^{(g_{1}g_{2})^{-1}}_{P(z)}\bigr)^{o}_{n_{m-1},
k_{m-1}}(u_{m-1})
\cdots \bigl(Y^{(g_{1}g_{2})^{-1}}_{P(z)}\bigr)^{o}_{n_{1}, k_{1}}
(u_{1})\lambda\bigr)\nn
 \qquad\times
\bigl(w_{1}\otimes Y^{g_{2}}(u_{m}, z_{m})Y^{g_{2}}_{n_{m+1},
k_{m+1}}(u_{m}) \cdots Y^{g_{2}}_{n_{m+l},
k_{m+l}}(u_{m+l})w_{2}\bigr)
\end{gather*}
is absolutely convergent on the region $|z|>|z_{m}|>0$
and its sum is equal to
\begin{gather}
f^{e}\bigl(z_{m}; u_{m}, w_{1}, Y^{g_{2}}_{n_{m+1},
k_{m+1}}(u_{m}) \cdots Y^{g_{2}}_{n_{m+l},
k_{m+l}}(u_{m+l})w_{2}; \nn
 \qquad
\bigl(Y^{(g_{1}g_{2})^{-1}}_{P(z)}\bigr)^{o}_{n_{m-1},
k_{m-1}}(u_{m-1})
\cdots \bigl(Y^{(g_{1}g_{2})^{-1}}_{P(z)}\bigr)^{o}_{n_{1}, k_{1}}
(u_{1})\lambda\bigr)\label{COM-inv-2}
\end{gather}
on the region $|z|>|z_{m}|>0$, $-\frac{3\pi}{2}<\arg (z_{m}-z)-\arg z
<-\frac{\pi}{2}$ and to
\begin{gather*}
f^{b_{z_{m}, z}^{-1}}\bigl(z_{m}; u_{m}, w_{1}, Y^{g_{2}}_{n_{m+1},
k_{m+1}}(u_{m}) \cdots Y^{g_{2}}_{n_{m+l},
k_{m+l}}(u_{m+l})w_{2}; \nn
 \qquad
\bigl(Y^{(g_{1}g_{2})^{-1}}_{P(z)}\bigr)^{o}_{n_{m-1},
k_{m-1}}(u_{m-1})
\cdots \bigl(Y^{(g_{1}g_{2})^{-1}}_{P(z)}\bigr)^{o}_{n_{1}, k_{1}}
(u_{1})\lambda\bigr)
\end{gather*}
on the region $|z|>|z_{m}|>0$, $\frac{\pi}{2}<\arg (z_{m}-z)-\arg z
<\frac{3\pi}{2}$.
By definition,
\begin{gather}
\bigl(\bigl(Y^{(g_{1}g_{2})^{-1}}_{P(z)}\bigr)^{o}(u_{m}, z_{m})
\bigl(Y^{(g_{1}g_{2})^{-1}}_{P(z)}\bigr)^{o}_{n_{m-1},
k_{m-1}}(u_{m-1})
\cdots \bigl(Y^{(g_{1}g_{2})^{-1}}_{P(z)}\bigr)^{o}_{n_{1}, k_{1}}
(u_{1})\lambda\bigr)\nn
 \qquad\times
\bigl(w_{1}\otimes Y^{g_{2}}_{n_{m+1},
k_{m+1}}(u_{m}) \cdots Y^{g_{2}}_{n_{m+l},
k_{m+l}}(u_{m+l})w_{2}\bigr)\label{COM-inv-2.1}
\end{gather}
is absolutely convergent on the region $|z_{m}|>|z|$
and its sum is equal to~\eqref{COM-inv-2}
on the re\-gion~${|z_{m}|>|z|}$, $|{\arg (z_{m}-z)-\arg z_{m}}|<\frac{\pi}{2}$.

But we know that
\begin{gather}
\sum_{n_{1}, \dots, n_{m-1}, n_{m+1}, \dots, n_{m+l}\in \C}
\sum_{k_{1},
\dots, k_{m-1}, k_{m+1},
\dots, k_{m+l}=0}^{K}f^{e}\bigl(z_{m}; u_{m}, w_{1}, Y^{g_{2}}_{n_{m+1},
k_{m+1}}(u_{m})\cdots
 \nn
\qquad\phantom{=\dots}{}
\times Y^{g_{2}}_{n_{m+l},
k_{m+l}}(u_{m+l})w_{2};
\bigl(Y^{(g_{1}g_{2})^{-1}}_{P(z)}\bigr)^{o}_{n_{m-1},
k_{m-1}}(u_{m-1})
\cdots \bigl(Y^{(g_{1}g_{2})^{-1}}_{P(z)}\bigr)^{o}_{n_{1}, k_{1}}
(u_{1})\lambda\bigr)\nn
\qquad\phantom{=\dots}{}\times
z_{1}^{-n_{1}-1}\cdots z_{m-1}^{-n_{m-1}-1}
z_{m+1}^{-n_{1}-1}\cdots z_{m+l}^{-n_{m+l}-1}\nn
\qquad\phantom{=\dots}{}\times
 (\log z_{1})^{k_{1}}\cdots (\log z_{m-1})^{k_{m-1}}
 (\log z_{m+1})^{k_{m+1}}\cdots (\log z_{m+l})^{k_{m+l}}
 \nn
\qquad
=\sum_{n_{1}, \dots, n_{m-1}, n_{m+1}, \dots, n_{m+l}\in \C}
\sum_{k_{1},
\dots, k_{m-1}, k_{m+1},
\dots, k_{m+l}=0}^{K}\bigl(\bigl(Y^{(g_{1}g_{2})^{-1}}_{P(z)}\bigr)^{o}_{n_{m-1},
k_{m-1}}(u_{m-1})\cdots\nn
\qquad\phantom{=}{}
\times \bigl(Y^{(g_{1}g_{2})^{-1}}_{P(z)}\bigr)^{o}_{n_{1}, k_{1}}
(u_{1})\lambda\bigr)
\bigl(w_{1}\otimes Y^{g_{2}}(u_{m}, z_{m})Y^{g_{2}}_{n_{m+1},
k_{m+1}}(u_{m}) \cdots Y^{g_{2}}_{n_{m+l},
k_{m+l}}(u_{m+l})w_{2}\bigr)\nn
 \qquad =
z_{1}^{-n_{1}-1}\cdots z_{m-1}^{-n_{m-1}-1}
z_{m+1}^{-n_{1}-1}\cdots z_{m+l}^{-n_{m+l}-1}\nn
\qquad\phantom{=}{}\times
 (\log z_{1})^{k_{1}}\cdots (\log z_{m-1})^{k_{m-1}}
 (\log z_{m+1})^{k_{m+1}}\cdots (\log z_{m+l})^{k_{m+l}}\label{COM-inv-3}
\end{gather}
as an iterated series of the multi-series~\eqref{COM-inv-1}
is absolutely convergent on the region
\[
|z_{1}|>\cdots>|z_{m-1}|>|z|
>|z_{m}|>\cdots>|z_{m+l}|>0
\]
 and its sum is equal to
\eqref{COM-inv-1.5}
on \smash{$\Omega_{m-1, 0, l+1}^{(1)}(z)$} or to~\eqref{COM-inv-1.8}
on \smash{$\Omega_{m-1, 0, l+1}^{(2)}(0, z)$}. We discuss~only the case that
\smash{$\Omega_{m-1, 0, l+1}^{(1)}(z)$} is nonempty. The case that~\smash{$\Omega_{m-1, 0, l+1}^{(2)}(0, z)$} is nonempty is similar and is omitted.
The sum of~\eqref{COM-inv-3} being equal to
\eqref{COM-inv-1.5}
on \smash{$\Omega_{m-1, 0, l+1}^{(1)}(z)$} means that the expansion of~\eqref{COM-inv-1.5}
on \smash{$\Omega_{m-1, 0, l+1}^{(1)}(z)$} can also be written
as the iterated series~\eqref{COM-inv-3}.
By the discussion above, the coefficients of the left-hand side
of~\eqref{COM-inv-3}
is equal to~\eqref{COM-inv-2.1} on the region
$|z_{m}|>|z|$, $|{\arg (z_{m}-z)-\arg z_{m}}|<\frac{\pi}{2}$.
So the expansion of~\eqref{COM-inv-1.5} on the region
\smash{$\Omega_{m, 0, l}^{(1)}(z)$} can be written as the iterated series
\begin{gather}
\sum_{n_{1}, \dots, n_{m-1}, n_{m+1}, \dots, n_{m+l}\in \C}
\sum_{k_{1},
\dots, k_{m-1}, k_{m+1},
\dots, k_{m+l}=0}^{K}
\bigl(\bigl(Y^{(g_{1}g_{2})^{-1}}_{P(z)}\bigr)^{o}(u_{m}, z_{m})\nn
\qquad\times
\bigl(Y^{(g_{1}g_{2})^{-1}}_{P(z)}\bigr)^{o}_{n_{m-1},
k_{m-1}}(u_{m-1})
\cdots \bigl(Y^{(g_{1}g_{2})^{-1}}_{P(z)}\bigr)^{o}_{n_{1}, k_{1}}
(u_{1})\lambda\bigr)\nn
\qquad\times
\bigl(w_{1}\otimes Y^{g_{2}}_{n_{m+1},
k_{m+1}}(u_{m}) \cdots Y^{g_{2}}_{n_{m+l},
k_{m+l}}(u_{m+l})w_{2}\bigr)\nn
\qquad\times
z_{1}^{-n_{1}-1}\cdots z_{m-1}^{-n_{m-1}-1}
z_{m+1}^{-n_{1}-1}\cdots z_{m+l}^{-n_{m+l}-1}\nn
\qquad\times
 (\log z_{1})^{k_{1}}\cdots (\log z_{m-1})^{k_{m-1}}
 (\log z_{m+1})^{k_{m+1}}\cdots (\log z_{m+l})^{k_{m+l}}.\label{COM-inv-4}
\end{gather}
Since the expansion of~\eqref{COM-inv-1.5} on the region
\smash{$\Omega_{m, 0, l}^{(1)}(z)$} must be absolutely convergent
as a~multisum, we see that the multi-series~\eqref{COM-inv-0}
corresponding to~\eqref{COM-inv-4} must be absolutely convergent
on the region \smash{$\Omega_{m, 0, l}^{(1)}(z)$} to~\eqref{COM-inv-1.5}.
Similarly, we can show that~\eqref{COM-inv-0} is absolutely convergent~on~the region \smash{$\Omega_{m, 0, l}^{(2)}(z)$} to~\eqref{COM-inv-1.7}.

This convergence result implies in particular
the absolute convergence
of~\eqref{COM-inv-0.1} on the re\-gion~$|z^{-1}|>
|z_{1}|>\cdots>|z_{l}|>0$.
Using the commutativity
for the twisted vertex operators~$Y^{g_{2}}$,
we see that for
$\sigma\in S_{m}$, the sums of~\eqref{COM-inv-0.1}
and
\eqref{COM-inv-0.2} are analytic extensions of each other.

For $u, v\in V$, since $\lambda$ satisfies the condition 1\,(a) in the
$P(z)$-compatibility condition, there exists $M\in \Z_{+}$ depending on
$u$ and $v$ only (see Remark~\ref{remark-M-i-j})
 such that
$z_{1}-z_{2}=0$ is not a~singularity of
$(z_{1}-z_{2})^{M}f^e_{2}(z_1, z_{2};u, v, w_{1}, w_{2};\lambda)$
and we have
\begin{equation}\label{COM-inv-5}
(z_{1}-z_{2})^{M}f^e_{2}(z_1, z_{2};u, v, w_{1}, w_{2};\lambda)
=(z_{1}-z_{2})^{M}f^e_{2}(z_2, z_{1};v, u, w_{1}, w_{2};\lambda).
\end{equation}
Since the expansion of $f^e_{2}(z_1, z_{2};u, v, w_{1}, w_{2};\lambda)$
on the region \smash{$\Omega^{(1)}_{2, 0, 0}(z)$} is
\[\bigl(\bigl(Y^{(g_{1}g_{2})^{-1}}_{P(z)}\bigr)^{o}(v, z_{2})
\bigl(Y^{(g_{1}g_{2})^{-1}}_{P(z)}\bigr)^{o}(u, z_{1})
\lambda\bigr)(w_{1}\otimes w_{2}),\]
we see that
\begin{gather}
(z_{1}-z_{2})^{M}
\bigl(\bigl(Y^{(g_{1}g_{2})^{-1}}_{P(z)}\bigr)^{o}(v, z_{2})
\bigl(Y^{(g_{1}g_{2})^{-1}}_{P(z)}\bigr)^{o}(u, z_{1})
\lambda\bigr)(w_{1}\otimes w_{2})\nn
\qquad=\sum_{i=0}^{M}\sum_{n\in \C}\sum_{k=0}^{K_{1}}
\sum_{n_{1}\in \C}\sum_{k_{1}=0}^{K_{1}}
\binom{M}{i}\bigl(\bigl(Y^{(g_{1}g_{2})^{-1}}_{P(z)}
\bigr)^{o}_{n_{1}, k_{1}}(v)
\bigl(Y^{(g_{1}g_{2})^{-1}}_{P(z)}\bigr)^{o}_{n, k}(u)
\lambda\bigr)(w_{1}\otimes w_{2})\cdot\nn
 \qquad\phantom{=}{}\times{\rm e}^{(n+M-i)\log z_{1}}
(\log z_{1})^{k}{\rm e}^{(n_{1}+i)\log z_{2}}
(\log z_{2})^{k_{1}}\label{com-9}
\end{gather}
must be convergent absolutely to
the left-hand side of~\eqref{COM-inv-5}
on \smash{$\Omega^{(1)}_{2, 0, 0}(z)$} (the region given by
$|z_{1}|>|z_{2}|>|z|$, $|{\arg (z_{1}-z)-\arg z_{1}}|<\frac{\pi}{2}$,
$|{\arg (z_{2}-z)-\arg z_{2}}|<\frac{\pi}{2}$).
On the other hand, $(z_{1}-z_{2})^{M}f_{2}^{e}(z_1, z_{2};u,
v, w_{1}, w_{2};\lambda)$ is analytic at $z_{1}-z_{2}=0$.
So~\eqref{com-9} is in fact absolutely convergent
to the left-hand side of~\eqref{COM-inv-5}
on the region
$|z_{1}|, |z_{2}|>|z|$, $|{\arg (z_{1}-z)-\arg z_{1}}|<\frac{\pi}{2}$,
$|{\arg (z_{2}-z)-\arg z_{2}}|<\frac{\pi}{2}$. Thus
\begin{equation}\label{com-10}
(z_{1}-z_{2})^{M}
\bigl(\bigl(Y^{(g_{1}g_{2})^{-1}}_{P(z)}\bigr)^{o}(u, z_{1})
\bigl(Y^{(g_{1}g_{2})^{-1}}_{P(z)}\bigr)^{o}(v, z_{2})
\lambda\bigr)(w_{1}\otimes w_{2})
\end{equation}
is absolutely convergent to the right-hand side of~\eqref{COM-inv-5}
also on the region
$|z_{1}|, |z_{2}|>|z|$, $|{\arg (z_{1}-z)-\arg z_{1}}|<\frac{\pi}{2}$,
$|{\arg (z_{2}-z)-\arg z_{2}}|<\frac{\pi}{2}$. From
\eqref{COM-inv-5}, \eqref{com-9} and~\eqref{com-10},
we obtain
\begin{align*}
&(z_{1}-z_{2})^{M}
\bigl(\bigl(Y^{(g_{1}g_{2})^{-1}}_{P(z)}\bigr)^{o}(v, z_{2})
\bigl(Y^{(g_{1}g_{2})^{-1}}_{P(z)}\bigr)^{o}(u, z_{1})
\lambda\bigr)(w_{1}\otimes w_{2})\nn
& \qquad =(z_{1}-z_{2})^{M}
\bigl(\bigl(Y^{(g_{1}g_{2})^{-1}}_{P(z)}\bigr)^{o}(u, z_{1})
\bigl(Y^{(g_{1}g_{2})^{-1}}_{P(z)}\bigr)^{o}(v, z_{2})
\lambda\bigr)(w_{1}\otimes w_{2})
\end{align*}
on the region $|z_{1}|, |z_{2}|>|z|$,
$|{\arg (z_{1}-z)-\arg z_{1}}|<\frac{\pi}{2}$,
$|{\arg (z_{2}-z)-\arg z_{2}}|<\frac{\pi}{2}$
for $\lambda\in \operatorname{COM}$, $u_{1}, v\in V$, $w_{1}\in W_{1}$
and $w_{2}\in W_{2}$,
which, by the definition of
\smash{$\bigl(Y^{(g_{1}g_{2})^{-1}}_{P(z)}\bigr)^{o}(u, x)$} for $u\in V$
above, is equivalent to~\eqref{wk-comm}.
\end{proof}

Since $\operatorname{COMP}\subset \operatorname{COM}$,
by Proposition~\ref{COM-inv},
\smash{$Y^{(g_{1}g_{2})^{-1}}_{P(z)}(u, x)\lambda$} for $u\in V$
and $\lambda\in \operatorname{COMP}$ is in
$\operatorname{COM}\{x\}[\log x]$. We now prove the following
stronger result.

\begin{prop}\label{COMP-inv}
The space $\operatorname{COMP}$ is invariant under the action
of the components of the twisted vertex operators
\smash{$Y^{(g_{1}g_{2})^{-1}}_{P(z)}(u, x)$} for $u\in V$.
\end{prop}
\begin{proof}
Let $\lambda\in \operatorname{COMP}$. Then $\lambda$
satisfies 2\,(a) in the $P(z)$-compatibility condition. We need only show that
for $v\in V$, $n\in \C$, $k=0, \dots, K$,
\smash{$\bigl(Y^{(g_{1}g_{2})^{-1}}_{P(z)}\bigr)^{o}_{n, k}(v)\lambda$}
also satisfies 2\,(a) in the $P(z)$-compatibility condition.

By the definition of \smash{$\bigl(Y^{(g_{1}g_{2})^{-1}}_{P(z)}\bigr)^{o}
(u, z_{1})$}
for $u\in V$ and the $P(z)$-compatibility condition for $\lambda$,
\[\bigl(\bigl(Y^{(g_{1}g_{2})^{-1}}_{P(z)}\bigr)^{o}(u,
z_{1})\lambda\bigr)\bigl(w_{1}\otimes Y^{g_{2}}_{n, k}(v)w_{2}\bigr)\]
and
\smash{$\lambda\bigl(Y^{g_{1}}(u, z_{1}-z)w_{1}\otimes Y^{g_{2}}_{n, k}(v)w_{2}\bigr)$}
for $v\in V$, $n\in \C$, $k=0, \dots, K$,
$w_{1}\in W_{1}$ and $w_{2}\in W_{2}$
are absolutely convergent to
$f_{1}^{e}\bigl(z_{1}; u, w_{1}, Y^{g_{2}}_{n, k}(v)w_{2}; \lambda\bigr)$
on the region $|z_{1}|>|z|$,
$|{\arg (z_{1}-z)}-\arg z_{1}|<\frac{\pi}{2}$ and
the region $|z|>|z_{1}-z|>0$, $|{\arg z_{1}-\arg z}|<\frac{\pi}{2}$,
respectively.
By Proposition~\ref{COM-comm},
\begin{gather*}
\bigl(\bigl(Y^{(g_{1}g_{2})^{-1}}_{P(z)}\bigr)^{o}(u,
z_{1})\lambda\bigr)(w_{1}\otimes
Y^{g_{2}}(v, z_{2})w_{2})\nn
\qquad=\sum_{n\in \C}\sum_{k=0}^{K}
\bigl(\bigl(Y^{(g_{1}g_{2})^{-1}}_{P(z)}\bigr)^{o}(u,
z_{1})\lambda\bigr)\bigl(w_{1}\otimes Y^{g_{2}}_{n, k}(v)w_{2}\bigr)
{\rm e}^{n\log z_{2}}(\log z_{2})^{k}
\end{gather*}
is absolutely convergent on the region $|z_{1}|>|z|
>|z_{2}|>0$ and its sum is equal to
\begin{equation}\label{COMP-inv-1.5}
f_{2}^{e}(z_{1}, z_{2};
u, v; w_{1}, w_{2}; \lambda)
\end{equation}
on $\Omega_{1, 0, 1}^{(1)}(z)$
and to
\begin{equation}\label{COMP-inv-1.7}
f_{2}^{b_{z_{2}, z}^{-1}}
(z_{1}, z_{2};
u, v; w_{1}, w_{2}; \lambda)
\end{equation}
on $\Omega_{1, 0, 1}^{(2)}(z)$.
Thus
\begin{align}
&\sum_{n\in \C}\sum_{k=0}^{K}
\lambda\bigl(Y^{g_{1}}(u, z_{1}-z)w_{1}\otimes Y^{g_{2}}_{n, k}(v)w_{2}\bigr)
{\rm e}^{n\log z_{2}}(\log z_{2})^{k}\nn
& \qquad =
\lambda\bigl(Y^{g_{1}}(u, z_{1}-z)w_{1}\otimes Y^{g_{2}}(v, z_{2})w_{2}\bigr)\label{COMP-inv-0.5}
\end{align}
is absolutely convergent on the region $|z_{1}|>|z|>|z_{1}-z|+|z_{2}|>0$
and its sum is equal to~\eqref{COMP-inv-1.5} on the region
$|z_{1}|>|z|>|z_{1}-z|+|z_{2}|>0$,
$|{\arg (z_{1}-z)-\arg z_{1}}|<\frac{\pi}{2}$,
$|{\arg z_{1}-\arg z}|<\frac{\pi}{2}$, $-\frac{3\pi}{2}< \arg (z_{2}-z)-\arg z<
-\frac{\pi}{2}$ and is equal to~\eqref{COMP-inv-1.7}
on the region $|z_{1}|>|z|>|z_{1}-z|+|z_{2}|>0$,
$|{\arg (z_{1}-z)-\arg z_{1}}|<\frac{\pi}{2}$,
$|{\arg z_{1}-\arg z}|<\frac{\pi}{2}$, $\frac{\pi}{2} <\arg (z_{2}-z)-\arg z<
\frac{3\pi}{2}$. On the other hand, since $(z_{1}-z, z_{2})=(0, 0)$
is a~regular singular point of~\eqref{COMP-inv-1.5} and~\eqref{COMP-inv-1.7},
we can expand them on the regions
$|z|>|z_{2}|, |z_{1}-z|>0$ to obtain a~series of the same form
as~\eqref{COMP-inv-0.5}. Thus we see that~\eqref{COMP-inv-0.5}
must be absolutely convergent on the region
$|z|>|z_{2}|, |z_{1}-z|>0$ and its sum is equal to
\eqref{COMP-inv-1.5} and~\eqref{COMP-inv-1.7} on
the regions $|z|>|z_{1}-z|+|z_{2}|>0$,
$|{\arg z_{1}-\arg z}|<\frac{\pi}{2}$,
$-\frac{3\pi}{2}< \arg (z_{2}-z)-\arg z<
-\frac{\pi}{2}$ and on the region $|z|>|z_{1}-z|+|z_{2}|>0$,
$|{\arg z_{1}-\arg z}|<\frac{\pi}{2}$, $\frac{\pi}{2} <\arg (z_{2}-z)-\arg z<
\frac{3\pi}{2}$, respectively.

The right-hand side of~\eqref{COMP-inv-0.5} is equal to
\begin{align}\label{COMP-inv-0.59}
\sum_{n\in \C}\sum_{k=0}^{K}
\lambda\bigl(Y^{g_{1}}_{n, k}(u)w_{1}\otimes Y^{g_{2}}(v, z_{2})w_{2}\bigr)
{\rm e}^{n\log (z_{1}-z)}(\log (z_{1}-z))^{k}.
\end{align}
We know that the series
\smash{$\lambda\bigl(Y^{g_{1}}_{n, k}(u)w_{1}\otimes Y^{g_{2}}(v, z_{2})w_{2}\bigr)$}
is absolutely convergent on the region~$|z|>|z_{2}|>0$ and its sum
is equal to
\smash{$f_{1}^{e}\bigl(z_{2}; v, Y^{g_{1}}_{n, k}(u)w_{1}, w_{2};\lambda\bigr)$}
on the region $|z|>|z_{2}|>0$, $-\frac{3\pi}{2}
<\arg (z_{2}-z)-\arg z_{2}<
-\frac{\pi}{2}$ and to
\smash{$f_{1}^{b_{z_{2}, z}}\bigl(z_{2}; v,
Y^{g_{1}}_{n, k}(u)w_{1}, w_{2};\lambda\bigr)$}
on the region $|z|>|z_{2}|>0$, $\frac{\pi}{2}<\arg (z_{2}-z)-\arg z_{2}<
\frac{3\pi}{2}$. We also know that the series
\[\bigl(\bigl(Y^{(g_{1}g_{2})^{-1}}_{P(z)}\bigr)^{o}
(v, z_{2})\lambda\bigr)\bigl(Y^{g_{1}}_{n, k}(u)w_{1}\otimes w_{2}\bigr)\]
is absolutely convergent to
\smash{$f_{1}^{e}\bigl(z_{2}; v, Y^{g_{1}}_{n, k}(u)w_{1}, w_{2};\lambda\bigr)$}
on the region $|z_{2}|>|z|>0$, $|{\arg (z_{2}-z)}-\arg z_{2}|<\frac{\pi}{2}$.
From this discussion, \eqref{COMP-inv-0.59} and the convergence
of~\eqref{COMP-inv-0.5}, we see that
\begin{align}
&\sum_{n\in \C}\sum_{k=0}^{K}
\bigl(\bigl(Y^{(g_{1}g_{2})^{-1}}_{P(z)}\bigr)^{o}
(v, z_{2})\lambda\bigr)(Y^{g_{1}}_{n, k}(u)w_{1}\otimes w_{2})
{\rm e}^{n\log (z_{1}-z)}(\log (z_{1}-z))^{k}\nn
&\qquad=\bigl(\bigl(Y^{(g_{1}g_{2})^{-1}}_{P(z)}\bigr)^{o}
(v, z_{2})\lambda\bigr)\bigl(Y^{g_{1}}(u, z_{1}-z)w_{1}\otimes w_{2}\bigr)\label{COMP-inv-0.6}
\end{align}
is absolutely convergent on the region $|z_{2}|>|z|>|z_{1}-z|>0$
and its sum is equal to~\eqref{COMP-inv-1.5} and~\eqref{COMP-inv-1.7} on the region
$|z_{2}|>|z|>|z_{1}-z|>0$, $|{\arg (z_{2}-z)-\arg z_{2}}|<\frac{\pi}{2}$,
$|{\arg z_{1}-\arg z}|<\frac{\pi}{2}$ and on the region
$|z_{2}|>|z|>|z_{1}-z|>0$, $|{\arg (z_{2}-z)-\arg z_{2}}|<\frac{\pi}{2}$,
$\frac{\pi}{2} <\arg (z_{2}-z)-\arg z<
\frac{3\pi}{2}$, respectively.

On the other hand, by Proposition~\ref{COM-comm},
\begin{align}\label{COMP-inv-0.7}
\bigl(\bigl(Y^{(g_{1}g_{2})^{-1}}_{P(z)}\bigr)^{o}(u,
z_{1})\bigl(Y^{(g_{1}g_{2})^{-1}}_{P(z)}\bigr)^{o}(v,
z_{2})\lambda\bigr)(w_{1}\otimes
w_{2})
\end{align}
is absolutely convergent on the region $|z_{2}|>|z_{1}|>|z|$
and its sum is equal to~\eqref{COMP-inv-1.5} on the region
$|z_{2}|>|z_{1}|>|z|$, $|{\arg (z_{1}-z)-\arg z_{1}}|<\frac{\pi}{2}$,
$|{\arg (z_{2}-z)-\arg z_{2}}|<\frac{\pi}{2}$. Taking the coefficients of
\smash{$\bigl(Y^{(g_{1}g_{2})^{-1}}_{P(z)}(v, z_{2})\bigr)^{o}$}
in both~\eqref{COMP-inv-0.6}
and~\eqref{COMP-inv-0.7}, we see that
\[\bigl(\bigl(Y^{(g_{1}g_{2})^{-1}}_{P(z)}\bigr)^{o}_{n, k}
(v)\lambda\bigr)\bigl(Y^{g_{1}}(u, z_{1}-z)w_{1}\otimes w_{2}\bigr)\]
and
\[\bigl(\bigl(Y^{(g_{1}g_{2})^{-1}}_{P(z)}\bigr)^{o}(u,
z_{1})\bigl(Y^{(g_{1}g_{2})^{-1}}_{P(z)}\bigr)^{o}_{n, k}(v)\lambda\bigr)(w_{1}\otimes
w_{2})\]
are absolutely convergent to
the coefficients of \smash{${\rm e}^{n\log z_{2}}(\log z_{2})^{k}$}
in the expansion of~\eqref{COMP-inv-1.5} near the singularity $z_{2}=\infty$
on the region $|z|>|z_{1}-z|>0$,
$|{\arg z_{1}-\arg z}|<\frac{\pi}{2}$ and
on the region ${|z_{1}|>|z|}$, $|{\arg (z_{1}-z)-\arg z_{1}}|<\frac{\pi}{2}$,
respectively. This is equivalent to 2\,(a) in the~$P(z)$\nobreakdash-com\-pat\-ibility
condition for \smash{$\bigl(Y^{(g_{1}g_{2})^{-1}}_{P(z)}\bigr)^{o}_{n, k}(v)
\lambda$}.
\end{proof}

For a~$P(z)$-intertwining map $I$
of type $\binom{W_{3}}{W_{1}W_{2}}$ and an element
$w_{3}'$ of $W_{3}'$, the element $\lambda_{I, w_{3}'}
\in \operatorname{COM}$ also has the following property.

\begin{prop}\label{int-gr-res-cond}
Consider the subspace \smash{$W_{\lambda_{I, w_{3}'}}$}
of $\operatorname{COM}$ obtained by applying the
coefficients of~the vertex operators
\smash{$Y^{(g_{1}g_{2})^{-1}}_{P(z)}(u, x)$} for all $u\in V$ to
\smash{$\lambda_{I, w_{3}'}$}. Then \smash{$W_{\lambda_{I, w_{3}'}}$} equipped with
\smash{$Y^{(g_{1}g_{2})^{-1}}_{P(z)}$}~is a~generalized
$(g_{1}g_{2})^{-1}$-twisted $V$-module. Moreover,
\smash{$W_{\lambda_{I, w_{3}'}}$} is a~generalized
$(g_{1}g_{2})^{-1}$-twisted $V$-submodule of the generalized
$(g_{1}g_{2})^{-1}$-twisted $V$-module \smash{$W_{1}\hboxtr_{P(z)}W_{2}$}.
\end{prop}
\begin{proof}
We need only verify that the map
given by $w_{3}'\mapsto \lambda_{I, w_{3}'}$
commutes with the vertex operators on $W_{3}'$ and $\operatorname{COM}$.
This verification is straightforward and is omitted.
\end{proof}

Note that the conformal element $\omega$ is in the fixed point
subalgebra $V^{\langle g_{1}, g_{2}\rangle}$ of $V$ under the
group $\langle g_{1}, g_{2}\rangle$ generated by $g_{1}$ and $g_{2}$.
Also note that $W_{1}$ and $W_{2}$ are $V^{\langle g_{1}, g_{2}\rangle}$-modules.
Then~we~can apply the construction and results in~\cite{HLZ4}
to $V^{\langle g_{1}, g_{2}\rangle}$-modules $W_{1}$ and $W_{2}$. In~particular,
\cite[equation~(5.85)]{HLZ4} give a~vertex operator \smash{$Y'_{P(z)}(u, x)$} for $u\in \!V^{\langle g_{1}, g_{2}\rangle}$
acting on~\smash{$(W_{1}\otimes W_{2})^{*}$}. It~is clear that on $\operatorname{COM}$,
\smash{$Y'_{P(z)}(u, x)=Y^{(g_{1}g_{2})^{-1}}_{P(z)}(u, x)$} for $u\in V^{\langle g_{1}, g_{2}\rangle}$.
In particular, we have the vertex operator \smash{$Y'_{P(z)}(\omega, x)$} acting
on \smash{$(W_{1}\otimes W_{2})^{*}$} and is equal to \smash{$Y^{(g_{1}g_{2})^{-1}}_{P(z)}(\omega, x)$}
on $\operatorname{COMP}$. Taking the coefficient of $x^{-2}$ in \smash{$Y'_{P(z)}(\omega, x)$},
we obtain an operator \smash{$L'_{P(z)}(0)$} on \smash{$(W_{1}\otimes W_{2})^{*}$}, which is
equal to the coefficient of $x^{-2}$ in \smash{$Y^{(g_{1}g_{2})^{-1}}_{P(z)}(\omega, x)$} on
$\operatorname{COMP}$.

\begin{prop}
For $v\in V$, the $L(0)$-commutator formula
\begin{equation}\label{L-0-commu}
\big[L'_{P(z)}(0), Y^{(g_{1}g_{2})^{-1}}_{P(z)}(v, x)\big]
=x\frac{\rm d}{{\rm d}x}Y^{(g_{1}g_{2})^{-1}}_{P(z)}(v, x)+Y^{(g_{1}g_{2})^{-1}}_{P(z)}(L_{V}(0)v, x)
\end{equation}
holds on $\operatorname{COM}$.
\end{prop}
\begin{proof}
We have the $L(0)$-commutator formula
\[
\big[L_{W_{2}}(0), Y^{g_{2}}(v, x)\big]=
x\frac{\rm d}{{\rm d}x}Y^{g_{2}}(v, x)+Y^{g_{2}}(L_{V}(0)v, x).\]
For $\lambda\in \operatorname{COM}$, by the definition of
$L'_{P(z)}(0)$ \cite[formula (5.110)]{HLZ4}, we have
\[
(L'_{P(z)}(0)\lambda)(w_{1}\otimes w_{2})
=\lambda((zL_{W_{1}}(-1)+L_{W_{1}}(0))w_{1}\otimes w_{2})
+\lambda(w_{1}\otimes L_{W_{2}}(0)w_{2})\]
for $w_{1}\in W_{1}$ and $w_{2}\in W_{2}$.
Also, the series
$\lambda\bigl(w_{1}\otimes Y^{g_{2}}(v, z_{1})w_{2}\bigr)$ is absolutely convergent
on the region $|z|>|z_{1}|>0$ and can be analytically extended
to a~multivalued analytic function~${f_{1}(z_{1}; v, w_{1}, w_{2}; \lambda)}$
on $M^{1}(0, z)$
with a~preferred single-valued branch $f_{1}^{e}(z_{1}; v, w_{1}, w_{2}; \lambda)$
on the region $M^{1}_{0}(0, z)$.
Moreover, it is convergent to $f_{1}^{e}(z_{1}; v, w_{1}, w_{2}; \lambda)$
on the region $|z|>|z_{1}|>0$, $-\frac{3\pi}{2}<\arg (z_{1}-z)-\arg z<-\frac{\pi}{2}$
and to \smash{$f_{1}^{b^{-1}_{z_{1}, z}}(z_{1}; v, w_{1}, w_{2}; \lambda)$}
on the region $|z|>|z_{1}|>0$, $\frac{\pi}{2}<\arg (z_{1}-z)-\arg z<\frac{3\pi}{2}$.

By the definition of \smash{$Y^{(g_{1}g_{2})^{-1}}_{P(z)}$}, we have
\smash{$\bigl(\bigl(Y^{(g_{1}g_{2})^{-1}}_{P(z)}\bigr)^{o}(v, z_{1})\lambda\bigr)(w_{1}\otimes w_{2})
=f_{1}^{e}(z_{1}; v, w_{1}, w_{2}; \lambda)$}
on the region $|z_{1}|>|z|$, \smash{$|{\arg (z_{1}-z)-\arg z_{1}}|<\frac{\pi}{2}$}.
Then we have
\begin{align*}
&\bigl(L_{P(z)}'(0)\bigl(Y^{(g_{1}g_{2})^{-1}}_{P(z)}\bigr)^{o}(v, z_{1})\lambda\bigr)(w_{1}\otimes w_{2})\nn
&\qquad=\bigl(\bigl(Y^{(g_{1}g_{2})^{-1}}_{P(z)}\bigr)^{o}(v, z_{1})\lambda\bigr)
\bigl((zL_{W_{1}}(-1)+L_{W_{1}}(0))w_{1}\otimes w_{2}\bigr)\nn
& \phantom{=}{}\qquad +\bigl(\bigl(Y^{(g_{1}g_{2})^{-1}}_{P(z)}\bigr)^{o}(v, z_{1})\lambda\bigr)(w_{1}\otimes L_{W_{2}}(0)w_{2})\nn
& \qquad =f_{1}^{e}(z_{1}; v, (zL_{W_{1}}(-1)+L_{W_{1}}(0))w_{1}, w_{2}; \lambda)
+f_{1}^{e}(z_{1}; v, w_{1}, L_{W_{2}}(0)w_{2}; \lambda)
\end{align*}
on the region $|z_{1}|>|z|$, $|{\arg (z_{1}-z)-\arg z_{1}}|<\frac{\pi}{2}$.
Also we have
\[
\bigl(\bigl(Y^{(g_{1}g_{2})^{-1}}_{P(z)}\bigr)^{o}(v, z_{1})L_{P(z)}'(0)\lambda\bigr)(w_{1}\otimes w_{2})
=f_{1}^{e}(z_{1}; v, w_{1}, w_{2}; L_{P(z)}'(0)\lambda).\]
But on the region $|z|>|z_{1}|>0$, $-\frac{3\pi}{2}<\arg (z_{1}-z)-\arg z<-\frac{\pi}{2}$,
we have
\begin{align*}
&f_{1}^{e}(z_{1}; v, w_{1}, w_{2}; L_{P(z)}'(0)\lambda)\nn
&\qquad=(L_{P(z)}'(0)\lambda)\bigl(w_{1}\otimes Y^{g_{2}}(v, z_{1})w_{2}\bigr)\nn
&\qquad=\lambda\bigl((zL_{W_{1}}(-1)+L_{W_{1}}(0))w_{1}\otimes Y^{g_{2}}(v, z_{1})w_{2}\bigr)
+\lambda\bigl(w_{1}\otimes L_{W_{2}}(0)Y^{g_{2}}(v, z_{1})w_{2}\bigr)\nn
&\qquad=\lambda\bigl((zL_{W_{1}}(-1)+L_{W_{1}}(0))w_{1}\otimes Y^{g_{2}}(v, z_{1})w_{2}\bigr)
+\lambda\bigl(w_{1}\otimes Y^{g_{2}}(v, z_{1})L_{W_{2}}(0)w_{2}\bigr)\nn
&\phantom{=}{} \qquad +\lambda\left(w_{1}\otimes z_{1}\frac{\rm d}{{\rm d}z_{1}}Y^{g_{2}}(v, z_{1})w_{2}\right)
+\lambda\bigl(w_{1}\otimes Y^{g_{2}}(L_{V}(0)v, z_{1})w_{2}\bigr)\nn
&\qquad=f_{1}^{e}(z_{1}; v, (zL_{W_{1}}(-1)+L_{W_{1}}(0))w_{1}, w_{2}; \lambda)
+f_{1}^{e}(z_{1}; v, w_{1}, L_{W_{2}}(0)w_{2}; \lambda)\nn
&\phantom{=}{} \qquad +z_{1}\frac{\rm d}{{\rm d}z_{1}}f_{1}^{e}(z_{1}; v, w_{1}, w_{2}; \lambda)
+f_{1}^{e}(z_{1}; L_{V}(0)v, w_{1}, w_{2}; \lambda).
\end{align*}
Then on the region $|z_{1}|>|z|$, $|{\arg (z_{1}-z)-\arg z_{1}}|<\frac{\pi}{2}$,
we have
\begin{align*}
&\bigl(\bigl(Y^{(g_{1}g_{2})^{-1}}_{P(z)}\bigr)^{o}(v, z_{1})L_{P(z)}'(0)\lambda\bigr)(w_{1}\otimes w_{2})\nn
&\qquad=f_{1}^{e}(z_{1}; v, w_{1}, w_{2}; L_{P(z)}'(0)\lambda)\nn
&\qquad=f_{1}^{e}(z_{1}; v, (zL_{W_{1}}(-1)+L_{W_{1}}(0))w_{1}, w_{2}; \lambda)
+f_{1}^{e}(z_{1}; v, w_{1}, L_{W_{2}}(0)w_{2}; \lambda)\nn
& \phantom{\qquad=}{} +z_{1}\frac{\rm d}{{\rm d}z_{1}}f_{1}^{e}(z_{1}; v, w_{1}, w_{2}; \lambda)
+f_{1}^{e}(z_{1}; L_{V}(0)v, w_{1}, w_{2}; \lambda)\nn
&\qquad=\bigl(\bigl(Y^{(g_{1}g_{2})^{-1}}_{P(z)}\bigr)^{o}(v, z_{1})\lambda\bigr)((zL_{W_{1}}(-1)+L_{W_{1}}(0))w_{1}
\otimes w_{2})\nn
& \phantom{\qquad=}{} +\bigl(\bigl(Y^{(g_{1}g_{2})^{-1}}_{P(z)}\bigr)^{o}(v, z_{1})\lambda\bigr)(w_{1}\otimes L_{W_{2}}(0)w_{2})
+z_{1}\frac{\rm d}{{\rm d}z_{1}}
\bigl(\bigl(Y^{(g_{1}g_{2})^{-1}}_{P(z)}\bigr)^{o}(v, z_{1})\lambda\bigr)(w_{1}\otimes w_{2})\nn
& \phantom{\qquad=}{}
+\bigl(\bigl(Y^{(g_{1}g_{2})^{-1}}_{P(z)}\bigr)^{o}(L_{V}(0)v, z_{1})\lambda\bigr)(w_{1}\otimes w_{2})\nn
&\qquad=\bigl(L_{P(z)}'(0)\bigl(Y^{(g_{1}g_{2})^{-1}}_{P(z)}\bigr)^{o}(v, z_{1})\lambda\bigr)(w_{1}
\otimes w_{2})+
\left(z_{1}\frac{\rm d}{{\rm d}z_{1}}\bigl(Y^{(g_{1}g_{2})^{-1}}_{P(z)}\bigr)^{o}(v, z_{1})\lambda\right)(w_{1}\otimes w_{2})
\nn
& \phantom{\qquad=}{} +\bigl(\bigl(Y^{(g_{1}g_{2})^{-1}}_{P(z)}\bigr)^{o}(L_{V}(0)v, z_{1})\lambda\bigr)(w_{1}\otimes w_{2}).
\end{align*}
Thus we obtain
\[
\big[\bigl(Y^{(g_{1}g_{2})^{-1}}_{P(z)}\bigr)^{o}(v, z_{1}), L_{P(z)}'(0)\big]
=z_{1}\frac{\rm d}{{\rm d}z_{1}}\bigl(Y^{(g_{1}g_{2})^{-1}}_{P(z)}\bigr)^{o}(v, z_{1})
+\bigl(Y^{(g_{1}g_{2})^{-1}}_{P(z)}\bigr)^{o}(L_{V}(0)v, z_{1}),\]
which is equivalent to~\eqref{L-0-commu}.
\end{proof}

Because~\eqref{L-0-commu} holds, the coefficient
\smash{$\bigl(Y^{(g_{1}g_{2})^{-1}}_{P(z)}\bigr)_{n, k}(u)$} of \smash{$x^{-n-1}(\log x)^{k}$}
in \smash{$Y^{(g_{1}g_{2})^{-1}}_{P(z)}(u, x)$} must have weight $\wt u-n-1$
when $u \in V$ is homogeneous. In particular, if $\lambda\in \operatorname{COM}$
is a~finite sum of generalized eigenvectors of $L_{P(z)}'(0)$,
the subspace $W_{\lambda}$ of $\operatorname{COM}$ obtained by
applying the~coefficients of the vertex operators
\smash{$Y^{(g_{1}g_{2})^{-1}}_{P(z)}(u, x)$} for all $u\in V$ to
$\lambda$ must be a~direct sum of generalized eigenspaces of $L_{P(z)}'(0)$, that is,
$W_{\lambda}=\amalg_{n\in\C}(W_{\lambda})_{[n]}$, where for $n\in \C$,
$(W_{\lambda})_{[n]}$ is the generalized eigenspace of $L_{P(z)}'(0)$
with eigenvalue $n$.

Proposition~\ref{int-gr-res-cond} together with Assumption~\ref{assum}
 in particular says that in the case that $\mathcal{C}$ is a~category of
$g$-twisted $V$-modules with $G$-actions, for \smash{$\lambda=\lambda_{I, w_{3}'}$},
there is a~$G$-action (in particular, an action of $(g_{1}g_{2})^{-1}$)
on \smash{$W_{\lambda_{I, w_{3}'}}$} such that for $n\in \C$,
\[
(W_{\lambda_{I, w_{3}'}})_{[n]}=\coprod_{\alpha\in \C/\Z}
(W_{\lambda_{I, w_{3}'}})_{[n]}^{[\alpha]},
\]
where for $\alpha\in \C/\Z$, \smash{$(W_{\lambda_{I, w_{3}'}})_{[n]}^{[\alpha]}$} is the
generalized eigenspace of
the action $(g_{1}g_{2})^{-1}$ with eigenvalues ${\rm e}^{2\pi \alpha {\rm i}}$ for $\alpha\in \C/\Z$ and
\[(g_{1}g_{2})^{-1}Y_{W_{\lambda_{I, w_{3}'}}}(v, x)w=Y_{W_{\lambda_{I, w_{3}'}}}\bigl((g_{1}g_{2})^{-1}v, x\bigr)(g_{1}g_{2})^{-1}w\]
for $v\in V$ and \smash{$w\in W_{\lambda_{I, w_{3}'}}$}.

Moreover, such an action of $(g_{1}g_{2})^{-1}$ can be extended to the action of $(g_{1}g_{2})^{-1}$
on $W_{1}\hboxtr_{P(z)}W_{2}$. It is possible that the action of $(g_{1}g_{2})^{-1}$ on $W_{1}\hboxtr_{P(z)}W_{2}$
can be further extended to subspaces of $\operatorname{COM}$
larger than on $W_{1}\hboxtr_{P(z)}W_{2}$ such that~\eqref{g1*g2-comp} holds. By Zorn's lemma, there exists a~maximal subspace
of $\operatorname{COM}$ containing $W_{1}\hboxtr_{P(z)}W_{2}$
with action of $(g_{1}g_{2})^{-1}$
such that the restriction of the action to $W_{1}\hboxtr_{P(z)}W_{2}$ is equal to the action of
$(g_{1}g_{2})^{-1}$ on $W_{1}\hboxtr_{P(z)}W_{2}$ and
\begin{equation}\label{g1*g2-comp}
(g_{1}g_{2})^{-1}Y^{(g_{1}g_{2})^{-1}}_{P(z)}(v, x)w=Y^{(g_{1}g_{2})^{-1}}_{P(z)}\bigl((g_{1}g_{2})^{-1}v, x\bigr)(g_{1}g_{2})^{-1}w
\end{equation}
holds for $v\in V$ and $w$ in this maximal subspace.

Using the vertex operators
\smash{$Y^{(g_{1}g_{2})^{-1}}_{P(z)}(u, x)$} for $u\in V$ and \smash{$L'_{P(z)}(0)$}
and motivated by Pro\-position~\ref{int-gr-res-cond} (especially, the discussion above on the
action of $(g_{1}g_{2})^{-1}$ on \smash{$W_{\lambda_{I, w_{3}'}}$}),
we also introduce the following condition for
$\lambda\in \operatorname{COM}$:

{\bf $\boldsymbol{P(z)}$-local-grading-restriction condition.}
\begin{itemize}
\item[(a)] The \textit{$P(z)$-grading condition}: $\lambda$ is a~(finite) sum
of generalized eigenvectors for the opera\-tor~$L'_{P(z)}(0)$.

\item[(b)]	Let $W_{\lambda}$ be the smallest subspace
of $(W_1\otimes W_2)^*$ containing $\lambda$
and stable under the \mbox{action} of all the coefficients of
the vertex operators \smash{$Y^{(g_{1}g_{2})^{-1}}_{P(z)}(u, x)$}
for all $u\in V$. Then
\smash{$\operatorname{dim}(W_{\lambda})_{[n]}<\infty$}, $ (W_{\lambda})_{[n]}=0$, for $\Re(n)$
sufficiently negative.
\end{itemize}

We denote the subspace of $P(z)$-local grading restricted functionals in
$\operatorname{COM}\subset
(W_1\otimes W_2)^*$~as ${\operatorname{LGR}_{P(z)}((W_1\otimes W_2)^*)}$,
or $\operatorname{LGR}$ for short. Clearly, the space LGR is closed
under the action~of~\smash{$Y^{(g_{1}g_{2})^{-1}}_{P(z)}(u, x)$}, $u\in V$.

\begin{thm}\label{comp-main}
For $\lambda$ satisfying the $P(z)$-compatibility
condition and the $P(z)$-local-grading-restriction condition,
the graded space $W_\lambda$ equipped with the twisted vertex operator map
\smash{$Y^{(g_{1}g_{2})^{-1}}_{P(z)}$} is a~grading-restricted
$(g_1g_2)^{-1}$-twisted generalized $V$-module. If for every element
$\lambda \!\in\! {(W_{1}\!\otimes\! W_{2})^{*}}\!$ satisfying both the $P(z)$-compatibility
condition and the $P(z)$-local grading restriction condition, $W_{\lambda}$ is a~$(g_{1}g_{2})^{-1}$-twisted generalized $V$-submodule of a~$(g_{1}g_{2})^{-1}$-twisted generalized $V$-module in $\operatorname{COMP}$ whose contragredient is
in $\mathcal{C}$
{\rm(}this holds in particular if $\mathcal{C}$ is the category of all grading-restricted
$g$-twisted generalized $V$-modules without $g$-actions for $g\in G$ or
 the category of all grading-restricted
$g$-twisted generalized $V$-modules with $G$-actions for $g\in G${\rm)},
then an
element $\lambda\in (W_{1}\otimes W_{2})^{*}$ is
in $W_{1}\hboxtr_{P(z)}W_{2}$ if and only if $\lambda$ satisfies
the $P(z)$-compatibility condition and the $P(z)$-local-grading-restriction
condition. In other words,
$
W_{1}\hboxtr_{P(z)}W_{2}=\operatorname{COMP}\cap\operatorname{LGR}$.
\end{thm}
\begin{proof}
The identity property
follows immediately from the definition of the twisted vertex oper\-ator
map \smash{$Y^{(g_{1}g_{2})^{-1}}_{P(z)}$}. The
$L(0)$-grading condition follows from
the $P(z)$-local-grading-restriction property.
The $L(-1)$-derivative property
follows from the
definition of \smash{$Y^{(g_{1}g_{2})^{-1}}_{P(z)}$} and the
$L(-1)$-derivative property of the twisted vertex operator
$Y^{g_{2}}$. We omit the details of the proofs of these properties.

We prove the equivariance property for
\smash{$Y^{(g_{1}g_{2})^{-1}}_{P(z)}$} now. It is equivalent to
\begin{align}\label{ascasc}
\bigl(\bigl(Y_{P(z)}^{(g_1g_2)^{-1}}\bigr)^o\bigr)^{p+1}
(g_1g_2u,z_{1})\tilde{\lambda}
=\bigl(\bigl(Y_{P(z)}^{(g_1g_2)^{-1}}\bigr)^o\bigr)^{p}(u,z_{1})
\tilde{\lambda}
\end{align}
for $u\in V$, $\tilde{\lambda}\in W_\lambda$.
By the definition of $\bigl(Y_{P(z)}^{(g_1g_2)^{-1}}\bigr)^o$,
we know that for $w_1\in W_1$, $w_2\in W_2$,
\[
\bigl(\bigl(Y_{P(z)}^{(g_1g_2)^{-1}}\bigr)^o(u,z)
\tilde{\lambda}\bigr)(w_{1}\otimes w_{2})
\]
is absolutely convergent on the region $|z_1|>|z|$,
$|{\arg(z_1-z)-\arg(z_1)}|<\frac{\pi}{2}$ to
$f_{1}^{e}\bigl(z_{1}; u, w_{1}, \allowbreak w_{2}; \tilde{\lambda}\bigr)$.
Let $b_{3}$ be the homotopy class containing a~loop given by a~circle centered at $0$
with radius larger than $|z|$ in the counterclockwise
direction on the complex $z_{1}$ plane. Then~\eqref{ascasc} is equivalent~to
\begin{equation}\label{ascasc-1}
f_{1}^{b_{3}}\bigl(z_{1}; g_{1}g_{2}u, w_{1}, w_{2}; \tilde{\lambda}\bigr)
=f_{1}^{e}\bigl(z_{1}; u, w_{1}, w_{2}; \tilde{\lambda}\bigr)
\end{equation}
for $u\in V$, $w_1\in W_1$, $w_2\in W_2$ and
$\tilde{\lambda}\in W_\lambda$.

By the equivariance property for $Y^{g_{2}}$ and the convergence of
$\tilde{\lambda}(w_{1}\otimes Y^{g_{2}}(v, z_{1}) w_{2})$
to the branch $f_{1}^{e}\bigl(z_{1}; v, w_{1}, w_{2}; \tilde{\lambda}\bigr)$ for $v\in V$
on the region $|z|>|z_{1}|>0$, $-\frac{3\pi}{2}
<\arg (z_{1}-z)-\arg z_{1}<-\frac{\pi}{2}$,
we obtain
\[f_{1}^{b_{z_{1}, 0}}\bigl(z_{1}; g_{2}v, w_{1}, w_{2}; \tilde{\lambda}\bigr)
=f_{1}^{e}\bigl(z_{1}; v, w_{1}, w_{2}; \tilde{\lambda}\bigr).\]
Using the equivariance property for $Y^{g_{2}}$ and the
convergence of
$\tilde{\lambda}\bigl(Y^{g_{1}}(v, z_{1}-z)w_{1}\otimes w_{2}\bigr)$,
we obtain similarly
\[f_{1}^{b_{z_{1}, z}}\bigl(z_{1}; g_{1}v, w_{1}, w_{2}; \tilde{\lambda}\bigr)
=f_{1}^{e}\bigl(z_{1}; v, w_{1}, w_{2}; \tilde{\lambda}\bigr).\]
Applying any homotopy class $b$ of loops
in the $z_{1}$ complex plane with $z$ and $0$ deleted to both sides of
this equality, we obtain
\[f_{1}^{b_{z_{1}, z}b}\bigl(z_{1}; g_{1}v, w_{1}, w_{2}; \tilde{\lambda}\bigr)
=f_{1}^{b}\bigl(z_{1}; v, w_{1}, w_{2}; \tilde{\lambda}\bigr)\]
for $v\in V$.
Then we have
\begin{align*}
f_{1}^{b_{z_{1}, z}b_{z_{1}, 0}}\bigl(z_{1};
g_{1}g_{2}u, w_{1}, w_{2}; \tilde{\lambda}\bigr)
&=f_{1}^{b_{z_{1}, 0}}(z_{1};
g_{2}u, w_{1}, w_{2}; \tilde{\lambda})=f_{1}^{e}\bigl(z_{1}; u, w_{1}, w_{2}; \tilde{\lambda}\bigr).
\end{align*}
But it is easy to see that $b_{z_{1}, z}b_{z_{1}, 0}=b_{3}$. Thus we
have proved~\eqref{ascasc-1}.

For $u, v\in V$,
$w\in W_{\lambda}$, $w'\in W_{\lambda}'$,
by Proposition~\ref{COM-comm},
we have
\begin{align}\label{comp-main-1}
&(x_1-x_2)^{M}\big\langle w',
Y^{(g_{1}g_{2})^{-1}}_{P(z)}(u, x_1)
Y^{(g_{1}g_{2})^{-1}}_{P(z)}(v, x_2)w\big\rangle\nn
& \qquad =(x_1-x_2)^{M}
\big\langle w',
Y^{(g_{1}g_{2})^{-1}}_{P(z)}(v, x_{2})
Y^{(g_{1}g_{2})^{-1}}_{P(z)}(u, x_{1})w\big\rangle,
\end{align}
where $M\in \Z_{+}$ depends only on $u$ and $v$.
Since $W_{\lambda}$ is lower bounded, by~\eqref{comp-main-1},
the left-hand side of~\eqref{comp-main-1}
has only finitely many terms in complex powers of $x_{1}$,
$x_{2}$ and integer powers of $\log x_{1}$, $\log x_{2}$.
Then
\[\big\langle w',
Y^{(g_{1}g_{2})^{-1}}_{P(z)}(u, x_1)
Y^{(g_{1}g_{2})^{-1}}_{P(z)}(v, x_2)w\big\rangle\]
is equal to this finite sum multiplied by
$(x_1-x_2)^{-M}$, which
is expanded in nonnegative powers of $x_{2}$.
Thus we have a~multivalued function of the form
\[f(z_1, z_2) = \sum_{i,
j, k, l = 0}^N a_{ijkl}z_1^{m_i}z_2^{n_j}(\operatorname{log}z_1)^k({\rm
log}z_2)^l(z_1 - z_2)^{-M},\]
where $m_{i}, n_{i}\in \C$ for $i=0, \dots, N$, with the
preferred branch
\[f^{e}(z_1, z_2)=\sum_{i,
j, k, l = 0}^N a_{ijkl}{\rm e}^{m_i\log z_{1}}{\rm e}^{n_j\log z_{2}}
(\log z_1)^k(\log z_2)^l(z_1 - z_2)^{-M}\]
such that
\[\big\langle w',
Y^{(g_{1}g_{2})^{-1}}_{P(z)}(u, z_1)
Y^{(g_{1}g_{2})^{-1}}_{P(z)}(v, z_2)w\big\rangle\]
is absolutely convergent on the region $|z_{1}|>|z_{2}|>0$
to $f^{e}(z_1, z_2)$.
From~\eqref{comp-main-1}, we also obtain the
commutativity, that is,
\[\big\langle w', Y^{(g_{1}g_{2})^{-1}}_{P(z)}(v, x_2)
Y^{(g_{1}g_{2})^{-1}}_{P(z)}(u, x_1)w\big\rangle\]
is absolutely convergent on the region $|z_{2}|>|z_{1}|>0$
to $f^{e}(z_1, z_2)$.

We now prove the associativity for \smash{$Y^{(g_{1}g_{2})^{-1}}_{P(z)}$}.
Since the associativity for \smash{$Y^{(g_{1}g_{2})^{-1}}_{P(z)}$} is equivalent to the associativity for \smash{$\bigl(Y^{(g_{1}g_{2})^{-1}}_{P(z)}\bigr)^{o}$}, we prove
this associativity. For $u, v\in V$,
$w_{1}\in W_{1}$, $w_{2}\in W_{2}$ and
$\tilde{\lambda}\in W_{\lambda}$,
by Proposition~\ref{COM-comm},
\begin{align}\label{comp-main-1.1}
\bigl(\bigl(Y^{(g_{1}g_{2})^{-1}}_{P(z)}\bigr)^{o}(v, z_2)
\bigl(Y^{(g_{1}g_{2})^{-1}}_{P(z)}\bigr)^{o}(u, z_1)\tilde{\lambda}\bigr)
(w_{1}\otimes w_{2})
\end{align}
is absolutely convergent
on the region $|z_{1}|>|z_{2}|>|z|$ and its sum is equal to
\begin{equation}\label{comp-main-1.6}
f_{2}^{e}\bigl(z_{1}, z_{2}; u, v, w_{1}, w_{2}; \tilde{\lambda}\bigr)
\end{equation}
on the region $|z_{1}|>|z_{2}|>|z|$, $|{\arg (z_{1}-z)-\arg z_{1}}|
<\frac{\pi}{2}$, $|{\arg (z_{2}-z)-\arg z_{2}}|
<\frac{\pi}{2}$. By the~definition of
\smash{$\bigl(Y^{(g_{1}g_{2})^{-1}}_{P(z)}\bigr)^{o}$},
\smash{$
\bigl(\bigl(Y^{(g_{1}g_{2})^{-1}}_{P(z)}\bigr)^{o}
((Y_{V})_{n}(u)v, z_2)\tilde{\lambda}\bigr)(w_{1}\otimes w_{2})$}
is absolutely convergent
on the region $|z_{2}|>|z|$ and its sum is equal to
\begin{equation}\label{comp-main-1.3}
f_{1}^{e}\bigl(z_{2}; (Y_{V})_{n}(u)v, w_{1}, w_{2}; \tilde{\lambda}\bigr)
\end{equation}
on the region $|z_{2}|>|z|$, $|{\arg (z_{2}-z)-\arg z_{2}}|<\frac{\pi}{2}$.
Also
$\tilde{\lambda}\bigl(w_{1}\otimes Y^{g_{2}}
((Y_{V})_{n}(u)v, z_2)w_{2}\bigr)$
is absolutely convergent on the region $|z|>|z_{2}|>0$ and its sum is equal to
\eqref{comp-main-1.3} on the region~${|z|>|z_{2}|>0}$,
$-\frac{3\pi}{2}<\arg (z_{2}-z)-\arg z<-\frac{\pi}{2}$ and is equal to
\[f_{1}^{b_{z_{2}, z}^{-1}}\bigl(z_{2}; (Y_{V})_{n}(u)v, w_{1}, w_{2};
\tilde{\lambda}\bigr)\]
on the region $|z|>|z_{2}|>0$,
$\frac{\pi}{2}<\arg (z_{2}-z)-\arg z<\frac{3\pi}{2}$.

By the $P(z)$-compatibility condition for
$\tilde{\lambda}$,
\smash{$\tilde{\lambda}\bigl(w_{1}\otimes Y^{g_{2}}(u, z_1)
Y^{g_{2}}(v, z_2)w_{2}\bigr)$}
is absolutely convergent to~\eqref{comp-main-1.6}
on the region $|z|>|z_{1}|>|z_{2}|>0$,
$-\frac{3\pi}{2}<\arg (z_{1}-z)-\arg z, \arg (z_{2}-z)-\arg z
<-\frac{\pi}{2}$ and to
\begin{equation}\label{comp-main-1.7}
f_{2}^{b_{z_{1}, z}^{-1}b_{z_{2}, z}^{-1}}
\bigl(z_{1}, z_{2}; u, v, w_{1}, w_{2}; \tilde{\lambda}\bigr)
\end{equation}
on the region $|z|>|z_{1}|>|z_{2}|>0$,
$\frac{\pi}{2}<\arg (z_{1}-z)-\arg z, \arg (z_{2}-z)-\arg z
<\frac{3\pi}{2}$, respectively.
By the associativity of the twisted vertex operator map
$Y^{g_{2}}$,
\begin{align}
&\sum_{k=0}^{K}\sum_{m\in \C}\bigg(\sum_{n\in \Z}
\tilde{\lambda}\bigl(w_{1}\otimes Y^{g_{2}}_{m, k}
((Y_{V})_{n}(u)v)w_{2}\bigr)(z_{1}-z_{2})^{-n-1}\bigg){\rm e}^{(-m-1)\log z_{2}}(\log z_{2})^{k}\nn
& \qquad =\sum_{k=0}^{K}\sum_{m\in \C}
\tilde{\lambda}\bigl(w_{1}\otimes Y^{g_{2}}_{m, k}
(Y_{V}(u, z_1-z_{2})v, z_2)w_{2}\bigr){\rm e}^{(-m-1)\log z_{2}}(\log z_{2})^{k}\nn
&\qquad=\tilde{\lambda}\bigl(w_{1}\otimes Y^{g_{2}}
(Y_{V}(u, z_1-z_{2})v, z_2)w_{2}\bigr) =\tilde{\lambda}\bigl(w_{1}\otimes Y^{g_{2}}(u, z_1)
Y^{g_{2}}(v, z_2)w_{2}\bigr)\label{comp-main-1.5}
\end{align}
on the region $|z|>|z_{1}|>|z_{2}|>|z_{1}-z_{2}|>0$,
$|{\arg z_{1}-\arg z_{2}}|<\frac{\pi}{2}$. But since $\tilde{\lambda}$ satisfies the $P(z)$-compatibility condition,
we know that the right-hand side of~\eqref{comp-main-1.5} is
equal to~\eqref{comp-main-1.6} on the re\-gion~${|z|>|z_{1}|>|z_{2}|>0}$,
$-\frac{3\pi}{2}<\arg (z_{1}-z)-\arg z, \arg (z_{2}-z)-\arg z
<-\frac{\pi}{2}$ and is equal to~\eqref{comp-main-1.7} on the region
$|z|>|z_{1}|>|z_{2}|>|z_{1}-z_{2}|>0$,
$-\frac{3\pi}{2}<\arg (z_{1}-z)-\arg z, \arg (z_{2}-z)-\arg z
<-\frac{\pi}{2}$. Moreover, by 1(b) of the $P(z)$-compatibility condition,
the component-isolated singularity at~${z_{2}=\infty}$, $z_{1}-z_{2}=0$ of the
multivalued analytic function \smash{$f\bigl(z_{1}, z_{2}; u, v, w_{1}, w_{2}; \tilde{\lambda}\bigr)$} is regular. Thus
it can be expanded as a~multiseries in powers of $z_{2}$, $z_{1}-z_{2}$ and nonnegative integral powers of
$\log z_{2}$ and $\log (z_{1}-z_{2})$ on the region $|z|>|z_{2}|>|z_{1}-z_{2}|>0$. By~\eqref{comp-main-1.5},
the single-valued branches~\eqref{comp-main-1.6} and~\eqref{comp-main-1.7} of
$f\bigl(z_{1}, z_{2}; u, v, w_{1}, w_{2}; \tilde{\lambda}\bigr)$ are expanded
as iterated series of the same form on the region $|z|>|z_{1}|>|z_{2}|>|z_{1}-z_{2}|>0$,
$|{\arg z_{1}-\arg z_{2}}|<\frac{\pi}{2}$. Thus the corresponding multiseries of this iterated series must be
the expansion of a~single-valued branch of $f\bigl(z_{1}, z_{2}; u, v, w_{1}, w_{2}; \tilde{\lambda}\bigr)$
and is thus absolutely convergent as a~multiseries on the region~${|z|>|z_{2}|>|z_{1}-z_{2}|>0}$.
In particular, the other iterated series
\begin{align}
&\sum_{n\in \Z}\left(\sum_{k=0}^{K}\sum_{m\in \C}
\tilde{\lambda}\bigg(w_{1}\otimes Y^{g_{2}}_{m, k}
((Y_{V})_{n}(u)v)w_{2}\bigg){\rm e}^{(-m-1)\log z_{2}}(\log z_{2})^{k}\right)(z_{1}-z_{2})^{-n-1}\nn
& \qquad =\sum_{n\in \Z}
\tilde{\lambda}(w_{1}\otimes Y^{g_{2}}
((Y_{V})_{n}(u)v, z_2)w_{2})(z_{1}-z_{2})^{-n-1}\label{comp-main-1.5.1}
\end{align}
of this multiseries is also absolutely convergent on the region $|z|>|z_{2}|>|z_{1}-z_{2}|>0$ to the same
single-valued analytic functions.
Thus we have proved that~\eqref{comp-main-1.5.1} is
absolutely convergent on the region $|z|>|z_{1}|>|z_{2}|>0$
and its sum is equal to~\eqref{comp-main-1.6}
 on the region
$|z|>|z_{1}|>|z_{2}|>0$,
$-\frac{3\pi}{2}<\arg (z_{1}-z)-\arg z, \arg (z_{2}-z)-\arg z
<-\frac{\pi}{2}$, $|{\arg z_{1}-\arg z_{2}}|<\frac{\pi}{2}$
and is equal to~\eqref{comp-main-1.7} on the region $|z|>|z_{1}|>|z_{2}|>0$,
$\frac{\pi}{2}<\arg (z_{1}-z)-\arg z, \arg (z_{2}-z)-\arg z
<\frac{3\pi}{2}$, $|{\arg z_{1}-\arg z_{2}}|<\frac{\pi}{2}$.
Then by the definition of \smash{$\bigl(Y^{(g_{1}g_{2})^{-1}}_{P(z)}\bigr)^{o}$},
\begin{align}\label{comp-main-2}
&\bigl(\bigl(Y^{(g_{1}g_{2})^{-1}}_{P(z)}\bigr)^{o}
(Y_{V}(u, z_1-z_{2})v, z_2)
\tilde{\lambda}\bigr)(w_{1}\otimes w_{2})\nn
&\qquad=\sum_{n\in \Z}
\bigl(\bigl(Y^{(g_{1}g_{2})^{-1}}_{P(z)}\bigr)^{o}
\bigl((Y_{V})^{n}(u)v, z_2\bigr)
\tilde{\lambda}\bigr)(w_{1}\otimes w_{2})(z_{1}-z_{2})^{-n-1}\nn
&\qquad=\sum_{n\in \Z}
f_{1}^{e}\bigl(z_{2}; (Y_{V})_{n}(u)v, w_{1}, w_{2}; \tilde{\lambda}\bigr)
(z_{1}-z_{2})^{-n-1}
\end{align}
is in fact the expansion of~\eqref{comp-main-1.6}
as a~Laurent series in $z_{1}-z_{2}$ near $z_{1}-z_{2}=0$
and then to expand
the coefficients as a~series in powers of $z_{2}$ and $\log z_{2}$
near $z_{2}=\infty$.
Thus we have shown that~\eqref{comp-main-1.1} and the left-hand side of~\eqref{comp-main-2}
are convergent on the region $|z_{1}|>|z_{2}|>|z|$
and $|z_{2}|>|z_{1}-z_{2}|, |z|$, respectively, and their sums are
equal to~\eqref{comp-main-1.6} on the region $|z_{1}|>|z_{2}|>|z|$,
$|{\arg (z_{1}-z)-\arg z_{1}}|, |{\arg (z_{2}-z)-\arg z_{2}}|<\frac{\pi}{2}$, $|z_{2}|>|z_{1}-z_{2}|, |z|$,
$|{\arg (z_{2}-z)-\arg z_{2}}|<\frac{\pi}{2}$,
$|{\arg z_{1}-\arg z_{2}}|<\frac{\pi}{2}$.

Since $\lambda$ satisfies the $P(z)$-local-grading-restriction condition,
\smash{$W_{\lambda}^{*}=\overline{W_{\lambda}'}$}. Then
we have a~linear map $I\colon W_{1}\otimes W_{2}\to W_{\lambda}^{*}=\overline{W_{\lambda}'}$
given by $\big\langle I(w_{1}\otimes w_{2}), \tilde{\lambda}\big\rangle=\tilde{\lambda}(w_{1}\otimes w_{2})$
for $w_{1}\in W_{1}$, $w_{2}\in W_{2}$ and $\tilde{\lambda}\in W_{\lambda}$. Then
\eqref{comp-main-1.1} is
absolutely convergent on the region $|z_{1}|>|z_{2}|>|z|$ means that the
series
\smash{$\big\langle I(w_{1}\otimes w_{2}), \bigl(Y^{(g_{1}g_{2})^{-1}}_{P(z)}\bigr)^{o}(v, z_2)
\bigl(Y^{(g_{1}g_{2})^{-1}}_{P(z)}\bigr)^{o}(u, z_1)\tilde{\lambda}\big\rangle$}
is absolutely convergent on the same region. Let \smash{$L_{P(z)}(0)$} be the adjoint operator of
\smash{$L_{P(z)}'(0)$}.
Then by~\eqref{L-0-commu}, we have
\begin{align}
&\big\langle x^{L_{P(z)}(0)} I(w_{1}\otimes w_{2}), \bigl(Y^{(g_{1}g_{2})^{-1}}_{P(z)}\bigr)^{o}(v, x_2)
\bigl(Y^{(g_{1}g_{2})^{-1}}_{P(z)}\bigr)^{o}(u, x_1)\tilde{\lambda}\big\rangle\nn
& \qquad =\big\langle I(w_{1}\otimes w_{2}), x^{L_{P(z)}'(0)} \bigl(Y^{(g_{1}g_{2})^{-1}}_{P(z)}\bigr)^{o}(v, x_2)
\bigl(Y^{(g_{1}g_{2})^{-1}}_{P(z)}\bigr)^{o}(u, x_1)\tilde{\lambda}\big\rangle\nn
& \qquad =\big\langle I(w_{1}\otimes w_{2}), \bigl(Y^{(g_{1}g_{2})^{-1}}_{P(z)}\bigr)^{o}\bigl(x^{-L_{V}(0)}v, x_2x^{-1}\bigr)
\bigl(Y^{(g_{1}g_{2})^{-1}}_{P(z)}\bigr)^{o}\nn
&\phantom{\qquad =}{}\times\bigl(x^{-L_{V}(0)}u, x_1x^{-1}\bigr)x^{L_{P(z)}'(0)}\tilde{\lambda}\big\rangle.\label{comp-main-2.1}
\end{align}
Substituting ${\rm e}^{n\log \xi}$, $\log \xi$, ${\rm e}^{n\log z_{1}}$, $\log z_{1}$,
${\rm e}^{n\log z_{2}}$, $\log z_{2}$ for $x^{n}$, $\log x$, $x_{1}^{n}$, $\log x_{1}$, $x_{2}^{n}$, $\log x_{2}$,
respectively, in~\eqref{comp-main-2.1}, we see that the left-hand side~\eqref{comp-main-2.1}
becomes
\begin{equation}\label{comp-main-2.2}
\bigl\langle {\rm e}^{(\log \xi)L_{P(z)}(0)} I(w_{1}\otimes w_{2}), \bigl(Y^{(g_{1}g_{2})^{-1}}_{P(z)}\bigr)^{o}(v, z_2)
\bigl(Y^{(g_{1}g_{2})^{-1}}_{P(z)}\bigr)^{o}(u, z_1)\tilde{\lambda}\bigr\rangle
\end{equation}
and the right-hand side becomes a~single-valued branch of a~multivalued analytic function
of~$\xi$,~$z_{1}$ and $z_{2}$
defined on the region $\bigl|z_{1}\xi^{-1}\bigr|>\bigl|z_{2}\xi^{-1}\bigr|>|z|$. The analytic function
obtained from the right-hand side of~\eqref{comp-main-2.1} can be expanded as an absolutely
convergent multiseries
in powers of $\xi$, $z_{1}$, $z_{2}$ and nonnegative integral powers of $\log \xi$, $\log z_{1}$,
$\log z_{2}$ on the region $0<|\xi|<1$, $|z_{1}|>|z_{2}|>|z|$. In particular,
the iterated series corresponding to this multiseries is also convergent absolutely.
Thus we see that
the coefficients of~\eqref{comp-main-2.2} as a~series in the variable $\xi$ is
an absolutely convergent series in the variable $z_{1}$ and $z_{2}$ on the region
$|z_{1}|>|z_{2}|>|z|$. In particular, for $n\in \C$,
\[
\big\langle \pi_{n}I(w_{1}\otimes w_{2}), \bigl(Y^{(g_{1}g_{2})^{-1}}_{P(z)}\bigr)^{o}(v, z_2)
\bigl(Y^{(g_{1}g_{2})^{-1}}_{P(z)}\bigr)^{o}(u, z_1)\tilde{\lambda}\big\rangle
\]
is absolutely convergent on the region $|z_{1}|>|z_{2}|>|z|$.
On the other hand, $\pi_{n}I(w_{1}\otimes w_{2})$ for~${n\in \C}$, $w_{1}\in W_{1}$ and $w_{2}\in W_{2}$
span $W_{\lambda}'$. In fact, if these elements do not span $W_{\lambda}'$, then there exists a~nonzero element
$\tilde{\lambda}\in W_{\lambda}$ such that \smash{$\bigl\langle I(w_{1}\otimes w_{2}), \tilde{\lambda}\bigr\rangle=0$}
for $w_{1}\in W_{1}$ and $w_{2}\in W_{2}$. But this means that $\tilde{\lambda}(w_{1}\otimes w_{2})=0$
for $w_{1}\in W_{1}$ and $w_{2}\in W_{2}$. Hence $\tilde{\lambda}=0$. Contradiction.
Thus we see that
\begin{equation}\label{comp-main-2.4}
\bigl\langle w', \bigl(Y^{(g_{1}g_{2})^{-1}}_{P(z)}\bigr)^{o}(v, z_2)
\bigl(Y^{(g_{1}g_{2})^{-1}}_{P(z)}\bigr)^{o}(u, z_1)\tilde{\lambda}\bigr\rangle
\end{equation}
is absolutely convergent on the region $|z_{1}|>|z_{2}|>|z|$ for $w'\in W_{\lambda}'$.
Since $W_{\lambda}$ and $W_{\lambda}'$ are lower bounded, the series~\eqref{comp-main-2.4}
must have only finitely many terms in negative real parts of powers of $z_{2}$
and finitely many terms in positive real parts of $z_{1}$. Then the absolute convergence of~\eqref{comp-main-2.4}
on the region $|z_{1}|>|z_{2}|>|z|$ means that it must be absolutely convergent on the region
$|z_{1}|>|z_{2}|>0$. Similarly, using
the left-hand side of~\eqref{comp-main-2}, we see that
\begin{equation}\label{comp-main-2.5}
\big\langle w', \bigl(Y^{(g_{1}g_{2})^{-1}}_{P(z)}\bigr)^{o}
(Y_{V}(u, z_1-z_{2})v, z_2)
\tilde{\lambda}\big\rangle
\end{equation}
 is absolutely convergent on the region $|z_{2}|>|z_{1}-z_{2}|>0$.
Since the sums of~\eqref{comp-main-1.1} and the left-hand side of~\eqref{comp-main-2}
are equal to the single-valued branches on the intersection of the corresponding regions,
we see that the sums of~\eqref{comp-main-2.4} and~\eqref{comp-main-2.5} are also
equal on the intersection of the corresponding regions.
Thus we have proved the associativity for \smash{$\bigl(Y^{(g_{1}g_{2})^{-1}}_{P(z)}\bigr)^{o}$},
which is equivalent to the associativity for \smash{$Y^{(g_{1}g_{2})^{-1}}_{P(z)}$}. This finishes the proof
that $W_\lambda$ equipped with~\smash{$Y^{(g_{1}g_{2})^{-1}}_{P(z)}$}
is a~grading-restricted generalized $(g_{1}g_{2})^{-1}$-twisted
$V$-module.

By Propositions~\ref{int-comp} and~\ref{int-gr-res-cond}, an element
of $W_{1}\hboxtr_{P(z)}W_{2}$ satisfies
the $P(z)$-compatibility condition and the $P(z)$-local-grading-restric\-tion
condition. We now assume that for an element $\lambda$ of $(W_{1}\otimes W_{2})^{*}$
satisfying the $P(z)$-compatibility condition and
the $P(z)$-local-grading-restric\-tion, $W_{\lambda}$ is a~$(g_{1}g_{2})^{-1}$-twisted generalized $V$-submodule of a~$(g_{1}g_{2})^{-1}$-twisted generalized $V$-mod\-ule~$W$ in $\operatorname{COMP}$ whose contragredient
$W'$ is in $\mathcal{C}$.
We want to prove \smash{$\lambda\in W_{1}\hboxtr_{P(z)}W_{2}$}.

Since $W^{*}$ is linearly isomorphic to
$\overline{W'}$, we shall identify
$W^{*}$ with $\overline{W'}$.
We define a~linear map~${I\colon W_{1}\otimes W_{2}
\to W^{*}}$
by
$\langle \mu, I(w_{1}\otimes w_{2})\rangle
=\mu(w_{1}\otimes w_{2})$
for $\mu\in W$, $w_{1}\in W_{1}$ and $w_{2}\in W_{2}$.
We define a~linear map
$
\Y_{I}\colon W_{1}\otimes W_{2}\to W'\{x\}[\log x]$, $
w_{1}\otimes w_{2}\mapsto \Y_{I}(w_{1}, x)w_{2}
$
by{\samepage
\[\Y_{I}(w_{1}, x)w_{2}=x^{L(0)}{\rm e}^{-(\log z)L(0)}
I\bigl(x^{-L(0)}{\rm e}^{(\log z)L(0)}w_{1}\otimes
x^{-L(0)}{\rm e}^{(\log z)L(0)}w_{2}\bigr)\]
for $w_{1}\in W_{1}$ and $w_{2}\in W_{2}$.}

From the definitions of $\Y_{I}$ and $L_{P(z)}'(0)$ on $\operatorname{COMP}\subset \operatorname{COM}$,
it is easy to see that $\Y_{I}$
satisfies the $L(-1)$-derivative property.

For $u_{1}, \dots, u_{k-1}\in V$, $w_{1}\in W_{1}$, $w_{2}\in W_{2}$
and $w\in W$, we have
\begin{gather}
\big\langle w, Y_{P(z)}^{(g_{1}g_{2})^{-1}}(u_{1}, z_{1})
\cdots Y_{P(z)}^{(g_{1}g_{2})^{-1}}(u_{k-1}, z_{k-1})
\Y_{I}(w_{1}, z_{k})w_{2}\big\rangle\nn
\qquad
=\bigl\langle w, Y_{P(z)}^{(g_{1}g_{2})^{-1}}(u_{1}, z_{1})
\cdots Y_{P(z)}^{(g_{1}g_{2})^{-1}}(u_{k-1}, z_{k-1})\nn
\phantom{\qquad=}{}\times
 {\rm e}^{(\log z_{k}-\log z)L(0)}
I\bigl({\rm e}^{-(\log z_{k}-\log z)L(0)}w_{1}\otimes
{\rm e}^{-(\log z_{k}-\log z)L(0)}w_{2}\bigr)\bigr\rangle\nn
\qquad=\bigl\langle {\rm e}^{(\log z_{k}-\log z)L(0)}w, Y_{P(z)}^{(g_{1}g_{2})^{-1}}\bigl(z_{k}^{-L(0)}z^{L(0)}u_{1},
{\rm e}^{-(\log z_{k}-\log z-\log z_{1})}\bigr)
 \nn
\phantom{\qquad=}{}\times
Y_{P(z)}^{(g_{1}g_{2})^{-1}}\bigl(z_{k}^{-L(0)}z^{L(0)}
u_{k-1}, {\rm e}^{-(\log z_{k}-\log z-\log z_{k-1})}\bigr)\nn
\phantom{\qquad=}{}\times
I\bigl({\rm e}^{-(\log z_{k}-\log z)L(0)}w_{1}\otimes
{\rm e}^{-(\log z_{k}-\log z)L(0)}w_{2}\bigr)\bigr\rangle\nn
\qquad=\bigl\langle \bigl(Y_{P(z)}^{(g_{1}g_{2})^{-1}}\bigr)^{o}
\bigl(z_{k}^{-L(0)}z^{L(0)}
u_{k-1}, {\rm e}^{-(\log z_{k}-\log z-\log z_{k-1})}\bigr)\cdots\nn
 \phantom{\qquad=}{}\times \bigl(Y_{P(z)}^{(g_{1}g_{2})^{-1}}\bigr)^{o}
\bigl(z_{k}^{-L(0)}z^{L(0)}u_{1},
{\rm e}^{-(\log z_{k}-\log z-\log z_{1})}\bigr)
{\rm e}^{(\log z_{k}-\log z)L(0)}w, \nn
\phantom{\qquad=\cdots}{}
I\bigl({\rm e}^{-(\log z_{k}-\log z)L(0)}w_{1}\otimes
{\rm e}^{-(\log z_{k}-\log z)L(0)}w_{2}\bigr)\bigr\rangle\nn
\qquad=\bigl(\bigl(Y_{P(z)}^{(g_{1}g_{2})^{-1}}\bigr)^{o}
\bigl(z_{k}^{-L(0)}z^{L(0)}
u_{k-1}, {\rm e}^{-(\log z_{k}-\log z-\log z_{k-1})}\bigr)\cdots\nn
\phantom{\qquad=}{} \times\bigl(Y_{P(z)}^{(g_{1}g_{2})^{-1}}\bigr)^{o}
\bigl(z_{k}^{-L(0)}z^{L(0)}u_{1},
{\rm e}^{-(\log z_{k}-\log z-\log z_{1})}\bigr)
{\rm e}^{(\log z_{k}-\log z)L(0)}w\bigr) \nn
\phantom{\qquad=}{}\times
\bigl({\rm e}^{-(\log z_{k}-\log z)L(0)}w_{1}\otimes
{\rm e}^{-(\log z_{k}-\log z)L(0)}w_{2}\bigr).\label{comp-main-3}
\end{gather}

By~\eqref{Pzaction}, the right-hand side of~\eqref{comp-main-3}
is equal to the series obtained by
expanding the function
\begin{gather*}
f_{k-1}^{e}\bigl(\xi_{1}, \dots, \xi_{k-1};
z_{k}^{-L(0)}z^{L(0)}u_{1}, \dots, z_{k}^{-L(0)}z^{L(0)}
u_{k-1}, \nn
 \qquad  {\rm e}^{-(\log z_{k}-\log z)L(0)}w_{1},
{\rm e}^{-(\log z_{k}-\log z)L(0)}w_{2}; {\rm e}^{(\log z_{k}-\log z)L(0)}w\bigr)
\end{gather*}
on the region $|\xi_{1}|>\cdots>|\xi_{k-1}|>|z|$
as series in powers of $\xi_{1}, \dots, \xi_{k-1}$ and nonnegative
integer powers of $\log \xi_{1}, \dots, \log \xi_{1}$ and then
substituting ${\rm e}^{-n(\log z_{k}-\log z-\log z_{1})}$, $\dots$,
${\rm e}^{-n(\log z_{k}-\log z-\log z_{k-1})}$ ($n\in \C$) for
${\rm e}^{n\log \xi_{1}}, \dots, {\rm e}^{n\log \xi_{k-1}}$ and
$\log z_{k}-\log z-\log z_{1}, \dots, \log z_{k}-\log z-\log z_{k-1}$
for $\log \xi_{1}, \dots,
\log \xi_{k-1}$, respectively.
Then the right-hand side of~\eqref{comp-main-3}
is absolutely convergent on the region
$\big|zz_{1}z_{k}^{-1}\big|>\cdots >\big|zz_{k-1}z_{k}^{-1}\big|>|z|$
or equivalently the region $|z_{1}|>\cdots>|z_{k-1}|>|z_{k}|>0$
and can be analytically extended to a~multivalued analytic function
on $M^{k-1}(0, z)$ with a~preferred branch.
By~\eqref{comp-main-3},
\[\big\langle w, Y_{P(z)}^{(g_{1}g_{2})^{-1}}(u_{1}, z_{1})
\cdots Y_{P(z)}^{(g_{1}g_{2})^{-1}}(u_{k-1}, z_{k-1})
\Y_{I}(w_{1}, z_{k})w_{2}\big\rangle\]
is absolutely convergent on the region
$|z_{1}|>\cdots>|z_{k-1}|>|z_{k}|>0$
and can be analytically extended to a~multivalued analytic function $g(z_{1}, \dots, z_{k})$
on $M^{k}(0, z)$ with a~preferred branch. Moreover, since a~component-isolated singularity of
the multivalued analytic function
\begin{align*}
&\prod_{1\le i\le j\le k-1}(\xi_{i}-\xi_{j})^{M_{ij}}f\bigl(\xi_{1}, \dots, \xi_{k-1};
z_{k}^{-L(0)}z^{L(0)}u_{1}, \dots, z_{k}^{-L(0)}z^{L(0)}
u_{k-1}, \nn
& \qquad {\rm e}^{-(\log z_{k}-\log z)L(0)}w_{1},
{\rm e}^{-(\log z_{k}-\log z)L(0)}w_{2}; {\rm e}^{(\log z_{k}-\log z)L(0)}w\bigr)
\end{align*}
of the form $(\xi_{1}, \dots, \xi_{k-1})-\beta=\delta$ for $\beta\in \C^{k-1}$
and $\delta\in \{0, \infty\}^{k-1}$ is regular, we see that a~component-isolated singularity
of the multivalued analytic function $g(z_{1}, \dots, z_{k})$ of the form
$\bigl(zz_{1}z_{k}^{-1}, \dots, zz_{k-1}z_{k}^{-1}\bigr)-\beta=\delta$ for $\beta\in \C^{k-1}$
and $\delta\in \{0, \infty\}^{k-1}$ is regular. But this is the same as saying that for fixed $z_{k}$, a~component-isolated singularity
of \smash{$\prod_{1\le i\le j\le k-1}(z_{i}-z_{j})^{M_{ij}}g(z_{1}, \dots, z_{k})$} of the form $(z_{1}, \dots, z_{k-1})-\beta'=\delta'$ for $\beta'\in \C^{k-1}$
and $\delta'\in \{0, \infty\}^{k-1}$.

We now prove the duality property for $\Y_{I}$. By~\eqref{comp-main-3},
for $w\in W$, $v\in V$, $w_{1}\in W_{1}$ and $w_{2}\in W_{2}$,
\begin{align}
&\big\langle w, Y_{P(z)}^{(g_{1}g_{2})^{-1}}(v, z_{1})
\Y_{I}(w_{1}, z_{2})w_{2}\big\rangle\nn
&\qquad=\bigl(\bigl(Y_{P(z)}^{(g_{1}g_{2})^{-1}}\bigr)^{o}
\bigl(z_{2}^{-L(0)}z^{L(0)}
v, {\rm e}^{-(\log z_{2}-\log z-\log z_{1})}\bigr)
{\rm e}^{(\log z_{2}-\log z)L(0)}w\bigr) \nn
& \phantom{\qquad=}{}\times
\bigl({\rm e}^{-(\log z_{2}-\log z)L(0)}w_{1}\otimes
{\rm e}^{-(\log z_{2}-\log z)L(0)}w_{2}\bigr),\label{comp-main-3.1}
\end{align}
which is absolutely convergent on the region $\big|z_{2}^{-1}zz_{1}\big|
>|z|$ or equivalently $|z_{1}|>|z_{2}|>0$ and its sum is equal to
\begin{gather}
f_{1}^{e}\bigl({\rm e}^{-(\log z_{2}-\log z-\log z_{1})};
z_{2}^{-L(0)}z^{L(0)}v,
{\rm e}^{-(\log z_{2}-\log z)L(0)}w_{1},
{\rm e}^{-(\log z_{2}-\log z)L(0)}w_{2};\nn
\qquad {\rm e}^{(\log z_{2}-\log z)L(0)}w\bigr)\label{comp-main-4}
\end{gather}
on the region $\big|z_{2}^{-1}zz_{1}\big|
>|z|$, $\big|{\arg \bigl(z_{2}^{-1}zz_{1}-z\bigr)-\arg z_{2}^{-1}zz_{1}}\big|<\frac{\pi}{2}$.
For $z_{1}, z_{2}\in \R_{+}^{2}$ satis\-fying
$z_{1}>z_{2}>0$, we have $\big|z_{2}^{-1}zz_{1}\big|=z_{1}z_{2}^{-1}|z|
>|z|$ and $\big|{\arg \bigl(z_{2}^{-1}zz_{1}-z\bigr)-\arg z_{2}^{-1}zz_{1}}\big|=\big|{\arg \bigl(z_{1}z_{2}^{-1}-1\bigr)z-\arg \bigl(z_{1}z_{2}^{-1}\bigr)z}\bigr|
|{\arg z-\arg z}|=0<\frac{\pi}{2}$. Hence the sum of the left-hand side of~\eqref{comp-main-3.1}
for such $z_{1}$ and $z_{2}$ is equal to~\eqref{comp-main-4}. But for an~absolutely convergent series of the form
of the left-hand side of~\eqref{comp-main-3.1}, its sum must be equal to one single valued branch
on the region~${|z_{1}|>|z_{2}|>0}$, $|{\arg (z_{1}-z_{2})-\arg z_{1}}|<\frac{\pi}{2}$. Thus the sum of
the left-hand side of~\eqref{comp-main-3.1} is equal to~\eqref{comp-main-4}
on the region $|z_{1}|>|z_{2}|>0$, $|{\arg (z_{1}-z_{2})-\arg z_{1}}|<\frac{\pi}{2}$.

By definition,
\begin{align}
&\big\langle w, \Y_{I}(w_{1}, z_{2})Y^{g_{2}}(v, z_{1})w_{2}\big\rangle\nn
&\qquad=\big\langle {\rm e}^{(\log z_{2}-\log z)L(0)}w,
\Y_{I}\bigl({\rm e}^{-(\log z_{2}-\log z)L(0)}w_{1}, z\bigr)\nn
& \phantom{\qquad=}{}\times
Y^{g_{2}}\bigl(z_{2}^{-L(0)}z^{L(0)}v, {\rm e}^{-(\log z_{2}-\log z-\log z_{1})}\bigr)
{\rm e}^{-(\log z_{2}-\log z)L(0)}w_{2}\big\rangle\nn
& \qquad =\bigl({\rm e}^{(\log z_{2}-\log z)L(0)}w\bigr)
\bigl({\rm e}^{-(\log z_{2}-\log z)L(0)}w_{1}\nn
& \phantom{\qquad=}{}\otimes
Y^{g_{2}}(z_{2}^{-L(0)}z^{L(0)}v, {\rm e}^{-(\log z_{2}-\log z-\log z_{1})})
{\rm e}^{-(\log z_{2}-\log z)L(0)}w_{2}\bigr)\label{comp-main-4.1}
\end{align}
is absolutely convergent on the region $|z|>\big|z_{2}^{-1}zz_{1}\big|>0$
or equivalently $|z_{2}|>|z_{1}|>0$ and its sum is equal to~\eqref{comp-main-4} on the region
$|z_{2}|>|z_{1}|>0$, $-\frac{3\pi}{2}<\arg \bigl(z_{2}^{-1}zz_{1}-z\bigr)
-\arg z<-\frac{\pi}{2}$. Using the same argument as above, we can show that
the sum of the left-hand side of~\eqref{comp-main-4.1} must be equal to
\eqref{comp-main-4} on the region
$|z_{2}|>|z_{1}|>0$, $-\frac{3\pi}{2}<\arg (z_{1}-z_{2})
-\arg z_{2}<-\frac{\pi}{2}$.

By 2\,(a) in the compatibility condition, we see that
\begin{gather}
\big\langle w, \Y_{I}\bigl(Y^{g_{1}}(v, z_{1}-z_{2})w_{1}, z_{2}\bigr)w_{2}\big\rangle\nn
\qquad=\big\langle {\rm e}^{(\log z_{2}-\log z)L(0)}w, \nn
\phantom{\qquad=}{}\times
 \Y_{I}\bigl(Y^{g_{1}}\bigl(z_{2}^{-L(0)}z^{L(0)}v,
z_{2}^{-1}zz_{1}-z\bigr){\rm e}^{-(\log z_{2}-\log z)L(0)}w_{1},
z\bigr){\rm e}^{-(\log z_{2}-\log z)L(0)}w_{2}\big\rangle\nn
\qquad=\bigl({\rm e}^{(\log z_{2}-\log z)L(0)}w\bigr)
\bigl(Y^{g_{1}}\bigl(z_{2}^{-L(0)}z^{L(0)}v,
z_{2}^{-1}zz_{1}-z\bigr){\rm e}^{-(\log z_{2}-\log z)L(0)}w_{1}\nn
\phantom{\qquad=}{}
\otimes {\rm e}^{-(\log z_{2}-\log z)L(0)}w_{2}\bigr),\label{comp-main-4.2}
\end{gather}
which is absolutely convergent on the region $|z|>{\big|z_{2}^{-1}zz_{1}-z\big|}>0$
or equivalently $|z_{2}|>{|z_{1}-z_{2}|}>0$ and its sum is equal
to~\eqref{comp-main-4} on the region
$|z|>\big|z_{2}^{-1}zz_{1}-z\big|>0$, $\big|{\arg z_{2}^{-1}zz_{1}}
-\arg z\big|<\frac{\pi}{2}$. Using the same argument as above, we can show that
the sum of the left-hand side of~\eqref{comp-main-4.2} must be equal to
\eqref{comp-main-4} on the region
$|z_{2}|>|z_{1}-z_{2}|>0$, $|{\arg z_{1}
-\arg z_{2}}|<\frac{\pi}{2}$. Thus the duality property for $\Y_{I}$
is proved and $\Y_{I}$ is a~twisted intertwining operator of type \smash{$\binom{W'}{W_{1}W_{2}}$}.
Then $I$ is a~twisted $P(z)$-intertwining map of the same type.

Now we have
$
\lambda(w_{1}\otimes w_{2})
=\langle \lambda, I(w_{1}\otimes w_{2})\rangle
=\lambda_{I, \lambda}(w_{1}\otimes w_{2})
$
for $w_{1}\in W_{1}$ and $w_{2}\in W_{2}$.
In particular, $\lambda=\lambda_{I, \lambda}$.
Since $W'$ is in $\mathcal{C}$, by definition, $\lambda_{I, \lambda}\in W_{1}\hboxtr_{P(z)}W_{2}$.
\end{proof}

\appendix

\section{A convergence lemma}

Let $A$ be a~finite subset of $\C/\Z$, $R_{\mu}\in \mu$ for
$\mu\in A$, $D$ a~subset of
$\cup_{\mu\in A}(R_{\mu}-\N)$, $\Delta\in \C$ and~$a_{n,j,i}\in \C$
for $n\in D$, $j=0, \dots, M$ and $i=1, \dots, N$. Consider
the triple series
\begin{align}
&\sum_{n\in D}\sum_{j=0}^{M}\sum_{i=0}^{N}a_{n,j,i}
{\rm e}^{(-\Delta+n+1)\log z_1}(\log z_1)^j{\rm e}^{(-n-1)
\log z_2}(\log z_2)^i\label{kkl}
\end{align}
for $z_{1}, z_{2}\in \C^{\times}$.

For any $z_1$, $z_2$ satisfying $|z_1|>|z_2|$,
 $|{\arg(z_1-z_2)-\arg z_1}|<\frac{\pi}{2}$,
we have
\begin{gather}
{\rm e}^{\alpha \log (z_1-z_2)}=\sum_{k\in \N}\binom{\alpha }{k}(-1)^k
{\rm e}^{(\alpha-k)\log z_1}z_2^k,\label{99}\\
\log (z_1-z_2)=\log z_1+\sum_{k\in \Z_{+}}\frac{(-1)^{k}}{k}
z_1^{-k}z_2^k,\label{99-1}
\end{gather}
and for $z_{2}\ne 0$, we have
\smash{$
{\rm e}^{\alpha l_{q_2}(-z_2)}={\rm e}^{\alpha \pi\i}
{\rm e}^{\alpha \log z_2}$}, \smash{$
l_{q_2}(-z_2)=\log z_2+\pi\i$},
where $\alpha\in\mathbb{C}$ and
 $q_2=0, 1$ if $\arg z_2<\pi$, $\arg z_2\geq\pi$, respectively.
 Note that in our notations, $\log z=l_{0}(z)$ for $z\in \C^{\times}$.

For $n\in D$, $j=0,\dots,M$,
$k\in\Z_{\geq0}$, $s=0,\dots,j$,
define $b_{n,j,k,s}\in\C$ as the coefficients of the following
formal power series expansion
\begin{align}
(x+y)^{-\Delta+n+1}\log(x+y)^j=\sum_{k\in\N}\sum_{s=0}^jb_{n,j,k,s}x^{-\Delta+n+1-k}y^k(\log x)^s,\label{b-defi}
\end{align}
where $x$ and $y$ are formal variables.
From~\eqref{99}, \eqref{99-1}
and~\eqref{b-defi}, when $|z_1|>|z_2|$ and
$|{\arg(z_1-z_2)-\arg z_1}|<\frac{\pi}{2}$, we have the following expansion:
\begin{gather}
{\rm e}^{(-\Delta+n+1)\log (z_{1}-z_{2})}(\log (z_1-z_2))^j\nonumber\\
\qquad
=\sum_{k\in\N}\sum_{s=0}^j(-1)^kb_{n,j,k,s}{\rm e}^{(-\Delta+n+1-k)
\log z_1}z_2^k(\log z_1)^s.\label{b-anly}
\end{gather}

Now consider
\begin{align}
\sum_{n\in D}\sum_{j=0}^{M}\sum_{i=0}^{N}
a_{n,j,i}{\rm e}^{(-\Delta+n+1)\log (z_1-z_{2})}
(\log (z_1-z_{2}))^j
{\rm e}^{(-n-1)\log (-z_2)}(\log (-z_2))^i.\label{88}
\end{align}
Using~\eqref{b-anly}, we can further expand each term in
the right-hand side of~\eqref{88} so that the right-hand side of~\eqref{88}
becomes the iterated sum
\begin{gather}
\sum_{n\in D}
\sum_{j=0}^{M}\sum_{i=0}^{N}a_{n,j,i}
\left(\sum_{k\in\N}\sum_{s=0}^j(-1)^kb_{n,j,k,s}
{\rm e}^{(-\Delta+n+1-k)\log z_1}z_2^k(\log z_1)^s\right)\nn
\qquad \times
{\rm e}^{(-n-1)\log (-z_2)}(\log (-z_2))^i\nn
 \phantom{\qquad \times}{} =\sum_{n\in D}
\sum_{j=0}^{M}\sum_{i=0}^{N}
\left(\sum_{k\in\N}\sum_{s=0}^j a_{n,j,i}(-1)^kb_{n,j,k,s}
{\rm e}^{(-\Delta+n+1-k)\log z_1}z_2^k(\log z_1)^s\right.\nn
 \phantom{\qquad \times =}{} \times\left.
{\rm e}^{(-n-1)\log (-z_2)}(\log (-z_2))^i\right),\label{iterate-sum-1}
\end{gather}
where the inner sum is absolutely convergent on the region
$|z_{1}|>|z_{2}|>0$.

We are interested in the convergence of the multisum
\begin{gather}
 \sum_{n\in D}\sum_{j=0}^{M}\sum_{i=0}^{N}\sum_{k\in\N}
\sum_{s=0}^j
a_{n,j,i}
(-1)^k b_{n,j,k,s}{\rm e}^{(-\Delta+n+1-k)
\log z_1}z_2^k\nonumber\\
\qquad\times
(\log z_1)^s
{\rm e}^{(-n-1)\log (-z_2)}(\log (-z_2))^i\label{iterate-sum-5.9}
\end{gather}
and the corresponding series
\begin{gather}
 \sum_{m\in D-\N}\sum_{s=0}^M
\sum_{i=0}^{N}
 \Bigg(\sum_{j=s}^{M}\sum_{\substack{n-k=m\\n\in D, k\in\N}}
a_{n,j,i}b_{n,j,k,s}\Bigg){\rm e}^{(-\Delta+m+1)
\log z_1}\nn
\qquad\times
(\log z_1)^s
{\rm e}^{(-m-1)\log (-z_2)}(\log (-z_2))^i\label{iterate-sum-6}
\end{gather}
in powers of $z_{1}$ and $z_{2}$ and nonnegative integral powers of
$\log z_{1} $ and $\log z_{2}$.

\begin{lem}\label{convergence}
Assume that
the triple series~\eqref{kkl} and the series obtained from
\eqref{kkl} by taking derivatives of each term in~\eqref{kkl}
with respect to $z_{1}$ and $z_{2}$
are absolutely convergent on the region given by
$|z_{1}|>|z_{2}|>0$.
Then the multisum~\eqref{iterate-sum-5.9}
is absolutely convergent on the region~${|z_1|>2|z_2|>0}$.
Assume in addition that
\eqref{kkl} is convergent on the region
$|z_{1}|>|z_{2}|>0$, $|{\arg (z_{1}-z_{2})-\arg z_{1}}|<\frac{\pi}{2}$ to a~single-valued analytic branch $f^{e}(z_{1}, z_{2})$ on $M^{2}_{0}$
of a~maximally extended multivalued analytic function on $M^{2}$ such that
$f^{e}(z_{1}-z_{2}, -z_{2})$ has no singular point in the
region $|z_{1}|>|z_{2}|>0$. Then the iterated sum
\eqref{iterate-sum-6} is absolutely
convergent on the region~${|z_1|>|z_2|>0}$,
$|{\arg z_1-\arg (z_1-z_{2})}|<\frac{\pi}{2}$ to
$f^{e}(z_{1}-z_{2}, -z_{2})$.
\end{lem}
\begin{proof}
Let $\tilde{n}=R_\mu-n$. Then the sum \smash{$\sum_{n\in D}$} in
\eqref{kkl} can be written as the same as
\smash{$\sum_{\mu\in A}\sum_{\tilde{n}\in\N}$}. So
the series~\eqref{kkl} can be written as
\begin{gather*}
\sum_{\mu\in A}\sum_{\tilde{n}\in\N}\sum_{j=0}^{M}\sum_{i=0}^{N}a_{R_\mu-\tilde{n},j,i}{\rm e}^{(-\Delta+R_\mu-\tilde{n}+1)\log z_1}(\log z_1)^j
{\rm e}^{(-R_\mu+\tilde{n}-1)\log z_2}(\log z_2)^i.
\end{gather*}
By assumption, the series
\begin{align}
&\sum_{\mu\in A}\sum_{\tilde{n}\in\N}\sum_{j=0}^{M}\sum_{i=0}^{N}
a_{R_\mu-\tilde{n},j,i}\frac{\partial^{p+q}}{\partial z_{1}^{p}\partial z_{2}^{q}}
{\rm e}^{(-\Delta+R_\mu-\tilde{n}+1)\log z_1}(\log z_1)^j
{\rm e}^{(-R_\mu+\tilde{n}-1)\log z_2}(\log z_2)^i.\!\!\!\label{1}
\end{align}
for $p, q\in \N$ is absolutely convergent on the region $|z_{1}|>|z_{2}|>0$.
For any $r>1$, consider
\begin{gather}
\sum_{\mu\in A}\sum_{\tilde{n}\in\N}\sum_{j=0}^{M}\sum_{i=0}^{N}a_{R_\mu-\tilde{n},j,i}
\frac{\partial^{p+q}}{\partial z_{1}^{p}\partial z_{2}^{q}} {\rm e}^{(-\Delta+R_\mu-\tilde{n}+1)\log z_1}(\log z_1)^j\nn
\qquad\times
{\rm e}^{(-R_\mu+\tilde{n}-1)\log z_2}(\log z_2)^i r^{\tilde{n}}.\label{1-1}
\end{gather}
Then the absolute convergence of~\eqref{1} on the
region $|z_{1}|>|z_{2}|>0$ is equivalent to the absolute convergence of
 \eqref{1-1} on the region $|z_{1}|>r|z_{2}|>0$.
 But the absolute convergence of~\eqref{1-1} on the region $|z_{1}|>r|z_{2}|>0$ is in turn equivalent to
 the absolute convergence of
\begin{equation}\label{prod-power-series}
\sum_{\tilde{n}\in\N}a_{R_\mu-\tilde{n},j,i}\frac{\partial^{p+q}}{\partial z_{1}^{p}\partial z_{2}^{q}}
\left(\frac{z_{2}}{z_{1}}\right)^{\tilde{n}}r^{\tilde{n}}
\end{equation}
on the same region $|z_{1}|>r|z_{2}|>0$.
But as a~power series in $\frac{z_{2}}{z_{1}}$, \eqref{prod-power-series} with $p, q=0$
has a~removable singularity at $\frac{z_{2}}{z_{1}}=0$.
Then we see that~\eqref{prod-power-series} with $p, q=0$ is absolutely convergent
on the region \smash{$|\frac{z_{2}}{z_{1}}|<\frac{1}{r}$} and is hence also uniformly
convergent on the closed region \smash{$|\frac{z_{2}}{z_{1}}|\le r_{2}$}
for any positive \smash{$r_{2}<\frac{1}{r}$}. In particular, by the property of
power series, \eqref{prod-power-series} for $p, q\in \N$ is also
absolutely convergent
on the region \smash{$|\frac{z_{2}}{z_{1}}|<\frac{1}{r}$} and is uniformly
convergent on the closed region~\smash{$|\frac{z_{2}}{z_{1}}|\le r_{2}$}.
Thus for such $r_{2}$,
\eqref{1-1} is a~(finite) linear combination of the uniformly convergent
series~\eqref{prod-power-series} on the region
\smash{$|\frac{z_{2}}{z_{1}}|\le r_{2}$}
with analytic functions of $z_{1}$ and $z_{2}$ on $M_{0}^{2}$
as coefficients.

Substituting $z_{1}-z_{2}$ and $-z_{2}$ for $z_{1}$ and $z_{2}$, respectively,
in~\eqref{1-1}, we see that
\begin{gather}
\sum_{\mu\in A}\sum_{\tilde{n}\in\N}\sum_{j=0}^{M}\sum_{i=0}^{N}
a_{n,j,i}(-1)^{q}\frac{\partial^{p+q}}{\partial z_{1}^{p}\partial z_{2}^{q}} {\rm e}^{(-\Delta+n+1)\log (z_1-z_{2})}
(\log (z_1-z_{2}))^j\nn
\qquad\times
{\rm e}^{(-n-1)\log (-z_2)}(\log (-z_2))^ir^{\tilde{n}}\label{88-1}
\end{gather}
is absolutely convergent on the region
$|z_1-z_{2}|>r|z_2|>0$ and is a~(finite) linear combination of
uniformly convergent series on the region
\smash{$|\frac{z_{2}}{z_{1}-z_{2}}|\le r_{2}$} with analytic coefficients on $M_{0}^{2}$
for any positive \smash{$r_{2}<\frac{1}{r}$}.
Now we expand each term in~\eqref{88-1}
using~\eqref{b-anly} and then passing the derivative~\smash{$\frac{\partial^{p+q}}{\partial z_{1}^{p}\partial z_{2}^{q}} $}
to each term since the right-hand side of~\eqref{b-anly} is essentially a~power series in~$\frac{z_{2}}{z_{1}}$
 to obtain the iterated sum
\begin{gather}
\sum_{\mu\in A}\sum_{\tilde{n}\in\N}
\sum_{j=0}^{M}\sum_{i=0}^{N}a_{R_\mu-\tilde{n},j,i}\frac{\partial^{p+q}}{\partial z_{1}^{p}\partial z_{2}^{q}}
\left(\left(\sum_{k\in\N}\sum_{s=0}^j(-1)^kb_{R_\mu-\tilde{n},j,k,s}\right.\right.
\nn
\phantom{=}{} \qquad\left. \left. \times{\rm e}^{(-\Delta+R_\mu-\tilde{n}+1-k)\log z_1}z_2^k(\log z_1)^s\right)
{\rm e}^{(-R_\mu+\tilde{n}-1)\log (-z_2)}(\log (-z_2))^ir^{\tilde{n}}\right)
\nn
\qquad=\sum_{\mu\in A}\sum_{\tilde{n}\in\N}
\sum_{j=0}^{M}\sum_{i=0}^{N}
\left(\sum_{k\in\N}\sum_{s=0}^j
a_{R_\mu-\tilde{n},j,i}
(-1)^k b_{R_\mu-\tilde{n},j,k,s}\frac{\partial^{p+q}}{\partial z_{1}^{p}\partial z_{2}^{q}}{\rm e}^{(-\Delta+R_\mu-\tilde{n}+1-k)
\log z_1}\right.
\nn
\phantom{=}{} \qquad\times
\left.
z_2^k(\log z_1)^s
{\rm e}^{(-R_\mu+\tilde{n}-1)\log (-z_2)}(\log (-z_2))^ir^{\tilde{n}}\right),\label{iterate-sum-1-1}
\end{gather}
where the inner sum is absolutely convergent on the region
$|z_{1}|>|z_{2}|>0$.
Then as a~subseries
of~\eqref{iterate-sum-1-1} divided by a~suitable function of the form
${\rm e}^{\lambda\log z_1}{\rm e}^{\mu\log (-z_2)}(\log (-z_2))^\nu$ for some
${\lambda, \mu \in \C}$, $\nu\in \N$, the series
\begin{align}\label{iterate-sum-2}
&\sum_{\tilde{n}\in\N}
\left(\sum_{k\in\N}\sum_{s=0}^j (-1)^{\tilde{n}+k}
a_{R_\mu-\tilde{n},j,i}
 b_{R_\mu-\tilde{n},j,k,s}\frac{\partial^{p+q}}{\partial z_{1}^{p}\partial z_{2}^{q}}\left(\frac{z_{2}}{z_{1}}\right)^{\tilde{n}+k}
(\log z_1)^s\right)r^{\tilde{n}}\nn
&\qquad=\sum_{\tilde{n}\in\N}
\left(\sum_{\tilde{k}\in\N}\sum_{s=0}^j (-1)^{\tilde{k}}
a_{R_\mu-\tilde{n},j,i}
 b_{R_\mu-\tilde{n},j,\tilde{k}-\tilde{n},s}\frac{\partial^{p+q}}{\partial z_{1}^{p}\partial z_{2}^{q}}
 \left(\frac{z_{2}}{z_{1}}\right)^{\tilde{k}}
(\log z_1)^s\right)r^{\tilde{n}},
\end{align}
where \smash{$b_{R_\mu-\tilde{n},j,\tilde{k}-\tilde{n},s}$} is defined to be $0$
when $\tilde{k}<\tilde{n}$,
is also absolutely convergent on the region~${|z_{1}-z_{2}|>r|z_{2}}|$ and is uniformly convergent on the closed region
\smash{$\bigl|\frac{z_{2}}{z_{1}-z_{2}}\bigr|\le r_{2}$} for any positive $r_{2}<\frac{1}{r}$,
where the inner sum is absolutely convergent in the region
$|z_{1}|>|z_{2}|$. In the region~${|z_{1}|>(1+r)|z_{2}|>0}$,
we have $|z_{1}-z_{2}|\ge |z_{1}|-|z_{2}|>r|z_{2}|>0$.
Then~\eqref{iterate-sum-2} is absolutely convergent on the region
$|z_{1}|>(1+r)|z_{2}|$ or $|\frac{z_{2}}{z_{1}}|<\frac{1}{1+r}$
with the inner sum absolutely convergent in the region
$|z_{1}|>|z_{2}|$ or $|\frac{z_{2}}{z_{1}}|<1$.
Note that for fixed $\frac{z_{2}}{z_{1}}$, $z_{1}$ can be any
complex number and thus $\log z_{1}$ can also be any nonzero complex number whose
imaginary part is in~$[0, 2\pi)$. Moreover,
\eqref{iterate-sum-2} is convergent absolutely on the same region
when $\log z_{1}$ is replaced by any choice of branches of the logarithm of $z_1$.
This means that $\zeta_{1}=\frac{z_{2}}{z_{1}}$ and
$\zeta_{2}=\log z_{1}$ can be viewed as
independent variables and, after the change of variables from $z_{1}$, $z_{2}$ to $\zeta_{1}$ and $\zeta_{2}$,
\begin{gather}
\sum_{\tilde{n}\in\N}
\Biggl(\sum_{\tilde{k}\in \N}\sum_{s=0}^j (-1)^{\tilde{k}}
a_{R_\mu-\tilde{n},j,i}
 b_{R_\mu-\tilde{n},j,\tilde{k}-\tilde{n},s}\left(-\zeta_{1}{\rm e}^{\zeta_{2}}
\frac{\partial}{\partial \zeta_{1}}+{\rm e}^{-\zeta_{2}}\frac{\partial}{\partial \zeta_{2}}\right)^{p} \nn
\qquad\times
\left({\rm e}^{-\zeta_{2}}\frac{\partial}{\partial \zeta_{1}}\right)^{q}\zeta_{1}^{\tilde{k}}
\zeta_{2}^s\Biggr)\zeta_{3}^{\tilde{n}}\label{iterate-sum-3.1}
\end{gather}
is absolutely convergent on the region given by
$|\zeta_{1}|<\frac{1}{1+r}$, $\zeta_{2}\in \C$ and $|\zeta_{3}|< r$. From
\eqref{iterate-sum-3.1}, we see that
\begin{align}\label{iterate-sum-3}
\sum_{\tilde{n}\in\N}
\left(\sum_{\tilde{k}\in \N}\sum_{s=0}^j (-1)^{\tilde{k}}
a_{R_\mu-\tilde{n},j,i}
 b_{R_\mu-\tilde{n},j,\tilde{k}-\tilde{n},s}
\frac{\partial^{\tilde{p}+\tilde{q}}}{\partial \zeta_{1}^{\tilde{p}}\partial \zeta_{2}^{\tilde{q}}}
\zeta_{1}^{\tilde{k}}\zeta_{2}^s\right)\zeta_{3}^{\tilde{n}}
\end{align}
for $\tilde{p}, \tilde{q}\in \N$ is absolutely convergent on the region given by
$|\zeta_{1}|<\frac{1}{1+r}$, $\zeta_{2}\in \C$ and $|\zeta_{3}|< r$.

On the other hand, on the region \smash{$|\frac{z_{2}}{z_{1}}|\le \frac{r_{2}}{1
+r_{2}}$} for positive $r_{2}<\frac{1}{r}$,
\[
\left|\frac{z_{2}}{z_{1}-z_{2}}\right|\le \frac{|z_{2}|}{|z_{1}|-|z_{2}|}
\le \frac{|z_{2}|}{\frac{1+r_{2}}{r_{2}}|z_{2}|-|z_{2}|}=r_{2}.
\]
Then for $\frac{z_{2}}{z_{1}}\in \C$ satisfying
$|\frac{z_{2}}{z_{1}}|\le \frac{r_{2}}{1
+r_{2}}$ for positive $r_{2}<\frac{1}{r}$ and $\log z_{1}\in \C$,
\eqref{iterate-sum-2} is uniformly convergent with the inner sum
also uniformly convergent on the same closed region.
Thus~\eqref{iterate-sum-3} is uniformly convergent on the region
$|\zeta_{1}|\le \frac{r_{2}}{1
+r_{2}}$, $\zeta_{2}\in \C$ and $|\zeta_{3}|\le r$
with the inner sum also uniformly convergent on the same closed region.
We have shown
that the series
\eqref{iterate-sum-3} for~${\tilde{p}, \tilde{q}\in \N}$ obtained by taking derivatives of each term in
the series~\eqref{iterate-sum-3} with $\tilde{p}, \tilde{q}=0$
with respect to $\zeta_{1}$, $\zeta_{2}$ and $\zeta_{3}$
is uniformly convergent
on the region
$|\zeta_{1}|\le \frac{r_{2}}{1
+r_{2}}$, $\zeta_{2}\in \C$ and $|\zeta_{3}|\le r$.
In particular, the derivatives of~\eqref{iterate-sum-3} with $\tilde{p}, \tilde{q}=0$ with respect
to $\zeta_{1}$, $\zeta_{2}$ and $\zeta_{3}$ exist and
is equal to the sum of the series obtained by taking the
corresponding derivatives of each term in~\eqref{iterate-sum-3} with $\tilde{p}, \tilde{q}=0$
on the region $|\zeta_{1}|< \frac{r_{2}}{1
+r_{2}}$, $\zeta_{2}\in \C$ and $|\zeta_{3}|<r$.
Then the sum of~\eqref{iterate-sum-3} with $\tilde{p}, \tilde{q}=0$
is an analytic function of $\zeta_{1}$, $\zeta_{2}$ and $\zeta_{3}$
on the same open region.

For $r>1$, let $\zeta_{1}$, $\zeta_{2}$, $\zeta_{3}$ be
complex numbers
satisfying \smash{$|\zeta_{1}|< \frac{1}{(1+r)r}$},
$\zeta_{2}\in \C$ and ${|\zeta_{3}|<r}$. Then we have
\smash{$\frac{|\zeta_{1}|}{1-|\zeta_{1}|}<(1+r)|\zeta_{1}|<\frac{1}{r}$}.
We choose $r_{2}$ such that \smash{$\frac{|\zeta_{1}|}{1-|\zeta_{1}|}<r_{2}<
(1+r)|\zeta_{1}|$}. From \smash{$\frac{|\zeta_{1}|}{1-|\zeta_{1}|}<r_{2}$},
we obtain \smash{$|\zeta_{1}|< \frac{r_{2}}{1+r_{2}}$}.
We also have \smash{$0<r_{2}< (1+r)|\zeta_{1}|<\frac{1}{r}$}.
Now $\zeta_{1}$, $\zeta_{3}$ satisfy
\smash{$|\zeta_{1}|\le \frac{r_{2}}{1+r_{2}}$} and $|\zeta_{3}|<r$.
This means that the sum of~\eqref{iterate-sum-3} with $\tilde{p}, \tilde{q}=0$ is analytic
and the derivatives can be
calculated term by term at $\zeta_{1}$, $\zeta_{2}$, $\zeta_{3}$.
So the sum of~\eqref{iterate-sum-3} with $\tilde{p}, \tilde{q}=0$ is analytic
and the derivatives can be
calculated term by term
on the polydisc \smash{$|\zeta_{1}|< \frac{1}{(1+r)r}$},
$\zeta_{2}\in \C$ and $|\zeta_{3}|<r$ for any $r>1$.
Since analytic functions on polydiscs
can be expanded as power series, the sum of~\eqref{iterate-sum-3} with $\tilde{p}, \tilde{q}=0$
can be expanded as a~power series
in $\zeta_{1}$, $\zeta_{2}$ and~$\zeta_{3}$ and the
coefficients of the power series expansion can be obtained using its
derivatives with respect to $\zeta_{1}$, $\zeta_{2}$ and~$\zeta_{3}$
evaluated at $\zeta_{1}=\zeta_{2}=\zeta_{3}=0$.
By taking the derivatives term by term,
we see that the coefficients of the power series expansion of this analytic function
are equal to the coefficients~$(-1)^{\tilde{k}}
a_{R_\mu-\tilde{n},j,i}
 b_{R_\mu-\tilde{n},j,\tilde{k}-\tilde{n},s}$ of
 the iterated series~\eqref{iterate-sum-3} with $\tilde{p}, \tilde{q}=0$.
 Since the power series expansion of an analytic function
 is an absolutely convergent multisum, we see that
 the triple series
\begin{equation}\label{iterate-sum-3.5}
\sum_{\tilde{n}\in\N}
\sum_{\tilde{k}\in \N}\sum_{s=0}^j (-1)^{\tilde{k}}
a_{R_\mu-\tilde{n},j,i}
 b_{R_\mu-\tilde{n},j,\tilde{k}-\tilde{n},s}\zeta_{1}^{\tilde{k}}
\zeta_{2}^s\zeta_{3}^{\tilde{n}}
\end{equation}
is absolutely convergent on the region \smash{$|\zeta_{1}|< \frac{1}{(1+r)r}$},
$\zeta_{2}\in \C$ and $|\zeta_{3}|<r$.

In particular, on the region $|z_{1}|>(1+r)r|z_{2}|$, taking
$\zeta_{1}=\frac{z_{2}}{z_{1}}$,
$\zeta_{2}=\log z_{1}$ and $\zeta_{3}=1$ in~\eqref{iterate-sum-3.5},
we see that the triple series
\begin{align}\label{iterate-sum-3.7}
\sum_{\tilde{n}\in\N}\sum_{\tilde{k}\in\N}\sum_{s=0}^j
(-1)^{\tilde{k}}
a_{R_\mu-\tilde{n},j,i}
 b_{R_\mu-\tilde{n},j,\tilde{k}-\tilde{n},s}
 \left(\frac{z_{2}}{z_{1}}\right)^{\tilde{k}}
(\log z_{1})^s
\end{align}
is absolutely convergent. Since $r$ is an arbitrary real number
satisfying $r>1$, \eqref{iterate-sum-3.7} is in fact absolutely convergent
on the region $|z_{1}|>2|z_{2}|$. Multiplying
${\rm e}^{(-\Delta+R_\mu+1)
\log z_1}{\rm e}^{(-R_\mu-1)\log (-z_2)}\allowbreak\times(\log (-z_2))^i$
to~\eqref{iterate-sum-3.7} and summing over
$A$, $j=0, \dots, M$ and $i=0, \dots, N$,
we see that the multiseries
{\samepage\begin{gather}
\sum_{\mu\in A}\sum_{\tilde{n}\in\N}
\sum_{j=0}^{M}\sum_{i=0}^{N}\sum_{\tilde{k}\in\N}
\sum_{s=0}^j(-1)^{\tilde{k}}
a_{R_\mu-\tilde{n},j,i}
 b_{R_\mu-\tilde{n},j,\tilde{k}-\tilde{n},s}{\rm e}^{(-\Delta+R_\mu+1)
\log z_1}\nn
 \qquad \times {\rm e}^{(-R_\mu-1)\log (-z_2)}(\log (-z_2))^i
 \left(\frac{z_{2}}{z_{1}}\right)^{\tilde{k}}
(\log z_{1})^s\nn
\phantom{ \qquad \times}{}=\sum_{n\in D}\sum_{j=0}^{M}\sum_{i=0}^{N}\sum_{k\in\N}
\sum_{s=0}^j
a_{n,j,i}
(-1)^k b_{n,j,k,s}{\rm e}^{(-\Delta+n+1-k)
\log z_1}z_2^k
(\log z_1)^s\nn
\phantom{ \qquad \times=}{}\times
{\rm e}^{(-n-1)\log (-z_2)}(\log (-z_2))^i\label{iterate-sum-3.8}
\end{gather}
is absolutely convergent on the region $|z_{1}|>2|z_{2}|>0$,
proving the first part of the lemma.}

In the case that the additional assumption holds,
from the proof above and the additional assumption,
the multisum~\eqref{iterate-sum-3.8}, which is equal to
the iterated series in the right-hand side of~\eqref{iterate-sum-1}, is
absolutely convergent on the region $|z_1|>2|z_2|>0$,
${|{\arg z_1-\arg (z_1-z_{2})}|<\frac{\pi}{2}}$ to~${f^{e}(z_{1}-z_{2}, -z_{2})}$. In particular,
\eqref{iterate-sum-6} as an iterated sum of~\eqref{iterate-sum-3.8}
is also absolutely convergent on the region $|z_1|>2|z_2|>0$,
$|{\arg z_1-\arg (z_1-z_{2})}|<\frac{\pi}{2}$ to
$f^{e}(z_{1}-z_{2}, -z_{2})$.
Since there is no singular point
of $f^{e}(z_{1}-z_{2}, -z_{2})$ in the region
$|z_{1}|>|z_{2}|>0$,
\eqref{iterate-sum-6} must also be absolutely convergent when
$|z_{1}|>|z_{2}|>0$ and is thus absolutely
convergent on the region $|z_{1}|>|z_{2}|>0$,
$|{\arg z_1-\arg (z_1-z_{2})}|<\frac{\pi}{2}$ to
$f^{e}(z_{1}-z_{2}, -z_{2})$.
\end{proof}

\begin{rema}
There is a~subtlety about the convergence regions for
\eqref{iterate-sum-3.8} and~\eqref{iterate-sum-6}.
Note that in general the multisum~\eqref{iterate-sum-3.8}
might not be absolutely convergent on the larger region~${|z_{1}|>|z_{2}|>0}$. This is because~\eqref{iterate-sum-3.8}
is not a~series in powers of $z_{1}$ and $z_{2}$ and
nonnegative integral powers of
$\log z_{1} $ and $\log z_{2}$. Even if we restore the variable $\zeta_{3}$
in the proof of the lemma above to obtain a~series
in powers of $z_{1}$, $z_{2}$, $\zeta_{3}$ and
nonnegative integral powers of
$\log z_{1} $ and $\log z_{2}$, since we do not have the assumption
that the sum of this series can be analytically extended to a~region containing $|z_{1}|>|z_{2}|>0$, this series
might not be absolutely convergent on any region
containing $|z_{1}|>|z_{2}|>0$.
\end{rema}

\begin{rema}\label{convergence-multi}
 Lemma~\ref{convergence} can be generalized to the case of more than two variables.
Consider the series
\begin{align}
&\sum_{\substack{n_{1}, \dots, n_{k}\in D\\n_{1}+\cdots +n_{k}=\Delta-k}}
\sum_{j_{1}, \dots, j_{k-1}=0}^{M}\sum_{i=0}^{N}a_{n_{1}, \dots, n_{k},j_{1},\dots, j_{k-1},i}
{\rm e}^{(-n_{1}-1)\log z_1}(\log z_1)^{j_{1}}\cdots \nn
& \qquad
\times {\rm e}^{(-n_{k-1}-1)\log z_l}(\log z_l)^{j_{k-1}} {\rm e}^{(-n_{k}-1)
\log z_{k}}(\log z_{k})^i\label{gen-kkl}
\end{align}
for $z_{1}, \dots, z_{k}\in \C^{\times}$. Replacing $z_{i}$ for $i=1, \dots, k-1$ and $z_{k}$
by $z_{i}-z_{k}$ and $-z_{k}$, respectively, we obtain the following series:
\begin{align}
&\sum_{\substack{n_{1}, \dots, n_{k}\in D\\n_{1}+\cdots +n_{k}=\Delta-k}}
\sum_{j_{1}, \dots, j_{k-1}=0}^{M}\sum_{i=0}^{N}a_{n_{1}, \dots, n_{k},j_{1},\dots, j_{k-1},i}
{\rm e}^{(-n_{1}-1)\log (z_1-z_{k})}\log (z_1-z_{k})^{j_{1}}\cdots\nn
& \qquad
 \times {\rm e}^{(-n_{k-1}-1)\log (z_{k-1}-z_{k})}(\log (z_{k-1}-z_{k}))^{j_{k-1}}{\rm e}^{(-n_{k}-1)
\log (-z_{k})}(\log (-z_{k}))^i.\label{gen-88}
\end{align}
Expanding \smash{${\rm e}^{(-n_{1}-1)\log (z_1-z_{k})}\log (z_1-z_{k})^{j_{1}}, \dots,
 {\rm e}^{(-n_{k-1}-1)\log (z_{k-1}-z_{k})}(\log (z_{k-1}-z_{k}))^{j_{l}}$} in~\eqref{gen-88} using~\eqref{b-anly}
with $-\Delta+n+1$ replaced by $-n_{1}-1, \dots, -n_{l}-1$, respectively, we
obtain an iterated series. This iterated series also has the corresponding multisum and
other iterated series. A~generalization of Lemma~\ref{convergence} says that
if~\eqref{gen-kkl} is absolutely convergent on the region~${|z_{1}|>\cdots > |z_{k}|>0}$
and its sum is equal to a~suitable single-valued branch of a~multivalued analytic function
on a~suitable subregion,
then the iterated series obtained by taking the sums in~\eqref{gen-88} first and then
taking the sums given by the expansions of ${\rm e}^{(-n_{1}-1)\log (z_1-z_{k})}(\log (z_1-z_{k})^{j_{1}}, \dots,
 {\rm e}^{(-n_{k-1}-1)\log (z_{k-1}-z_{k})}(\log (z_{k-1}-z_{k}))^{j_{l}}$ is absolutely
convergent on the region $|z_{1}|>\cdots >|z_{k}|>0$ and its sum is equal to a~given
single-valued branch on a~suitable subregion. Here we omit the details of this generalization and its proof.
\end{rema}

\subsection*{Acknowledgements}
We are grateful to the referees for their
detailed comments, questions, and suggestions on an early version of this paper.


\pdfbookmark[1]{References}{ref}
\LastPageEnding

\end{document}